\documentclass[reqno,a4paper]{amsart}
\usepackage{amssymb}
\usepackage[fontsize=12pt]{scrextend}
\usepackage{blindtext}
\usepackage[maxbibnames=99]{biblatex}
\usepackage{amsmath,amsthm,amsfonts,subfigure,hyperref,etoolbox,color,bm,soul,comment,breqn,mathtools,url}
\usepackage{tikz}
\usepackage{verbatim}
\usetikzlibrary{arrows,calc,patterns,shapes,decorations.markings}

\newcommand{\sgn}{\textrm{sgn}}
\newcommand{\ep}{\varepsilon}
\theoremstyle{plain}
\newtheorem{definition}{Definition}[section]
\newtheorem{theorem}[definition]{Theorem}
\newtheorem{lemma}[definition]{Lemma}

\newtheorem{proposition}[definition]{Proposition}

\newtheorem{corollary}[definition]{Corollary}
\newtheorem{remark}[definition]{Remark}

\allowdisplaybreaks

\usepackage[utf8]{inputenc}
\usepackage{fancyhdr,lipsum}

\usepackage[margin=1in, headheight=13pt]{geometry}

\DeclarePairedDelimiter\ceil{\lceil}{\rceil}

\usepackage{biblatex}
\AtEndEnvironment{definition}{\null\hfill\qedsymbol}%
\AtEndEnvironment{theorem}{\null\hfill\qedsymbol}%

\makeatletter
\def\l@subsection{\@tocline{2}{0pt}{3pc}{5pc}{}}
\def\l@subsubsection{\@tocline{2}{0pt}{6pc}{5pc}{}}
\makeatother

\numberwithin{figure}{section}
\numberwithin{equation}{section}

\begin{document}
\begin{titlepage}
\title[The Mathieu Differential Equation and Generalizations to Infinite Fractafolds]{The Mathieu Differential Equation and Generalizations to Infinite Fractafolds}
\author[Shiping Cao]{Shiping Cao$^1$}
\email{sc2873@cornell.edu}
\author[Anthony Coniglio]{Anthony Coniglio$^2$}
\email{coniglio@iu.edu}
\author[Xueyan Niu]{Xueyan Niu$^3$}
\email{xyniu@connect.hku.hk}
\author[Richard Rand]{Richard Rand$^4$}
\email{rhr2@cornell.edu}
\author[Robert Strichartz]{Robert Strichartz$^5$}
\email{str@math.cornell.edu}

\begin{abstract}
One of the well-studied equations in the theory of ODEs is the Mathieu differential equation. A common approach for obtaining solutions is to seek solutions via Fourier series by converting the equation into an infinite system of linear equations for the Fourier coefficients. We study the asymptotic behavior of these Fourier coefficients and discuss the ways in which to numerically approximate solutions. We present both theoretical and numerical results pertaining to the stability of the Mathieu differential equation and the properties of solutions.  Further, based on the idea of using Fourier series, we provide a method in which the Mathieu differential equation can be generalized to be defined on the infinite Sierpinski gasket. We discuss the stability of solutions to this fractal differential equation and describe further results concerning properties and behavior of these solutions.
\end{abstract}

\footnotetext[1]{Cornell University, mathematics graduate student (sc2873@cornell.edu)}
\footnotetext[2]{Indiana University Bloomington, class of 2019 (coniglio@iu.edu)}
\footnotetext[3]{University of Hong Kong, class of 2019 (xyniu@connect.hku.hk)}
\footnotetext[4]{Professor of Mathematics, Cornell University (rhr2@cornell.edu)}
\footnotetext[5]{Professor of Mathematics, Cornell University (str@math.cornell.edu)}

\end{titlepage}

\maketitle

\tableofcontents

\newpage

\section{Introduction}

The Mathieu differential equation, as will be defined in Section 2, takes the form $\frac{d^2u}{dt^2}+\left(\delta+\varepsilon\cos t\right)u=0$, where $\delta,\varepsilon\in\mathbb{R}$ are fixed, and $u:\mathbb{R}\to\mathbb{R}$. Named after French mathematician Émile Léonard Mathieu (1835-1890), the origin of the Mathieu differential equation stems from real-world phenomena. For example, it describes the motion of a pendulum subject to a periodic driving force. See \cite{stoker} for more details.

One topic pertaining to the Mathieu differential equation which has been much researched is the stability of solutions. See \cite{josephbarry,hochstadt,levykeller,loud,verticalconvergence}. Readers may also read \cite{stoker,randlecturenotes} for a brief introduction to this area.  Most of the existing literature concerns the parameter space, i.e. the space of $\delta$-$\ep$ pairs, whose choices can drastically alter the behavior of solutions. This paper presents research results and new phenomena both on the parameter space and on solutions.

On the other hand, another area of mathematics which has been actively researched in recent years is analysis on fractals, based on J. Kigami's construction of Laplacians on post-critically finite self-similar sets. See \cite{kigami1,kigami2,kigami,differentialequationsonfractals}.  It is of interest to define and study the Mathieu differential equation on fractal domains, and one suitable domain is the infinite Sierpinski gasket developed by R.S. Strichartz \cite{blowup}. Existing research on the infinite Sierpinski gasket and on other `fractafolds' be found in \cite{fractalfold1,fractalfold5,fractalfold4,fractalfold3,fractalfold2,teplyaev}. In particular, we will use the important result concerning the spectrum of the Laplacian on the infinite Sierpinski gasket by A. Teplyaev \cite{teplyaev}. We will explore a way to generalize the Mathieu differential equation to be defined on the infinite Sierpinski gasket.

This paper is organized as follows. In Section 2, we provide relevant background on the Mathieu differential equation, give definitions which will be used throughout the remainder of the paper, and describe some of the methods used to study the relationship that the values of $\delta$ and $\varepsilon$ have to the solutions of the equation; in doing so we will discuss `transition curves' and the `truncation method' used to study solutions to the Mathieu differential equation. In Section 3 we present a number of proofs for theorems pertinent to the Mathieu differential equation and provide justification for the methods presented in Section 2. Some of our methods and proofs are modified versions from \cite{asai,Agil}. In Section 4 we provide a discussion of computational results in studying the Mathieu differential equation, including results related to the asymptotic behavior of transition curves and the convergence of solutions. In Section 5 we give an overview of the Sierpinski gasket ($SG$) and provide definitions of various terms in fractal analysis, such as the fractal Laplacian and the infinite Sierpinski gasket ($SG_\infty$), to be used in the remaining sections. In Section 6 we extend the content of Sections 2, 3, and 4 to the fractal setting by explaining how the `Mathieu differential equation defined on the real line' can be generalized to a `Mathieu differential equation defined on the infinite Sierpinski gasket.' We describe how solutions to this generalization can be studied by considering solutions expanded as a linear combination of eigenfunctions of the fractal Laplacian. We also describe how the `truncation method' on the line can be used to study the Mathieu differential equation on $SG_\infty$. In Section 7 we provide computational results and observations about the shape and asymptotic behavior of the transition curves for the Mathieu differential equation on $SG_\infty$ and also about the behavior of solutions. In Section 8 we describe further research that can be done on the Mathieu differential equation and its fractal generalizations. Section 9, the Appendix, describes an alternate approach to the `truncation method' which involves partitioning Fourier coefficients into various equivalence classes.

A website for this research has been created at \mdseries\urlstyle{same}\url{http://pi.math.cornell.edu/~aac254/}.  We invite the reader to visit this website, as it contains plots, graphs, data, and other information gathered from the research.

\section{Definitions and Methods}
In this section we give formal definitions pertaining to the Mathieu differential equation on the real line and then describe the main methods used to study it. \begin{comment}\st{Theoretical justification for the use of these methods will be given in Section \ref{truncation}. These methods will be generalized to the fractafold setting in Section 6.}\end{comment}

\subsection{Background \& Definitions}

We begin by defining the \textit{Mathieu differential equation} on the real line.\\

\begin{definition}[Mathieu Differential Equation]
\label{mathieuequationdefinition}
The \underline{Mathieu differential equation} (on the real line) is defined as
\begin{equation}
\label{mathieuequationequation}
\frac{d^2u}{dt^2}+\left(\delta+\varepsilon\cos t\right)u=0,
\end{equation}
where $\delta,\varepsilon\in\mathbb{R}$ are fixed and $u:\mathbb{R}\to\mathbb{R}$ is unknown.
\end{definition}

%The Mathieu differential equation (often simply referred to as 'the Mathieu differential equation') was named for French mathematician Émile Léonard Mathieu (1835-1890). The Mathieu differential equation has a number of real-world applications, such as to the motion of a pendulum (see chapter 5 of [3] for a further discussion). 

Henceforth, the abbreviation `MDE' will often be used to denote `Mathieu differential equation.' 
Also, unless otherwise specified, the term `Mathieu differential equation' (and its abbreviation), when used in Sections 2, 3, and 4, will always refer to the MDE defined on the \textit{real line}, as opposed to the later sections which discuss the MDE as defined on a \textit{fractal} domain.

The particular values of $\delta$ and $\varepsilon$ chosen can drastically alter the corresponding solutions to the MDE (see \cite{stoker}, for example). With this in mind, we make the following definition.\\

\begin{definition}[Stable and Unstable $\delta$-$\varepsilon$ Pairs]
\label{stableandunstablepairsdefinition}

(a). An ordered pair $(\delta,\varepsilon)\in\mathbb{R}^2$ is a \underline{stable} pair of values if every solution to the corresponding MDE remains bounded for all $t\in\mathbb{R}.$

(b). An ordered pair $(\delta,\varepsilon)\in\mathbb{R}^2$ is an \underline{unstable} pair of values if there exists a solution to the corresponding MDE which is unbounded.\\
\end{definition}

\begin{figure}[h]
  \centering
  \includegraphics[scale=0.10]{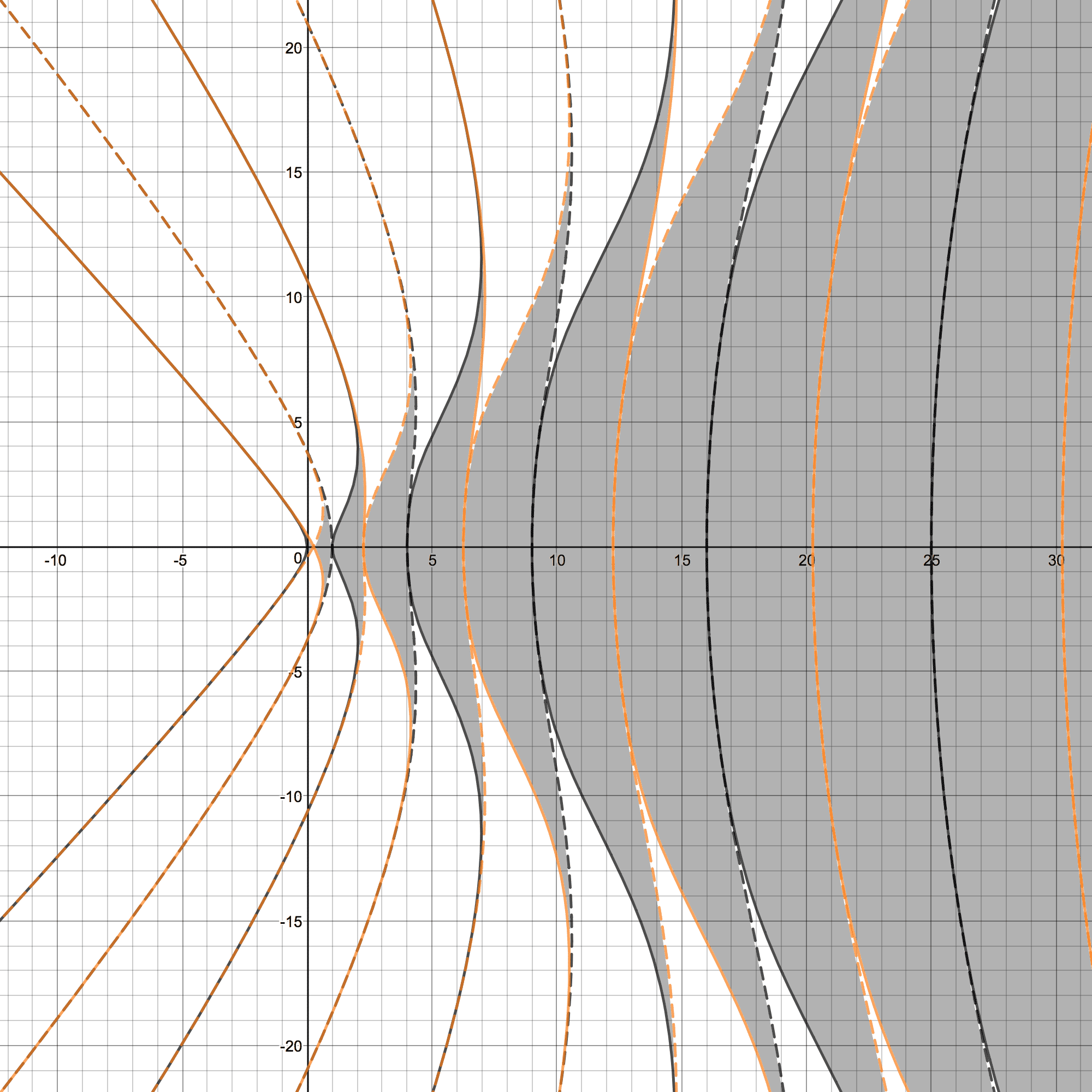}
  \caption{Stable and Unstable Regions of the $\delta$-$\varepsilon$ plane.}\label{deltaepsilonplanefigure}
\end{figure}

In order to understand which values of $\delta$ and $\varepsilon$ correspond to stable pairs and which correspond to unstable pairs, it is useful to examine the $\delta$-$\varepsilon$ plane given in Figure \ref{deltaepsilonplanefigure}. In this figure, the horizontal axis is the $\delta$-axis, and the vertical axis is the $\varepsilon$-axis. The gray shaded region of the plane corresponds to $(\delta,\varepsilon)$ pairs which are stable and is called the \textit{stable region}; the white region of the plane corresponds to $(\delta,\varepsilon)$ pairs which are unstable and is called the \textit{unstable region}. The solid curves and dashed curves which are either orange or black are called the \textit{transition curves} and form the boundary between the stable region and the unstable region. The curves which share the same color (either orange or black) and the same format (either solid or dashed) share certain properties in common which will be discussed in Section 2.2.

\label{sec:introduction}

There exists a systematic method for determining the stable and unstable regions shown in Figure \ref{deltaepsilonplanefigure}, which we describe in Section 2.2. A useful result regarding this method is given in \cite{stoker} as follows:\\

\begin{theorem}
\label{transitioncurvetheorem}
A pair $(\delta,\varepsilon)$ lies on a transition curve for the MDE if and only if there exists a corresponding nontrivial solution $u:\mathbb{R}\to\mathbb{R}$ which is periodic with period $2\pi$ or $4\pi$.
\end{theorem}

Theorem \ref{transitioncurvetheorem} motivates us to study the solutions to the MDE with Fourier series.

\begin{remark}
In this paper, we say $f$ is periodic of period $2N\pi$, or $f$ is $2N\pi$-periodic, if $N$ is the smallest positive integer such that $f$ admits the Fourier expansion
\[f(t)=\sum_{j=0}^{\infty} a_j\cos \left(\frac{j}{N}t\right)+\sum_{j=1}^{\infty} b_j\sin\left(\frac{j}{N}t\right).\]
\begin{comment}
\st{With this, if $f$ is not a trigonometric polynomial, $2N\pi$ is just the smallest period of $f$.} Of course, $f$ admits a Fourier expansion with $N$ replaced by $M$ such that $N|M$, \st{but a lot of term vanishes \textcolor{blue}{[?]}}.
\end{comment}
\end{remark}

%Thus, in order to determine the stable and unstable regions of plane, it suffices to Theorem \ref{transitioncurvetheorem} motivates us to study the solutions to the Mathieu differential equation with Fourier series.

%In this plot of the $\delta$-$\varepsilon$ plane, the horizontal axis is the $\delta$-axis, and the vertical axis is the $\varepsilon$-axis. The shaded region represents regions of the $\delta$-$\varepsilon$ plane corresponding to stable pairs of values, and the white regions correspond to unstable pairs of values. The curves that separate these two regions are called the \textit{transition curves}.

\begin{comment}
\begin{itemize}
\item Each transition curve crosses the $\delta$-axis exactly once, and does so at a point of the form \textcolor{green}{$(\frac{n^2}{4},0),$} with $n\in\mathbb{N}_{\geq0}$.

\item To each $\delta_n,$ where $n\geq1,$ there are exactly two transition curves that cross the $\delta$-axis through $(\frac{n^2}{4},0).$

\item The transition curves are symmetric about the $\delta$-axis.

\item Each transition curve asymptotes to have slope -1 as $\varepsilon\to\infty$ (and, hence, by the aforementioned symmetry of the transition curves, asymptotes to have slope +1 as $\varepsilon\to-\infty$). \textcolor{green}{See [,] for better estimation.}

\item No two transition curves intersect, except possibly at points of the form \textcolor{green}{$(\frac{n^2}{4},0).$} See []

\item Points on the transition curves where $\varepsilon=0$ are stable, whereas all other points on the transition curves are unstable. See []
\end{itemize}
\end{comment}

\subsection{Fourier Expansion and Matrix Form}

Let us suppose we have a periodic solution $u$ with period $2\pi$ or $4\pi$ to the Mathieu differential equation. In this case, we can write the Fourier series expansion of $u$ as
$$u(t)=\sum_{j=0}^\infty a_j\cos\left(\frac{j}{2} t\right)+ \sum_{j=1}^\infty b_j\sin\left(\frac{j}{2} t\right).$$
%If we assume $u$ is $2\pi$-periodic, then we would set $L=2\pi$, and if we assume $u$ is $4\pi$-periodic, then we would set $L=4\pi$. For the sake of generality, we let $L=2\pi N,$ $N\in\mathbb{N}$, in which case our Fourier expansion for $u$ becomes

%\begin{equation}
%   u(t)=\sum_{k=0}^\infty a_k\cos\left(\frac{k}{N}t\right)+ %\sum_{k=1}^\infty b_k\sin\left(\frac{k}{N}t\right). 
%\end{equation}

If we plug this Fourier series for $u$ into the MDE, rearrange terms, and use some trigonometric identities, we find that the function $u$ given above solves the MDE if and only if the two infinite systems of linear homogeneous equations shown below for the cosine and sine coefficients are satisfied. See \cite{randlecturenotes} and \cite{stoker} for more details.

\[
\text{cosine coefficients} \left\{
\begin{matrix}
    \delta a_0+\frac{\varepsilon}{2}a_2=0,\\
    \left(\delta-\left(\frac12\right)^2\right)a_1+\frac{\varepsilon}{2}\left(a_{3}+a_{1}\right)=0,\\
    (\delta-1)a_2+\ep a_0+\frac{\ep}{2}a_4=0,\\
\left(\delta-\left(\frac{j}2\right)^2\right)a_j+\frac{\varepsilon}{2}\left(a_{j-2}+a_{j+2}\right)=0,\quad(j\geq 3),\\ \\
 \end{matrix}\right.
\]

and

\[
\text{sine coefficients} \left\{
\begin{matrix}
    \left(\delta-\left(\frac12\right)^2\right)b_1+\frac{\varepsilon}{2}(b_3-b_1)=0,\\
    (\delta-1)b_2+\frac{\varepsilon}{2}b_4=0,\\
    \left(\delta-\left(\frac{j}2\right)^2\right)b_j+\frac{\varepsilon}{2}(b_{j-2}+b_{j+2})=0,\quad(j\geq 3).\\
 \end{matrix}
 \right.
\]

To solve these systems of equations, we consider putting them into matrix form. Readers can check that the above two systems of equations are equivalent to the following four equations in matrix form collectively, which correspond to cosine or sine coefficients with even or odd indexes. All the matrices are tridiagonal.

%We present the most important matrices when $N$ equals to 2, which correspond to the $2\pi$ or $4\pi$ periodic solutions on the transition curves. (\delta-\lambda_1)

For the cosine coefficients with even indexes, we have

$$
\underbrace{\begin{pmatrix}
        \delta & \frac{\varepsilon}{2} &\\ \varepsilon & \delta-1^2 & \frac{\varepsilon}{2} &\\
        & \frac{\varepsilon}{2} & \delta-2^2 & \frac{\varepsilon}{2}\\
        & & \frac{\varepsilon}{2} & \delta-3^2 & \frac{\varepsilon}{2}\\
        & & & \ddots & \ddots & \ddots
    \end{pmatrix}}_A
    \underbrace{\begin{pmatrix}
    a_0 \\ a_2 \\ a_4 \\ a_6 \\ \vdots
    \end{pmatrix}}_{\textbf{a}_{even}}
    =\underbrace{\begin{pmatrix}
    0 \\ 0 \\ 0 \\ 0 \\ \vdots
    \end{pmatrix}}_{\textbf{0}}.
    \\
$$

For the sine coefficients with even indexes, we have

$$\underbrace{\begin{pmatrix}
        \delta-1^2 & \frac{\varepsilon}{2} &\\ \frac{\varepsilon}{2} & \delta-2^2 & \frac{\varepsilon}{2} &\\
        & \frac{\varepsilon}{2} & \delta-3^2 & \frac{\varepsilon}{2}\\
        & & \ddots & \ddots & \ddots
    \end{pmatrix}}_{B}
    \underbrace{\begin{pmatrix}
    b_2 \\ b_4 \\ b_6 \\ \vdots
    \end{pmatrix}}_{\textbf{b}_{even}}
    =\begin{pmatrix}
    0 \\ 0 \\ 0 \\ \vdots
    \end{pmatrix}.$$

For cosine coefficients with odd indexes, we have

$$
\underbrace{\begin{pmatrix}
        \delta-\frac{1}{4}+\frac{\varepsilon}{2} & \frac{\varepsilon}{2} &\\ \varepsilon & \delta-\frac{9}{4} & \frac{\varepsilon}{2} &\\
        & \frac{\varepsilon}{2} & \delta-\frac{25}{4} & \frac{\varepsilon}{2}\\
        & & \frac{\varepsilon}{2} & \delta-\frac{49}{4} & \frac{\varepsilon}{2}\\
        & & & \ddots & \ddots & \ddots
    \end{pmatrix}}_C
    \underbrace{\begin{pmatrix}
    a_1 \\ a_3 \\ a_5 \\ a_7 \\ \vdots
    \end{pmatrix}}_{\textbf{a}_{odd}}
    =\begin{pmatrix}
    0 \\ 0 \\ 0 \\ 0 \\ \vdots
    \end{pmatrix}.
    \\
$$

For sine coefficients with odd indexes, we have

$$
\underbrace{\begin{pmatrix}
        \delta-\frac{1}{4}-\frac{\varepsilon}{2} & \frac{\varepsilon}{2} &\\ \varepsilon & \delta-\frac{9}{4} & \frac{\varepsilon}{2} &\\
        & \frac{\varepsilon}{2} & \delta-\frac{25}{4} & \frac{\varepsilon}{2}\\
        & & \frac{\varepsilon}{2} & \delta-\frac{49}{4} & \frac{\varepsilon}{2}\\
        & & & \ddots & \ddots & \ddots
    \end{pmatrix}}_D
    \underbrace{\begin{pmatrix}
    b_1 \\ b_3 \\ b_5 \\ b_7 \\ \vdots
    \end{pmatrix}}_{\textbf{b}_{odd}}
    =\begin{pmatrix}
    0 \\ 0 \\ 0 \\ 0 \\ \vdots
    \end{pmatrix}.
    \\
$$

We can solve the above four matrix equations separately. If either of the first two equations has a nontrivial solution in $\ell^2$, then the MDE has a $2\pi$-periodic solution, since
$$u(t):=\sum_{j=0,2,4,6,...} a_j\cos\left(\frac{j}{2} t\right)+\sum_{j=2,4,6,...} b_j\sin\left(\frac{j}{2} t\right)$$ solves the MDE and is $2\pi$-periodic, comparing with Remark 2.4.
\begin{comment}
\begin{align*}
    u(t)&=\sum_{j=0,2,4,6,...} a_j\cos\left(\frac{j}{2} t\right)+\sum_{j=2,4,6,...} b_j\sin\left(\frac{j}{2} t\right)\\
    &=\sum_{j=0,2,4,6,...} a_j\cos\left(\frac{j}{2} t+\pi j\right)+\sum_{j=2,4,6,...} b_j\sin\left(\frac{j}{2} t+\pi j\right)\\
    &=\sum_{j=0,2,4,6,...} a_j\cos\left(\frac{j}{2} (t+2\pi )\right)+\sum_{j=2,4,6,...} b_j\sin\left(\frac{j}{2} (t+2\pi )\right)\\
    &=u(t+2\pi)
\end{align*}
\end{comment}
Similarly, if either the third or fourth equations has a nontrivial solution in $\ell^2$, then the MDE has a $4\pi$-periodic solution, since
$$u(t):=\sum_{j=1,3,5,...} a_j\cos\left(\frac{j}{2} t\right)+\sum_{j=1,3,5,...} b_j\sin\left(\frac{j}{2} t\right)$$ solves the MDE and is $4\pi$-periodic, also comparing with Remark 2.4. Interested readers can also read the appendix for equations in matrix form for periodic solutions with larger periods, say $2N\pi$ with $N$ an arbitrary positive integer. 

%\st{which means $\textbf{a}_{even}$ or $\textbf{b}_{even}$ is not 0, then we get a $2\pi$-periodic solution to the MDE. Similarly, if the third and fourth equations have a nontrivial solution in $\ell^2$, which means $\textbf{a}_{odd}$ or $\textbf{b}_{odd}$ is not \textbf{0}, then we get a $4\pi$-periodic solution to the MDE.}
Now we return our discussion to Figure \ref{deltaepsilonplanefigure}, where we plot the stable and unstable regions. According to Theorem \ref{transitioncurvetheorem}, the transition curves consist of $(\delta,\ep)$ pairs with nontrivial $2\pi$- or $4\pi$-periodic solutions. Equivalently, by the above discussion, $(\delta,\ep)$ lies on a transition curve if and only if at least one of the four equations in matrix form discussed above has a nontrivial solution in $\ell^2$. Since we develop four different equations in matrix form, we use different colors (orange and black) and line formats (solid and dashed) in Figure \ref{deltaepsilonplanefigure} to distinguish between the four matrix equations:
\begin{itemize}
    \item If $(\delta,\ep)$ falls on a black solid line, then the equation $Ax=0$ has a nontrivial solution.
    \item If $(\delta,\ep)$ falls on a black dashed line, then the equation $Bx=0$ has a nontrivial solution.
    \item If $(\delta,\ep)$ falls on an orange solid line, then the equation $Cx=0$ has a nontrivial solution.
    \item If $(\delta,\ep)$ falls on an orange dashed line, then the equation $Dx=0$ has a nontrivial solution.
\end{itemize}

%In particular, \st{the matrices $A$ and $B$ correspond to solutions of period $2\pi$, while the other two matrices correspond to solutions of period $4\pi$.} nontrivial solutions with $\textbf{a}_{odd}=\textbf{b}_{odd}=\textbf{0}$ are $2\pi$ periodic, and nontrivial solutions with $\textbf{a}_{even}=\textbf{b}_{even}=\textbf{0}$ are $4\pi$ periodic. 

\subsection{Truncation Method}

In finite-dimensional linear algebra, a homogeneous matrix system of equations $Mx=0$, where $M$ is an $m\times m$ matrix, has a nontrivial solution if and only if the determinant of $M$ is zero. Although our four systems of matrix equations above are infinite, we hope to use techniques from finite-dimensional linear algebra to study the properties of these infinite systems.

For each infinite matrix, fix $m\in\mathbb{N}$ and consider its $m\times m$ leading principal submatrix. Below we give the truncated matrices $A_m$ and $B_m$ of $A$ and $B,$ respectively, as examples. $C_m$ and $D_m$ are defined similarly.

$$
A_m=\begin{pmatrix}
        \delta & \frac{\varepsilon}{2} &\\ \varepsilon & \delta-1^2 & \frac{\varepsilon}{2} &\\
        & \frac{\varepsilon}{2} & \delta-2^2 & \frac{\varepsilon}{2}\\
        & & \frac{\varepsilon}{2} & \delta-3^2 & \frac{\varepsilon}{2}\\
        & & & \ddots & \ddots & \ddots\\
        & & & & \frac{\varepsilon}{2} & \delta-(m-1)^2
    \end{pmatrix}
$$

$$
B_m=\begin{pmatrix}
        \delta-1^2 & \frac{\varepsilon}{2} &\\ \frac{\varepsilon}{2} & \delta-2^2 & \frac{\varepsilon}{2} &\\
        & \frac{\varepsilon}{2} & \delta-3^2 & \frac{\varepsilon}{2}\\
        & & \frac{\varepsilon}{2} & \delta-4^2 & \frac{\varepsilon}{2}\\
        & & & \ddots & \ddots & \ddots\\
        & & & & \frac{\varepsilon}{2} & \delta-m^2
    \end{pmatrix}
$$
Taking the determinant of each truncated matrix yields an algebraic expression involving $\delta$ and $\varepsilon$. Setting each expression equal to zero, we can then plot each equation in the $\delta$-$\varepsilon$ plane and obtain a set of algebraic curves in the variables $\delta$ and $\varepsilon$. If we choose $m$ sufficiently large, the $\delta$-$\varepsilon$ curves derived from the truncated matrices will be very close to the true transition curves corresponding to the infinite matrices. This statement is made precise by Ikebe et al. in \cite{asai}, and we present their results in Section 3.2.

\begin{comment}This \textgreen{proposition} does hold and was proved by Ikebe et al. in \cite{asai} for \textcolor{blue}{infinite} symmetric matrices, which we will show in \textcolor{blue}{Theorem 3.6} the next section.\end{comment}

\subsection{Backward Recursion}
In addition to studying the parameter space, it is also of interest to investigate the set of periodic solutions to the MDE. Motivated by the discussions in Section 2.2, we will study the Fourier coefficients of periodic solutions. For the computational results presented the paper, when numerically computing Fourier coefficients of solutions, instead of computing $a_2$, $a_4$, $a_6$,... successively given an initial value $a_0$ (and similarly for the other three matrix equations) using the recursion relation on Page 4, we will compute Fourier coefficients by setting initial value $a_{n_0}$ for some large index $n_0$ and compute $a_{n_0-2}$, $a_{n_0-4}$, $a_{n_0-6}$,... successively. The former method is called the `forward recursion method', whereas the latter method is called the `backward recursion method.' We choose to use the backward recursion method instead of the more commonly-used forward recursion method due to the instability of the forward recursion method, which will be discussed in Section 3.1. Also see \cite{Agil} for more details on the backward recursion method.

\begin{comment}
In addition to the transition curves, it is also of interest to investigate the set of solutions for the Mathieu differential equation. In particular, we will study the Fourier series solutions for the MDE. Before the end of this section, we point out that we will not use the forward recursion method in our computation for Fourier coefficients, which means we will not start our computation from $a_0$. Instead, we will set the initial value for a large index, and compute the coefficients with smaller indices by using the recursion relation backward. Related theorems will be stated in the next section.
\end{comment}

\begin{comment}
\begin{figure}
  \centering
%%%\includegraphics[scale = .35]{DesmosIntervalPlot}
\caption{Regions of stable pairs (shaded regions), unstable pairs (white regions), and the transition curves (source: \ref{})}
\end{figure}
\end{comment}

\begin{comment}
Recall that, for any given $\delta$ and $\varepsilon,$ specifying $a_0$ and $a_1$ uniquely determines the values of $a_k$ via equations \ref{}. Similarly, for any given $\delta$ and $\varepsilon$, specifying $b_1$ and $b_2$ uniquely determines the values of $b_k$ for all $k\geq3$ via equations \ref{}. Therefore, it is natural to wonder how the limiting behavior of the sequences $\{a_k\}_{k=0}^\infty$ and $\{b_k\}_{k=1}^\infty$, as $k\to\infty$, depends on $\delta$ and $\varepsilon$. We have the following theorem, which can be found in \ref{}, originally proven by Henri Poincar\'e.
\end{comment}

\section{Three-Term Recurrence Relations and Truncations of the Infinite Matrices}
\label{truncation}

\begin{comment}
\st{The ways in which we study the MDE on the line are twofold. First, we truncate each infinite matrix and take its determinant to obtain a relation between $\delta$ and $\varepsilon$ for which the MDE has a nontrivial periodic solution. Secondly, we compute the corresponding Fourier coefficients and study the behavior of the solution.}

\textcolor{blue}{[I think we should delete the above paragraph, because the following paragraph is somewhat similar to it and already serves as an introductory paragraph itself to Section 3]}
\end{comment}

In this section, we will provide theoretical foundation for the methods we use. In doing so we will discuss the asymptotic convergence of Fourier coefficients and provide an error estimate for the truncation method we use for approximating the transition curves. (We will extend these techniques to the fractal MDE in later sections.)

In this section, we work in generality by considering matrices of the form
\begin{equation}\label{eqn3}
\begin{pmatrix}
\delta-\lambda_1-\gamma\varepsilon & \alpha_1\varepsilon\\
\beta_1\varepsilon & \delta-\lambda_2 & \alpha_2\varepsilon\\
& \beta_2\varepsilon & \delta-\lambda_3 & \alpha_3\varepsilon\\
& & \beta_3\varepsilon & \delta-\lambda_4 & \ddots\\
& & &\ddots & \ddots
\end{pmatrix},
\end{equation}
where $\alpha_j,\beta_j,\lambda_j\in\mathbb{R}$ for $j\geq1$, and $\gamma\in\mathbb{R}$; further, we assume that $\{\alpha_j\}_{j=1}^\infty$ and $\{\beta_j\}_{j=1}^\infty$ are bounded sequences and that $\lim_{j\to\infty}\lambda_j=\infty$. Note that matrices $A,B,C,$ and $D$ as defined in Section 2.2 are special cases of (\ref{eqn3}). In addition, we always assume $\alpha_j,\beta_j\neq 0$, for all $j\geq 0$, in this section. We will consider equations of the form 
\begin{equation}\label{eqn1}\begin{pmatrix}
\delta-\lambda_1-\gamma\varepsilon & \alpha_1\varepsilon\\
\beta_1\varepsilon & \delta-\lambda_2 & \alpha_2\varepsilon\\
& \beta_2\varepsilon & \delta-\lambda_3 & \alpha_3\varepsilon\\
& & \beta_3\varepsilon & \delta-\lambda_4 & \ddots\\
& & &\ddots & \ddots 
\end{pmatrix}\begin{pmatrix}c_1\\c_2\\c_3\\c_4\\\vdots\end{pmatrix}=\begin{pmatrix}0\\0\\0\\0\\\vdots\end{pmatrix},\end{equation}
where $\bm{c}=(c_1,c_2,\cdots)^T\in\ell^2:=\left\{(x_1, x_2, ...)^T:\sum_{j=1}^\infty \left|x_j\right|^2<\infty\right\}$ is a real sequence.

\subsection{Asymptotic Behavior of $\{c_j\}_{j=1}^\infty$}
Three-term recurrence relations (TTRRs) and difference equations have been well-studied throughout the last century. Since (\ref{eqn1}) naturally gives a TTRR,
\begin{equation}\label{TTRRmathieu}
\begin{cases}
(\delta-\lambda_1-\gamma\varepsilon)c_1+\alpha_1\ep c_{2}=0\text{,}\\
\beta_{j-1}\ep c_{j-1}+(\delta-\lambda_j)c_j+\alpha_j\ep c_{j+1}=0\text{, for }j\geq 2,
\end{cases}
\end{equation}
we can use existing theorems to study the asymptotic behavior of the sequence $\{c_j\}_{j=1}^\infty$ as $j\to\infty$. Below we first present the well-known Poincar\'e Theorem (see \cite{Agil}), also called the Poincar\'e-Perron Theorem, which describes the asymptotic behavior of solutions to general TTRRs, and then we will see how the theorem applies to the equations above.\\

\begin{theorem} [Poincar\'e-Perron Theorem, \cite{Agil}] \label{TTRRPP}
For $j\geq1$, let $c_{j+1}+b_jc_j+a_jc_{j-1}=0$ with $\lim_{j\to\infty}b_j=b$ and $\lim_{j\to\infty}a_j=a$. Let $t_1$ and $t_2$ denote the zeros of the characteristic equation $t^2+bt+a=0$. Then, if $|t_1|\neq|t_2|$, the difference equation has two linearly independent solutions $\{x_n\}$ and $\{y_n\}$ satisfying 
\[\lim_{j\to\infty}\frac{x_j}{x_{j-1}}=t_1,\lim\limits_{j\to\infty}\frac{y_j}{y_{j-1}}=t_2.\]
If $|t_1|=|t_2|$, then $$\limsup_{j\to\infty}|c_j|^{\frac{1}{j}}=|t_1|$$ for any nontrivial solution $\{c_j\}$ to the recurrence relation.
\end{theorem}

With this theorem in mind, we consider the TTRR
\begin{equation}\label{eqn2}
f_jc_{j+1}+d_jc_j+g_{j-1}c_{j-1}=0 \quad\quad(j\geq1),
\end{equation}
where $\{f_j\}$ and $\{g_j\}$ neither of which contains $0$ as one of its terms, are uniformly bounded, and where $|d_j|\to\infty$ as $j\to\infty$. Notice that, with our assumptions, the solution of Equation \ref{eqn2} is uniquely determined by specifying the first two terms $c_1$ and $c_2$. We can derive the following proposition from the Poincar\'e-Perron Theorem, adapting its proof directly from the proof of the Poincar\'e-Perron Theorem in \cite{Agil}. \\
\begin{comment}
Notice that for equation, we actually will need to study the TTRR the TTRR we need to study is
\[\beta_{j-1}\varepsilon c_{j-1}+(\beta-\lambda_j)c_j+\alpha_j\varepsilon c_{j+1}=0,\]
which obviously satisfies the above condition is in the form of equation (3)
\end{comment}

\begin{proposition}\label{ttrrasym} Assume $\{f_j\}$ and $\{g_j\}$, neither of which is equal to the zero sequence, are uniformly bounded, and assume that $|d_j|\to \infty$ as $j\to\infty$. Then the TTRR (\ref{eqn2}) has two linearly independent solutions $\{x_j\}$ and $\{y_j\}$ satisfying
$$\frac{x_{j}}{x_{j-1}}\sim -\frac{g_{j-1}}{d_j},\quad \frac{y_{j}}{y_{j-1}}\sim -\frac{d_{j-1}}{f_{j-1}}.$$
\end{proposition}

\textit{Proof.} The proof is modified directly from the proof for the Poincar\'e-Perron Theorem in \cite{Agil}. For convenience, we assume $d_j\neq 0$ for any $j$, since otherwise we can consider the TTRR defined for $j\geq M$ for some large $M$. Let  $\hat{c}_j=\left(\prod_{i=1}^{j-1}\frac{f_i}{d_{i}}\right)c_j$. Then the TTRR becomes 
\[\hat{c}_{j+1}+\hat{c}_{j}+\frac{f_{j-1}g_{j-1}}{d_{j-1}d_j}\hat{c}_{j-1}=0.\]
By the Poincar\'e-Perron theorem, there exists a pair of linearly independent solutions of this new recurrence, $\{\hat{x}_j\}$ and $\{\hat{y}_j\}$, satisfying
\[\lim_{j\to\infty}\frac{\hat{x}_j}{\hat{x}_{j-1}}=0, \lim_{j\to\infty}\frac{\hat{y}_j}{\hat{y}_{j-1}}=-1.\]
For $j$ sufficiently large, $\hat{x}_j$ does not vanish. Indeed, if $\hat{x}_{j_0}=0$ for some $j_0\in\mathbb{N}$ then the TTRR implies $\frac{\hat{x}_{j_0+2}}{\hat{x}_{j_0+1}}=-1$, which cannot happen if $j_0$ were sufficiently large (since $\frac{\hat{x}_{j}}{\hat{x}_{j-1}}\to0$). Similarly, for $j$ sufficiently large $\hat{y}_j$ does not vanish, for if $\hat{y}_{j_0}=0$ for some $j_0\in\mathbb{N}$ then $\frac{\hat{y}_{j_0}}{\hat{y}_{j_0-1}}=0$, which cannot happen for $j_0$ sufficiently large (since $\frac{\hat{y}_j}{\hat{y}_{j-1}}\to-1$).

Let $H^{(1)}_j=\hat{x}_j/\hat{x}_{j-1}$ and $H^{(2)}_j=\hat{y}_j/\hat{y}_{j-1}$. Then 
\[H^{(i)}_{j+1}H^{(i)}_j+H^{(i)}_j+\frac{f_{j-1}g_{j-1}}{d_{j-1}d_j}=0,\quad i=1,2.\]
Multiplying each side by $H_j^{(l)},l\neq i$, we get 
\[\begin{cases}
  H^{(1)}_{j+1}H^{(1)}_jH^{(2)}_j+H^{(1)}_jH^{(2)}_j+\frac{f_{j-1}g_{j-1}}{d_{j-1}d_j}H^{(2)}_j=0,\\
  H^{(2)}_{j+1}H^{(2)}_jH^{(1)}_j+H^{(1)}_jH^{(2)}_j+\frac{f_{j-1}g_{j-1}}{d_{j-1}d_j}H^{(1)}_j=0.
\end{cases}\]
Subtracting one equation from the other, we obtain
\[H^{(1)}_jH^{(2)}_j\left(H^{(1)}_{j+1}-H^{(2)}_{j+1}\right)=\frac{f_{j-1}g_{j-1}}{d_{j-1}d_j}\left(H^{(1)}_j-H^{(2)}_j\right).\]
So, since
\[H_j^{(1)}\sim 0 \text{ and } H_j^{(2)}\sim -1,\] we have
\[H_j^{(1)}=\frac{f_{j-1}g_{j-1}}{d_{j-1}d_j}\frac{1}{H^{(2)}_j}\frac{H^{(1)}_j-H^{(2)}_j}{H^{(1)}_{j+1}-H^{(2)}_{j+1}}\sim -\frac{f_{j-1}g_{j-1}}{d_{j-1}d_j}.\]

Therefore,
$$\frac{x_j}{x_{j-1}}=\frac{\left(\prod_{i=1}^{j-1}\frac{d_i}{f_{i}}\right)\hat{x}_j}{\left(\prod_{i=1}^{j-2}\frac{d_i}{f_{i}}\right)\hat{x}_{j-1}}=H^{(1)}_j\frac{d_{j-1}}{f_{j-1}}\sim -\frac{g_{j-1}}{d_j},$$ and
$$\qquad\qquad\qquad\qquad\frac{y_j}{y_{j-1}}=\frac{\left(\prod_{i=1}^{j-1}\frac{d_i}{f_{i}}\right)\hat{y}_j}{\left(\prod_{i=1}^{j-2}\frac{d_i}{f_{i}}\right)\hat{y}_{j-1}}=H_j^{(2)}\frac{d_{j-1}}{f_{j-1}}\sim -\frac{d_{j-1}}{f_{j-1}}.\qquad\qquad\qquad\qquad\hfill\square$$

Let $\{z_j\}$ be an arbitrary solution to the TTRR \ref{eqn2}. We say $\{z_j\}$ is a \textit{minimal solution} of the TTRR if it obeys $\frac{z_{j}}{z_{j-1}}\sim -\frac{g_{j-1}}{d_j}$, and we say $\{z_j\}$ is a \textit{dominant solution} of the TTRR if it satisfies $\frac{z_{j}}{z_{j-1}}\sim -\frac{d_{j-1}}{f_{j-1}}$. Note that $\lim_{j\to\infty}z_j=0$ if $\{z_j\}$ is a minimal solution and $\lim_{j\to\infty}z_j=\infty$ if $\{z_j\}$ is a dominant solution.

We can apply Proposition \ref{ttrrasym} to study Equation 3.3, since Equation 3.3 consists of an equation for $c_1,c_2$ and a TTRR of the form \ref{eqn2}.  

\begin{comment}
$\{x_j\}_{j=1}^\infty$ and $\{y_j\}_{j=1}^\infty$, such that
$$\frac{x_{j}}{x_{j-1}}\sim -\frac{\alpha_{j-1}\varepsilon}{\beta-\lambda_j},\quad \frac{y_{j}}{y_{j-1}}\sim -\frac{\beta-\lambda_{j-1}}{\beta_{j-1}\varepsilon}.$$
\end{comment}

\begin{corollary}\label{coasym} Assume $\ep\neq 0$.

(a).  Equation 3.2 has a unique solution up to multiplication by a constant.

(b). Let $\bm{c}=(c_1,c_2,c_3,...)^T$ be a nontrivial solution to Equation (3.2). If $\bm{c}\in\ell^2$, then $c_j$ will decay with the asymptotic behavior $\frac{c_j}{c_{j-1}}\sim \frac{\beta_{j-1}\varepsilon}{\lambda_j-\delta}$. If $\bm{c}\notin\ell^2$, then $c_j$ will diverge with the asymptotic behavior $\frac{c_{j}}{c_{j-1}}\sim\frac{\lambda_{j-1}-\delta}{\alpha_{j-1}\varepsilon}$.
\end{corollary}
\textit{Proof.}  It is equivalent to studying Equation \ref{TTRRmathieu}. 

(a). Given $c_1\in\mathbb{R}$, we can solve $c_2,c_3,\cdots$ inductively, which gives us a solution to Equation \ref{TTRRmathieu}. In this way, any solution to Equation \ref{TTRRmathieu} is uniquely determined by $c_1$.

(b) Any solution $\textbf{c}=\{c_j\}_{j=1}^\infty$ can be written as a linear combination of $\{x_j\}$ and $\{y_j\}$ as given in in Proposition \ref{ttrrasym}, and hence there exist $p,q\in\mathbb{R}$ such that $c_j=px_j+qy_j$ for all $j$. Note that if $q\neq0$ then $\{z_j\}$ is a dominant solution, and if $q=0$ then $\{z_j\}$ is a minimal solution.

If $\{c_j\}\in\ell^2$ then $\lim_{j\to\infty}c_j=0$ and thus $\{c_j\}$ cannot be a dominant solution; thus, we must have $p=0$ and hence $\frac{c_j}{c_{j-1}}\sim \frac{\beta_{j-1}\varepsilon}{\lambda_j-\delta}$. On the other hand, if $\{c_j\}\notin\ell^2$ then $p\neq0$ and hence $\frac{c_{j}}{c_{j-1}}\sim\frac{\lambda_{j-1}-\delta}{\alpha_{j-1}\varepsilon}$.

%Let $\{z_j\}$ be an arbitrary solution to the TTRR \ref{eqn2}. Any two solutions $\{x_j\}$ and $\{y_j\}$ as in Proposition \ref{ttrrasym} form a basis of the set of solutions to Equation \ref{eqn2}, and hence there exist $p,q\in\mathbb{R}$ such that $z_j=px_j+qy_j$ for all $j$. We say $\{z_j\}$ is a \textit{minimal solution} of the TTRR if it obeys $\frac{z_{j}}{z_{j-1}}\sim -\frac{g_{j-1}}{d_j}$, and we say $\{z_j\}$ is a \textit{dominant solution} of the TTRR if it satisfies $\frac{z_{j}}{z_{j-1}}\sim -\frac{d_{j-1}}{f_{j-1}}$. If $q\neq0$ then $\{z_j\}$ is a dominant solution, and of $q=0$ then $\{z_j\}$ is a minimal solution.}

\hfill$\square$\\

Corollary \ref{coasym} shows that there is a sharp contrast between the two possible types of behavior that a solution to Equation (\ref{eqn1}) can have. Since $\{\alpha_j\}$ and $\{\beta_j\}$ are bounded sequences, when $(\delta,\varepsilon)$ are fixed and properly chosen any solution $\{c_j\}$ to (\ref{eqn1}) will converge to $0$ very rapidly; otherwise, all nontrivial solutions will tend to infinity very fast.

If $c_1\in\mathbb{R}$ is an initial value which corresponds to a minimal solution, Proposition \ref{ttrrasym} shows that numerically computing the solution with the forward recursion method is unstable, since a small error in computation will lead to a dominant solution, in which $\{c_j\}$ explodes. On the other hand, the backward recursion method can give a good approximation for a minimal solution. In \cite{Agil}, detailed discussions are given on this method, as well as a corresponding error estimate. We briefly state the result here.

\begin{proposition}[\cite{Agil}] Fix a large $N$, and set $c_{N+1}=0, c_{N}=1$. Compute $c_{N-1},c_{N-2},\cdots,c_1$ one-by-one backward with TTRR (\ref{eqn2}). If $\{x_k\}$ is an arbitrary minimal solution and $\{y_k\}$ is an arbitrary dominant solution, then $\{c_j\}$ satisfy
\[\frac{c_{k}}{c_{k-1}}\left( \frac{x_{k}}{x_{k-1}}\right)^{-1}-1=\frac{r_N}{r_{k-1}}\cdot\frac{1-r_k/r_{k-1}}{1-r_N/r_{k-1}},\]
where $r_k=\frac{x_k}{y_k}$ and $1\leq k\leq N$. %\st{Here, the definitions of $x_k$ and $y_k$ are as given in Proposition \ref{ttrrasym}.}
\end{proposition}

\subsection{The Truncation Method}
In this part, we introduce the truncation method for determining which $(\delta,\varepsilon)$ pairs yield nontrivial solutions of (\ref{eqn1}). We start from the observation 
%\begin{dmath*}
\begin{align*}
\begin{pmatrix}
\delta-\lambda_1-\gamma\varepsilon & \alpha_1\varepsilon\\
\beta_1\varepsilon & \delta-\lambda_2 & \alpha_2\varepsilon\\
& \beta_2\varepsilon & \delta-\lambda_3 & \alpha_3\varepsilon\\
& & \beta_3\varepsilon & \ddots & \ddots\\
& & &\ddots 
\end{pmatrix}&=\delta I-\begin{pmatrix}
\lambda_1+\gamma\varepsilon & -\alpha_1\varepsilon\\
-\beta_1\varepsilon & \lambda_2 & -\alpha_2\varepsilon\\
& -\beta_2\varepsilon & \lambda_3 & -\alpha_3\varepsilon\\
& & -\beta_3\varepsilon & \ddots & \ddots\\
& & &\ddots
\end{pmatrix}\\
&=\delta I-T(\varepsilon),
\end{align*}
%\end{dmath*}
where $I$ is the identity matrix and the shorthand notation $T(\varepsilon)$ is used in the last equality. Clearly, for a fixed $\varepsilon$, equation (\ref{eqn1}) has nontrivial solutions if and only if $\delta$ is an eigenvalue of $T(\varepsilon)$. 

The eigenvalue problem for infinite tridiagonal matrices has been widely studied. In \cite{asai}, there is a result on the error between eigenvalues of the truncated matrices and those of the corresponding infinite matrix. We state the result in Theorem \ref{symtrunc} below.

Here we consider the infinite symmetric tridiagonal matrix $T$ of the form
\[T=\begin{pmatrix}
d_1 & f_1\\
f_1 & d_2 & f_2\\
& f_2 & d_3 & f_3\\
& & f_3 & \ddots & \ddots\\
& & &\ddots 
\end{pmatrix},\]
where $d_n\to\infty$ as $n\to\infty$ and $\{f_n\}$ bounded. Denote its $n\times n$ truncation by $T_n$, i.e., 
\[T_n=\begin{pmatrix}
d_1 & f_1\\
f_1 & d_2 & f_2\\
& \ddots &\ddots & \ddots\\
& & f_{n-1}&d_n 
\end{pmatrix}.\]

\begin{theorem}\label{symtrunc} (\cite{asai}) Let $T$ and $T_n$ $(n\geq1)$ be given as above.\\
(a). $T$ has pure point spectrum.\\
(b). If $\delta$ is a given simple eigenvalue of $T$, then there exists, for each $n\in\mathbb{N}$, an eigenvalue $l_n$ of $T_n$ such that the sequence $\{l_n\}_{n=1}^\infty$ satisfies $l_n\to\delta$ as $n\to\infty$. For any such sequence the error is given by 
\[l_n-\delta=\frac{f_{n+1}c_nc_{n+1}}{\bm{c}^T\bm{c}}(1+o(1)),\]
where $\bm{c}=(c_1,c_2,\cdots)^T\in \ell^2$ is an eigenvector corresponding to $\delta$.
\end{theorem}

We can extend the above result to nonsymmetric tridiagonal matrices.\\

\begin{theorem} Let $T$ be an infinite tridiagonal matrix of the form
\[T=\begin{pmatrix}
d_1 & f_1\\
g_1 & d_2 & f_2\\
& g_2 & d_3 & f_3\\
& & g_3 & \ddots & \ddots\\
& & &\ddots 
\end{pmatrix},\]
where $d_n\to\infty$ as $n\to\infty$, and $\{f_n\}_{n=1}^\infty$ and $\{g_n\}_{n=1}^\infty$ are bounded, positive, and nonzero. Let $T_n$ be the $n$ $\times$ $n$ truncation of $T$. If $\delta$ is a given simple eigenvalue of $T$, then there exists, for each $n\in\mathbb{N}$, an eigenvalue $l_n$ of $T_n$ such that the sequence $\{l_n\}_{n=1}^\infty$ satisfies $l_n\to\delta$ as $n\to\infty$. For any such sequence the error is given by
$$l_n-\delta=\frac{\sqrt{f_{n+1}g_{n+1}}c_nc_{n+1}}{\sum_{j=1}^\infty \kappa_jc_j^2}(1+o(1)),$$
where $\kappa_j:=\sqrt{\prod_{i=1}^{j-1}\frac{g_i}{f_i}}$ and where $\bm{c}=(c_1,c_2,\cdots)^T\in \ell^2$ is an eigenvector corresponding to $\delta$.
\end{theorem}

\textit{Proof.} (a) Let 
\[T'=\begin{pmatrix}
d_1 & \sqrt{f_1g_1}\\
\sqrt{f_1g_1} & d_2 & \sqrt{f_2g_2}\\
& \sqrt{f_2g_2} & d_3 & \sqrt{f_3g_3}\\
& & \sqrt{f_3g_3} & \ddots & \ddots\\
& & &\ddots
\end{pmatrix}.\]
 We observe that
\[T=\begin{pmatrix}
\kappa_1 \\
& \kappa_2\\
& & \kappa_3\\
& & &\ddots
\end{pmatrix}\cdot T'\cdot \begin{pmatrix}
\kappa^{-1}_1 \\
& \kappa^{-1}_2\\
& & \kappa^{-1}_3\\
& & &\ddots
\end{pmatrix}.\]
It is clear that $T'$ is a self-adjoint operator from $\ell^2$ to $\ell^2$ with pure point spectrum by Theorem \ref{symtrunc}.

Let
\[\bm{c}'=diag(\kappa^{-1}_1,\kappa^{-1}_2,\cdots)\bm{c}=(\kappa^{-1}_1c_1,\kappa^{-1}_2c_2,\cdots)^T,\]
where $\bm{c}=(c_1,c_2,\cdots)^T$. Clearly, the eigenvalue equation $(\delta I-T)\bm{c}=0$ holds if and only if $(\delta I-T')\bm{c}'=0$ in a pointwise sense, where $I$ is the identity transformation.

In addition, from Corollary \ref{coasym} we deduce that there is a one-to-one correspondence between eigenvectors of $T$ and $T'$ in $\ell^2$. Indeed, if $\bm{c}$ is an eigenvector of $T$ in $\ell^2$, we have $\frac{c'_j}{c'_{j-1}}=\frac{\kappa^{-1}_jc_j}{\kappa^{-1}_{j-1}c_{j-1}}\sim-\frac{\sqrt{f_{j-1}g_{j-1}}}{d_j}$ because $\frac{c_j}{c_{j-1}}\sim -\frac{g_{j-1}}{d_j}$ by Corollary \ref{coasym}, which implies $\bm{c}'$ is an eigenvector of $T'$ in $\ell^2$, and for the same reason if $\bm{c}'\in\ell^2$ is an eigenvector of $T'$ then $\bm{c}=diag(\kappa_1,\kappa_2,\cdots)\bm{c'}$ is an eigenvector of $T$ in $\ell^2$.

Now, using Theorem \ref{symtrunc} we deduce that there is a sequence of eigenvalues of the truncated $T'_n$ that converges to the eigenvalue $\delta$ of $T'$ and hence of $T$, since $T$ and $T'$ have the same eigenvalues. In addition, noticing that, for each $n\in\mathbb{N}$, $T_n$ has the same eigenvalue as $T'_n$, we see that there is a sequence of eigenvalues of $T_n$ which converges to $\delta$. Lastly, the error estimate comes by applying Theorem 3.5 to $T'$ with the eigenvector $\bm{c}'=\text{diag}(\kappa^{-1}_1,\kappa^{-1}_2,\cdots)\bm{c}$.\hfill$\square$

\section{Observations and Analysis for the Mathieu Differential Equation on the Line}

The Mathieu differential equation has been an important topic in differential equations due to its numerous real-world applications. However, most of the existing work on the MDE focuses on theoretical analysis of the stable and unstable regions of the $\delta$-$\varepsilon$ plane, such as the asymptotic behavior of the transition curves, and there has not been much computational work done on the stability curves as well as on the solutions to the Mathieu differential equation. In this section, in addition to describing some of our theoretical results, we will explain our computational results on the intricate shape of the $\delta$-$\varepsilon$ curves and on the converging behavior of the solutions.

\subsection{The $\delta$-$\varepsilon$ Plot}

%Now that some theoretical background has been given in the previous section, we will now discuss some results.

%As discussed previously, the transition curves are determined through the following process. First, one expands the function $u$ in terms of a ($2\pi$-periodic or $4\pi$-periodic) Fourier expansion. Then, one plugs this into the Mathieu differential equation and finds two infinite systems of linear homogeneous equations that must be satisfied by the cosine and sine coefficients of the Fourier expansion. One then puts these equations into matrix form, takes the determinant of a truncated matrix, and sets the resulting expression equal to 0. This yields the implicit equation involving $\delta$ and $\varepsilon$ which produces the transition curves in the $\delta$-$\varepsilon$ plane and thus determines the regions of stability and the regions of instability. A figure is shown below.

In this part, we use the truncation method introduced in Section 2 and Section 3 to study the $\delta$-$\varepsilon$ curves in more detail.

%Figure \ref{dyadic} shows the $\delta$-$\varepsilon$ plane with transition curves, as well as with the stable and unstable regions. Here we truncate each infinite matrix to be size $20\times 20$.

%It can be shown that the transition curves for the Mathieu differential equation (resulting from $2\pi$- and $4\pi$-periodic expansion) share a number of interesting properties. We state below, without proof, some of these properties (see [3] for a further discussion):

%Hence, from this figure it is easy to see that, once one determines where the transition curves lie on the $\delta\varepsilon$-plane, one can determine the stability for all other $(\delta,\varepsilon)$ in the plane since they are, roughly speaking, separated into regions `between' the transition curves.

\begin{figure}[h]
  \centering
  \includegraphics[scale=0.15]{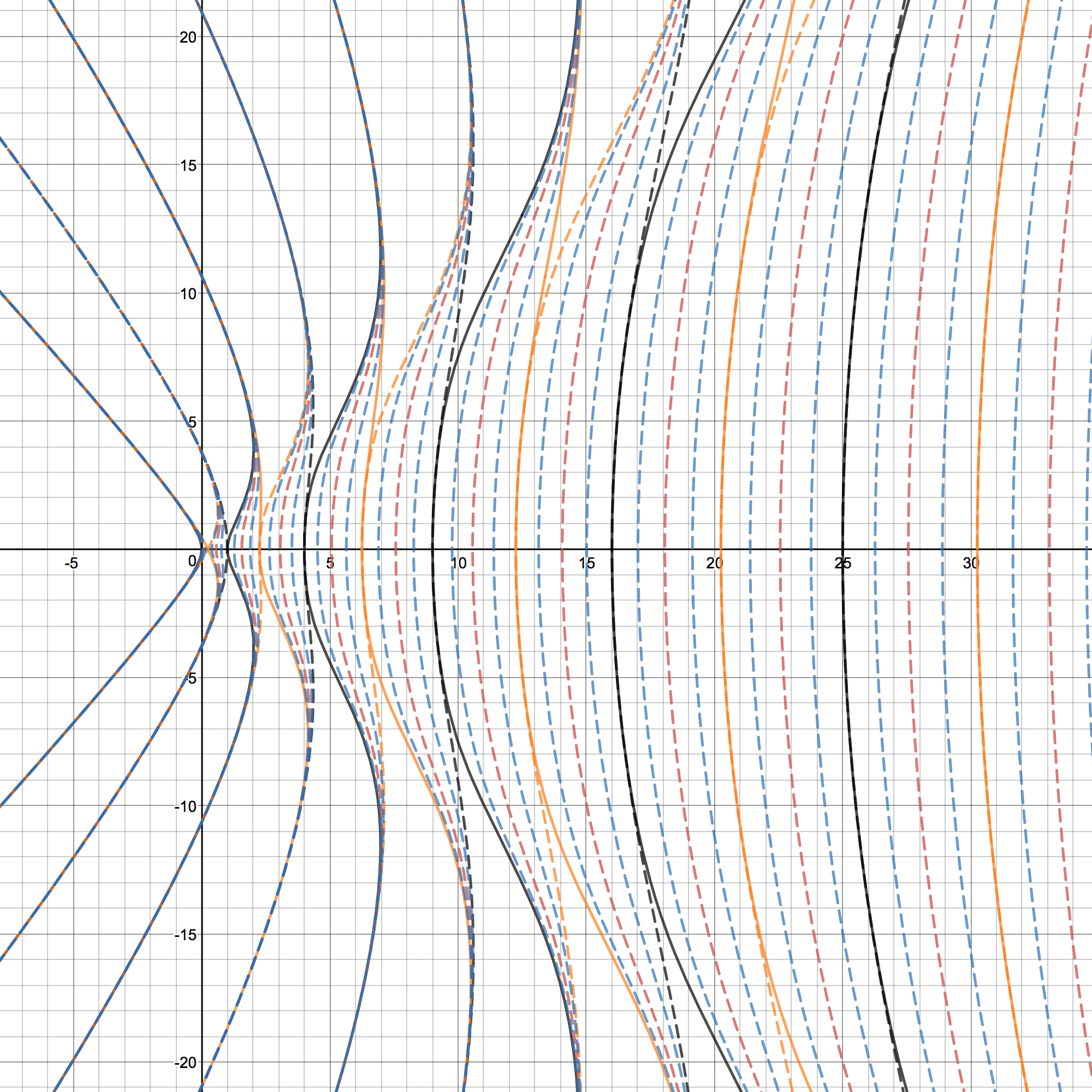}
  \caption{Curves corresponding to expansions of larger periods: curves with $2\pi$-periodic solutions (black solid and black dashed), curves with $4\pi$-periodic solutions (orange solid and orange dashed), curves with $8\pi$-periodic solutions (red dashed) and curves with $16\pi$-periodic solutions (blue dashed).}\label{dyadic}
\end{figure}

As discussed above, the transition curves are found via $2\pi$-periodic and $4\pi$-periodic Fourier expansions of solutions. However, one might wonder what the curves would be if the same process is undertaken for Fourier expansions of \textit{larger} periods, say periods $8\pi,16\pi,32\pi,...$. The answer is that curves corresponding to expansions of larger periods in fact "fill in" the stable regions whose boundary is obtained using the $2\pi$- and $4\pi$-periodic expansions. See Figure \ref{dyadic} for an illustration. \\

In fact, this ``fill in" property remains valid for even larger periods, and a proof is given below.

\begin{proposition}
Let $D$ be the set of $\delta$-$\varepsilon$ pairs such that the corresponding MDE has a solution of period $2^k\pi$ for some $k\geq 3$. Then $D$ is dense in the stable region.
\end{proposition}

\textit{Proof.} First, if $(\delta,\varepsilon)$ is a stable pair, then by the discussion given in Chapter VI of \cite{stoker} we can find two linearly independent solutions $u_1$ and $u_2$ of the MDE, such that 
\[\begin{pmatrix} u_1(2\pi)\\u'_1(2\pi)\end{pmatrix}=e^{i2\pi \theta}\begin{pmatrix}u_1(0)\\u'_1(0)\end{pmatrix},\quad \begin{pmatrix}u_2(2\pi)\\u'_2(2\pi)\end{pmatrix}=e^{-i2\pi \theta}\begin{pmatrix}u_2(0)\\u'_2(0)\end{pmatrix},\]
for some $\theta\in [0,\frac{1}{2}]$. $\theta$ is uniquely determined by the pair $(\delta,\ep)$, so we may view $\theta=\theta(\delta,\ep)$ as a function from the stable region to $[0,1/2]$.

In fact, by the Floquet Theorem, there exists a $2\times 2$ matrix $A$ depending only on $(\delta,\ep)$ such that, for any solution $u$ of the MDE,
\[\begin{pmatrix}
u(t+2\pi)\\u'(t+2\pi)
\end{pmatrix}=A\begin{pmatrix}
u(t)\\u'(t)
\end{pmatrix}.\]
So $(\delta,\ep)$ is a stable pair only if the two eigenvalues of $A$ both have norm $1$, which can be written as $e^{i2\pi\theta}$ and $e^{-i2\pi\theta}$. Readers can find details in \cite{stoker}. Clearly, $A$ depends continuously on the parameter $(\delta,\ep)$, which shows that $\theta=\theta(\delta,\ep)$ is a continuous function on the stable region.

Next, we fix $\theta_0\in[0,\frac{1}{2}]$ and $\ep\in\mathbb{R}$, and show that only countably many $\delta$'s can be found such that $\theta_0=\theta(\delta,\ep)$. Recall that for each such $\delta$ we have a solution $u_1$ to the corresponding MDE so that $e^{-i\theta_0 t}u_1(t)$ is a $2\pi$-periodic function. We have the Fourier series expansion $e^{-i\theta_0 t}u_1(t)=\sum_{j=-\infty}^\infty c_je^{ijt}$, which implies 
\[u_1(t)=e^{i\theta_0 t}\sum_{j=-\infty}^\infty c_je^{ijt}.\]
Plugging this into the Mathieu differential equation, we get the following equation for the coefficients, 
\[
\begin{pmatrix}
\quad\ddots & \ddots & \ddots\\
    & & \delta-(\theta_0-2)^2       & \frac{\varepsilon}{2} &  & \\ 
    & & \frac{\varepsilon}{2} & \delta-(\theta_0-1)^2 & \frac{\varepsilon}{2}   &   &  \\
           & & & \frac{\varepsilon}{2} & \delta-\theta_0^2 & \frac{\varepsilon}{2} &  &  \\
         & & & & \frac{\varepsilon}{2} & \delta-(\theta_0+1)^2 & \frac{\varepsilon}{2} &  &  \\
           & & & & & \ddots & \ddots &\ddots\\
\end{pmatrix}
\begin{pmatrix}
    \vdots \\
    c_{-2} \\
    c_{-1} \\
    c_0 \\
    c_1 \\
    c_2 \\
    \vdots \\
\end{pmatrix}
=\begin{pmatrix}
    \vdots \\
    0 \\
    0 \\
    0 \\
    0 \\
    0 \\
    \vdots \\
\end{pmatrix}.
\]
From this we can see that there are only countably many values of $\delta$ such that the corresponding equation has a nontrivial solution $\{c_j\}_{j=-\infty}^\infty \in \{\bm{c}\in\mathbb{C}^{\mathbb{Z}}|\sum_{j=-\infty}^\infty |c_j|^2<\infty\}$. 
%\textcolor{blue}{[What does $w$ stand for?]}\textcolor{green}{[doubly, just a notation..., dinstinct from $\ell^2$]} In fact, when $\delta=0$, the matrix as an operator $\ell^2_w\to \ell^2_w$ has a compact inverse. 

%On the other hand, we can also fix a stable pair $(\delta,\ep)$. In view of the eigenvalues of the transformation matrix $A$ and their uniqueness, $\theta(\delta,\ep)$ becomes a function from the stable region to $[0,1/2]$. Obviously, $\theta(\delta,\ep)$ is continuous on the stable region, since $A$ depends continuously on the parameters $(\delta,\ep)$. See \cite{stoker} for discussions on $A$.

Now, we prove the proposition by contradiction. Suppose $D$ is not dense in the stable region. Then, there is a ball in the $\delta$-$\ep$ plane where $\theta$ does not take any dyadic rational value, noticing that when $\theta$ is a dyadic rational the MDE has $2^k\pi$-periodic solutions $u_1,u_2$ for some $k$. This means $\theta$ equals a constant $\theta_0$ on that ball, since $\theta$ is a continuous function of $(\delta,\ep)$. However, for this fixed $\theta=\theta_0$, and a fixed $\ep$ in the ball, there are only countably many possible values for $\delta$, which can not fill in the ball. This gives a contradiction.\hfill$\square$\\

Next we will discuss the asymptotic behavior of the stable and unstable regions. Existing works including \cite{hochstadt}, \cite{levykeller}, and \cite{verticalconvergence} cover two important observations:

\begin{enumerate}
    \item First, as illustrated in Figure \ref{vertical}, the width of each stable band gets thinner and thinner as $|\varepsilon|$ tends to $\infty$. In \cite{verticalconvergence}, it is shown that the width of the $k$th stable band will decrease exponentially on $|\varepsilon|\geq k^2$ as $|\varepsilon|\to\infty$.
\begin{figure}[h]
  \centering
  \includegraphics[scale=0.37]{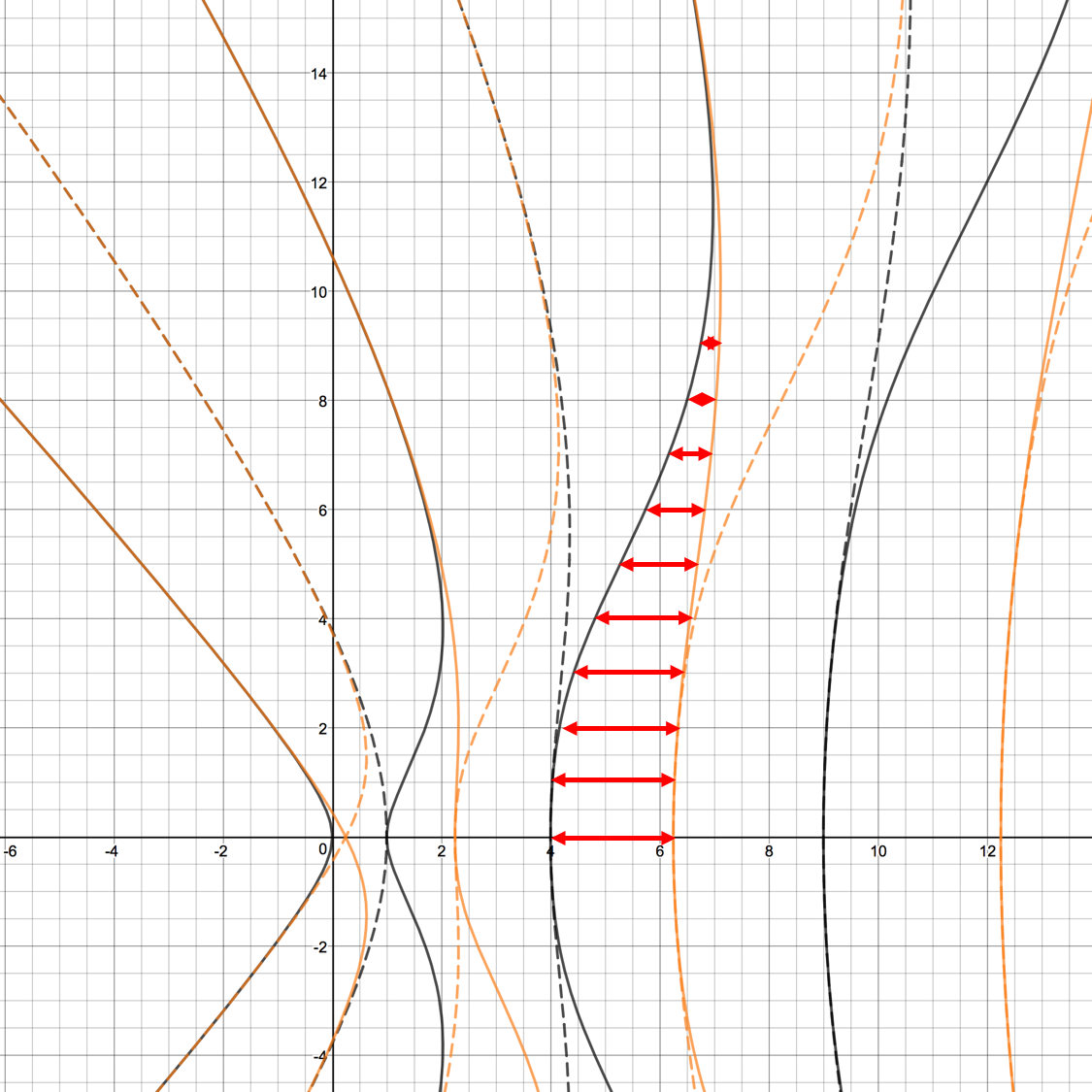}
  \caption{The width of the $5$th stable band.}\label{vertical}
\end{figure}
\item Secondly, for each fixed $\varepsilon$, the width of $k$th unstable band becomes smaller as $k$ increases. An estimate is given in \cite{josephbarry} and \cite{hochstadt} with different methods, both yielding
\[d_k=\frac{2|\varepsilon|^k}{2^m\big((k-1)!\big)^2}\left(1+o\left(\frac{|\ep|^2}{k^2}\right)\right),\]where $d_k$ is the width of the $k$th unstable band with $\varepsilon$ fixed.
\end{enumerate}

\begin{figure}
\centering
Fitting curve of the form: $\frac{a}{b+e^{cx}}$.\qquad\qquad\qquad\qquad\qquad\qquad\qquad\qquad\qquad\qquad\vspace{0.15cm}

\includegraphics[scale=0.45]{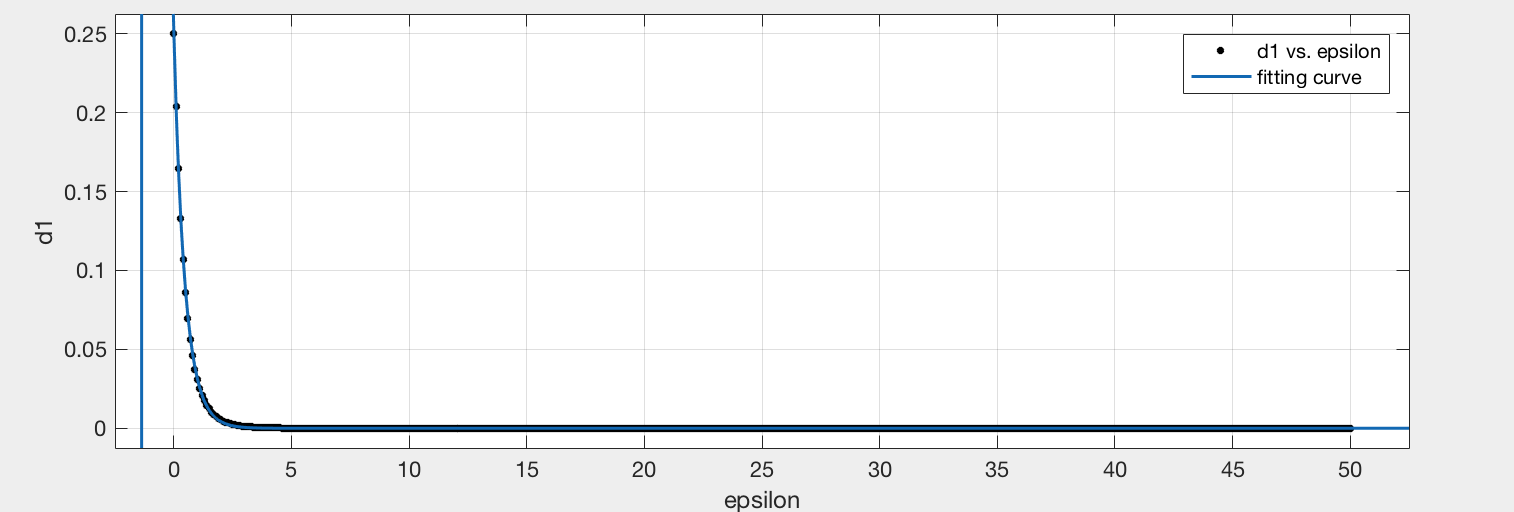}

$1$st stable band, $a=0.2278,b=-0.09431,c=1.994$\vspace{0.15cm}

\includegraphics[scale=0.45]{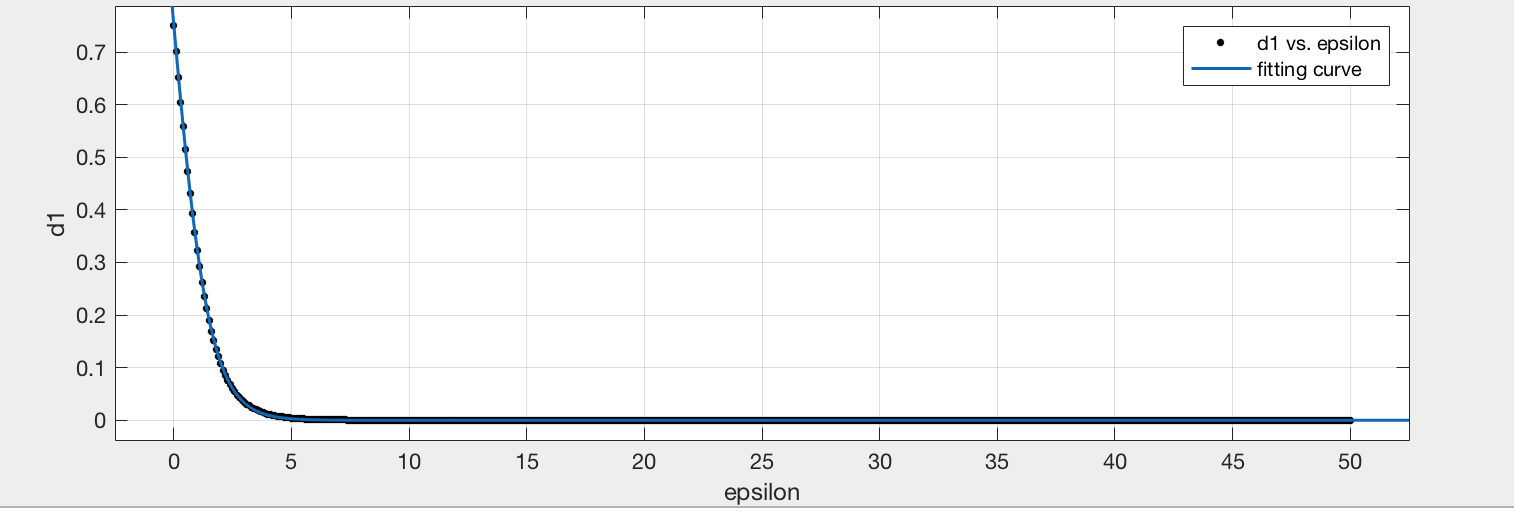}

$2$nd stable band, $a=1.390,b=0.8504,c=1.24$\vspace{0.15cm}

\includegraphics[scale=0.45]{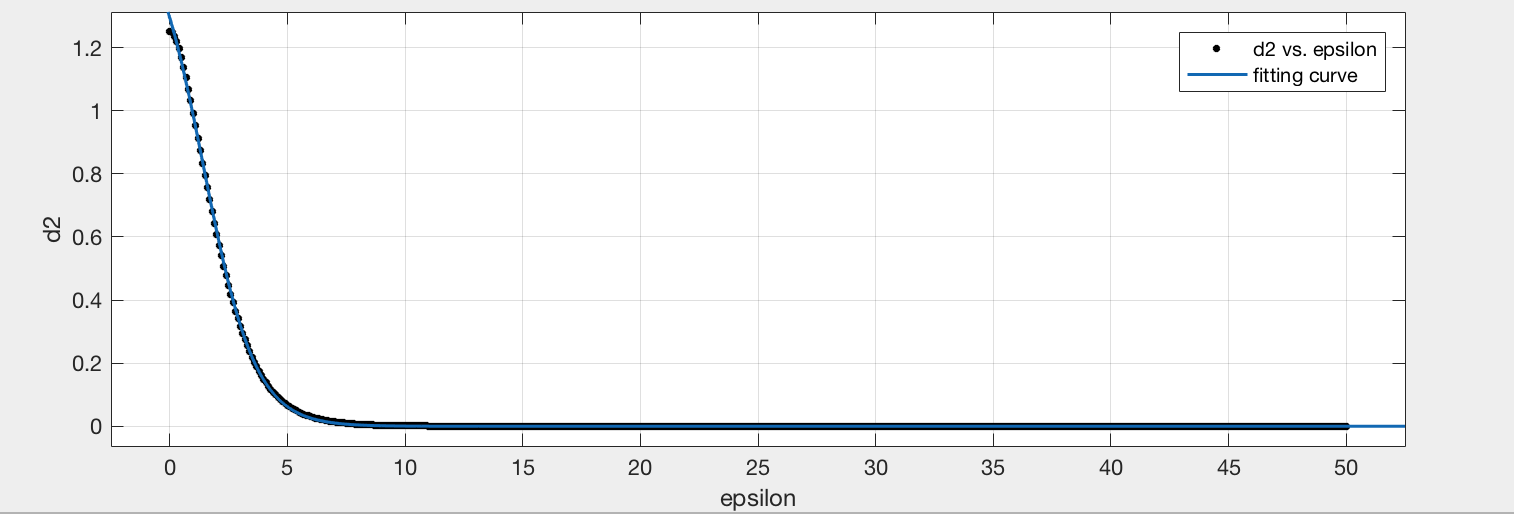}

$3$rd stable band, $a=6.045,b=3.674,c=0.9081$\vspace{0.15cm}

\includegraphics[scale=0.45]{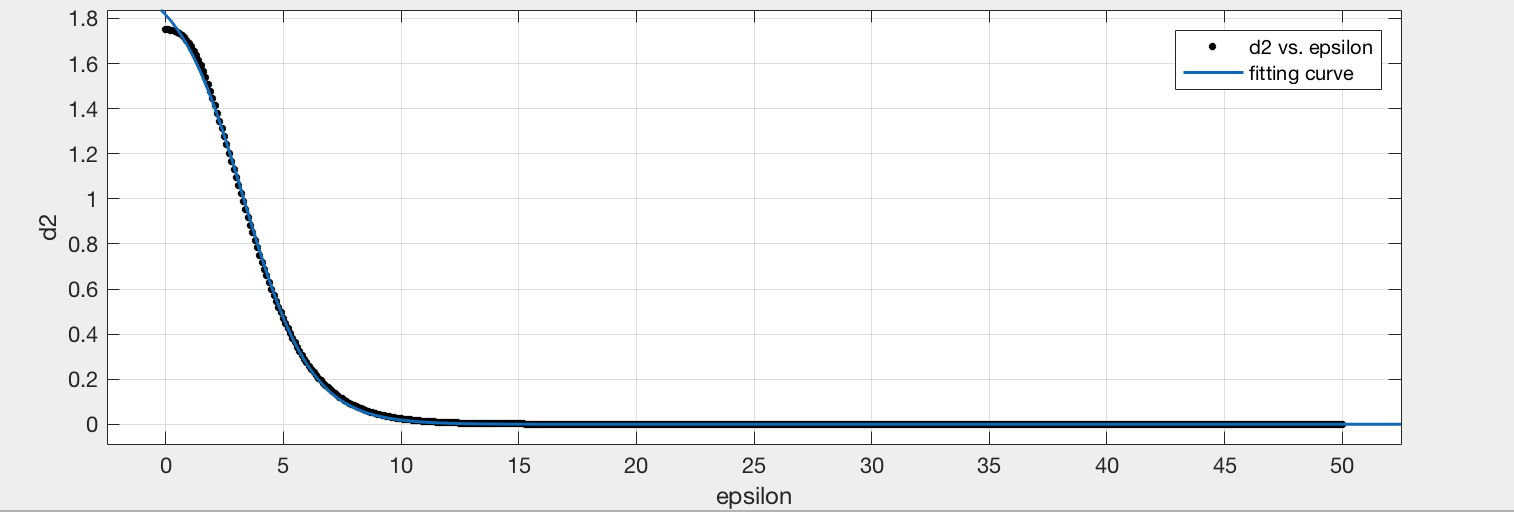}

$4$th stable band, $a=20.17,b=2.507,c=0.6972$\vspace{0.15cm}

\includegraphics[scale=0.45]{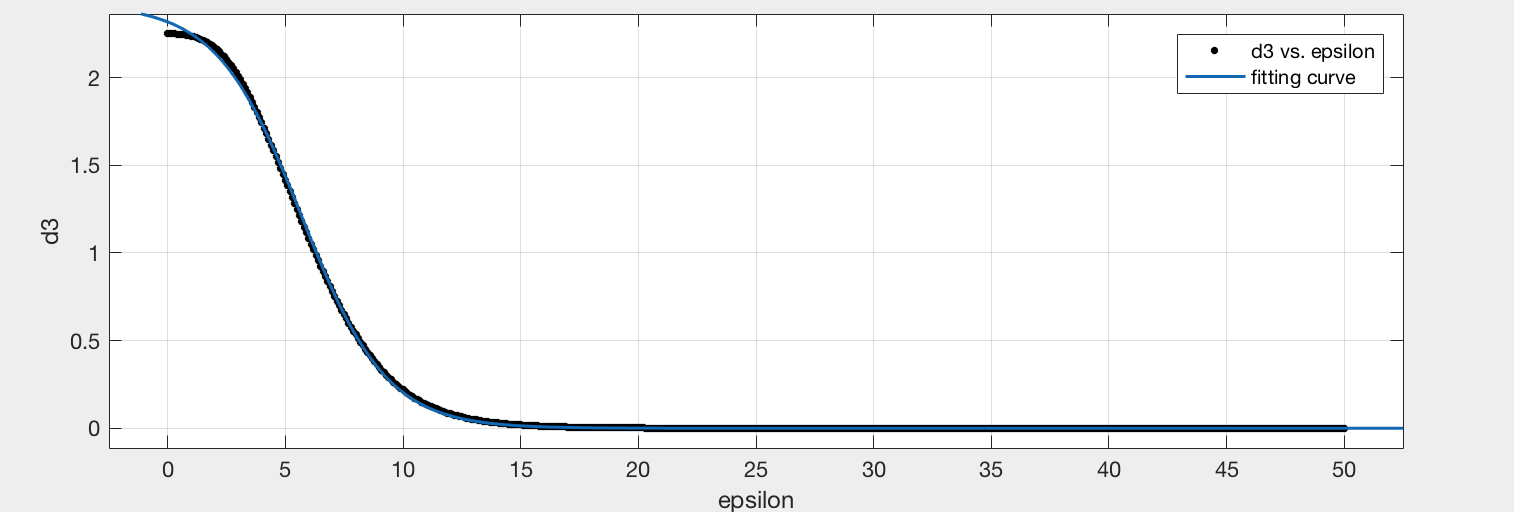}

$5$th stable curve, $a=56.06,b=23.20,c=0.5527$\vspace{0.15cm}

\caption{The width of first $5$ stable bands.}\label{wid1}
\end{figure}

\begin{figure}
\centering
Fitting curve of the form: $\frac{a}{b+e^{cx}}$.\qquad\qquad\qquad\qquad\qquad\qquad\qquad\qquad\qquad\qquad\vspace{0.1cm}

\includegraphics[scale=0.45]{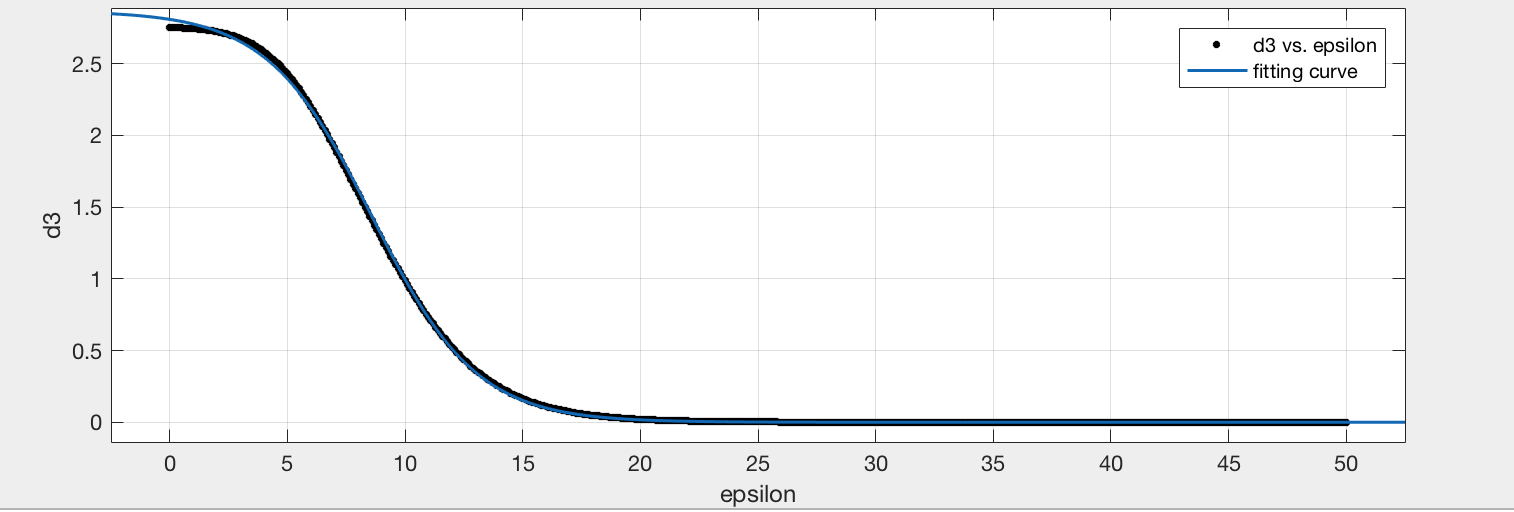}

$6$th stable band, $a=136.8,b=12.38,c=0.4496$\vspace{0.15cm}

\includegraphics[scale=0.45]{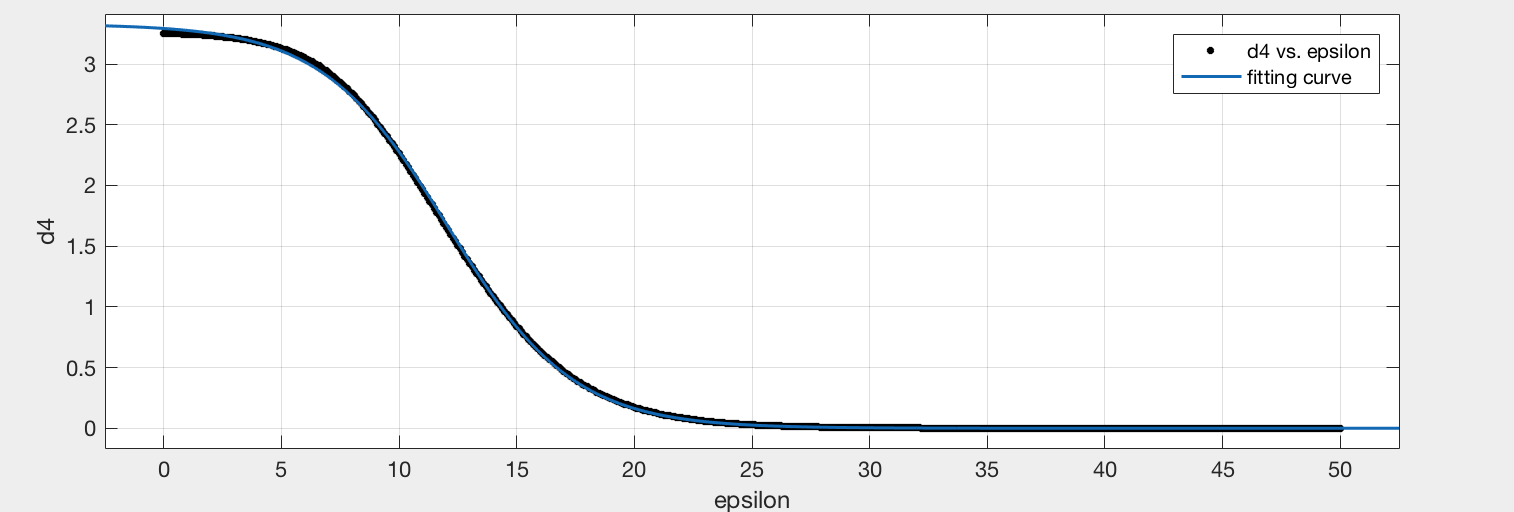}

$7$th stable band, $a=303.0,b=20.39,c=0.374$\vspace{0.15cm}

\includegraphics[scale=0.45]{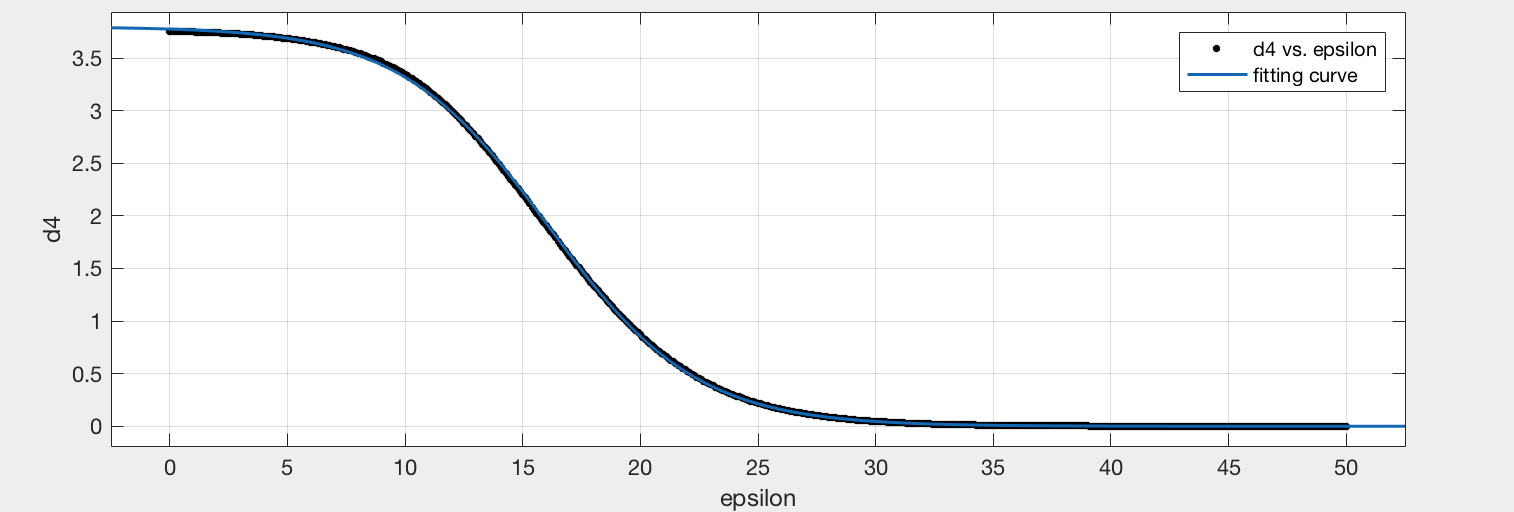}

$8$th stable band, $a=623.9,b=164.2,c=0.3168$\vspace{0.15cm}

\includegraphics[scale=0.45]{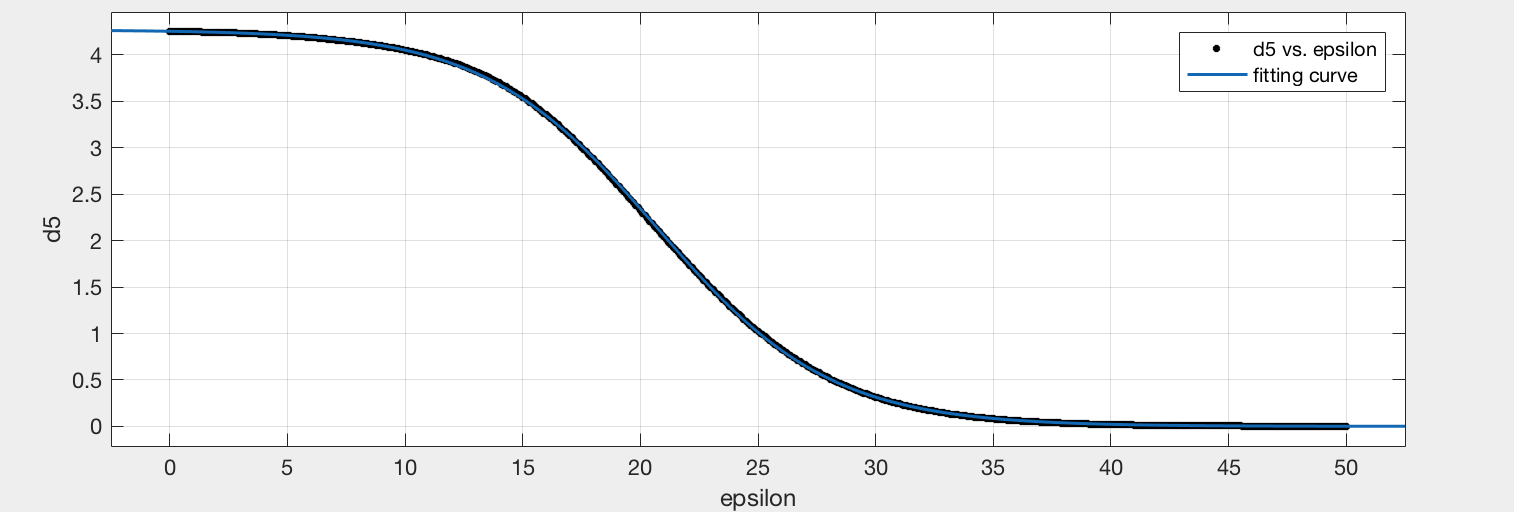}

$9$th stable band, $a=1212,b=283.9,c=0.2726$\vspace{0.15cm}

\includegraphics[scale=0.45]{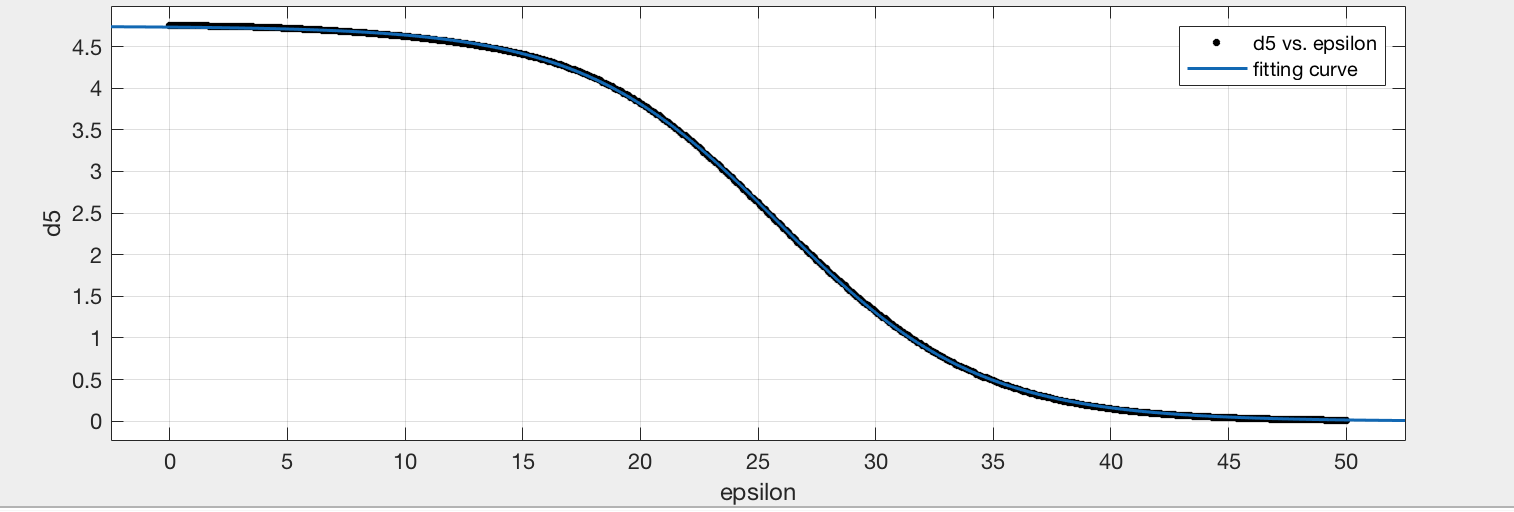}

$10$th stable curve, $a=2248,b=473.9,c=0.2377$\vspace{0.15cm}
\caption{The width of $6$th-$10$th stable bands.}\label{wid2}
\end{figure}

Since the study of the type of convergence in the second observation is comparatively complete, we focus on the type in the first observation, which concerns the width of the stable band as $|\varepsilon|\to\infty$. In particular, we study the first ten stable bands by computing the width of each band from $\varepsilon=0$ to $\varepsilon=50$ in increments of $0.1$. A graph for each stable band, plotting width vs epsilon, are shown below. Then we use the curve fitting toolbox in MATLAB to estimate the fitting curves between $\ep$ and width in each stable band. Please see Figure \ref{wid1} and Figure \ref{wid2} for the result and fitting curves.
%\textcolor{blue}{Below, we provide ten graphs, one for each of the first ten stable bands, which plot band width vs. $\varepsilon.$} In each of these plots, created using MATLAB, $\varepsilon$ ranges from $\varepsilon=0$ to $\varepsilon=50$ in increments of $0.1,$ and each uses a curve-fitting toolbox to draw fitting curves. Please see Figure} \ref{wid1} and Figure \ref{wid2}.}\\

%The estimate for the width of $k-th$ unstable interval is complete, \textcolor{blue}{and} we do not give further experiment. For the asymptotic behavior when $\varepsilon\to\infty$, we choose the first ten stable bands and compute the width of stable band from $\varepsilon=0$ to $\varepsilon=50$. Please see Figure \ref{wid1} and Figure \ref{wid2} for the result and fitting curves.\\

The last topic we will discuss before moving on to talk about solutions to the MDE is to approximate and simplify the irregular boundary between stable and unstable regions for practical use. We can see that most $(\delta,\varepsilon)$ pairs in the first and the fourth quadrants which lie below by the line $\varepsilon=\delta$ and above by the line $\varepsilon=-\delta$ are stable, while most $(\delta,\varepsilon)$ values in those quadrants outside that region form unstable pairs. For any $w>0$, let $R_w:=\{(\delta,\varepsilon)\in\mathbb{R}^2:0<\varepsilon<w \text{ and }-\delta<\varepsilon<\delta\}$. We are interested in determining the probability, for various values of $w$, that a $\delta$-$\varepsilon$ pair is stable, given that it lies in the triangle $R_w.$

\begin{figure}[htp] % not h only
\centering
\subfigure[$R_{w_1}$]{%
\includegraphics[width=0.4\textwidth]{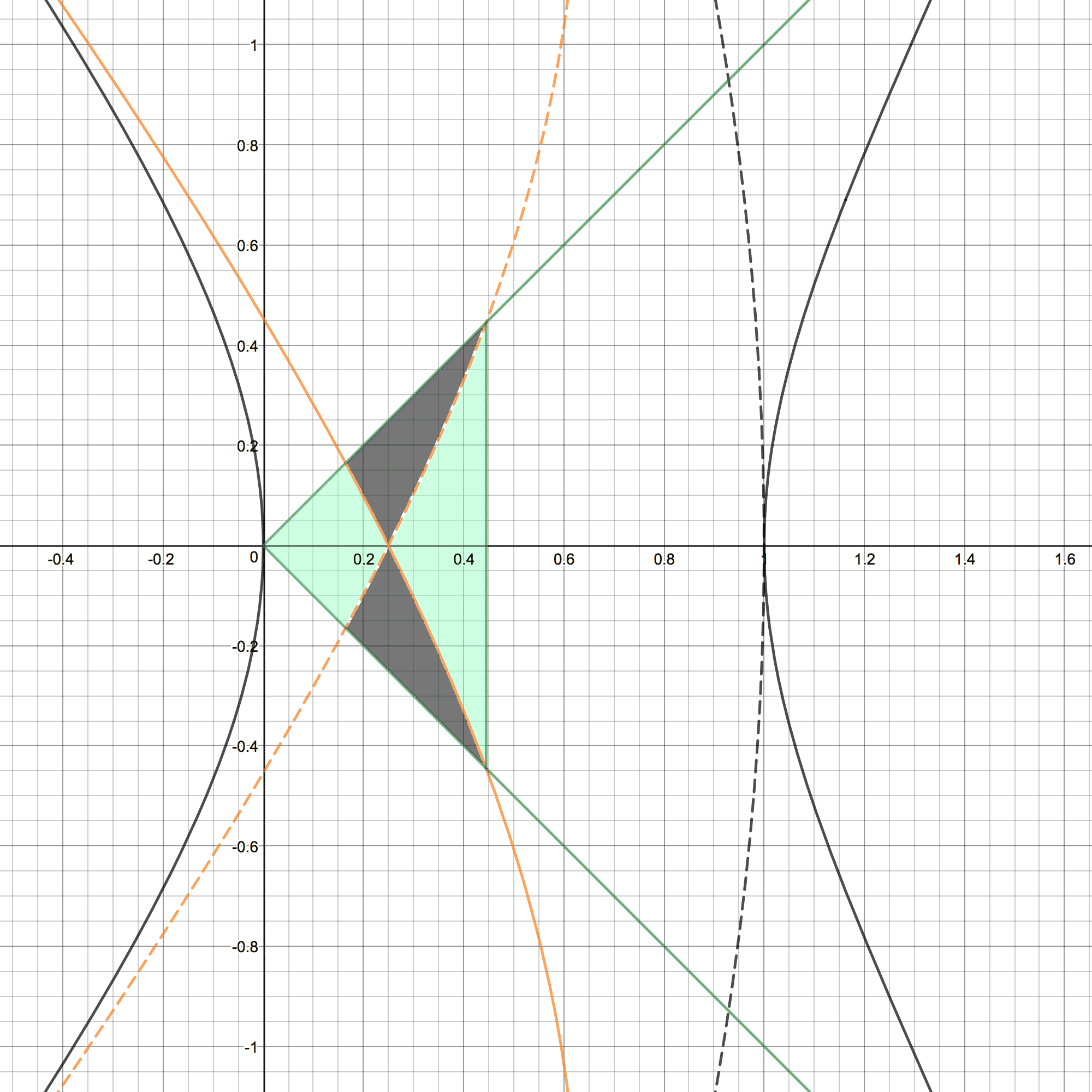}%
\label{fig:q_time_all}%
}\hfil
\subfigure[$R_{w_2}$]{%
\includegraphics[width=0.4\textwidth]{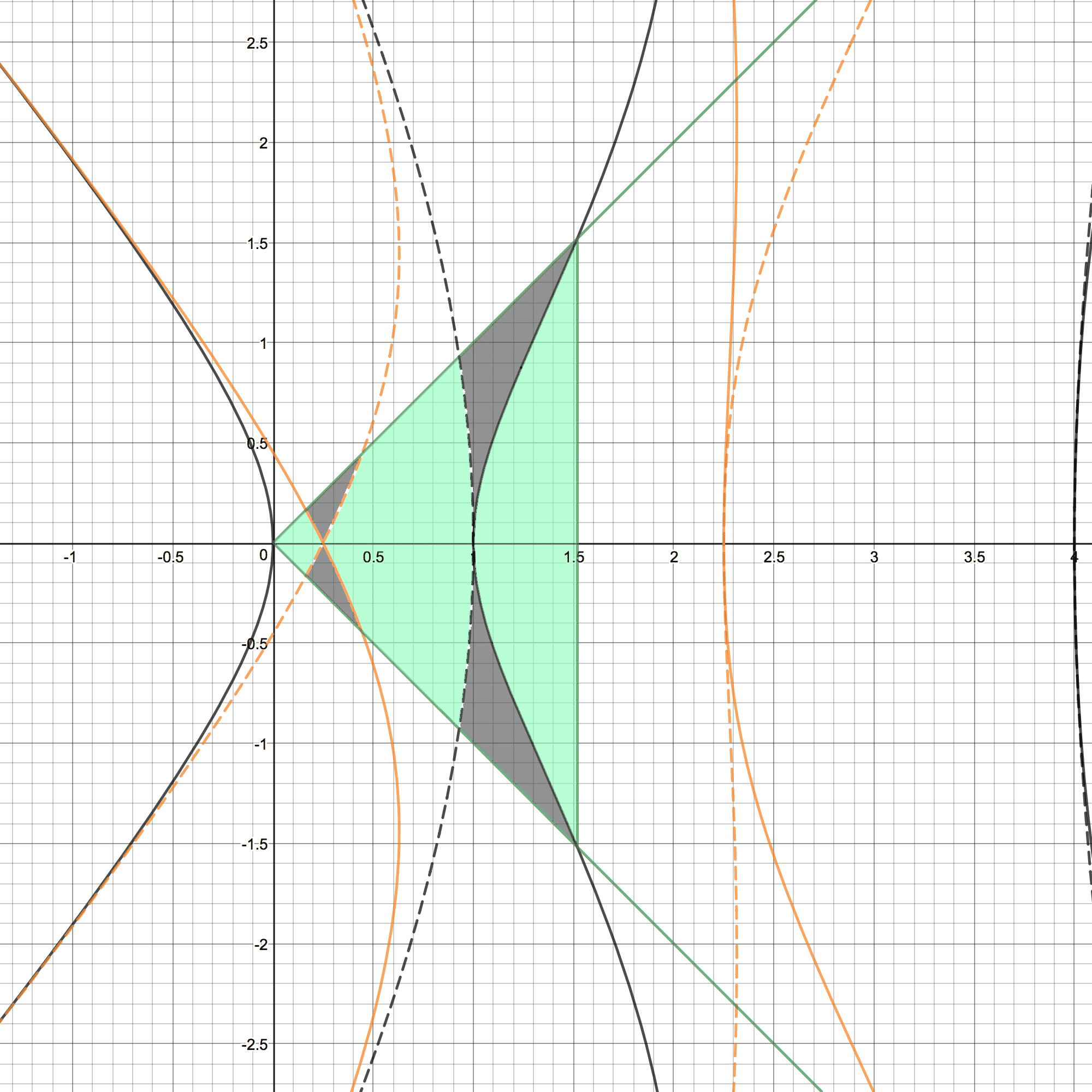}%
\label{fig:q_time_sat}%
}

\subfigure[$R_{w_3}$]{%
\includegraphics[width=0.4\textwidth]{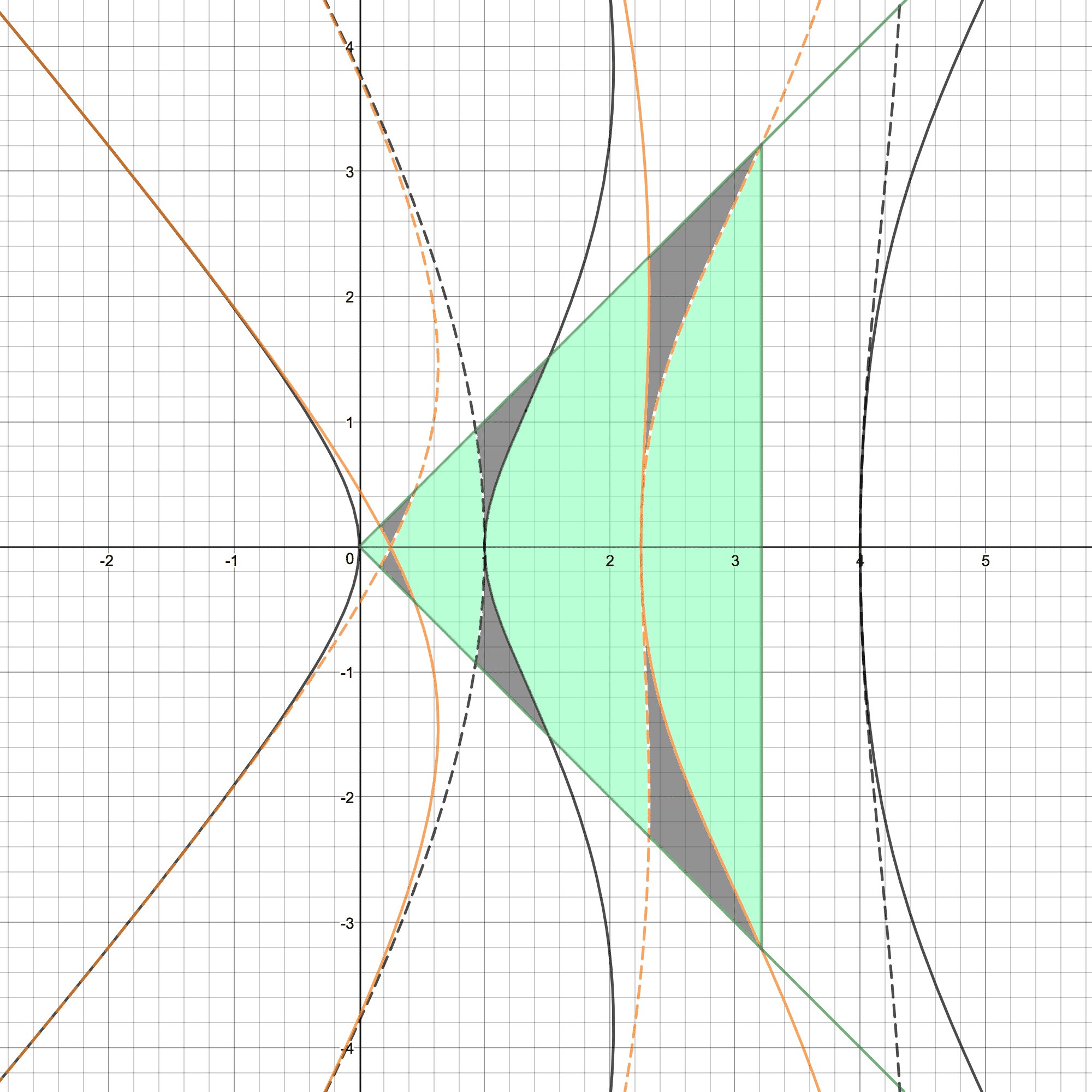}%
\label{fig:q_time_unsat}%
}\hfil
\subfigure[$R_{w_4}$]{%
\includegraphics[width=0.4\textwidth]{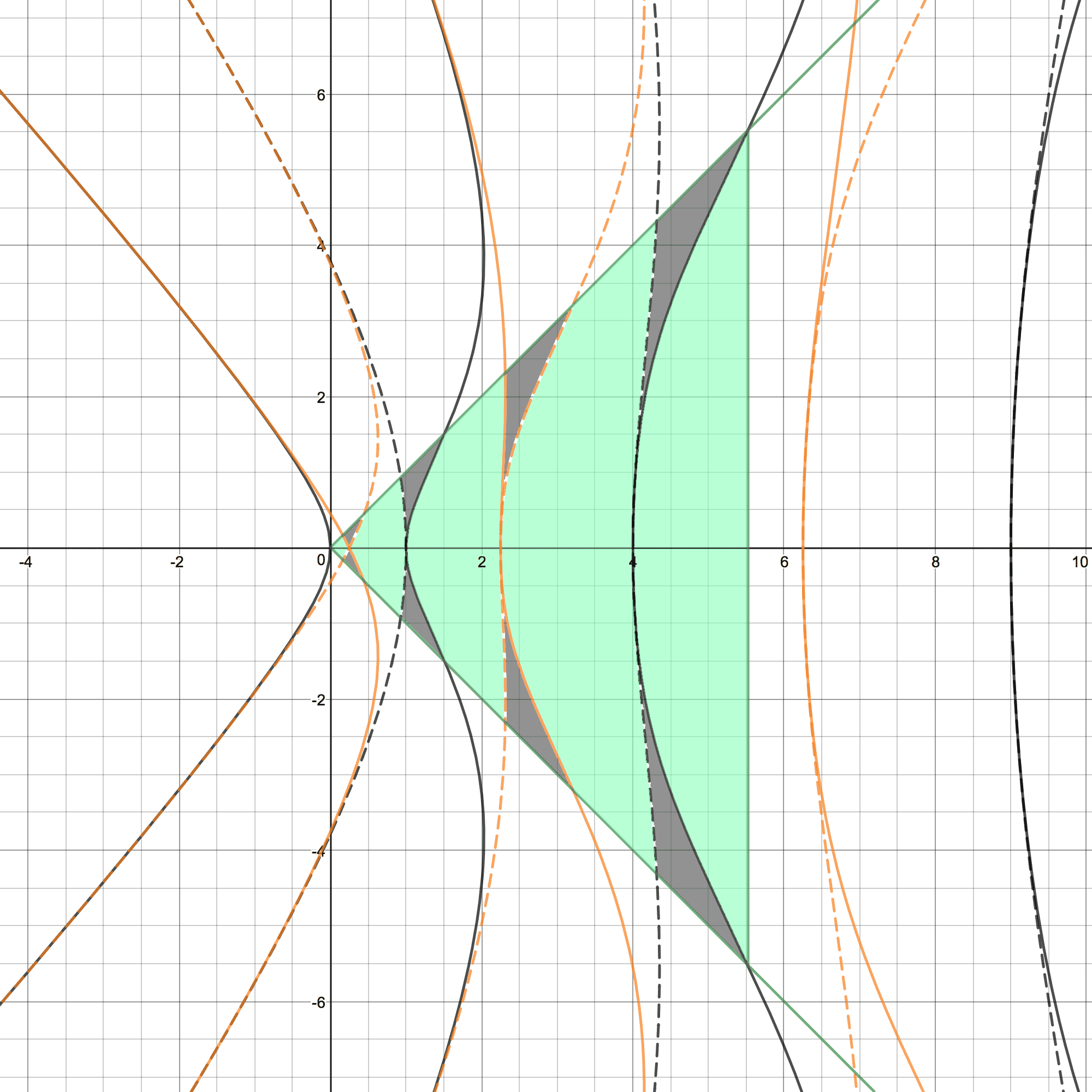}%
\label{fig:q_time_sample}%
}

\caption{The triangle area corresponding to $R_{w_1},R_{w_2},R_{w_3},$ and $R_{w_4}$. The green area is the stable region, and shaded area is the unstable region.}\label{triarea}
\end{figure}

So, we numerically compute the probabilities with different choices of $w$. Notice that the line $\varepsilon=\delta$ (as well as the line $\varepsilon=-\delta$) alternatively passes through stable and unstable regions. So, we define $w_i$ to be the $\delta$-coordinate of the point of intersection between the line $\varepsilon=\delta$ and the right boundary of the $i$-th unstable region; also, for each $i,$ we define $P_i$ to be the probability that a $\delta$-$\varepsilon$ pair is stable, given that it is in $R_{w_i}$. Readers can see figure \ref{triarea} for an illustration of how we divide the regions into triangles.

$P_i$ can also be considered as the probability of getting a stable $(\delta, \ep)$ pair within the triangular regions, which can characterize how well the irregular boundaries can be approximated by lines $\ep=\pm\delta$. We list the $P_i,1\leq i\leq 18$ in Table \ref{Prob}.

\begin{table}[htbp]
\caption{The probability $P_i$'s.}\label{Prob}
\centering
\begin{tabular}{l | l }
\hline
$i$ & $P_i$ \\
\hline
1 & 0.625056436387445\\
2 & 0.784428425813594\\
3 & 0.845139663995868\\
4 & 0.878143787704672\\
5 & 0.899154589232086\\
6 & 0.913800056779179\\
7 & 0.924632580597566\\
8 & 0.932988858616457\\
9 & 0.939640860147421\\
10 & 0.945067300763657\\
11 & 0.949581603650156\\
12 & 0.953397960289362\\
13 & 0.956667925443240\\
14 & 0.959502140424685\\
15 & 0.961982968671620\\
16 & 0.964174006346465\\
17 & 0.966133547589692\\
18 & 0.967239011961280\\
\hline
\end{tabular}
\end{table}

\begin{figure}[h]
\centering
\includegraphics[scale=0.7]{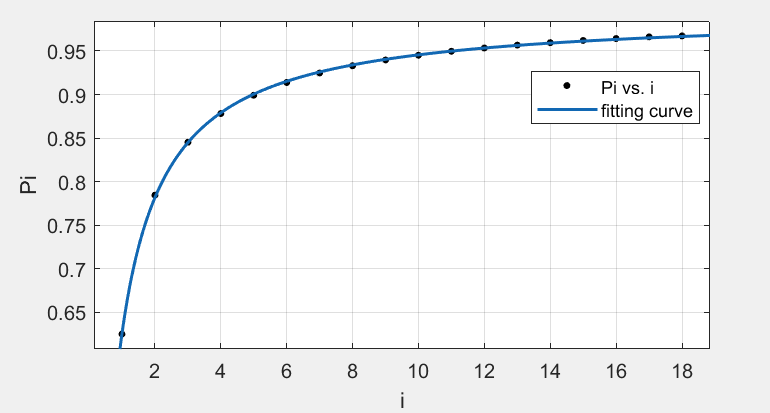}
\caption{Probabilities $P_i$ and the fitting curve.}
\label{prob}
\end{figure}

In addition, we also obtain a surprisingly good fitting curve of the form $P_i=\frac{ai+b}{i+c}$, with $a=0.9938,b=-0.1424,c=0.3608$. See figure \ref{prob} for details. \\

\subsection{Solutions of the Mathieu Differential Equation}
Now we discuss the periodic solutions for the Mathieu differential equation. In Corollary \ref{coasym}, we showed that the sequence $\{c_j\}$ of Fourier coefficients to periodic solutions corresponding to $\delta$-$\varepsilon$ pairs on the transition curves converge rapidly. Hence, given a $\delta$-$\varepsilon$ pair on a transition curve, we can calculate explicitly the coefficients of the corresponding solution with a finite truncation of the Fourier series and use those coefficients to plot, with very high accuracy, solutions to the MDE.

One question one might investigate is as follows. Suppose one fixes a transition curve and considers periodic solutions corresponding to various $(\delta, \varepsilon)$ points along that transition curve. How do properties of solutions change as $\varepsilon$ varies?

To be more precise, we need some new notation. Recall from Section 2.2 that, when solving $2\pi$- or $4\pi$- periodic solutions, we can rewrite the MDE as a system of linear equations for Fourier coefficents, which collectively can be rewritten as four independent homogeneous matrix equations, with corresponding matrices $A,B,C,D$. The transition curves consist of the $\delta$-$\ep$ pairs such that one of these matrices degenerates. In this way, a transition curve can be labeled with a matrix and an integer $k\geq 1$, so that for each $\delta$-$\ep$ pair on this transition curve, $\delta$ is the $k$th smallest real number such that the matrix degenerates with the fixed $\ep$. We want a more convenient way to label these transition curves, and to achieve this, we proceed as follows. Fix a transition curve, let $\ep=0$ on the transition curve, and solve the corresponding equation in matrix form. According to the discussion in Section 2.2, we can then construct a periodic solution to the MDE, and the solution is clearly an eigenfunction of $\frac{d^2}{dt^2}$. For example, if we consider the transition curve labeled by $A$ and $k=2$ we will get $\cos t$ as a solution with the above process. In this way, a transition curve can be labeled with an eigenfunction of $\frac{d^2}{dt^2}$. In light of this, we define $p(\varphi,\varepsilon)$ to be the point in the $\delta$-$\ep$ plane with prescribed $\ep$-coordinate which lies on the transition curve corresponding to the eigenfunction $\varphi$.

%\st{The periodic solution for a $\delta-\ep$ pair on a transition curve is unique up to multiplication by a constant, except at the intersection with the $\delta$ axis. Noticing that all these solutions take the same form of Fourier series expansion out of the following four possibilities $\sum_{j=0}^\infty c_j\cos jt$, $\sum_{j=1}^\infty c_j\sin jt$, $\sum_{j=0}^\infty c_j\cos \frac{2j+1}{2}t$ or $\sum_{j=0}^\infty c_j\sin \frac{2j+1}{2}t$, we can find a unique solution with the same kind of Fourier series expansion to the MDE up to multiplication as that corresponding to the intersection point of the transition curve and the $\delta$ axis, and this is immediately an eigenfunction of $-\frac{d^2}{dt^2}$. It is also easy to see that the solutions vary continuously as the pair moves along the curve. In this way, the transition curve is uniquely characterized by an eigenfunction $\varphi$ of $-\frac{d^2}{dt^2}$, and for convenience, we define $p(\varphi,\varepsilon)$ to be the point on the transition curve corresponding to $\varphi$ with prescribed $\ep$ coordinate.}

As shown in Section 2, the transition curves can be organized into four different classes, made up of points $\{p(\sin kt,\varepsilon)\},\{p(\cos kt,\varepsilon)\},\{p(\sin \frac{2k+1}{2}t,\varepsilon)\},\{p(\cos\frac{2k+1}{2}t,\varepsilon)\}$ with our new notation. It is natural to study them separately.

\begin{figure}[h]
\centering
\includegraphics[scale=0.3]{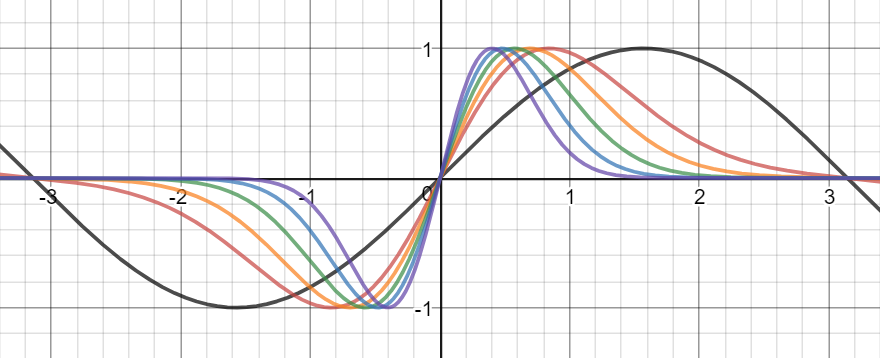}
\caption{Normalized solutions corresponding to $p(\sin t,0)$ (solid black), $p(\sin t,5)$ (red), $p(\sin t,10)$ (orange), $p(\sin t,20)$ (green), $p(\sin t,40)$ (blue), $p(\sin t,80)$ (purple).}
\label{linesolsconv1}
\end{figure}

\begin{figure}[h]
\centering
\includegraphics[scale=0.3]{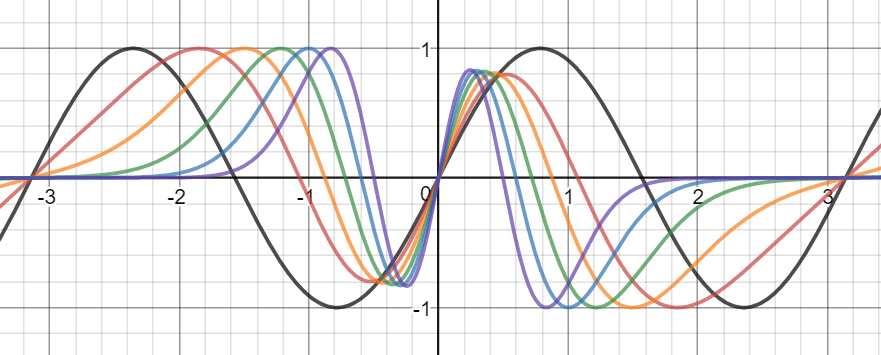}
\caption{Normalized solutions corresponding to $p(\sin 2t,0)$ (solid black), $p(\sin 2t,5)$ (red), $p(\sin 2t,10)$ (orange), $p(\sin 2t,20)$ (green), $p(\sin 2t,40)$ (blue), $p(\sin 2t,80)$ (purple).}
\label{linesolsconv2}
\end{figure}

\begin{figure}[h]
\centering
\includegraphics[scale=0.3]{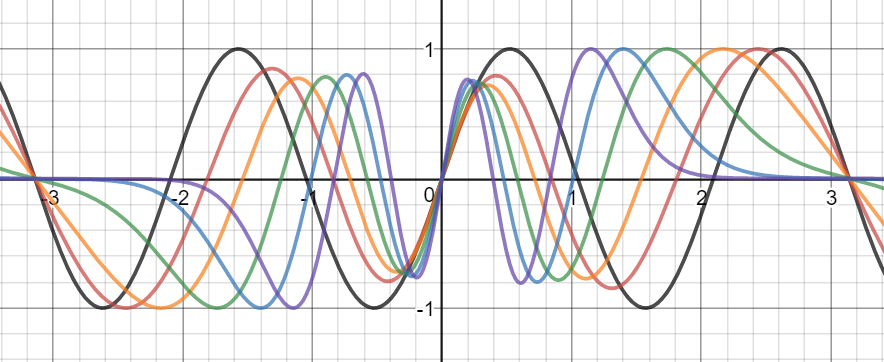}
\caption{Normalized solutions corresponding to $p(\sin 3t,0)$ (solid black), $p(\sin 3t,5)$ (red), $p(\sin 3t,10)$ (orange), $p(\sin 3t,20)$ (green), $p(\sin 3t,40)$ (blue), $p(\sin 3t,80)$ (purple).}
\label{linesolsconv3}
\end{figure}

First we consider the class $\{p(\sin kt,\varepsilon)\}$. Because of the symmetry of the transition curves across the $\delta$-axis, it suffices to only consider positive values of $\varepsilon$. In fact, if $u$ is a solution corresponding to the pair $(\ep,\delta)$, then $u(t+\pi)$ is a solution corresponding to the pair $(-\ep,\delta)$ since $-\cos t=\cos(t+\pi)$. The normalized solutions are plotted in Figures \ref{linesolsconv1}-\ref{linesolsconv3}, where by `normalized solution' we mean a solution $u$ for which $\max_{t\in\mathbb{R}}u(t)=1$. In plotting the graphs we have used the particular values $k=1,2,3$ and $\ep=5,10,20,40,80,160$. $15\times 15$ truncated matrices are used to compute the solutions. The horizontal axis is the $t$-axis, and the vertical axis is the $u$-axis. 

%\st{that we have scaled each solution vertically by a constant factor so as to make its maximum value equal to 1}

An interesting pattern can be observed from Figure \ref{linesolsconv1}.  \textit{The sequence of $t$-coordinates of these maximal points appears to approach $0$ monotonically as $\ep$ increases.} Figures \ref{linesolsconv2} and \ref{linesolsconv3} show the same behavior as well.

\begin{figure}[h]
\centering
\includegraphics[scale=0.33]{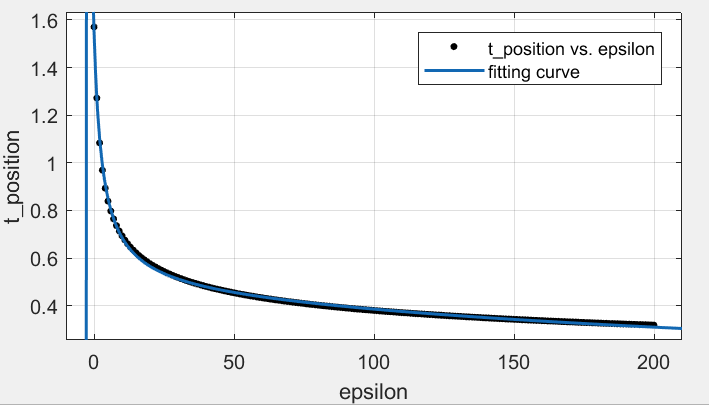}\\
$t$-position of maximal point on curve $p(\sin t,\ep)$, with $a=3.031,b=5.736,c=10.18,d=14.88$.\vspace{0.4cm}

\includegraphics[scale=0.33]{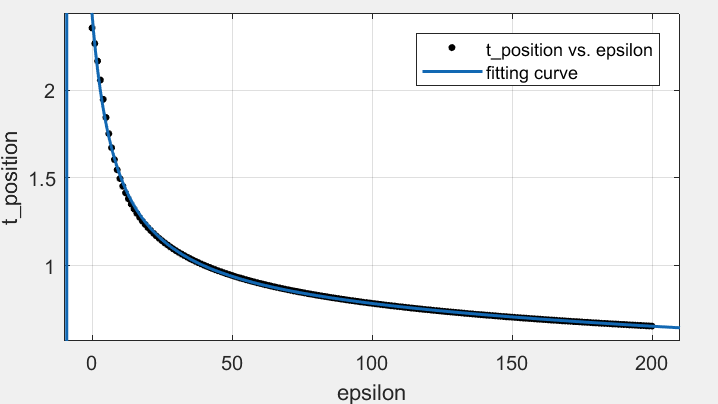}\\
$t$-position of maximal point on curve $p(\sin 2t,\ep)$, with $a=10.31,b=22.8,c=17.43,d=29.09$.\vspace{0.4cm}

\includegraphics[scale=0.33]{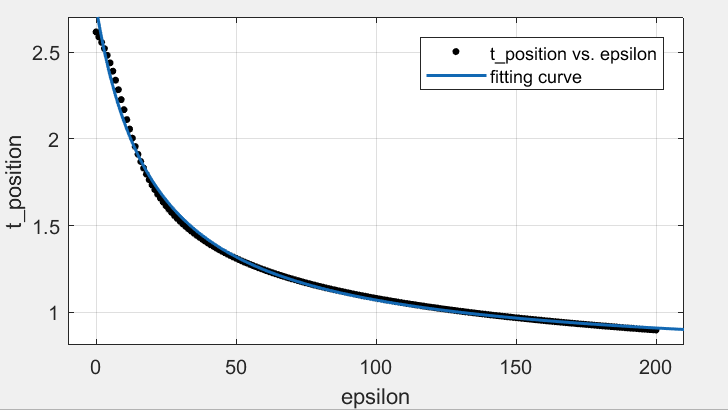}\\
$t$-position of maximal point on curve $p(\sin 3t,\ep)$, with $a=432.5,b=1345,c=607.2,d=1257$.
\caption{t-position of maximal points, with fitting curve $t=\frac{a\ep+b}{\ep^2+c\ep+d}$.}
\label{max_position_sin}
\end{figure}

We now use curves to fit a curve of the $t$-coordinate of these relative maxima in $[0,\pi]$ as a function of $\ep$. We still pick the same three transition curves, and take $\ep$ from $1$ to $200$, in increments by $1$. The relationship of $t$ coordinate of maximal points and $\ep$ on the three curves are shown in figure \ref{max_position_sin}, and the fitting curve is of the form $t=\frac{a\ep+b}{\ep^2+c\ep+d}$.
%In Figure \ref{max2_position_sin}, we have plotted the $t$-coordinate for the  shows three plots. \st{Each plot for each of the three transition curves described above }

\begin{figure}[h]
\centering
\includegraphics[scale=0.33]{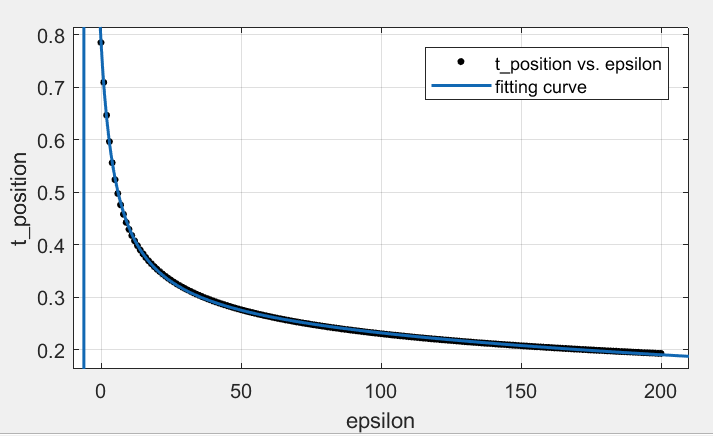}\\
$t$-position of the second maximal point on curve $p(\sin 2t,\ep)$, with $a=2.272,b=4.602,c=12.75,d=19.84$.\vspace{0.4cm}

\includegraphics[scale=0.33]{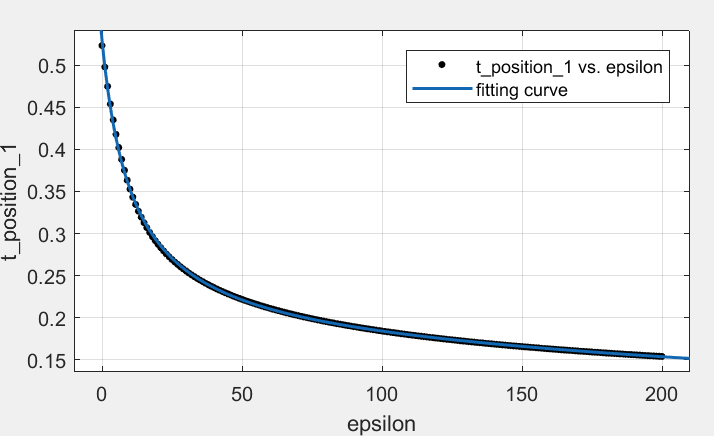}\\
$t$-position of the second maximal point on curve $p(\sin 3t,\ep)$, with $a=710.8.5,b=1894,c=5324,d=1.045\times 10^4$.
\caption{t-position of the second maximal points, with fitting curve $t=\frac{a\ep+b}{\ep^2+c\ep+d}$.}
\label{max2_position_sin}
\end{figure}

For the case of $p(\sin 2t,\ep)$ and $p(\sin 3t,\ep)$, does the same convergence behavior occur if one considers the sequence of the first minima when $t>0$? We also do the computation, and show the results in figure \ref{max2_position_sin}, where the same kind of fitting curve also works.

\begin{figure}[h]
\centering
\includegraphics[scale=0.33]{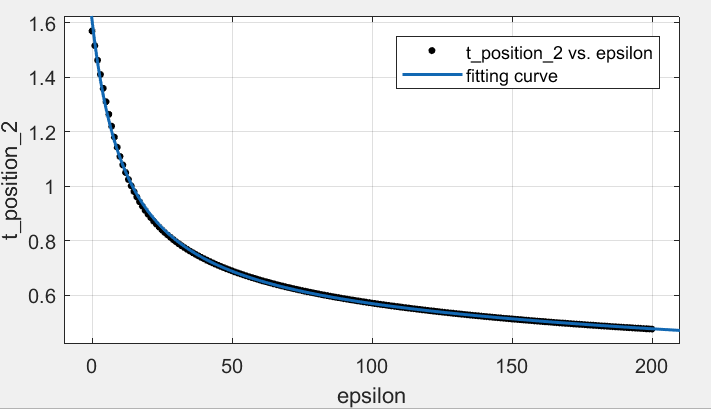}
\caption{t-position of the third maximal points on curve $p(\sin 3t,\ep)$, with fitting curve $t=\frac{a\ep+b}{\ep^2+c\ep+d}$, $a=9.567,b=23.07,c=22.86,d=40.46$.}
\label{max3_position_sin}
\end{figure}

In case of $p(\sin 3t,\ep)$, the same observation applies to the second maxima when $t>0$. See figure \ref{max3_position_sin} for details.

\begin{figure}[h]
\centering
\includegraphics[scale=0.33]{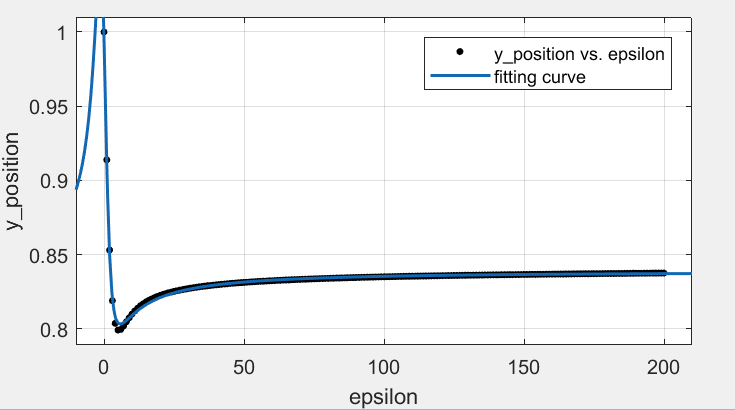}
\caption{u-position of the second maximal points on curve $p(\sin 2t,\ep)$, with fitting curve $t=\frac{a\ep^2+b\ep+c}{\ep^2+d\ep+e}$. $a=0.8393,b=-0.5638,c=5.566,d=-0.1973,e=5.571$.}
\label{max2_yposition_sin2t}
\end{figure}

One can also consider the behavior of the $u$-coordinate of these various sequences of local extrema. See figure \ref{max2_yposition_sin2t} for the behavior of $u$-coordinate of the first minimal points on the curve $p(\sin 2t,\ep)$, where we can fit the points well with a rational function of the form $\frac{a\ep^2+b\ep+c}{\ep^2+b\ep+c}$.

\begin{figure}[h]
\centering
\includegraphics[scale=0.33]{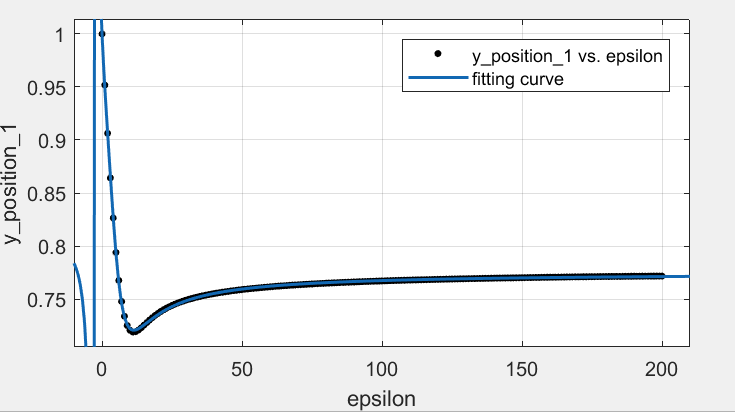}\\
$u$-position of the second maximal point on curve $p(\sin 3t,\ep)$, with $a=0.7784,b=-3.71,c=21.57,d=150,e=-3.899,f=30.57,g=149.9$.\vspace{0.4cm}

\includegraphics[scale=0.33]{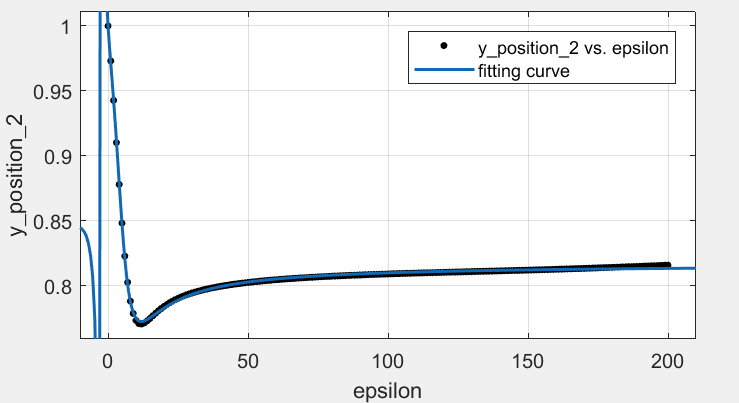}\\
$u$-position of the third maximal point on curve $p(\sin 3t,\ep)$, with $a=0.8166,b=-4.106,c=20.44,d=136.9,e=-4.234,f=25.11,g=136.8$.
\caption{u-position of the second and third maximal points corresponding to curve $p(\sin 3t,\ep)$, with fitting curve $u=\frac{a\ep^3+b\ep^2+c\ep+d}{\ep^3+e\ep^2+f\ep+g}$.}
\label{max_yposition_sin3t}
\end{figure}

For the curve $p(\sin 3t,\ep)$, there are three local extremes on each solution, and we also do the computation for the $u$-coordinate for the first minimal and the second maximal value. The fitting curve is a little more complicated, of the form $t=\frac{a\ep^3+b\ep^2+c\ep+d}{\ep^3+d\ep^2+e\ep+g}$. See figure \ref{max_yposition_sin3t}.

\begin{figure}[h]
\centering
\includegraphics[height=3.5cm]{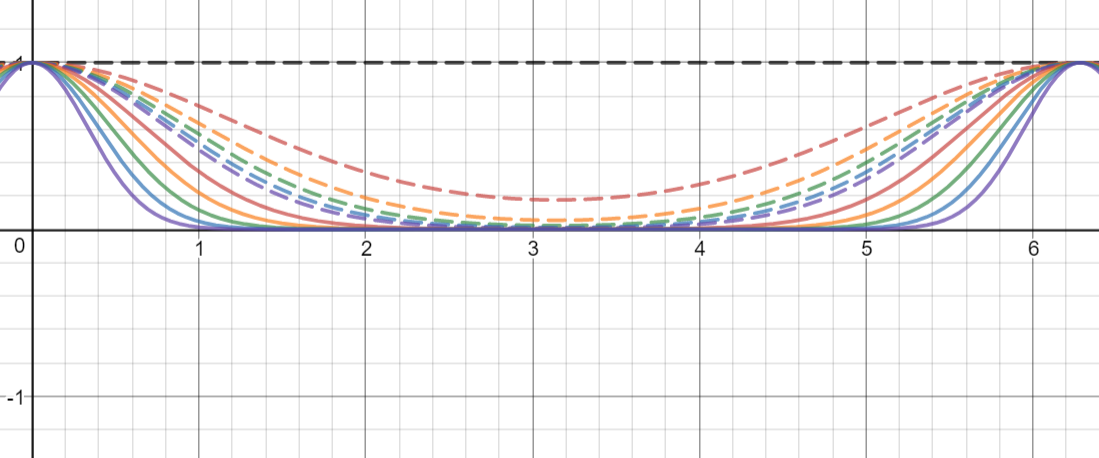}
\caption{Normalized solutions corresponding to $p(\cos 0t,\ep)$, with $\ep=0,1,2,3,4,5,10,20,40,80,160$.}
\label{solutioncos0t}
\end{figure}

\begin{figure}[h]
\centering
\includegraphics[height=3.5cm]{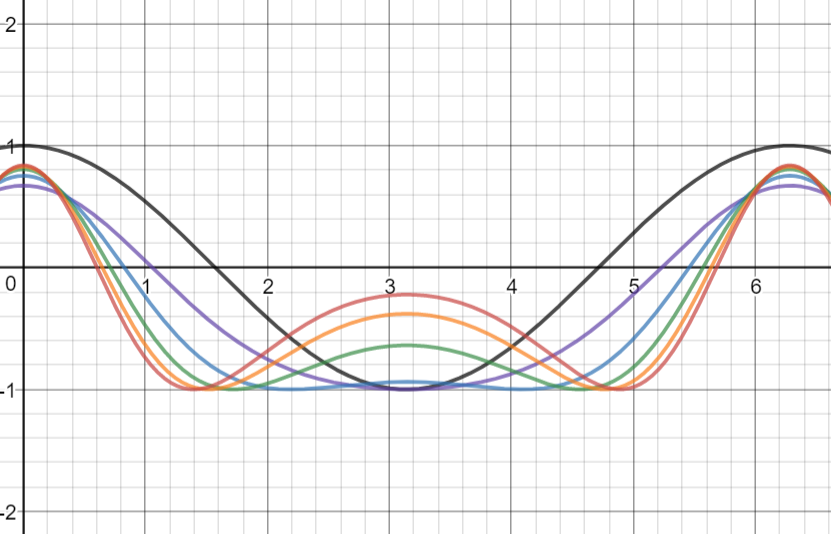}\quad
\includegraphics[height=3.5cm]{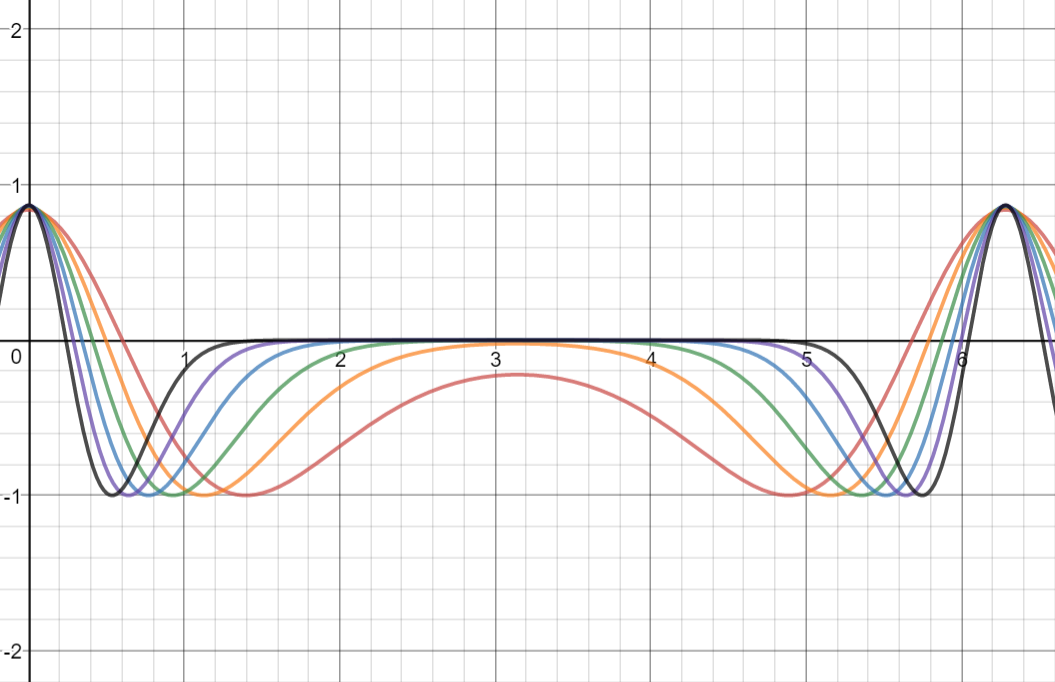}
\caption{Normalized solutions corresponding to $p(\cos t,\ep)$, with $\ep=0,1,2,3,4,5,10$ in the left graph, and $\ep=10,20,40,80,160$ in the right graph.}
\label{solutioncost}
\end{figure}

\begin{figure}[h]
\centering
\includegraphics[height=3.5cm]{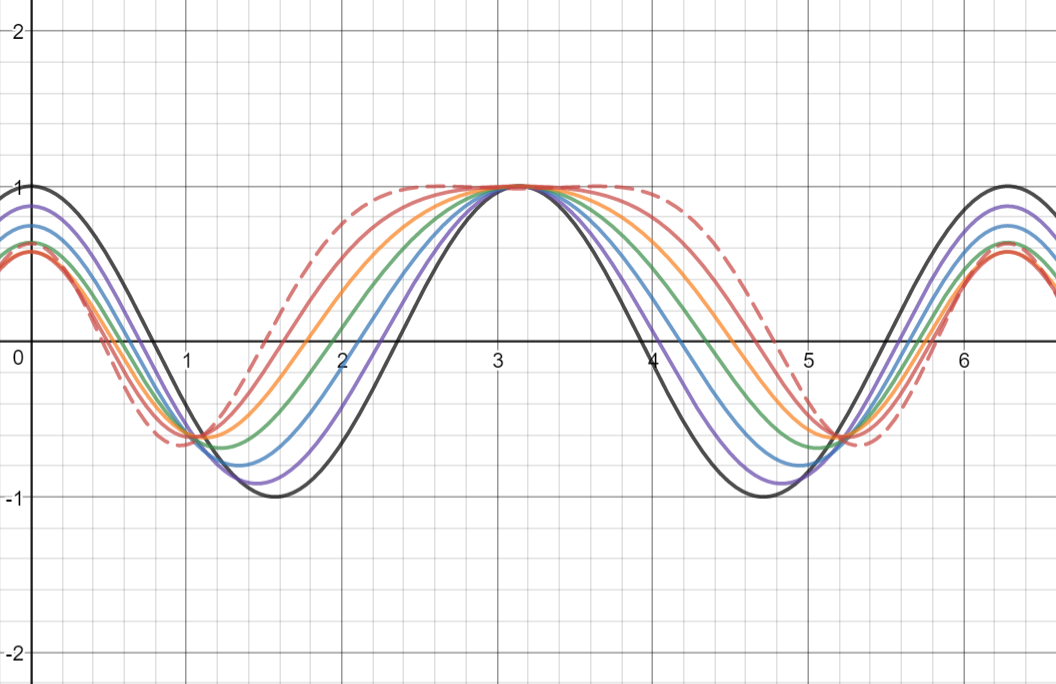}\quad
\includegraphics[height=3.5cm]{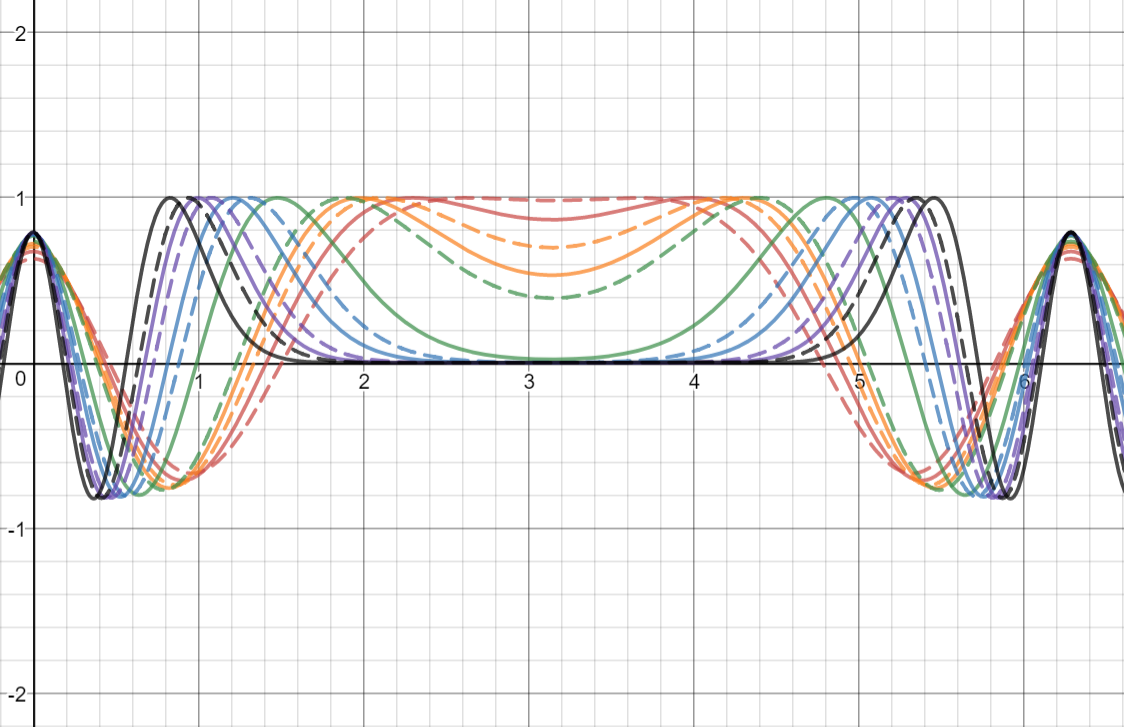}
\caption{Normalized solutions corresponding to $p(\cos 2t,\ep)$, with $\ep=0,1,2,3,4,5,6$ in the left, and $\ep=6,7,8,9,10,20,30,40,60,80,100,160$ in the right.}
\label{solutioncos2t}
\end{figure}

Next, we turn to the solutions on curves $p(\cos kt,\ep)$, and we do experiments for $k=0,1,2$. Please look at figure \ref{solutioncos0t},\ref{solutioncost},\ref{solutioncos2t} for the results.

\begin{figure}[h]
\centering
\includegraphics[height=3.5cm]{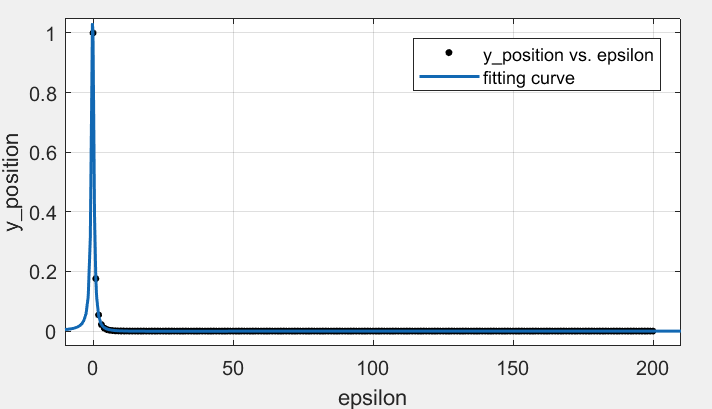}
\caption{The $u$ coordinate of minimum points of solutions for points $p(\cos 0,\ep)$, with $\ep=0,1,\cdots, 200$. Fitting curve $u=\frac{a\ep+b}{\ep^2+c\ep+d}$, with $a=-0.02171,b=0.2895,c=0.2289,d=0.2895$.}
\label{min_yposition_cos0}
\end{figure}

We can see that the shapes of these curves are different from those in figure \ref{linesolsconv1} \ref{linesolsconv2},\ref{linesolsconv3}. The first curve $p(\cos 0t,\ep)$ is of course very special, as it corresponds to the constant functions. It turns out the solution is always positive on this curve, as a consequence of Sturm's Theorem \cite{loud}. Also, we have the minimal values of the solutions at $\pi$, see figure \ref{min_yposition_cos0} for the minimal values for different choices of $\ep$. We use a fitting curve of the form $u=\frac{a\ep+b}{\ep^2+c\ep+d}$.

\begin{figure}[h]
\centering
\includegraphics[height=3.5cm]{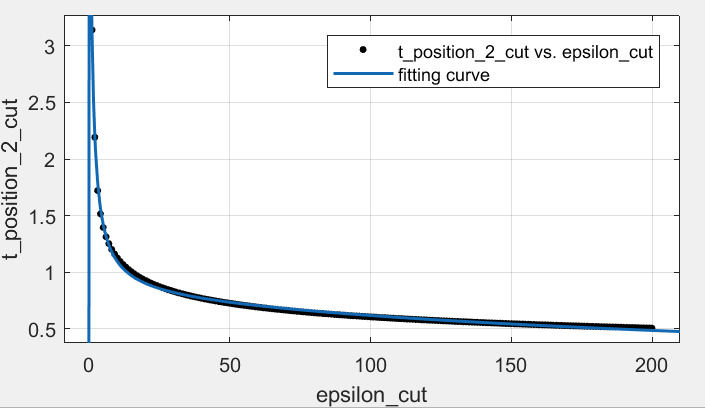}
\caption{The $t$ coordinate of minimum points of solutions for points $p(\cos t,\ep)$, with $\ep=2,\cdots, 200$. Fitting curve $y=\frac{a\ep+b}{\ep^2+c\ep+d}$, with $a=241,9,b=1284,c=310.8,d=177.1$.}
\label{min_tposition_cost}
\end{figure}

The solutions on the second curve is a little more complicated. As shown in figure \ref{solutioncost}, we can see that minimal value is achieved at $t=\pi$ in $[0,2\pi]$ for $\ep=0$ and $\ep=1$, but as $\ep$ becomes larger, the minimum point splits into two different minimum points. This is an interesting phenomenon, and we have an explanation for this in Theorem \ref{thm44}. We do computations on the $t$-coordinate of the minimum points as $\ep$ gets larger. Please see figure \ref{min_tposition_cost} for the data, and we use a fitting curve of the form $t=\frac{a\ep+b}{\ep^2+c\ep+d}$.

\begin{figure}[h]
\centering
\includegraphics[height=3.5cm]{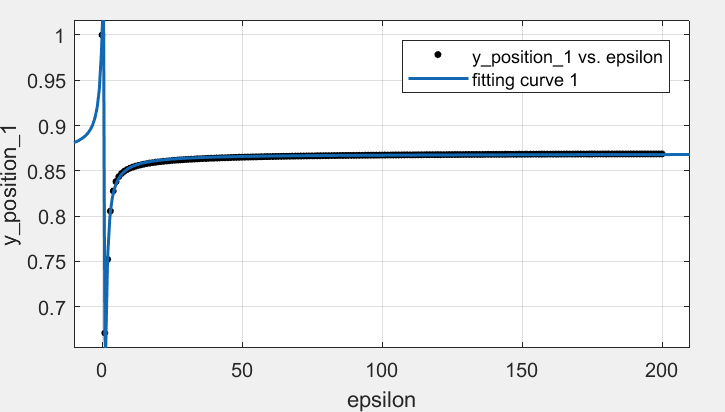}
\caption{The $u$ coordinate of minimum points of solutions for points $p(\cos t,\ep)$, with $\ep=0,1,2,\cdots, 200$. Fitting curve $u=\frac{a\ep^2+b\ep+c}{\ep^2+d\ep+e}$, with $a=0.8687,b=-1.631,c=0.8498,d=-1.72$ and $e=0.8497$.}
\label{max1_yposition_cost}
\end{figure}

\begin{figure}[h]
\centering
\includegraphics[height=3.5cm]{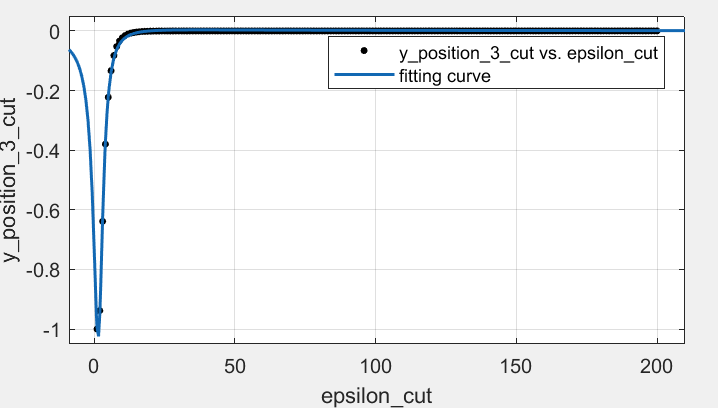}
\caption{The $u$ coordinate of minimum points of solutions for points $p(\cos t,\ep)$, with $\ep=0,1,2,\cdots, 200$. Fitting curve $u=\frac{a\ep+b}{\ep^2+c\ep+d}$, with $a=0.2495,b=-4.991,c=-3.037,d=6.722$.}
\label{max3_yposition_cost}
\end{figure}

Another thing we need to highlight here is that the maximal absolute value does not occur at $0$ or $\pi$, but at some other points. So, as usual, we normalize the solutions so that their minimum value are $-1$, and we compute values of solutions at the two other local peaks at $0$ and $\pi$, which is shown in figure \ref{max1_yposition_cost}. It is also interesting to see in figure \ref{max3_yposition_cost} that the local maximal value at $0$ drop rapidly near $\ep=0$, and then increase quickly.

\begin{figure}[h]
\centering
\includegraphics[height=3.5cm]{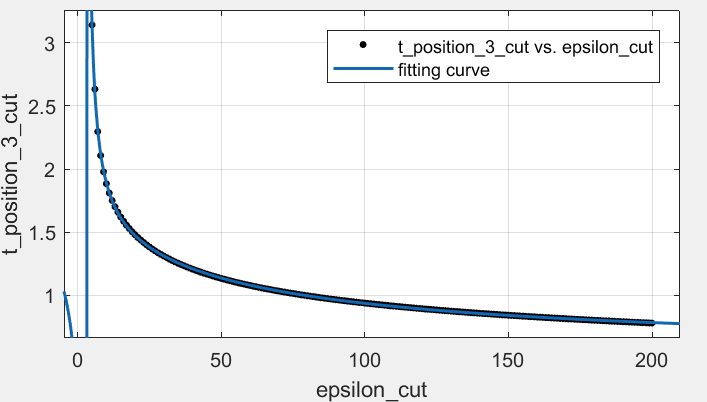}
\caption{The $t$ coordinate of maximum points of solutions for points $p(\cos 2t,\ep)$, with $\ep=5,6,\cdots, 200$. Fitting curve $t=\frac{a\ep^2+b\ep+c}{\ep^2+d\ep+e}$, with $a=0.5572,b=88.5,c=-3.009,d=55.35$ and $e=-157.7$.}
\label{peak3_tposition_cos2t}
\end{figure}

The solutions on the curve $p(\cos 3t,\ep)$ is even more complicated. Still, we observe the split of one peak. When $\ep=0,1,2,3,4,5$, we only see one maximum point, but as $\ep$ becomes larger, we see two maximum points. We also do computation on the $t$ coordinate of one of these splited maximum point. See figure \ref{peak3_tposition_cos2t} for the computation result, where we use fitting curves of the form $t=\frac{a\ep^2+b\ep+c}{\ep^2+d\ep+e}$.

We also do experiment on the t coordinate of the first minimum points and the $u$-coordinate of the peaks on the curve $p(\cos 2t,\ep)$.

\begin{figure}[h]
\centering
\includegraphics[height=3.5cm]{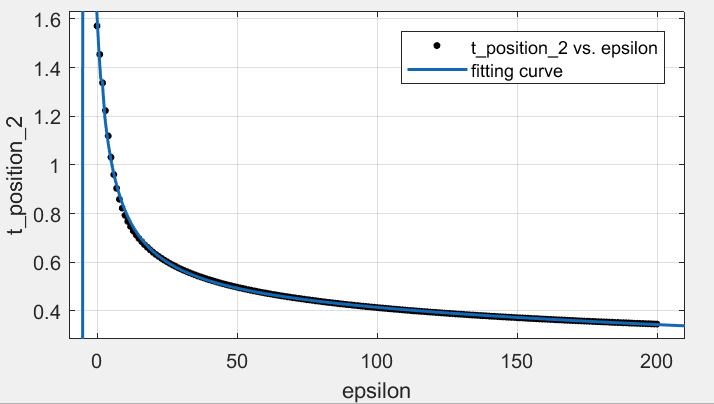}
\caption{The $t$ coordinate of the second peak of solutions for points $p(\cos 2t,\ep)$, with $\ep=5,6,\cdots, 200$. Fitting curve $t=\frac{a\ep^2+b\ep+c}{\ep^3+d\ep^2+e\ep+f}$, with $a=269.3,b=4745,c=-9715,d=641.3,e=-2040$ and $f=-6017$.}
\label{peak2_tposition_cos2t}
\end{figure}

%\textcolor{green}{As usual, we also do experiment on the $y$ coordinate of the peaks of solutions on the curve $p(\cos 2t,\ep)$.}

\begin{figure}[h]
\centering
\includegraphics[height=3.5cm]{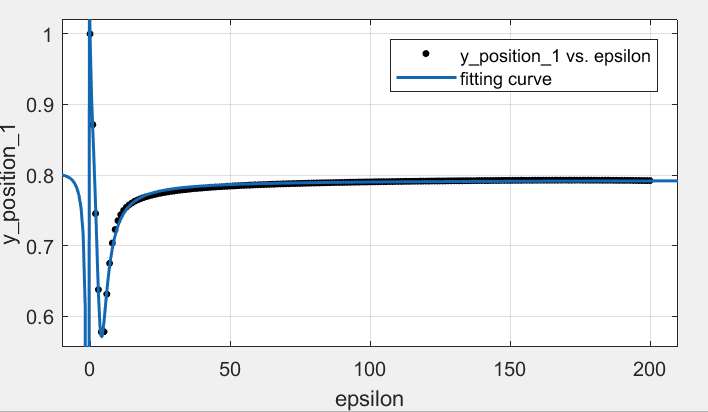}
\caption{The $t$ coordinate of the second peak of solutions for points $p(\cos 2t,\ep)$, with $\ep=1,2,\cdots, 200$. Fitting curve $t=\frac{a\ep^3+b\ep^2+c\ep+d}{\ep^3+e\ep^2+f\ep+g}$, with $a=0.7931,b=-4.642,c=6.596,d=13.45,e=--5.508,f=9.829$ and $g=13.44$.}
\label{peak1_yposition_cos2t}
\end{figure}

\begin{figure}[h]
\centering
\includegraphics[height=3.5cm]{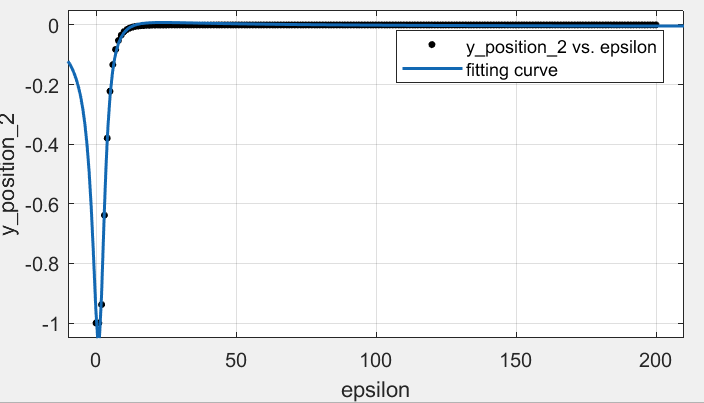}
\caption{The $t$ coordinate of the second peak of solutions for points $p(\cos 2t,\ep)$, with $\ep=1,2,\cdots, 200$. Fitting curve $y=\frac{a\ep^2+b\ep+c}{\ep^2+e\ep+f}$, with $a=-0.006875,b=0.7009,c=-8.537,d=-2.457,e=8.857$.}
\label{peak2_yposition_cos2t}
\end{figure}

\begin{figure}[h]
\centering
\includegraphics[height=3.5cm]{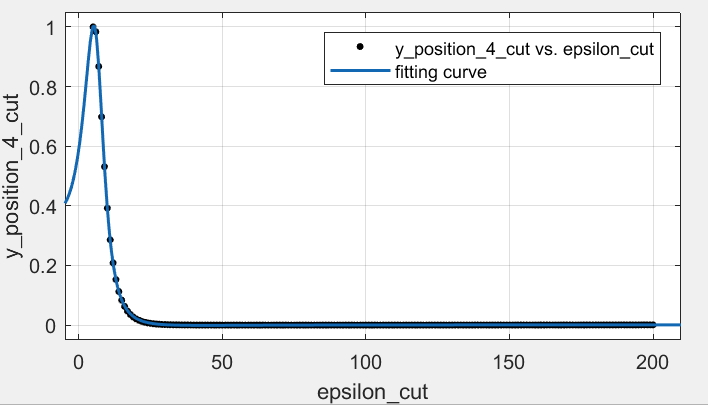}
\caption{The $t$ coordinate of the second peak of solutions for points $p(\cos 2t,\ep)$, with $\ep=5,6,\cdots, 200$. Fitting curve $y=\frac{a\ep^2+b\ep+c}{\ep^3+d\ep^2+e\ep+f}$, with $a=0.2268,b=-25.1,c=587.7,d=4.775,e=-156$ and $f=1002$.}
\label{peak4_yposition_cos2t}
\end{figure}

\begin{figure}[h]
\centering
\includegraphics[height=3.5cm]{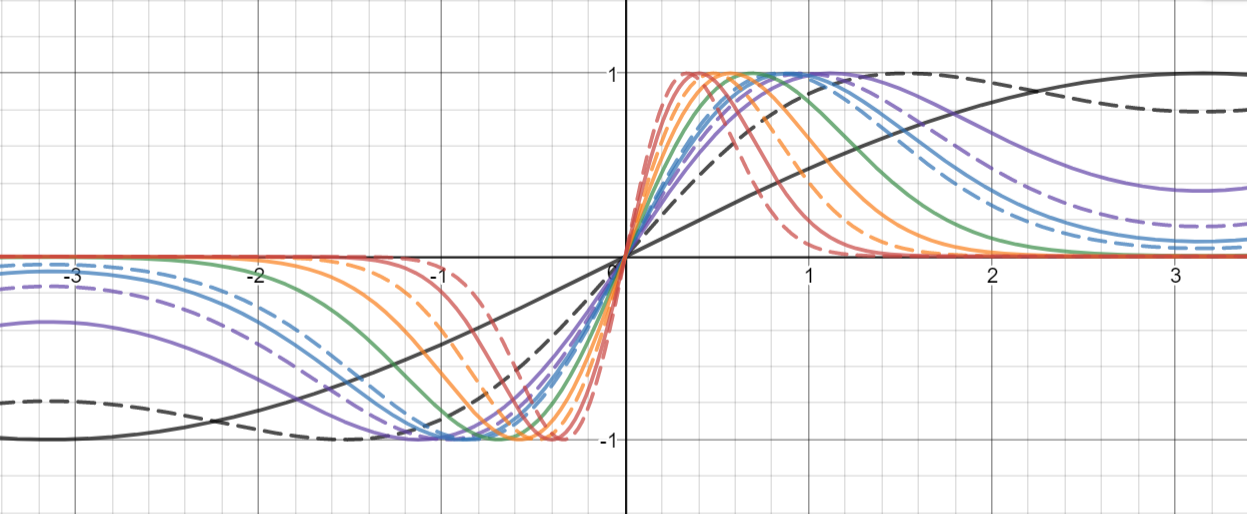}
\caption{Normalized solutions corresponding to $p(\sin \frac12t,\ep)$, with $\ep=0,1,2,3,4,5,10,20,40,80,160$.}
\label{solutionsin1/2t}
\end{figure}

\begin{figure}[h]
\centering
\includegraphics[height=3.5cm]{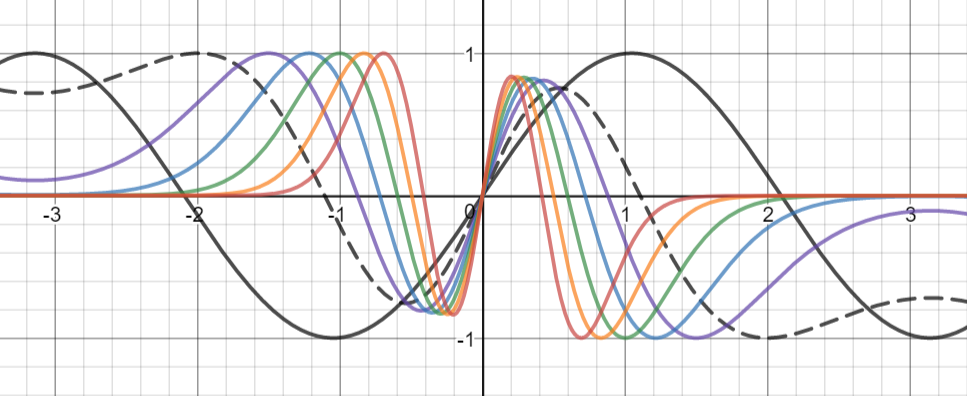}
\caption{Normalized solutions corresponding to $p(\sin\frac32t,\ep)$, with $\ep=0,5,10,20,40,80,160$.}
\label{solutionsin3/2t}
\end{figure}

\begin{figure}[h]
\centering
\includegraphics[height=3.5cm]{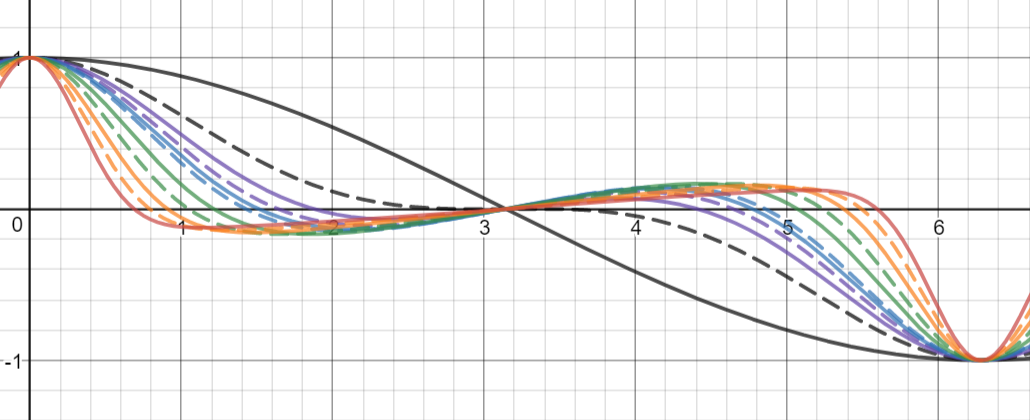}
\caption{Normalized solutions corresponding to $p(\cos \frac12t,\ep)$, with $\ep=0,1,2,3,4,5,10,20,40,80,160$.}
\label{solutioncos1/2t}
\end{figure}

\begin{figure}[h]
\centering
\includegraphics[height=3.5cm]{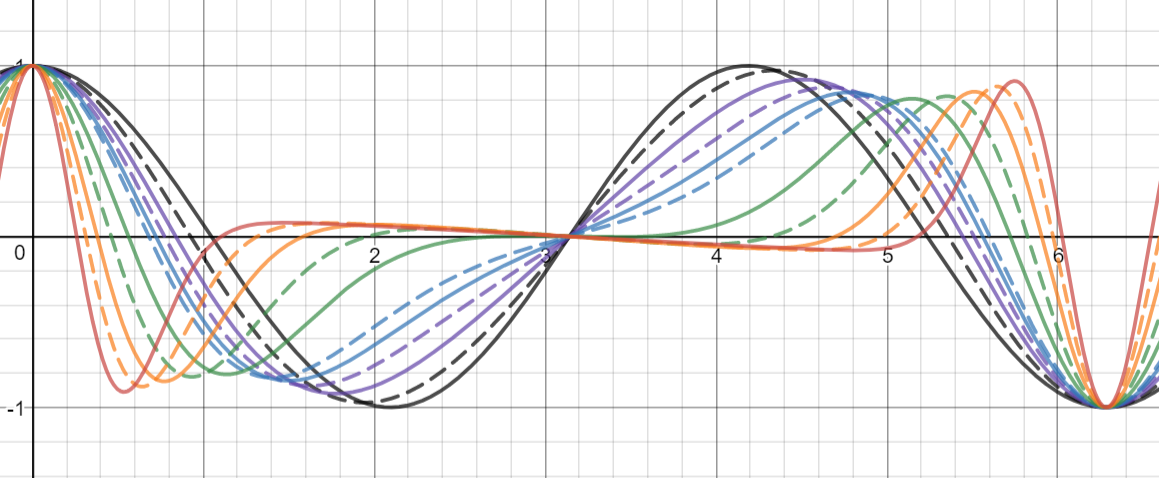}
\caption{Normalized solutions corresponding to $p(\cos \frac32t,\ep)$, with $\ep=0,1,2,3,4,5,10,20,40,80,160$.}
\label{solutioncos3/2t}
\end{figure}

\subsection{Explanation on the Behavior of the Solutions}
Before the end of this section, we give a proof on the convergence of the t-coordinates of the peaks. This interesting point is that the behaviors are closely related to the asymptotic behavior of the transition curves, i.e. how $\delta$ behaves as $\varepsilon$ increases to $\infty$. There are several works on this problem, and here we refer \cite{loud} for one version. Readers can find other dicussions in \cite{verticalconvergence}.

\begin{theorem} [W.S. Loud]\label{theorem42} For fixed $\ep>0$, let $\delta$ be the $k$th smallest value on the transition curves. Then we have
\[\delta=-\ep+(k-\frac{1}{2})\sqrt{2\ep}+O(\ep^{1/2}),\]
\textit{as} $\ep\to +\infty$.
\end{theorem}

Notice that in the above theorem, we fix a transition curve and let $\ep\to\infty$, which is exactly the way we study the behavior of the solutions. Now, we apply the theorem to derive some interesting facts that we have observed in our experiments. 

\begin{lemma}\label{lemma43} For fixed $\ep>0$, let $\delta$ be the $k$th smallest value on the transition curves, and $u$ be a nontrivial periodic solution for the corresponding MDE. Then, for large $\ep>0$, there is no local extremum of $u$ in $[c_\ep,\pi)$ where $c_\ep\sim \ep^{-1/4}$.
\end{lemma}

\textit{Proof.} Let $c_\ep$ be the zero of $\delta+\ep\cos t$ in $(0,\pi)$, where theorem \ref{theorem42} guarantee the existence of $c_\ep$ when $\ep$ is large. In addition, we have the estimate by using Taylor expansion of $\cos t$ locally at $0$, 
\[\lim_{\ep\to\infty} \frac{c_\ep}{\ep^{1/4}}=(2k-1)^{\frac{1}{2}}2^{\frac{3}{4}}.\]

\begin{comment}
-----------
We claim that $u(t)u'(t)<0$ for all $t\in(c_\ep,\pi)$. Suppose not. Then $u(x)u'(x)\geq0$ for some $x\in(c_\ep,\pi)$. We seek to derive a contradiction.

Note: Since $\sin t>0$ on $(0,\pi)$, $\frac{d}{dt}[\delta+\ep\cos t]=-\ep\sin t<0$ and thus $\delta+\ep\cos t$ is decreasing on this interval. Hence, $\delta+\ep\cos t$ can have only one zero in $(0,\pi)$ and hence has no zero in $(c_\ep,\pi)$. So, for $c_\ep<t<\pi$, $0=\delta+\ep\cos c_\ep>\delta+\ep\cos t$ and thus $0<-(\delta+\ep\cos t)$.

So, for each $t$ in $(0,\pi)$, it must be the case that $u(t)$ and $u''(t)$ are both positive, both negative, or both zero.

Case 1: $u(x)u'(x)>0$ and $u(x)>0$. Then $u'(x)>0$ and $u''(x)>0$. Now $u'(y)>0$ for all $x<y<\pi$. Suppose not. Then $u'(y)=0$ for some $y$.

Now, $u(t)$ increases on $(x,\pi)$; that is, $u'(t)>0$ on $(x,\pi)$. 

We claim that if $x<y<\pi$ then 

Case 2: $u(x)u'(x)>0$ and $u(x)<0$. Then $u'(x)<0.$ 

Case 3: $u(x)u'(x)=0$ and $u(x)=0$. Since we assume $u$ is a nontrivial solution, we must have $u'(x)\neq0$. Since $u'$ is continuous, there exists a small $h>0$ such that $|u'(t)|>0$ for all $t\in(x-h,x+h)$, so that $u$ is increasing on all of $(x-h,x+h)$ or is decreasing on all of $u$ and thus $u(t)\neq0$ on $(x,x+h)$.

Case 4: $u(x)u'(x)=0$ and $u'(x)=0$. Since we assume $u$ is a nontrivial solution, we must have $u(x)\neq0$.

\[\begin{cases}
\frac{d}{dt}u=u',\\
\frac{d}{dt}u'=(-\delta-\ep\cos t)u,
\end{cases}\]

----------\end{comment}
We will show $u'(t)\neq 0$ on $[c_\ep,\pi)$, which immediately implies the lemma.

First, we show there is no $t\in [c_\ep,\pi)$ such that $u(t)u'(t)>0$, by contradiction. Without loss of generality, we assume $u(t_0)>0$ and $u'(t_0)>0$ for some $t_0\in [c_\ep,\pi)$. Then we can show that $u(t)>0,u'(t)>0$ for any $t\in [t_0,\pi]$. If this is not true, we can take $t_1=\inf\{t\in [t_0,\pi]:u(t)\leq 0\text{ or }u'(t)\leq0\}$, then $u(t)$ and $u'(t)$ are both increasing on $[t_0,t_1]$ as 
\[\begin{cases}
\frac{d}{dt}u=u',\\
\frac{d}{dt}u'=-(\delta+\ep\cos t)u,
\end{cases}\]
where $-(\delta+\ep\cos t)>-(\delta+\ep\cos c_\ep)=0$ on $(c_\ep,\pi)$. As a consequence, $u(t_1)>0$ and $u'(t_1)>0$, which contradicts the definition of $t_1$, noticing that both $u,u'$ are continuous. As a result, we see that $u(t)>0,u'(t)>0$ for any $t\in [t_0,\pi]$. On the other hand, $u$ satisfies the boundary condition $u(\pi)=0$ or $u'(\pi)=0$, since $u$ takes a Fourier series expansion in one of the following forms (see \cite{stoker}): $\sum_{j=0}^\infty c_j\cos jt$, $\sum_{j=1}^\infty c_j\sin jt$, $\sum_{j=0}^\infty c_j\cos \frac{2j+1}{2}t$ or $\sum_{j=0}^\infty c_j\sin \frac{2j+1}{2}t$. This gives a contradiction.

Using the above observation, we can see  $u'(t)\neq 0$ for any $t\in [c_\ep, \pi)$. Otherwise, if $u'(t)=0$ for some $t\in [c_\ep, \pi)$, we should have $u(t)\neq 0$, so that $u(t+h)u'(t+h)>0$ for some small $h>0$. \hfill$\square$\\

\begin{theorem}\label{thm44}  (a). Consider a fixed transition curve characterized by the solution $u(t)=\sin kt$ or $u(t)=\cos (k+\frac{1}{2})t$. Fix a large $\ep>0$, and let $u$ be a nontrivial periodic solution. Then, we have all the local extremums of $u$ in $\bigcup_{n\in \mathbb{Z}}(-c_\ep+2n\pi,c_\ep+2n\pi)$, where $c_\ep$ tends to $0$ as $\ep\to \infty$. ($c_\ep$ depends on $k$)

(b). Consider a fixed transition curve characterized by the solution $u(t)=\cos kt$ or $u(t)=\sin (k+\frac{1}{2})t$. Fix a large $\ep>0$, and let $u$ be a nontrivial periodic solution. Then, we have all the local extremums of $u$ in $\pi\mathbb{Z} \cup\big(\bigcup_{n\in \mathbb{Z}}(-c_\ep+2n\pi,c_\ep+2n\pi)\big)$, where $c_\ep$ tends to $0$ as $\ep\to \infty$. ($c_\ep$ depends on $k$)

In addition, there is a critical value $\alpha\geq 0$ given by the equation
 $$\delta-\ep=0,$$
 where we view $\delta$ as a function of $\ep$ on a fixed transition curve. Then $\pi$ is a local maximum of $|u|$ if $0\leq \ep<\alpha$, and $\pi$ is a local minimum of $|u|$ is $\ep>\alpha$. When $\ep=\alpha$ and the curve is not characterized by $\cos 0t=1$, $\pi$ is a local maximum of $|u|$.
\end{theorem}

\textit{Proof.} (a) Lemma \ref{lemma43} shows that there are no local extremums in $[c_\ep,\pi)$. By symmetry, we do not have local extremums in $\bigcup_{n\in\mathbb{Z}} [c_\ep,2n\pi+\pi)\cup (2n\pi+\pi,c_\ep]$. It remains to show $(2k+1)\pi,k\in\mathbb{Z}$ are not local extremums in these cases. It is enough to show that $\pi$ is not a local extremum.

Note that for any case in part (a), a nontrivial $2\pi$- or $4\pi$-periodic solution $u$ corresponding to a point $(\delta,\ep)$ on a transition curve is written as a sum of sines or as a sum of cosines:
$u(t)=\sum_{k=0}^\infty a_k\sin(kt)$ or $u(t)=\sum_{k=0}^\infty a_k\cos(k+\frac{1}{2})t$.

It is easy to see that $u(\pi)=0$. We therefore conclude that $u'(\pi)\neq0$, since otherwise $u$ would be the trivial solution. As a result, $\pi$ is not a local extremum.
\begin{comment}
\st{(as a solution to the MDE is uniquely determined by specifying the value of $u$ and $u'$ at a particular point)}

\st{In the case where $u$ is a sum of cosines, we have that $u$ is an odd function around $t=\pi$, since it is a sum of functions which are  Since $u$ is a nontrivial solution, we therefore conclude that $u'(0)\neq0$, since otherwise $u$ would be the trivial solution (as a solution to the MDE is uniquely determined by specifying the value of $u$ and $u'$ at a particular point).}

\st{(a) is then an immediate consequence of lemma 4.3, since we can check that $\pi$ is not a local extremum since $u'(\pi)\neq 0$.}
\end{comment}

(b) We only need to understand the behavior of $u$ at $\pi$. The idea is essentially the same as the proof of Lemma 4.3, where we look at the sign of $\delta+\ep\cos t$. It is also clear that $u(\pi)\neq 0$ and $u'(\pi)=0$ in this case.

Notice that $u''(\pi)=-(\delta-\ep)u(\pi)$, and we set $u(\pi)>0$ by multipling a constant, so that we can conveniently look at the absolute value of $u$. If $\delta+\ep\cos \pi=\delta-\ep>0$, then $\pi$ is a local maximum of $|u|$; if $\delta-\ep<0$, then $\pi$ is a local minimum of $|u|$. For the case that $\ep=\delta$, we can see that $\delta-\ep\cos t>0$ on $(0,2\pi)\setminus \{\pi\}$ except for the case $\delta=\ep=0$, which only happens on the curve characterized by $\cos 0t$. So $\pi$ is still a local maximum of $|u|$. 

In addition, $\delta-\ep$ is a strictly decreasing function of $\ep$ on $(0,\infty)$, so we can find a critical point $\alpha$ such that $\delta-\ep>0$ on $(0,\alpha)$, and $\delta-\ep<0$ on $(\alpha,\infty)$. To show that $\delta-\ep$ is strictly decreasing, notice that $\delta$ is the $k$th smallest value such that the MDE has a $2\pi$- or $4\pi$-periodic solution, which is equivalent to say that $\delta$ is the $k$th smallest eigenvalue of the self-adjoint operator $-\Delta-\ep\cos t$ on $L^2(\mathbb{R}/4\pi\mathbb{Z})$. Let $\ep>\ep'>0$, and $\delta,\delta'$ be the corresponding eigenvalues. Then using Raleigh quotient, we get
\[\begin{aligned}
\delta'&=\max_{u\in M_{k,\ep'}} \frac{\langle(-\Delta-\ep'\cos t) u,u\rangle_{L^2(\mathbb{R}/4\pi\mathbb{Z})}}{\|u\|_{L^2(\mathbb{R}/4\pi\mathbb{Z})}}\\
&>\max_{u\in M_{k,\ep'}} \frac{\langle(-\Delta-\ep'\cos t-(\ep-\ep')(\cos t+1)) u,u\rangle_{L^2(\mathbb{R}/4\pi\mathbb{Z})}}{\|u\|_{L^2(\mathbb{R}/4\pi\mathbb{Z})}}\\
&=\max_{u\in M_{k,\ep'}} \frac{\langle(-\Delta-\ep\cos t) u,u\rangle_{L^2(\mathbb{R}/4\pi\mathbb{Z})}}{\|u\|_{L^2(\mathbb{R}/4\pi\mathbb{Z})}}-(\ep-\ep')\\
&>\inf_{M} \max_{u\in M} \frac{\langle(-\Delta-\ep\cos t) u,u\rangle_{L^2(\mathbb{R}/4\pi\mathbb{Z})}}{\|u\|_{L^2(\mathbb{R}/4\pi\mathbb{Z})}}-(\ep-\ep')=\delta-(\ep-\ep'),
\end{aligned}\]
where $M_{k,\ep'}$ is the space spanned by the smallest $k$ eigenfunctions of $-\Delta-\ep'\cos t$ and the infimum in the last line is taken over all the $k$ dimensional subspaces of $dom\Delta$ on $\mathbb{R}/4\pi\mathbb{Z}$.  \hfill$\square$\\

\textbf{Remark.} In Theorem 4.2, Lemma 4.3 and Theorem 4.4, we considered the case that $\ep\geq 0$. Since $-\cos t=\cos (t+\pi)$, for the case that $\ep\leq 0$, the MDE can be rewritten as 
$$\frac{d^2}{dt^2}\hat{u}+(\delta-\ep\cos t)\hat{u}=0,$$
with $\hat{u}(t)=u(t-\pi)$. We can still use Theorem 4.2, Lemma 4.3 and Theorem 4.4 to study solutions to the MDE with coefficient $\ep<0$, by applying a shift of $\pi$.

%\clearpage

\section{The Sierpinski Gasket and the Fractal Laplacian}
In this section we give preliminary definitions and results concerning analysis on fractals and provide a basic introduction to the `infinite' Sierpinski gasket will be given. This serves to set up the discussion of the generalization of the MDE to the fractal setting described in in Section 6.

\subsection{Sierpinski Gasket}

\begin{figure}[h]
 \centering
 \includegraphics[scale = .45]{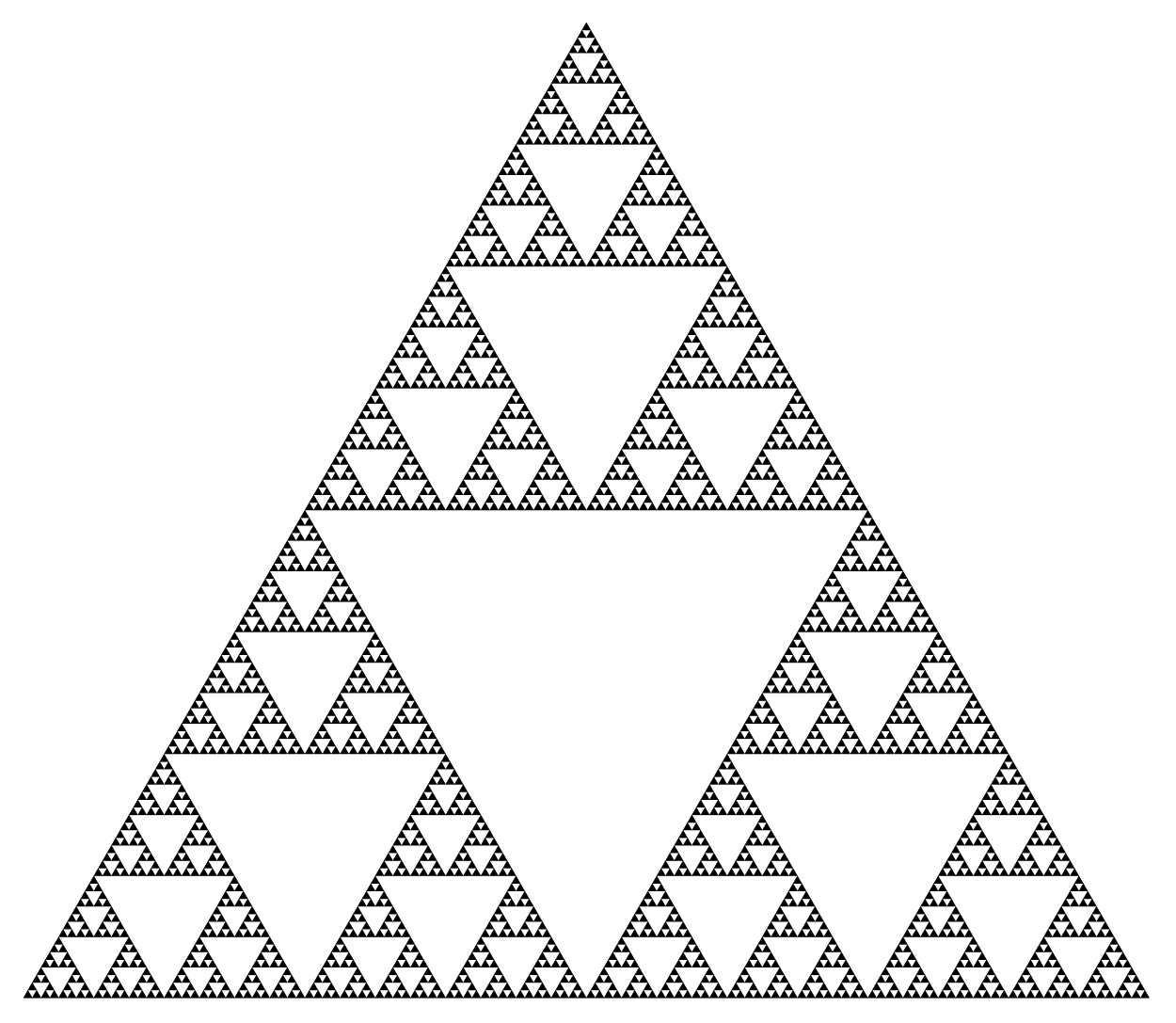}
 \caption{The Sierpinski Gasket}\label{SG}
\end{figure}

Consider the three contraction mappings $\{F_i:\mathbb{R}^2\to\mathbb{R}^2\}_{i=0,1,2}$ given by
\begin{equation*}
\begin{cases}
    F_0(x)=\frac{1}{2}x\\
    F_1(x)=\frac12x+\left(\frac12,0\right)\\
    F_2(x)=\frac12x+\left(\frac14,\frac{\sqrt{3}}{4}\right),
\end{cases}
\end{equation*}
where $x\in\mathbb{R}^2$.  Then $\{F_i\}_{i=0,1,2}$ form an `iterated function system' (see page 133 of \cite{falconer}). By Theorem 9.1 in \cite{falconer}, there exists a unique nonempty compact set $K\subset\mathbb{R}^2$ such that (see \cite{Hutchinson})
$$K=\bigcup_{i=0}^2F_i(K).$$ 
Then $K$ is defined to be the \textit{Sierpinski gasket}, often denoted $SG$. See Figure \ref{SG}.

\begin{figure}[h]\label{agra}
 \centering
 \includegraphics[scale = .35]{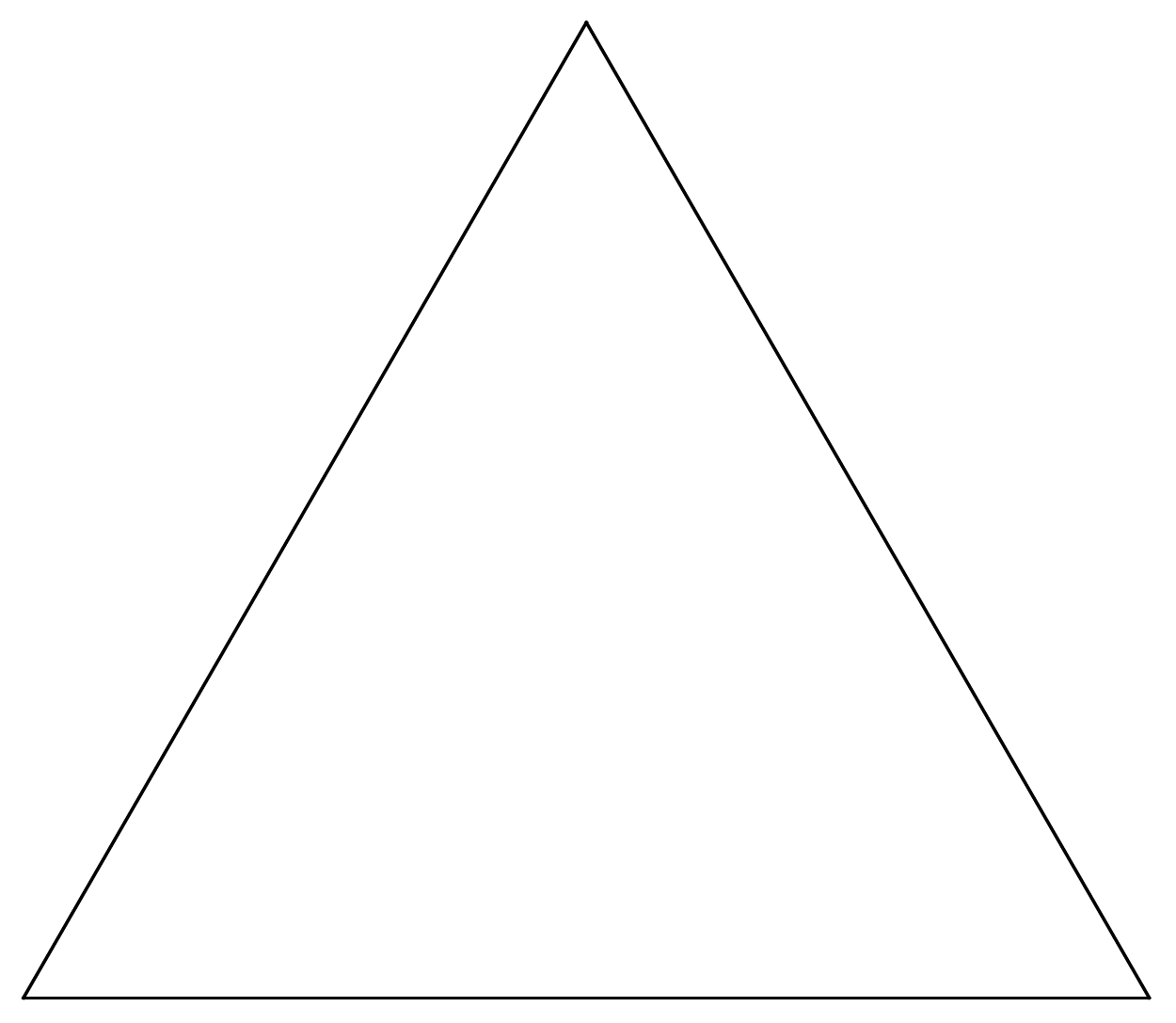}\qquad
 \includegraphics[scale = .35]{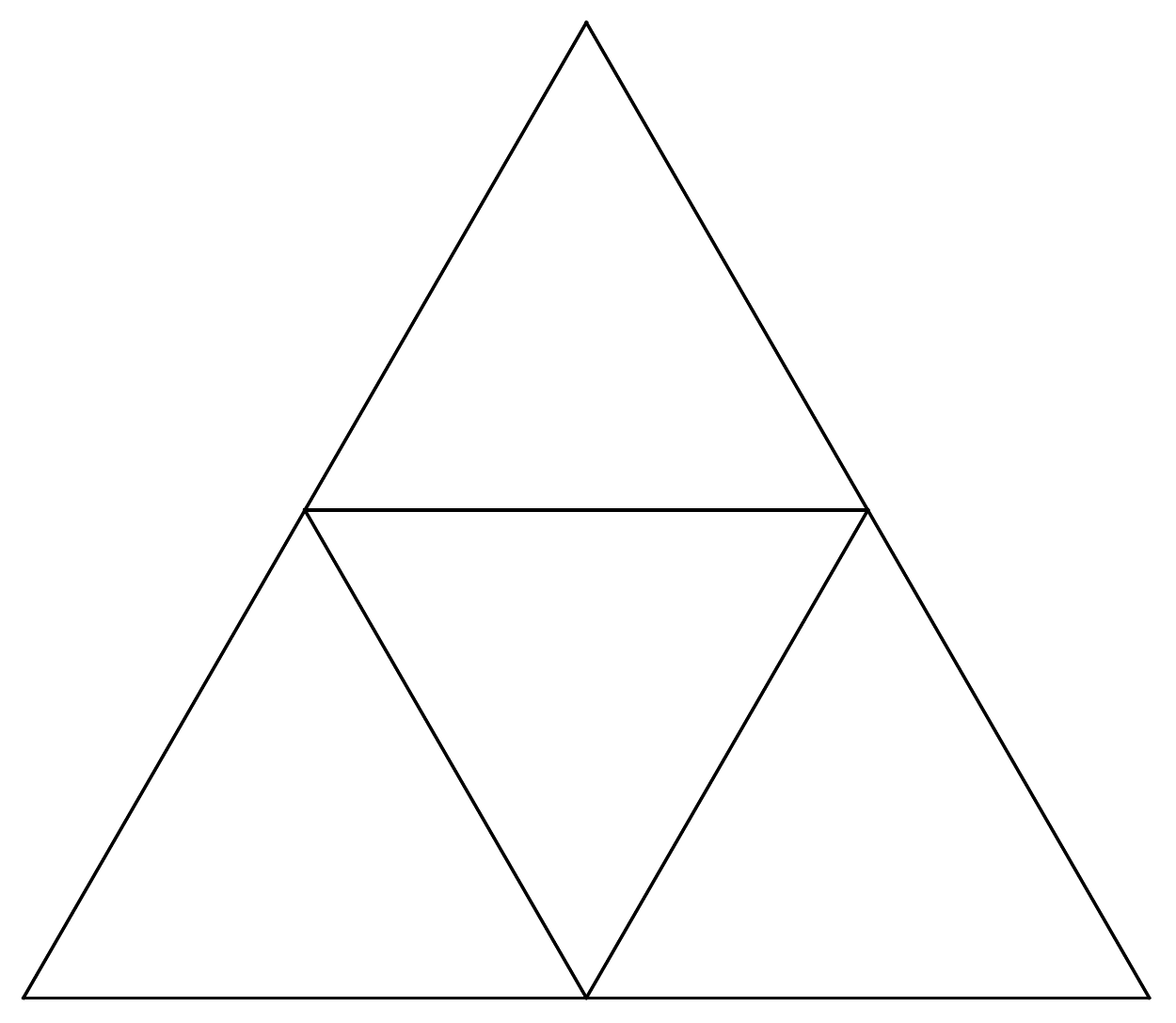}\qquad
 \includegraphics[scale = .35]{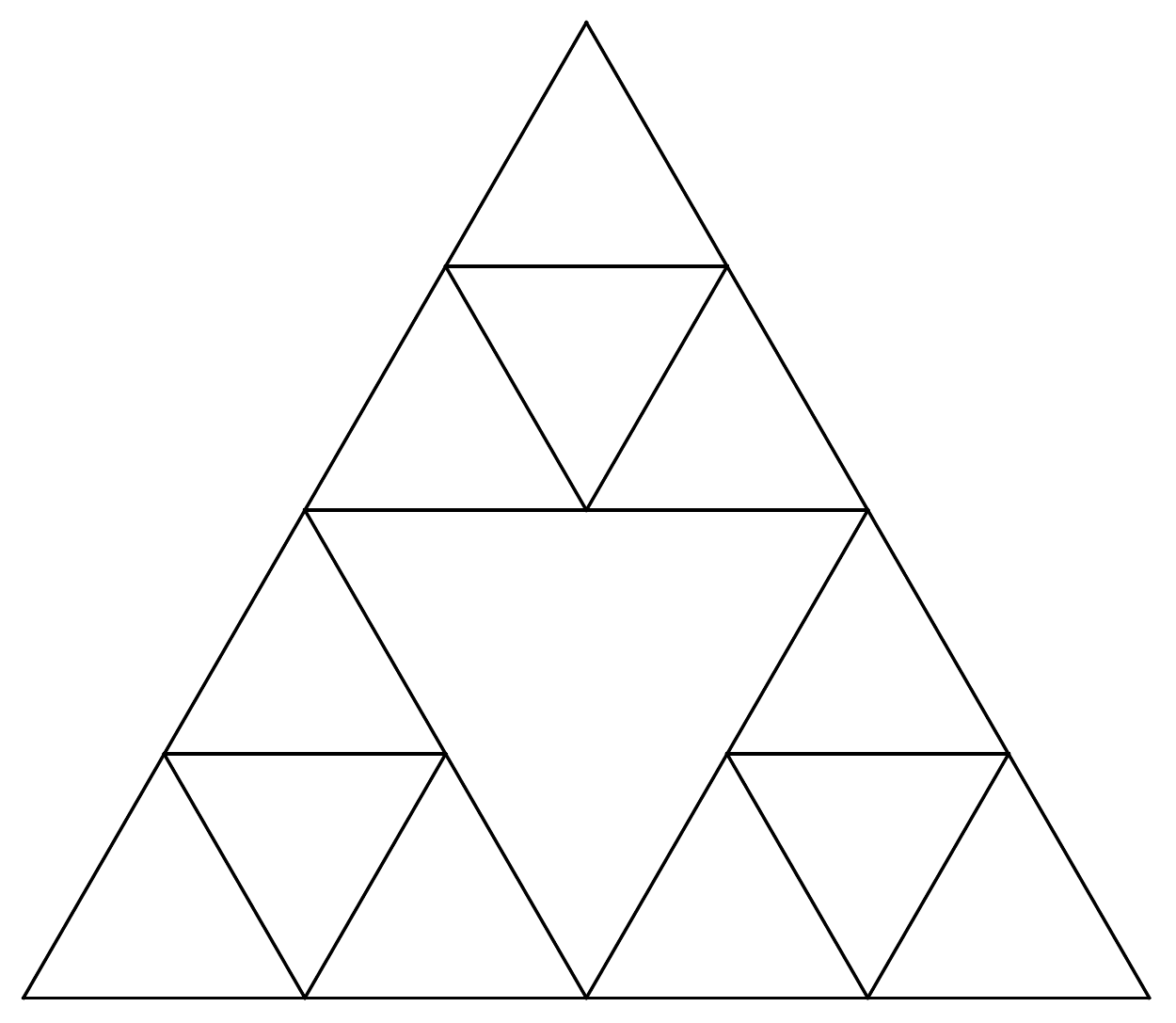}\qquad
 \caption{Approximating graphs $\Gamma_0,\Gamma_1,\Gamma_2$}
\end{figure}

In studying $SG$ it is useful to use its graph approximations, constructed as follows. Let $q_0, q_1,$ and $q_2$ be the unique fixed points of $F_0, F_1,$ and $F_2,$ respectively. Define a vertex set $V_0:=\{q_0,q_1,q_2\}\subset SG$. Then $F_i(SG)\cap F_j(SG)=F_iV_0\cap F_jV_0$,for any $i\neq j$. We refer to $V_0$ as the \textit{boundary} of $SG$. We further define vertex sets $V_n\subset SG$ for $n\geq1$ inductively by  $V_n:=\bigcup_{i=0}^2F_i(V_{n-1})$ and let $V_*:=\bigcup_{n=0}^\infty V_n$ be the set of all vertices. Note that $V_*$ is dense in $SG$. For an $m$-tuple $w=\left(w_1,w_2,\dotsb,w_m\right)$, where $w_j\in\{0,1,2\}$ for each $w_j\text{ } (1\leq j\leq m)$, we define $F_w$ by $$F_w:=F_{w_1}\circ F_{w_2}\circ\dotsb\circ F_{w_m}$$ and say that $w$ is a \textit{word} of length $\left|w\right|=m$. With this, an edge relation $\sim_m$ on $V_m$ can be introduced as follows: for $x,y\in V_m$ we say $x\sim_m y$ if and only if there exists a word $w$ of length $m$ and unequal indices $i,j\in\{0,1,2\}$ such that $x=F_w(q_i)$ and $y=F_w(q_j)$. This relation on $V_m$ gives a sequence of graphs $\Gamma_m$ approximating SG, with the vertex set $V_m$ and the edge set $E_m=\{\{x,y\}|x\sim_m y\}$. See Figure \ref{agra} for $\Gamma_0,\Gamma_1,\Gamma_2$.

\subsection{Fractal Laplacian}

Now we are ready to define the \textit{Laplacian} on $SG$. Suppose $u:SG\to\mathbb{R}.$ We define the level-$m$ \textit{discrete Laplacian} $\Delta_m$ by

$$
\Delta_mu(x):=\sum_{y\sim_mx}{\big(u(y)-u(x)\big),\quad\quad x\in V_m\setminus V_0}.
$$
Then we define the \textit{continuous Laplacian} $\Delta$ by
$$
\Delta u(x)=\frac32\lim_{m\to\infty}5^m\Delta_mu(x),\quad\quad x\in V_*\setminus V_0.$$
If the limit above converges uniformly on $V_*\setminus V_0$ to a continuous function, we say $u\in\text{dom} \Delta$ . In this case, we extend $\Delta u$ to all of $SG$, including points not in $V_*\setminus V_0$, by continuity (recall that $V_*\setminus V_0$ is dense in $SG$). The continuous Laplacian on SG is the analog of the usual `second-order derivative' on the line.

We shall mention now the the following proposition derived by O. Ben-Bassat, R.S. Strichartz,and A. Teplyaev in \cite{Oben}, which will be of interest in the next section.\\

\begin{theorem}[O. Ben-Bassat, R.S. Strichartz \& A. Teplyaev]\label{notindomain}
    Let $u$ be a nonconstant function in dom$\Delta.$ Then $u^2$ is \textit{not} in dom$\Delta$.
\end{theorem}

A function $u$ satisfying $-\Delta u=\lambda u$ for some number $\lambda$ is called an \textit{eigenfunction} of $\Delta$ with \textit{eigenvalue} $\lambda$.

If a function $u:SG\to\mathbb{R}$ satisfies $u|_{V_0}=0,$ then we say that $u$ satisfies the \textit{Dirichlet boundary condition}, and the eigenvalue problem
 \[\begin{cases}
      -\Delta u=\lambda u,\\
      u|_{V_0}=0
\end{cases}\]
is called the \textit{Dirichlet eigenvalue problem}. A function $u$ satisfying both of these equations is called a \textit{Dirichlet eigenfunction} of $\Delta.$

Similarly, we have a notion of a `Neumann condition' as follows. Define the \textit{normal derivative} $\partial_n(q_i)$ for $i\in\{0,1,2\}$ by
$$\partial_nu(q_i):=\lim_{m\to\infty}\left(\frac53\right)^m\big(2u(q_i)-u(F_i^m(q_{i+1}))-u(F_i^m(q_{i-1}))\big),$$
where we have identified indices modulo 3.

If a function $u:SG\to\mathbb{R}$ satisfies $\partial u(q_i)=0$ for $i=0,1$ and $2,$ then we say that $u$ satisfies the \textit{Neumann boundary condition}, and the eigenvalue problem
\[\begin{cases}
      -\Delta u=\lambda u,\\
      \partial u(q_i)=0\text{, for } i=0,1,2
\end{cases}\]
is called the \textit{Neumann eigenvalue problem}. A function $u$ satisfying both of these equations is called a \textit{Neumann eigenfunction} of $\Delta.$

Dirichlet and Neumann eigenfunctions on $SG$ are the analog of sine and cosine functions on the line.

\subsection{Spectral Decimation}
A method for explicitly computing all possible eigenvalues and eigenfunctions of $\Delta$ was introduced in \cite{decimation} using a process called \textit{spectral decimation}. Below we briefly discuss some results from spectral decimation we will use. Readers can find detailed discussion on spectral decimation in \cite{decimation} and \cite{differentialequationsonfractals}.\\

\begin{proposition}
Suppose $\lambda_m\neq 2,5,6$, and $\lambda_{m-1}$ is given by
\begin{equation}\label{decimation}
\lambda_{m-1}=\lambda_m(5-\lambda_m).
\end{equation}
(a) If $u$ is a $\lambda_{m-1}$-eigenfunction of $\Delta_{m-1}$ on $V_{m-1}$, then it can be uniquely extended to be a $\lambda_m$-eigenfunction of $\Delta_m$ defined on $V_m$.

(b) Conversely, if $u$ is a $\lambda_m$-eigenfunction of $\Delta_m$ on $V_m$, then $u|_{V_{m-1}}$ is a $\lambda_{m-1}$-eigenfunction of $\Delta_{m-1}$ on $V_{m-1}$.
\end{proposition}

If we want to extend an eigenfunction of $\Delta_{m-1}$ with eigenvalue $\lambda_{m-1}$ to an eigenfunction of $\Delta_m$ using the proposition above, we have two choices, except when $\lambda_{m-1}=6$ (as we will see below), in which to extend the eigenfunction:
$$\lambda_m=\frac{5\pm\sqrt{25-4\lambda_{m-1}}}{2}.$$ For convenience, define the functions $\psi_+(x)$ and $\psi_-(x)$ by
$$\psi_-(x)=\frac{5-\sqrt{25-4x}}{2},\quad\psi_+(x)=\frac{5+\sqrt{25-4x}}{2}.$$

The numbers $2,5,6$ are called \textit{forbidden eigenvalues}, and it turns out that each Dirichlet eigenfunction of $\Delta$ comes from a 2-, 5-, or 6-eigenfunction of $\Delta_{m_0}$ for some $m_0\geq 0$, while all the Neumann eigenfunctions come from 5- or 6-eigenfunctions of $\Delta_{m_0}$ for some $m_0\geq 0$. If $u$ is a Dirichlet or Neumann eigenfunction, we call $m_0$ the \textit{generation of birth} and $u|_{V_{m_0}}$ the \textit{initial function}.

Suppose $u$ is an eigenfunction of $\Delta$ with eigenvalue $\lambda$ arising from initial eigenvalue $\lambda_{m_0}.$ Then we say that $\lambda$ is a \textit{2-series} (resp., \textit{5-series} or \textit{6-series}) eigenvalue if $\lambda_{m_0}=2$ (resp., $\lambda_{m_0}=5$ or 6).

With a fixed generation of birth $m_0$ and initial function $f$, we can extend the function level-by-level. We can first fix a sequence $\{\varepsilon_m\}_{m=1}^\infty$ of $\pm$, with only finitely many $-$, and then let $\lambda_{m}=\psi_{\varepsilon_{m-m_0}}(\lambda_{m-1})$ for $m>m_0$ inductively. Then the function is extended to be an eigenfunction of $\Delta$, with the corresponding eigenvalue
\[\lambda:=\frac{3}{2}\lim_{m\to\infty}5^m\lambda_m.\]
In fact, all the eigenfunctions with a given generation of birth and initial function can be generated by the above recipe. Also, if the initial eigenvalue is $6$, we can only choose $\varepsilon_1=+1$, as $\psi_-(6)=2$ is a forbidden eigenvalue. 

Now, suppose we fix a generation of birth $m_0$, fix an initial function on $V_{m_0}$, and let $$E=\big\{e=\{e_j\}_{j=1}^N:N\in\mathbb{N} \text{ and } e_j\in\{+,-\},e_1=+\big\}\cup\emptyset,$$ where $\emptyset$ denotes the empty sequence. Define $\psi_e=\psi_{e_1}\psi_{e_2}\cdots \psi_{e_{|e|}}$, where $|e|$ is the length of $e$. In particular, $\psi_\emptyset(x)=x$.

Let $\Psi(x)=\frac{3}{2}\lim_{l\to\infty}5^{l}\psi^l_-(x)$. Then, we can deduce the following possibilities:

1. If $\lambda_{m_0}=2$ or $5$, then all the possible eigenvalues of $\Delta$ having generation of birth $m_0$ are given by
$$\begin{aligned}
\lambda_e&=\frac{3}{2}\lim_{l\to\infty}5^{|e|+l+m_0}\psi^l_{-}\psi_e(\lambda_{m_0})\\&=5^{|e|+m_0}\Psi(\psi_e(\lambda_{m_0})), \quad e\in E.
\end{aligned}$$

2. If $\lambda_{m_0}=6$, then all the possible eigenvalues of $\Delta$ having generation of birth $m_0$ are given by
$$\begin{aligned}
\lambda_e&=\frac{3}{2}\lim_{l\to\infty}5^{|e|+l+m_0+1}\psi^l_{-}\psi_e(\psi_+(6))\\&=5^{|e|+m_0+1}\Psi(\psi_e(3)), \quad e\in E.
\end{aligned}$$
Here we remark that when we use the notation $\lambda_e$, we always assume that we have a fixed generation of birth $m_0$ and initial eigenvalue $\lambda_{m_0}$.

There is a method, given by \cite{fde}, in which to arrange in increasing order the set of eigenvalues arising a fixed generation of birth and initial eigenvalue. The idea is to translate each finite sequence in $E$ into a binary number. The process is as follows.

Given $e\in E$ of length $|e|=n$, let $d(e)$ be the integer with binary (base-2) expansion $d(e)=\sum_{i=1}^{n} 2^{n-i} d_i(e)$, where 
\[
d_i(e)=\begin{cases} 1, &\text{ if }i=1,\\1-d_{i-1}(e),&\text{ if }i\neq 1\text{ and }e_i=+,\\
d_{i-1}(e),&\text{ if }i\neq 1\text{ and }e_i=-.\end{cases}
\]
In addition, we set $d(\emptyset)=0$. Then $\lambda_e$ is the $\big(d(e)+1\big)$-th smallest eigenvalue. For example, if $e=(e_1,e_2,e_3,e_4,e_5,e_6,e_7)=(+,-,+,+,+,-,-)$, then $d(e)=d(+,-,+,+,+,-,-)=1101000_2=104,$ where $1101000_2$ is written in base-2. Thus, $\lambda_e$ corresponding to the sequence $e$ is the 105th smallest eigenvalue. Also, note that, in particular, $\lambda_{e=\emptyset}$ is the smallest eigenvalue, $\lambda_{e=(+)}$ is the 2nd smallest eigenvalue, and $\lambda_{e=(+,+)}$ is the 3rd smallest eigenvalue. We can, of course, reverse this process so that, given an in integer $d$ we can find the $e\in E$ corresponding to the $d$-th smallest eigenvalue.

Another fact is that the sequence of eigenvalues $\{\lambda_{n,m_0}\}_{n=1}^\infty$, where $\lambda_{1,m_0}<\lambda_{2,m_0}<\lambda_{3,m_0}<\lambda_{4,m_0}<...$, corresponding to a fixed generation of birth $m_0$ and a fixed initial eigenfunction grow according to the power law $n^{\log 5/\log 2}$. In addition, if we define the eigenvalue counting function $\rho_{m_0}:\mathbb{R}_+\to\mathbb{N}$ to be 
\[\rho_{m_0}(x)=\#\{e\in E:\lambda_{e,m_0}\leq x\},\]
then we have the following proposition concerning the asymptotic behavior of $\rho_{m_0}$.\\

\begin{proposition}
For a fixed generation of birth $m_0$ and fixed initial function on $V_{m_0}$, there exists a $\log 5$-periodic continuous function $g(t)$ such that 
\[\lim_{x\to\infty} \left(\frac{\rho_{m_0}(x)}{x^{\log2/\log5}}-g(\log x)\right)=0.\]
\end{proposition} 

\textit{Proof.} To avoid the high multiplicity of eigenfunctions corresponding to a same eigenvalue, we fix an initial eigenfunction instead of just fixing an initial eigenvalue. In our setting, the eigenfunction is unique for each eigenvalue, so we only need to count the number of eigenvalues.

For convenience, we prove the proposition for generation of birth $m_0$ and initial eigenvalue $5$. For initial eigenvalue $2$ and $6$, the arguments are essentially the same. We will show that $\frac{\rho_{m_0}(5^nx)}{2^nx^{\log2/\log5}}$ converges to $g(\log x)$ uniformly on some interval $[e^c,5e^c)$ as $n\to \infty$. First, we have the following observations.

Observation 1: For eigenvalues of generation of birth $m_0$ and initial eigenvalue $5$,  if $n$ is fixed then we have 
$$5^{-|e|}\lambda_{e}\in \cup_{|e'|=n+1,e'\in E}5^{m_0}\Psi\psi_{e'}([0,5])=\cup_{l\in \{-,+\}^n}5^{m_0}\Psi\psi_+\psi_l([0,5]),$$
for all $|e|>n$ where $\psi_l=\psi_{l_1}\circ\psi_{l_2}\cdots \psi_{l_n}$ for each word $l=(l_1,l_2,\cdots,l_n)\in\{-,+\}^n=\big\{l=\{l_j\}_{j=1}^n:l_j=+\text{ or }-\big\}$.

Notice that $5^{-{|e|}}\lambda_{e}\in 5^{m_0}\Psi\circ \psi_e([0,5])$, and the fact $\psi_-([0,5])\subset [0,5],\psi_+([0,5])\subset [0,5]$, we have $5^{-{|e|}}\lambda_{e}\in 5^{m_0}\Psi\circ \psi_{e_1e_2\cdots e_{n+1}}([0,5])$. The observation follows.

For each $l\in\bigcup_{n=0}^\infty\{-,+\}^n$, we denote the endpoints of $5^{m_0}\Psi\psi_+\psi_l([0,5])$ by $a_l$ and $b_l$, i.e., $5^{m_0}\Psi\psi_+\psi_l([0,5])=[a_l,b_l]$. By some easy computation, we can get the following observation.

Observation 2: Assume $l=(l_1,l_2,\cdots,l_n)$. Then $$\rho_{m_0}(5^kb_l)=2^{k-1}+2^{k-n-1}\big(1+d(+l)-2^n\big)=2^{k-n-1}(1+d(+l)),\text{ if }k>n.$$

To show observation 2, we need to consider two cases. First, for any $e$ with $|e|\leq k-1$, we have $\lambda_e<5^{k}b_l$, which counts for $2^{k-1}$ eigenvalues. Second, we have $2^{k-n-1}$ eigenvalues $\lambda_e$ in each interval of the form $5^k5^{m_0}\Psi\psi_+\psi_{l'}([0,5])$, where $l'\in\{+,-\}^n$. In fact, for $|e|=k$,  $\lambda_e\in 5^k5^{m_0}\Psi\psi_+\psi_{l'}([0,5])$ if and only if $(e_1,e_2,\cdots,e_{n+1})=(+,l'_1,l'_2,\cdots,l'_n)$, and we have $2$ free choices for each $e_{j},n+2\leq j\leq k$. There are $d(+l)+1-2^n$ intervals in $[0,b_l]$, noticing that $b_{l'}\leq b_{l}$ if and only if $d(+l')\leq d(+l)$. Combining the above facts, we get the second term $2^{k-n-1}\big(1+d(+l)-2^n\big)$.

Now, fix $n$ and consider $l\in \{-,+\}^n$. It is easy to see that $\frac{\rho_{m_0}(5^kb_l)}{2^kb_l^{\log 2/\log 5}}$ converges as $k\to\infty$, since it is a constant for $k\geq n+1$. We denote the limit $g(\log b_l)$. 

Next, we look at general $x$. Note that we can find some constant $c$ such that $$\log(5^{m_0}\Psi(\psi_+([0,5])))\subset [c,c+\log 5),$$ since $5^{m_0-1}\Psi\psi_+(5)<5^{m_0}\Psi\psi_+(0)$. We want to show that $\frac{\rho_{m_0}(5^nx)}{2^nx^{\log2/\log5}}$ converges to some function $g(\log x)$ uniformly on $c\leq\log x\leq c+\log5$ as $n\to \infty$.

In fact, if we fix $n$ and look at $b_l\leq x\leq b_{l'}$ for some $l,l'\in \{+,-\}^n$ such that $d(+l')=d(+l)+1$, we have 
\[\begin{aligned}
&g(\log b_l)\left(\frac{b_l}{x}\right)^{\log2/\log5}=\frac{\rho(5^kb_l)}{2^k{b_l}^{\log2/\log5}}\cdot \left(\frac{b_l}{x}\right)^{\log2/\log5}=\frac{\rho(5^kb_l)}{2^kx^{\log2/\log5}}\\&\leq\frac{\rho(5^kx)}{2^kx^{\log2/\log5}}\leq\frac{\rho(5^kb_{l})+2^{k-n-1}}{2^kx^{\log2/\log5}}=g(\log b_l)\left(\frac{b_l}{x}\right)^{\log2/\log5}+x^{-\log2/\log5}2^{-n-1},
\end{aligned}\]
for any $k\geq n+1$. In addition, if $e^c\leq x\leq\min_{l\in\{+,-\}^n}a_l$, we have
\[ \frac{\rho_{m_0}(5^kx)}{2^kx^{\log2/\log5}}= \frac{\rho_{m_0}(5^{k-1}b_l)}{2^kx^{\log2/\log5}}=\frac{1}{2}g(\log b_l)\left(\frac{b_l}{x}\right)^{\log2/\log5},\]
where $l=(-,-,-,\cdots)\cdots\in \{+,-\}^n$. Similarly, we have for $\max_{l\in\{+,-\}^n}b_l\leq x\leq 5e^c$, 
\[\frac{\rho_{m_0}(5^kx)}{2^kx^{\log2/\log5}}=\frac{\rho_{m_0}(5^kb_l)}{2^kx^{\log2/\log5}}=g(\log b_l)\left(\frac{b_l}{x}\right)^{\log2/\log5},\]
where $l=(+,-,-,\cdots)\in \{+,-\}^n$.

The above discussions shows that  $\frac{\rho_{m_0}(5^nx)}{2^nx^{\log2/\log5}}$ converges to some function $g(\log x)$ uniformly on $c\leq\log x\leq c+\log5$ as $n\to \infty$. Also we can easily see that $g$ is continusous with $g(c)=g(\log 5)$ from the estimates. In fact, for each $x$, we can find a small neighbourhood such that for $n>k$ and any $y$ in the neighborhood, $\left|\frac{\rho(5^ky)}{2^ky^{\log2/\log5}}-\frac{\rho(5^kx)}{2^kx^{\log2/\log5}}\right|<2^{-n}x^{-\log2/\log5}$. The estimate obviously holds for the limit function. 

We can extend $g$ to be periodic on $\mathbb{R}$, and the theorem follows immediately.\hfill$\square$

%\section{-decay rate (i.e. convergence rate) of $c_j$'s and explosion rate of $c_j$'s (include explicit asymptotic formulas)}

%\section{-above bullet point implies that truncations are justified}

%\section{-curve plots in delta-epsilon plane}

%\section{Observations of delta-epsilon curves for the line}

%\subsection{-the delta-epsilon curves for 4pi, 8pi, 16pi, etc expansion "fill in" the stable regions of the curves. It suffic/es to use only $N=2^k$ (not other values of $N$). The reasoning is that, since the diadic rationals $\{\frac{1}{2^k}:k\in\mathbb{N}\}$ is dense in the real line, plotting all curves with period $2^k$ will constitute a dense set of curves between the $2\pi$ and $4\pi$ transition curves. Taking the closure of the union of these curves gives the stable region.}

\subsection{Infinite Sierpinski Gasket}

In the last part of this section, we introduce the \textit{infinite Sierpinski gasket} (\textit{$SG_\infty$}). It is a particular example of fractal blow-ups introduced in \cite{blowup} by R. S. Strichartz. 

Recall that the Sierpinski gasket is defined by the self-similar identity, $SG=\bigcup_{i=0}^2 F_i(SG)$, where each $F_i$ is a contraction mapping $\mathbb{R}^2\to\mathbb{R}^2$ of contraction ratio $\frac{1}{2}$ for $i=0,1,2$, as defined earlier in this section. The infinite Sierpinski gasket is constructed as follows. \\

\begin{definition}
Suppose a sequence $\mathcal{K}=\{k_n\}_{n\geq 1}$, $k_n\in\{0,1,2\}$, is fixed. Define $SG_M=F^{-1}_{\mathcal{K},M}SG$, where $F_{\mathcal{K},M}=F_{k_1}F_{k_2}\cdots F_{k_M}$. Then the infinite Sierpinski Gasket $SG_\infty$ is defined by $SG_\infty=\cup_{M\geq 1}SG_M$.
\end{definition} 

The Laplacian $\Delta_\infty$ on $SG_\infty$ can be defined locally with graph approximation in a same way as on $SG$. In \cite{teplyaev}, a Sierpinski lattice was introduced to describe the infinite graphs that approximate $SG_\infty$. Define 
\[V_{(m)}:=\bigcup_{M=1}^\infty F^{-1}_{\mathcal{K},M}V_{M+m},\]
and say $x\sim_{(m)}y$ if $F_{\mathcal{K},M}(x)\sim_{m+M}F_{\mathcal{K},M}(y)$ for some $M$. Then the resulting infinite graph is called a Sierpinski lattice. We can still define the discrete Laplacian on the lattices by 
$$
\Delta_{(m)}u(x):=\sum_{y\sim_{(m)}x}{\big(u(y)-u(x)\big),\quad\quad x\in V_{(m)}}.
$$
Then the \textit{continuous Laplacian} $\Delta$ is defined by
$$
\Delta _\infty u(x)=\frac32\lim_{m\to\infty}5^m\Delta_{(m)}u(x),\quad\quad\forall x\in V_*\setminus V_0.$$

One of the most important results on $SG_\infty$ was A. Teplyaev's theorem (see \cite{teplyaev}) below showing that the Laplacian $\Delta_\infty$ on $SG_\infty$ has pure point spectrum, which means the eigenfunctions of the Laplacian form a complete set.

\begin{theorem}[A. Teplyaev] The Laplacian $\Delta_\infty$ is self-adjoint in $L^2(SG_\infty,\mu)$, where $\mu$ is the Hausdoff measure on $SG_\infty$. The spectrum of $\Delta_\infty$ is pure point (i.e., the eigenfunctions of $\Delta_\infty$ form a basis of  $L^2(SG_\infty,\mu)$) and each eigenvalue has infinite multiplicity. The set of eigenfunctions with compact support is complete in $L^2(SG_\infty,\mu)$.
\end{theorem}

As a result of the theorem, spectral decimation still works on $SG_\infty$. Each of the eigenfunctions of $\Delta_\infty$ is an extension of an eigenfunction of $\Delta_{(m_0)}$ with eigenvalue $5$ or $6$ by spectral decimation. The only difference here is that the generation of birth $m_0$ takes values in $\mathbb{Z}$ instead of $\mathbb{N}$. All the results concerning eigenvalues from a same generation of birth and initial function in the previous section, including Proposition 5.3, still hold on $SG_\infty$.

\begin{comment}
\textcolor{blue}{[Somewhere in Section 5 we need to define what a `fracafold' is]}
\end{comment}

\section{Extending the Mathieu Differential Equation to Infinite Fractafolds}

Now we are ready to discuss how we will define the Mathieu differential equation on an infinite fractafold.

\subsection{Defining the Fractal MDE}

Recall that the MDE, defined on the real line, is given by
$$\frac{d^2u}{dt^2}+\left(\delta+\varepsilon\cos t\right)u=0,$$ where $u$ is a function from $\mathbb{R}$ to $\mathbb{R}$.

The first questions we wish to address are ``What should the fractal space be that replaces the line?'' and ``what should `periodic function' mean?''. We choose to consider the infinite Sierpinski gasket $SG_\infty$ to be our domain. By ``periodic function,'' we mean a function on $SG_\infty$ which is identical on all the copies of $SG$ of the same size. In particular, if we are given a function $u$ on $SG$ with the boundary conditions

\begin{equation}
\label{cases}
\begin{cases}
    \partial_nu(q_l)=0\quad\text{for }l=0,1,2&\\
    u(q_0)=u(q_1)=u(q_2),&
\end{cases}
\end{equation}
we can get a periodic function on $SG_\infty$ by translating the function to other copies.

But what about the differential equation? The first step in finding a fractal analog of the MDE defined on the line is to replace $\frac{d^2}{dt^2}$ with the fractal Laplacian $\Delta$, since $\Delta$ is the analog of the second-derivative operator.

Now, what to do with the $\varepsilon\left(\cos x\right) u$ term? Recall from Theorem \ref{notindomain} above that the multiplication of two nonconstant functions in dom$\Delta$ may result in a function which is not in dom$\Delta$. Thus, we cannot simply replace $\cos x$ by a function in dom$\Delta$, and so we must figure out a suitable analog of multiplication by cosine.

Recall that, in the line case, we sought solutions $u$ in the form of a Fourier expansion in terms of cosines and sines. Note, however, that \textit{cosines and sines on the line are Neumann and Dirichlet eigenfunctions, respectively, of $\frac{d^2}{dt^2}$.} Hence, we will adopt a form of Mathieu's equation which is compatible with functions $u$ that can be written as a linear combination of Neumann eigenfunctions, motivated by Equation (\ref{cases}). 

We choose to only consider functions which have Neumann eigenfunction expansions of the form
\begin{equation}\label{sumofphi}
u(x)=\sum_{j=0}^\infty c_j\varphi_j(x)
\end{equation}
where each $\varphi_j$ is a Neumann eigenfunction function defined as follows. Fix a generation of birth, a series (5-series or 6-series), and an initial eigenfunction $\varphi$ such that $\Delta_{m_0}\varphi(x)=\lambda_{m_0}\varphi(x)$ for all $x\in V_{m_0}$ and $\varphi(q_0)=\varphi(q_1)=\varphi(q_2)$. With the spectral decimation algorithm introduced in Section 5.3, we obtain a set of Neumann eigenfunctions $\varphi_i$ of $\Delta$ extended from $\varphi$, and we write $\lambda_i$ for the eigenvalue corresponding to $\varphi_i$. We still take the order $\lambda_1<\lambda_2<\lambda_3<\cdots$ as in Section 5.3. Readers can also find more details on Neumann eigenfunctions on $SG$ in \cite{differentialequationsonfractals}. The reason we fix a common initial Neumann eigenfunction is that, if we do not fix such an initial eigenfunction and instead consider the set of all Neumann eigenfunctions of $\Delta$, then we cannot order their eigenvalues in a discrete way as above.

To this end, suppose we have a function $u$ on $SG$ which can be written as a linear combination of Neumann eigenfunctions as in Equation \ref{sumofphi}, where the $\varphi_j$ ($j\geq0$) satisfy the definition in the previous paragraph. Then, for any $\varphi_j$ with $j\geq2$ we define multiplication by cosine as follows:
$$\left(cos t\right)\varphi_j:=\frac12\varphi_{j-1}+\frac12\varphi_{j+1}.$$
The motivation for this definition comes from the fact that, for the line case, the Neumann eigenfunction $\cos (jt)$ obeys the following trigonometric property when multiplied by $\cos t$:
$$\cos t\cos \left(jt\right)=\frac12\cos(j-1)t+\frac12\cos(j+1)t.$$

As for $j\geq 1$ we consider two possibilities, each of which will be described in Section 6.2. In addition, we will consider another two variant versions in the following subsection. 

\begin{comment}

In this case, we have

\begin{dmath*}
\left(\delta+\varepsilon\cos x\right)\left(\sum_{j=0}^\infty c_j\cos\left(jx\right)\right)=\sum_{j=0}^\infty\delta c_j\cos(jx)+\varepsilon c_0\cos x+\frac12\sum_{j=1}^\infty\varepsilon c_j\cos(j+1)x+\frac12\sum_{j=1}^\infty\varepsilon c_j\cos(j-1)x=\left(\delta c_0+\frac12\varepsilon c_1\right)+\sum_{j=1}^\infty\left(\delta c_j+\frac12\varepsilon c_{j-1}+\frac12\varepsilon c_{j+1}\right)\cos(jx)
\end{dmath*}

So, if we want the fractal analog, we should start with a list of a basis of eigenfunctions $\varphi_j$ satisfying the two conditions in Equation \ref{eqn61}. \end{comment}

Before moving onto Section 6.2, we talk about one particular case, in which

    \begin{equation}
    \begin{cases}
    (\cos t)\varphi_0=\varphi_1\\
    (\cos t)\varphi_1=\frac{1}{2}\varphi_0+\frac{1}{2}\varphi_2.\\
    \end{cases}
    \end{equation}

Here, $\varphi_0$ and $\varphi_1$ play an analogous role to $\cos(0t)$ and $\cos(1t)$, respectively, in the line case, since

    \begin{equation}
    \begin{cases}
    (\cos t)\cos(0j)=(\cos t)\cdot1=\cos t=\cos(1t)\\
    (\cos t)\cos (1t)=(\cos t)(\cos t)=\cos^2 t=\frac12+\frac12\cos(2t)=\frac12\cos(0t)+\frac12\cos(2t).\\
    \end{cases}
    \end{equation}

\begin{comment}
\st{where we consider $\varphi_j$ as analogue of $\cos(jt)$. In this case, }

$$(\cos t)\varphi_0=\varphi_1, (\cos t)\varphi_1=\frac{1}{2}\varphi_0+\frac{1}{2}\varphi_2.$$
\end{comment}

Using our developments thus far, the fractal Mathieu differential equation would say
\begin{dmath*}
(\delta c_0+\frac{1}{2}\ep c_1)+(-\lambda_1c_1+\ep c_0+\delta c_1+\frac{1}{2}\ep c_2)\varphi_1+\sum_{j\geq2}\left(-\lambda_jc_j+\frac12\varepsilon c_{j-1}+\delta c_j+\frac12\varepsilon c_{j+1}\right)\varphi_j=0,
\end{dmath*}
where $\lambda_j$ is the eigenvalue of $\Delta$ corresponding to eigenfunction $\varphi_j$. Since the set $\{\varphi_j\}$ of eigenfunctions is linearly independent, we must have 
    \begin{equation}
    \begin{cases}
    \delta c_0+\frac{1}{2}\ep c_1=0\\
    -\lambda_1c_1+\ep c_0+\delta c_1+\frac{1}{2}\ep c_2=0\\
    -\lambda_jc_j+\frac12\varepsilon c_{j-1}+\delta c_j+\frac12\varepsilon c_{j+1},\text{ }j\geq2.\\
    \end{cases}
    \end{equation}

\begin{comment}
$$-\lambda_jc_j+\frac12\varepsilon c_{j-1}+\delta c_j+\frac12\varepsilon c_{j+1}=0$$
for all $j\geq1$.
\end{comment}

Putting these equations into matrix form we obtain

\begin{equation}
    \begin{pmatrix}
        \delta & \frac{\varepsilon}{2} &\\ \varepsilon & \delta-\lambda_1 & \frac{\varepsilon}{2} &\\
        & \frac{\varepsilon}{2} & \delta-\lambda_2 & \frac{\varepsilon}{2}\\
        & & \frac{\varepsilon}{2} & \delta-\lambda_3 & \frac{\varepsilon}{2}\\
        & & & \ddots & \ddots & \ddots
    \end{pmatrix}
    \begin{pmatrix}
    c_0 \\ c_1 \\ c_2 \\ c_3 \\ \vdots
    \end{pmatrix}
    =\begin{pmatrix}
    0 \\ 0 \\ 0 \\ 0 \\ \vdots
    \end{pmatrix}.
    \\
    \end{equation}

Note that this matrix takes the form of Equation \ref{eqn3}, which reproduce below:

\[\begin{pmatrix}
\delta-\lambda_1-\gamma\varepsilon & \alpha_1\varepsilon\\
\beta_1\varepsilon & \delta-\lambda_2 & \alpha_2\varepsilon\\
& \beta_2\varepsilon & \delta-\lambda_3 & \alpha_3\varepsilon\\
& & \beta_3\varepsilon & \delta-\lambda_4 & \ddots\\
& & &\ddots & \ddots
\end{pmatrix}.\]

%Notice that we fix a generation of birth $m_0$, fix a series (5-series or 6-series), and consider the set of Neumann eigenvalues that arise from these conditions. They are discrete and can be ordered $\lambda_1<\lambda_2<\lambda_3<\cdots$. Set $\lambda_0=0$. These are the eigenvalues we use in the matrix above. 

\subsection{Variants of the Mathieu Differential Equation on the Line}

We will consider 4 different `versions', $M_1$, $M_2$, $M_3$, and $M_4$, of the coefficient matrix  for the fractal MDE $M_ix=0$ ($i=1,2,3,4$):
\begin{itemize}

\item Version 1:
\[M_1:=\begin{pmatrix}
\delta-\lambda_0 & \frac12\varepsilon\\
\varepsilon & \delta-\lambda_1 & \frac12\varepsilon\\
& \frac12\varepsilon & \delta-\lambda_2 & \frac12\varepsilon\\
& & \frac12\varepsilon & \delta-\lambda_3 & \ddots\\
& & &\ddots & \ddots
\end{pmatrix}.\]

This matrix is reminiscent of the cosine matrices for the line case. Note that the first term in the second row is $\varepsilon$, not $\frac12\varepsilon$.

The recursion relation for the coefficients $c_j$ becomes
    \begin{equation}
    \begin{cases}
    (\delta-\lambda_0)c_0+\frac12\varepsilon c_1=0,\\
    \varepsilon c_0+(\delta-\lambda_1)c_1+\frac12\varepsilon c_2=0\\
    \frac12\varepsilon c_{j-1}+(\delta-\lambda_j)c_j+\frac12\varepsilon c_{j+1}=0,\quad(j\geq2).\\
    \end{cases}
    \end{equation}

\item Version 2:
\[M_2:=\begin{pmatrix}
\delta-\lambda_1 & \frac12\varepsilon\\
\frac12\varepsilon & \delta-\lambda_2 & \frac12\varepsilon\\
& \frac12\varepsilon & \delta-\lambda_3 & \frac12\varepsilon\\
& & \frac12\varepsilon & \delta-\lambda_4 & \ddots\\
& & &\ddots & \ddots
\end{pmatrix}.\]

This matrix is reminiscent of the sine matrices for the line case. Note that the eigenvalues start from $\lambda_1$ instead of $\lambda_0$.

The recursion relation for the coefficients $c_j$ becomes
    \begin{equation}
    \begin{cases}
    (\delta-\lambda_1)c_1+\frac12\varepsilon c_2=0,\\
    \frac12\varepsilon c_{j-1}+(\delta-\lambda_j)c_j+\frac12\varepsilon c_{j+1}=0\quad(j\geq2).\\
    \end{cases}
    \end{equation}

\item Versions 3 and 4:

For Version 3 and for Version 4 we take the following approach.

We now consider a variant of the Mathieu differential equation, given by

$$\Delta u+\left(\delta+\varepsilon A\right)u=0,$$

where $A$ is the operator analogous to multiplication by cosine. For $j\geq 2$, define
$A\cdot\varphi_j:=\alpha_{j-1}\varphi_{j-1}+\beta_j\varphi_{j+1},$ with $\alpha_j$ and $\beta_j$ satisfying
  \begin{equation}
    \label{systemforalphabeta}
    \begin{cases}
    \alpha_{j-1}+\beta_j=1,\\
    \alpha_{j-1}\sqrt{\lambda_{j-1}}+\beta_{j}\sqrt{\lambda_{j+1}}=\sqrt{\lambda_j},
    \end{cases}
    \end{equation}
for $j\geq 1$. One can solve system \ref{systemforalphabeta} to find the solutions of $\alpha_j$ and $\beta_j$ 
    \begin{equation}
    \label{solutionsforalphabeta}
    \begin{cases}
    \alpha_{j}=\frac{\sqrt{\lambda_{j+2}}-\sqrt{\lambda_{j+1}}}{\sqrt{\lambda_{j+2}}-\sqrt{\lambda_{j}}},\text{ }(j\geq 0),\\
    \\
    \beta_j=\frac{\sqrt{\lambda_{j}}-\sqrt{\lambda_{j-1}}}{\sqrt{\lambda_{j+1}}-\sqrt{\lambda_{j-1}}},\text{ }(j\geq 1).
    \end{cases}
    \end{equation}
 Note that $0<\alpha_j<1$ and $0<\beta_j<1$ for all $j$, and hence the sequences $\{\alpha_j\}$ and $\{\beta_j\}$ are each uniformly bounded.

The motivation for the setup above is that, if we plug in $\lambda_j=j^2$, which corresponds to the eigenvalues on the line, Equation \ref{solutionsforalphabeta} yields $\alpha_j=\beta_j=\frac{1}{2}$. For $j=0$ and $j=1$, we still consider two cases, which give us Version 3 and Version 4 as follows.

We let Version 3 be

\[M_3:=\begin{pmatrix}
\delta-\lambda_0 & \alpha_0\varepsilon\\
\varepsilon & \delta-\lambda_1 & \alpha_1\varepsilon\\
& \beta_1\varepsilon & \delta-\lambda_2 & \alpha_2\varepsilon\\
& & \beta_2\varepsilon & \delta-\lambda_3 & \ddots\\
& & &\ddots & \ddots
\end{pmatrix}\]

where $\alpha_j$ and $\beta_j$ are as given in Equation \ref{solutionsforalphabeta}.

This matrix is reminiscent of the cosine matrices for the line case. Note that the first term in the second row has an extra factor of $2$.

The recursion relation for the coefficients $c_j$ becomes
    \begin{equation}
    \begin{cases}
    (\delta-\lambda_0)c_0+\alpha_0\varepsilon c_1=0,\\
    \varepsilon c_0+(\delta-\lambda_1)c_1+\alpha_1\varepsilon c_2=0\\
    \beta_{j-1}\varepsilon c_{j-1}+(\delta-\lambda_j)c_j+\alpha_j\varepsilon c_{j+1}=0,\quad(j\geq2).\\
    \end{cases}
    \end{equation}

We let Version 4 be

\[M_4:=\begin{pmatrix}
\delta-\lambda_1 & \alpha_1\varepsilon\\
\beta_1\varepsilon & \delta-\lambda_2 & \alpha_2\varepsilon\\
& \beta_2\varepsilon & \delta-\lambda_3 & \alpha_3\varepsilon\\
& & \beta_3\varepsilon & \delta-\lambda_4 & \ddots\\
& & &\ddots & \ddots
\end{pmatrix}\]

where $\alpha_j$ and $\beta_j$ are as given in Equation \ref{solutionsforalphabeta}.

This matrix is reminiscent of the sine matrices for the line case. Again, note that the eigenvalues start from $\lambda_1$ instead of $\lambda_0$.

The recursion relation for the coefficients $c_j$ becomes 
    \begin{equation}
    \begin{cases}
    (\delta-\lambda_1)c_1+\alpha_1\varepsilon c_2=0,\\
    \beta_{j-1}\varepsilon c_{j-1}+(\delta-\lambda_j)c_j+\alpha_j\varepsilon c_{j+1}=0,\quad(j\geq2).\\
    \end{cases}
    \end{equation}

\end{itemize}

\section{Ovservations and Analysis for the Fractal Mathieu Differential Equation}

In this section we parallel our results presented in Section 4 by giving a discussion of the asymptotic behavior of the transition curves (defined below) for the fractal MDE and of the convergence of solutions. We also describe a phenomenon wherein the $SG_\infty$ transition curves form a prominent `diamond' pattern.

\subsection{The $\delta$-$\varepsilon$ Plot}

%Now that some theoretical background has been given in the previous section, we will now discuss some results.

%As discussed previously, the transition curves are determined through the following process. First, one expands the function $u$ in terms of a ($2\pi$-periodic or $4\pi$-periodic) Fourier expansion. Then, one plugs this into the Mathieu differential equation and finds two infinite systems of linear homogeneous equations that must be satisfied by the cosine and sine coefficients of the Fourier expansion. One then puts these equations into matrix form, takes the determinant of a truncated matrix, and sets the resulting expression equal to 0. This yields the implicit equation involving $\delta$ and $\varepsilon$ which produces the transition curves in the $\delta$-$\varepsilon$ plane and thus determines the regions of stability and the regions of instability. A figure is shown below.

\begin{figure}[!h]
\centering
\includegraphics[height=3cm]{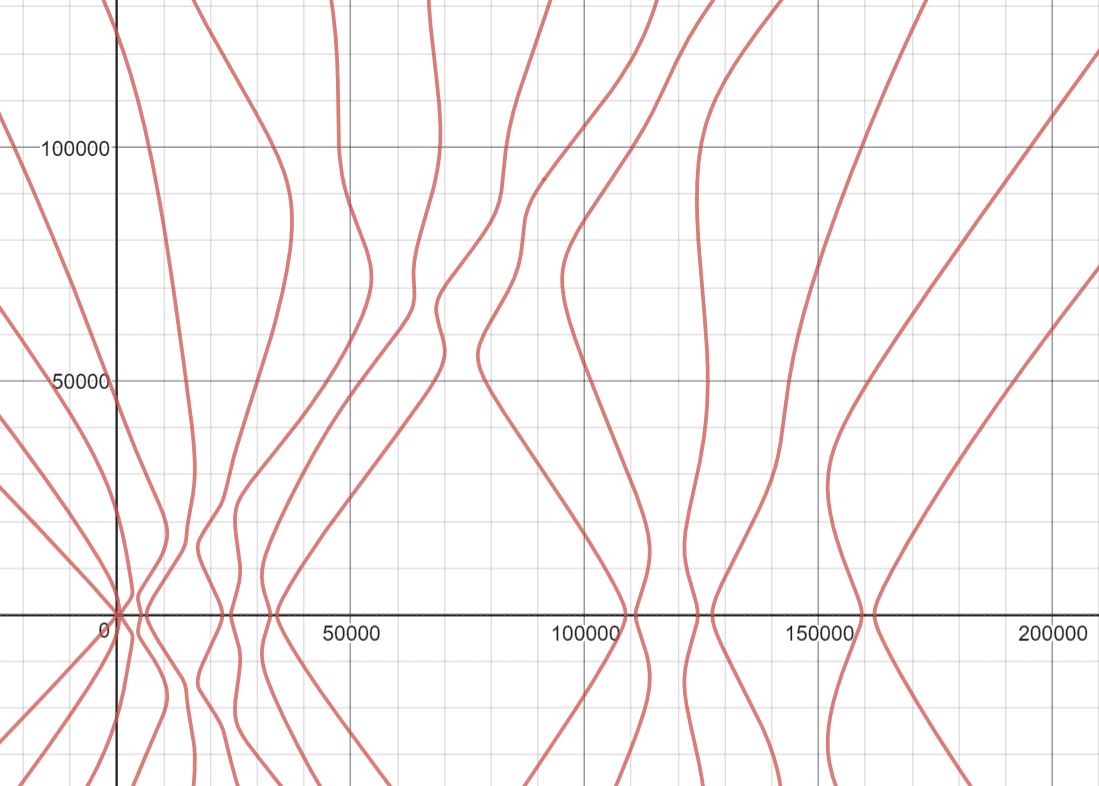}\quad
\includegraphics[height=3cm]{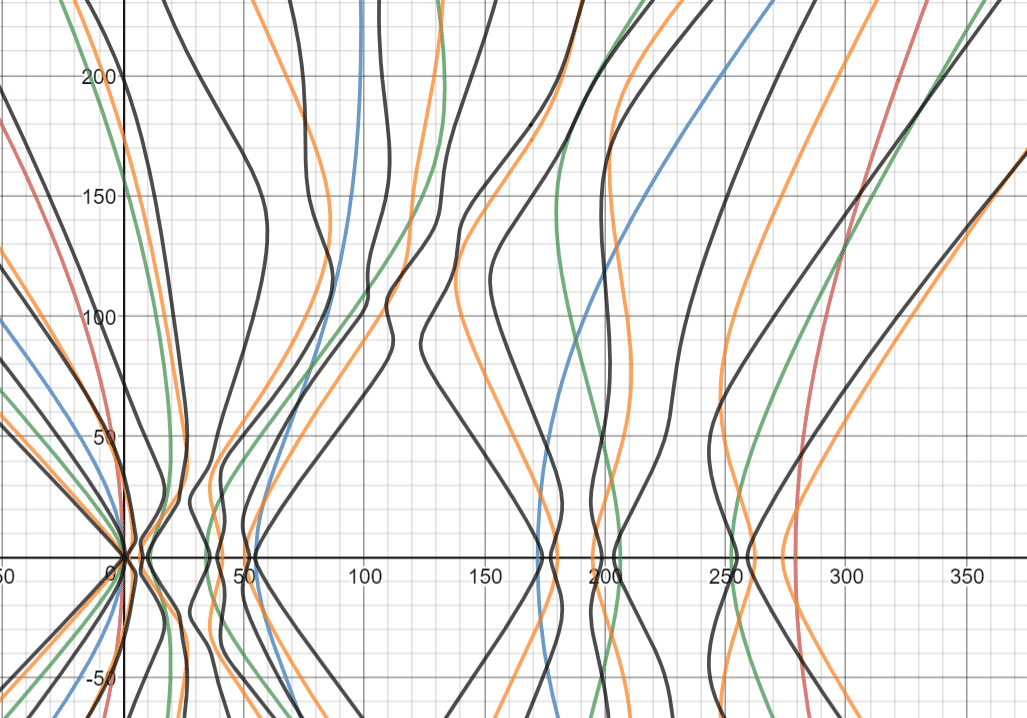}
\caption{Transition curves for version 1, 5 series. The left picture shows transition curves of a single matrix corresponding to generation of birth $0$. The right picture shows transition curves of generation of birth 0 (red), -1 (blue), -2 (green), -3 (orange), -4 (black).}
\label{coscurves-half-5series}
\end{figure}

\begin{figure}[!h]
\centering
\includegraphics[height=3cm]{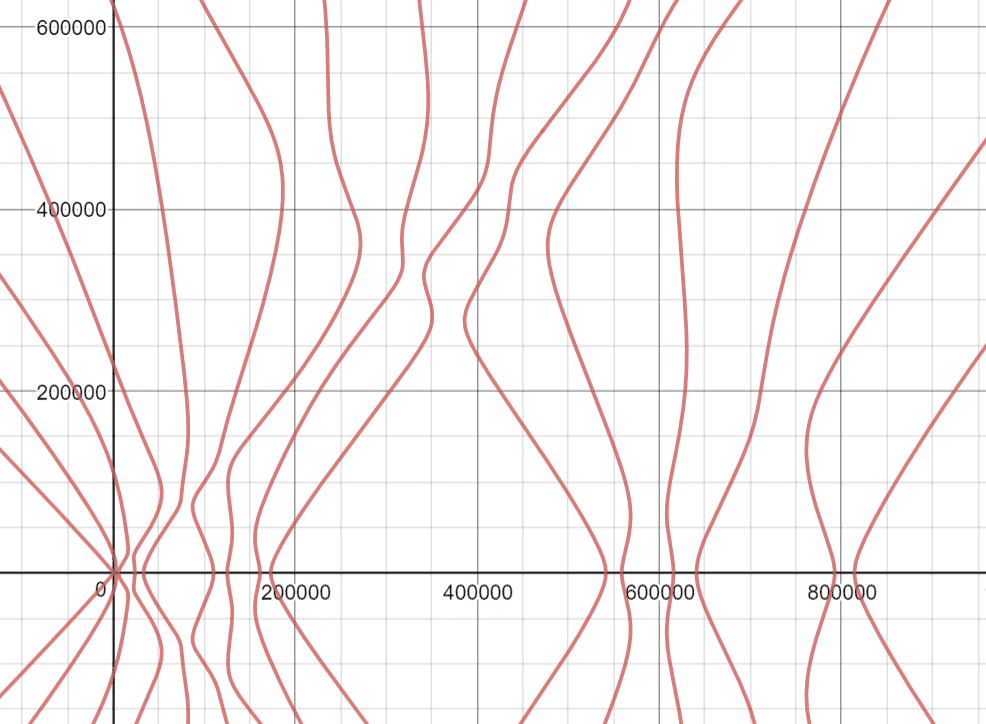}\quad
\includegraphics[height=3cm]{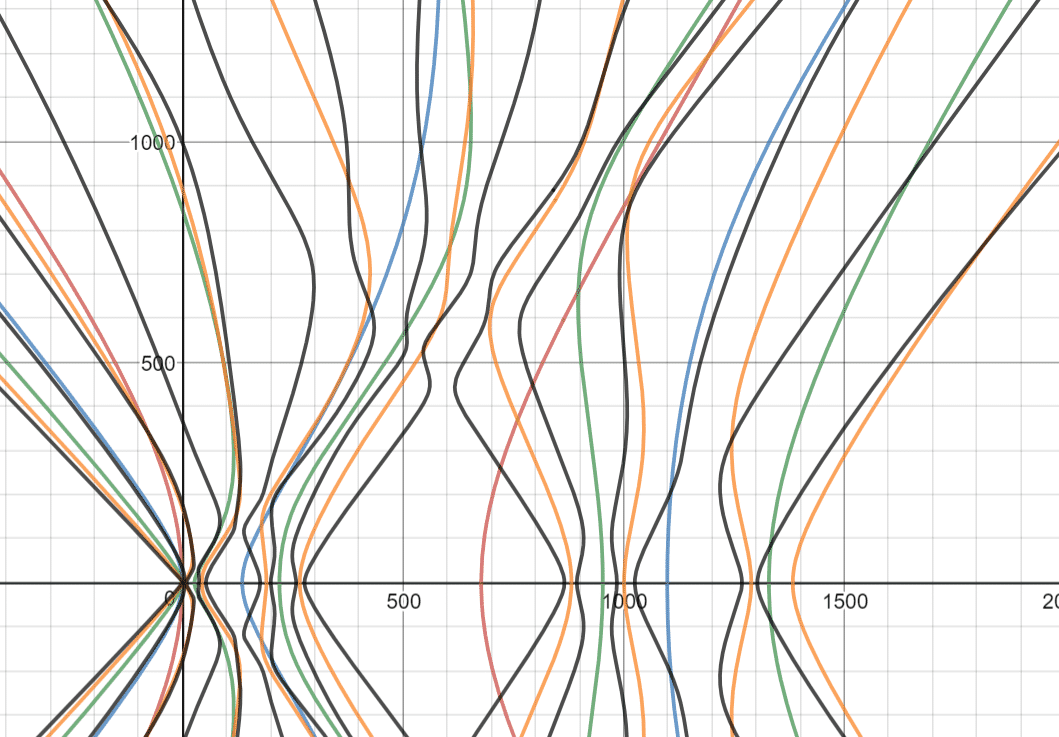}
\caption{Transition curves for version 1, 6 series. The left picture shows transition curves of a single matrix corresponding to generation of birth $0$. The right picture shows transition curves of generation of birth 0 (red), -1 (blue), -2 (green), -3 (orange), -4 (black).}
\label{coscurves_half_6series}
\end{figure}

\begin{figure}[!h]
\centering
\includegraphics[height=3cm]{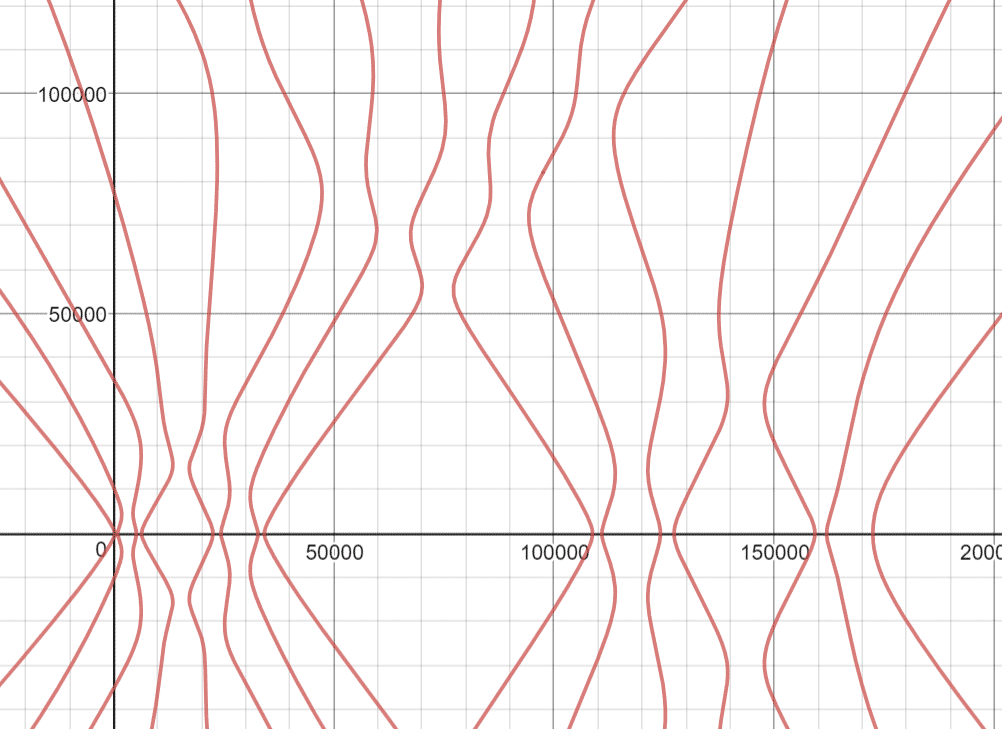}\quad
\includegraphics[height=3cm]{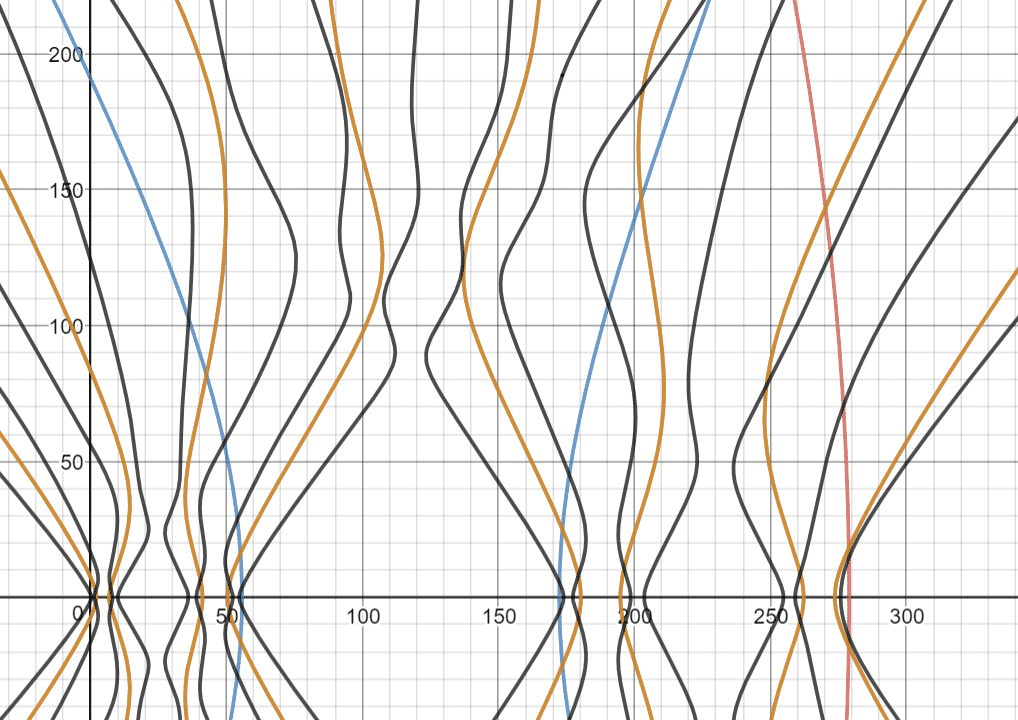}
\caption{Transition curves for version 2, 5 series. The left picture shows transition curves of single matrix corresponding to generation of birth $0$. The right picture shows transition curves of generation of birth 0 (red), -1 (blue), -2 (green), -3 (orange), -4 (black).}
\label{sincurves-half-5series}
\end{figure}

\begin{figure}[!h]
\centering
\includegraphics[height=3cm]{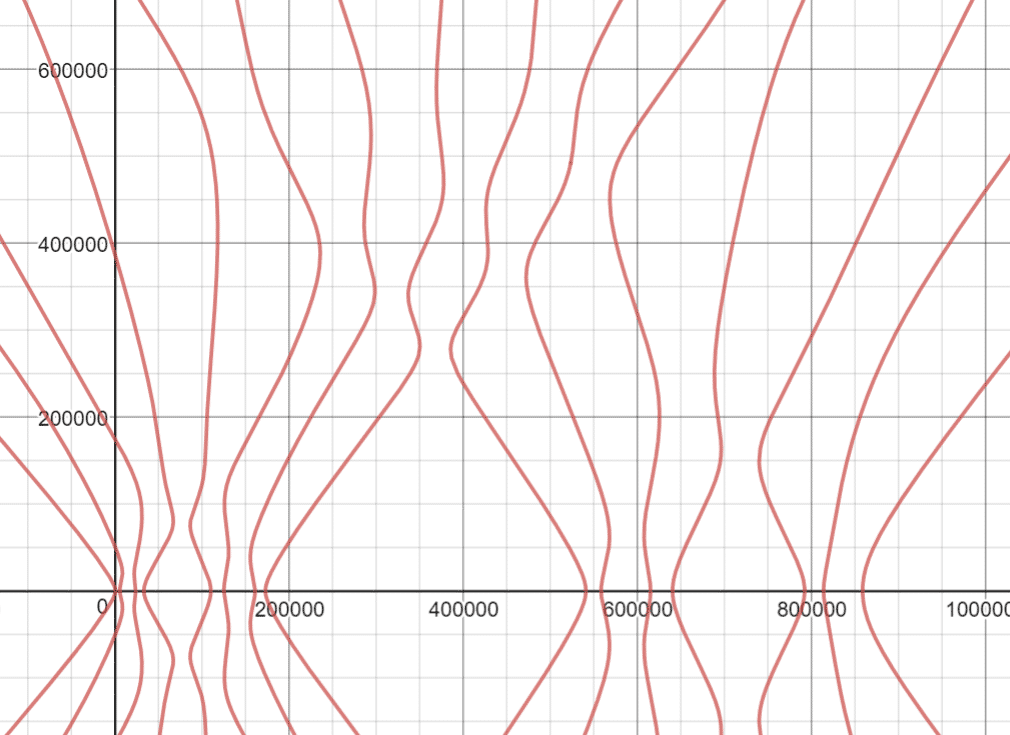}\quad
\includegraphics[height=3cm]{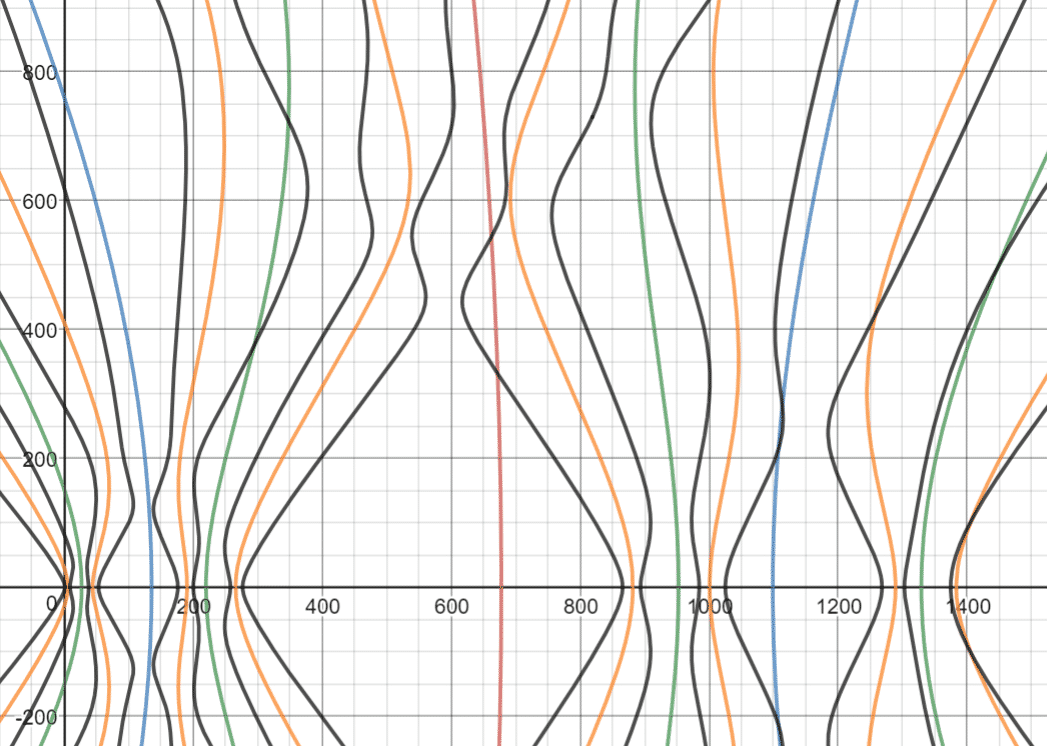}
\caption{Transition curves for 
version 2, 6 series. The left picture shows transition curves of single matrix corresponding to generation of birth $0$. The right picture shows transition curves of generation of birth 0 (red), -1 (blue), -2 (green), -3 (orange), -4 (black).}
\label{sincurves_half_6series}
\end{figure}

\begin{figure}[!h]
\centering
\includegraphics[height=2.8cm]{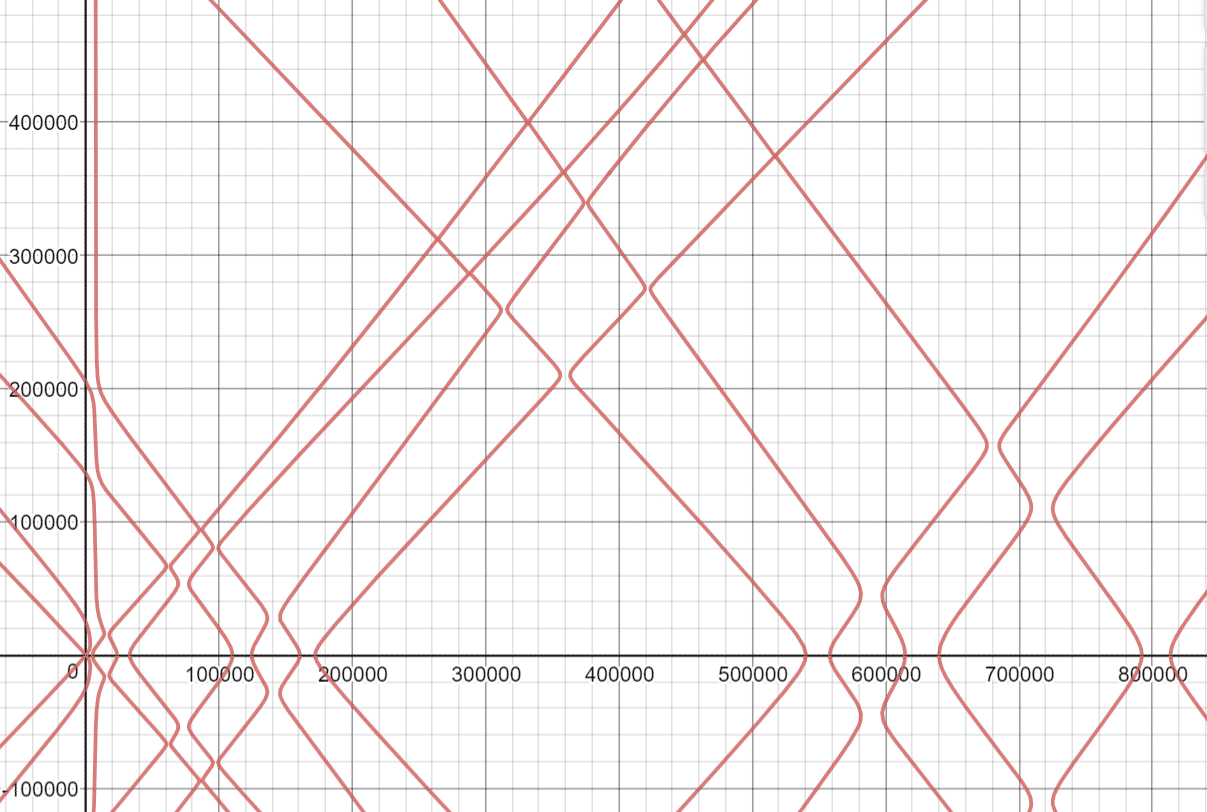}\quad
\includegraphics[height=2.8cm]{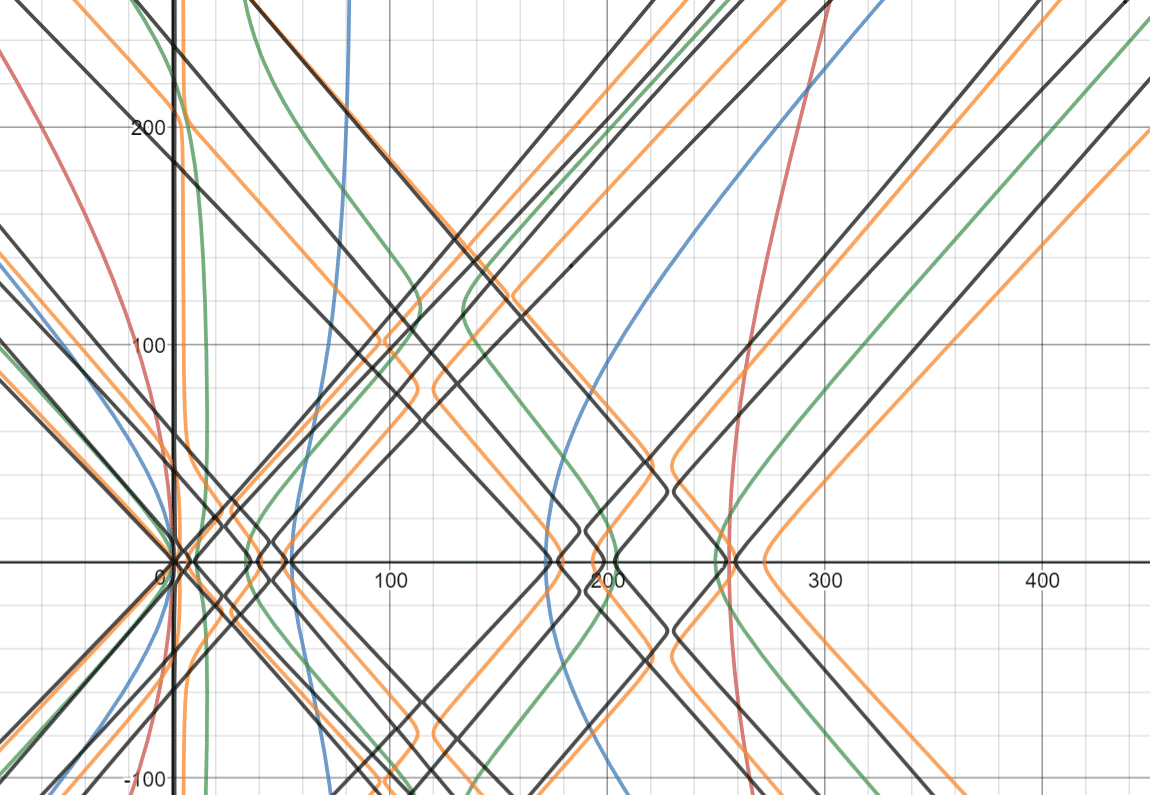}
\caption{Transition curves for version 3, 5 series. The left picture shows transition curves of single matrix corresponding to generation of birth $0$. The right picture shows transition curves of generation of birth 0 (red), -1 (blue), -2 (green), -3 (orange), -4 (black).}
\label{coscurves-squareroot-5series}
\end{figure}

\begin{figure}[!h]
\centering
\includegraphics[height=2.8cm]{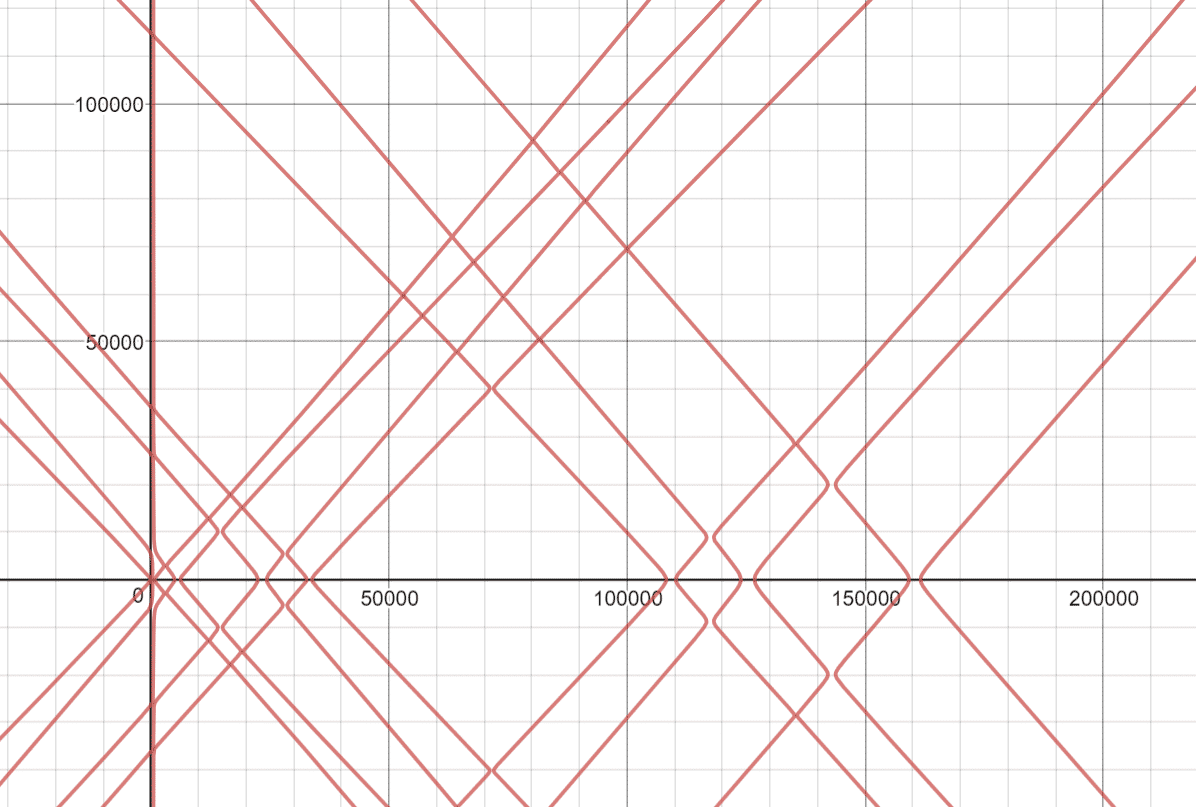}\quad
\includegraphics[height=2.8cm]{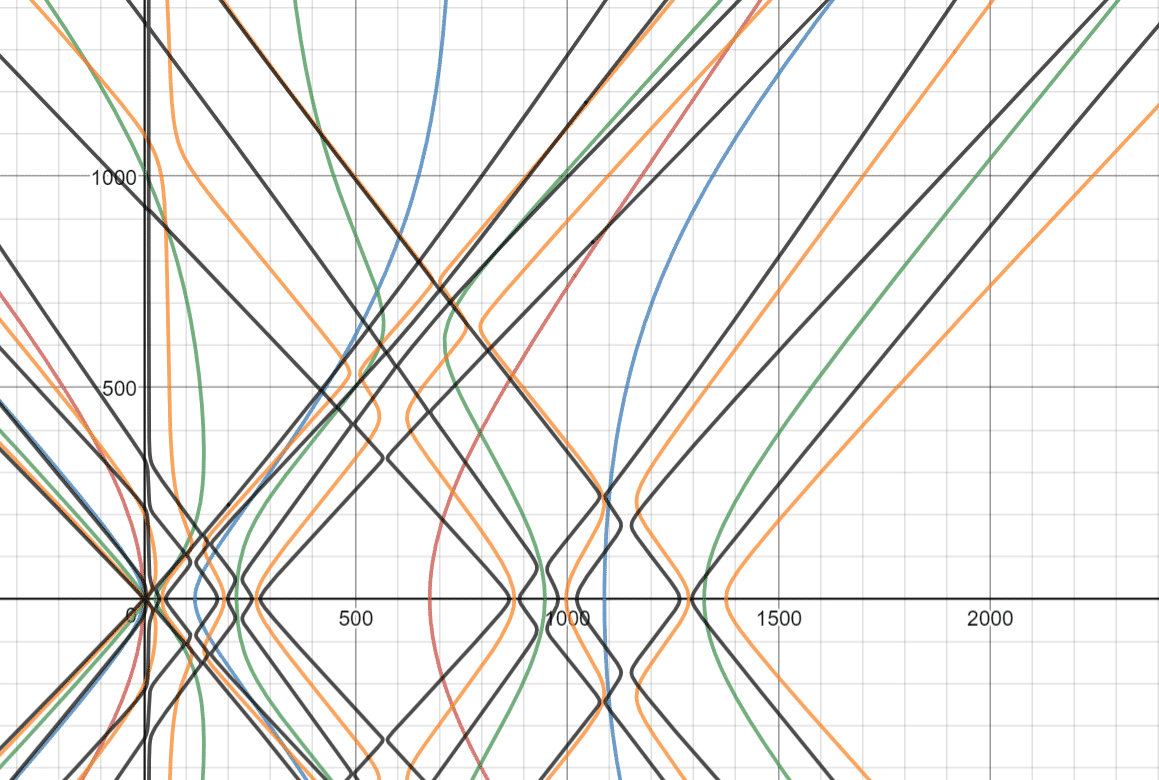}
\caption{Transition curves for version 3, 6 series. The left picture shows transition curves of single matrix corresponding to generation of birth $0$. The right picture shows transition curves of generation of birth 0 (red), -1 (blue), -2 (green), -3 (orange), -4 (black).}
\label{coscurves-squareroot-6series}
\end{figure}

\begin{figure}[!h]
\centering
\includegraphics[height=2.8cm]{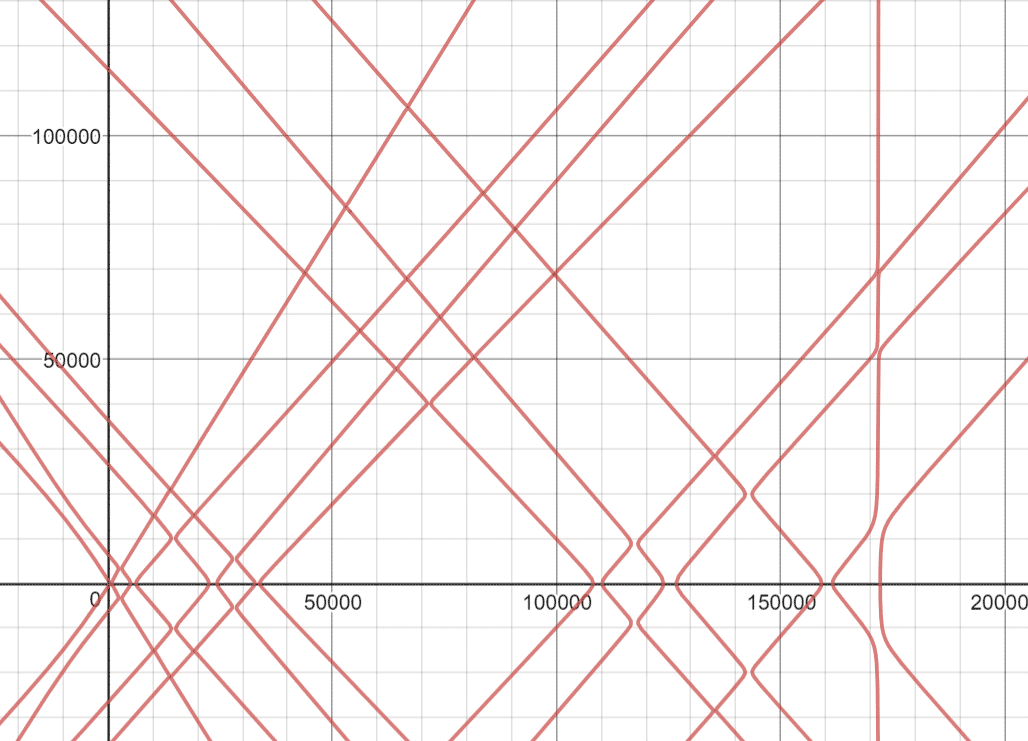}\quad
\includegraphics[height=2.8cm]{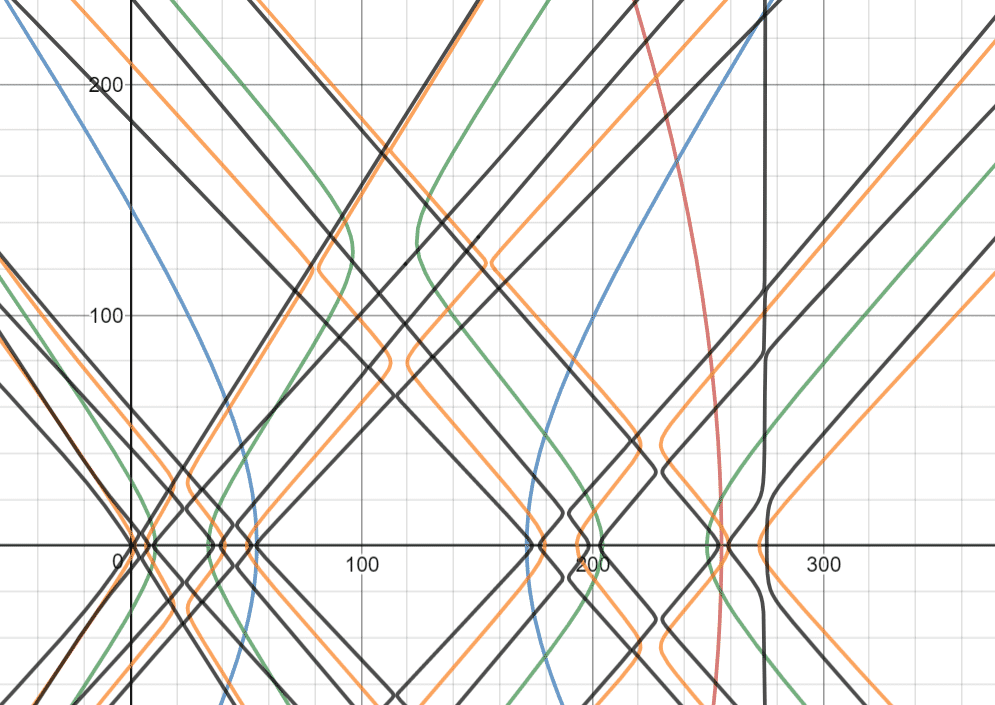}
\caption{Transition curves for version 4, 5 series. The left picture shows transition curves of single matrix corresponding to generation of birth $0$. The right picture shows transition curves of generation of birth 0 (red), -1 (blue), -2 (green), -3 (orange), -4 (black).}
\label{sincurves-squareroot-5series}
\end{figure}

\begin{figure}[!h]
\centering
\includegraphics[height=2.8cm]{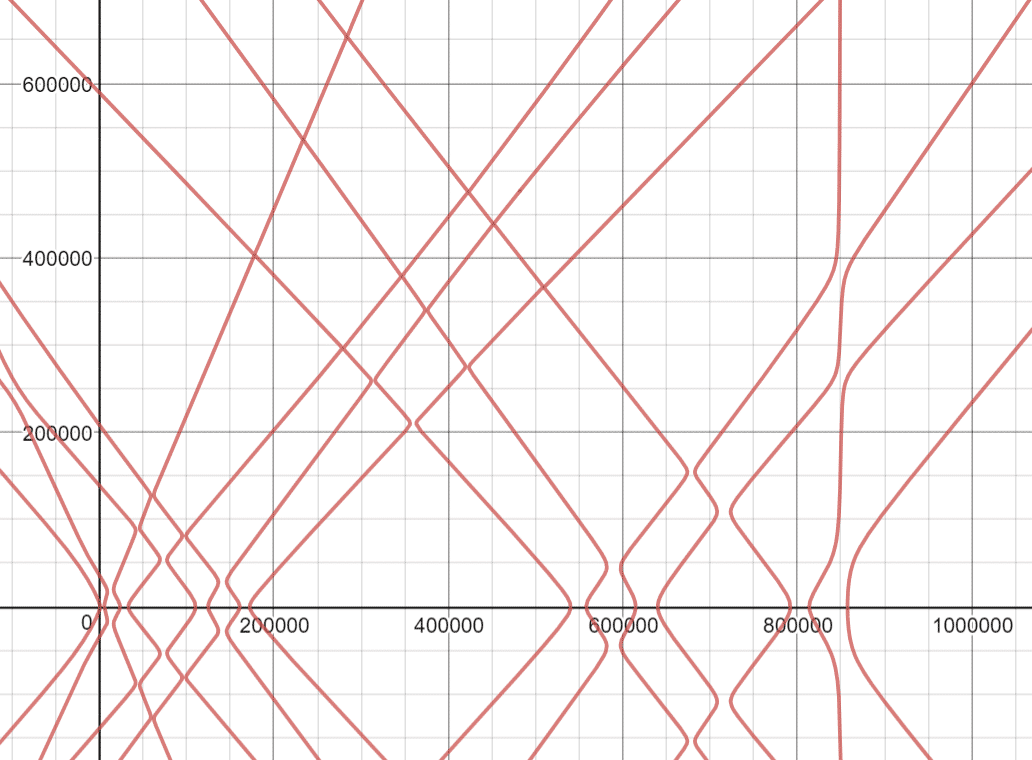}\quad
\includegraphics[height=2.8cm]{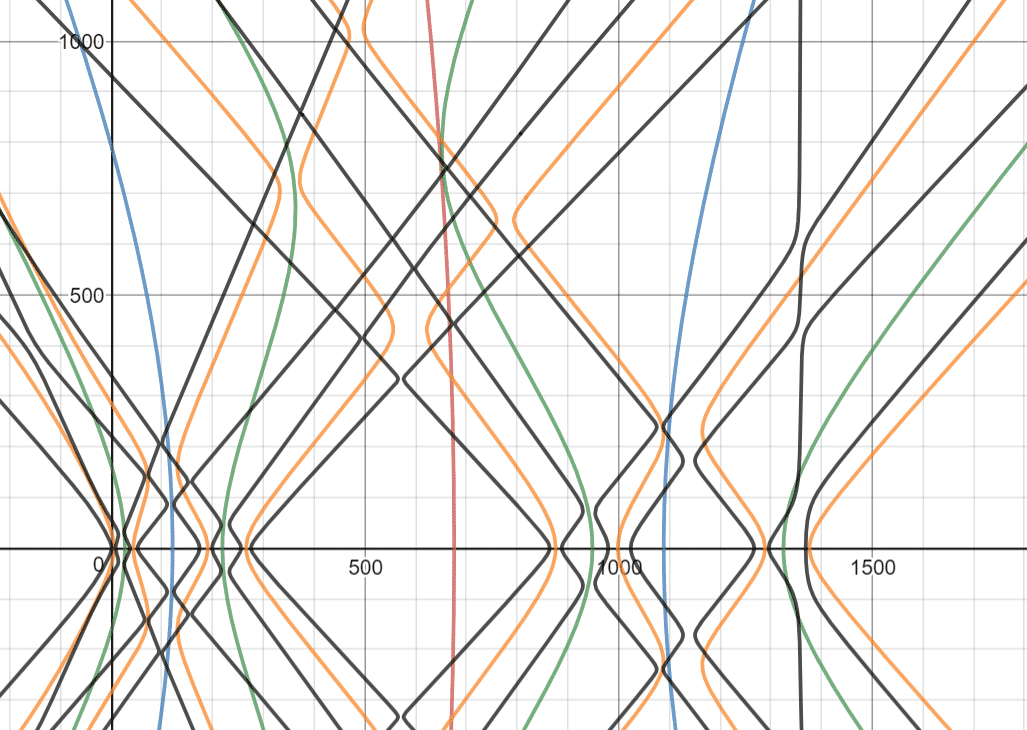}
\caption{Transition curves for version 4, 6 series. The left picture shows transition curves of single matrix corresponding to generation of birth $0$. The right picture shows transition curves of generation of birth 0 (red), -1 (blue), -2 (green), -3 (orange), -4 (black).}
\label{sincurves-squareroot-6series}
\end{figure}

In Section 6, we introduced 4 different Versions $M_ix=0$, ($i=1,2,3,4$) of the fractal MDE in matrix form on $SG_\infty$. If one fixes one of the four Versions of the MDE as above, then the points $(\delta,\ep)$ in the $\delta$-$\ep$ plane for which the $M_ix=0$ has a nontrivial solution in $\ell^2$ make up the \textit{transition curves} for that Version of the MDE. In Figures \ref{coscurves-half-5series}-\ref{sincurves-squareroot-6series}, we show the transition curves for each of the four Versions.

%Version 1: 

%Version 2:

%It can be shown that the transition curves for the Mathieu differential equation (resulting from $2\pi$- and $4\pi$-periodic expansion) share a number of interesting properties. We state below, without proof, some of these properties (see [3] for a further discussion):

%Hence, from this figure it is easy to see that, once one determines where the transition curves lie on the $\delta\varepsilon$-plane, one can determine the stability for all other $(\delta,\varepsilon)$ in the plane since they are, roughly speaking, separated into regions `between' the transition curves.

%Version 3:

%Version 4:
\clearpage{}

%\clearpage{}

As can be viewed from Figures \ref{coscurves-half-5series}-\ref{sincurves-squareroot-6series}, in each plot there is a `diamond' pattern formed between adjacent transition curves. That is, the curves seem to `bounce' off each other before parting in separate directions. This pattern does not appear in the transition curves for the Mathieu differential equation on the line. This is one way in which the transition curves for $SG_\infty$ are different than those for the line.

We give a short proof here concerning the asymptotic behavior of the transition curves for Versions 1 and 2.

\begin{proposition}\label{prop71} For version 1 and version 2, fix one transition curve. The ratio $\frac{\delta}{\varepsilon}$ converges to $-1$ if $\varepsilon\to\infty$. The ratio $\frac{\delta}{\varepsilon}$ converge to $1$ if $\varepsilon \to-\infty$.
\end{proposition}

\textit{Proof.} In fact, $M_1$ and $M_2$ are of the form $A+\varepsilon\cdot B$, where $A$ is a diagonal matrix, and $B$ has one of the following two forms,
\[B=\begin{bmatrix}
0 &\frac{1}{2}\\
\frac{1}{2}& 0 &\frac{1}{2}\\
& \frac{1}{2} & 0 &\frac{1}{2}\\
& & \ddots&\ddots
\end{bmatrix}\text{, or }
B=\begin{bmatrix}
0 &\frac{1}{2}\\
1 & 0 &\frac{1}{2}\\
& \frac{1}{2} & 0 &\frac{1}{2}\\
& & \ddots&\ddots\end{bmatrix}.\]
The first matrix is real symmetric with continuous spectrum $[-1,1]$. To see this, we only need observe that $B=U\cdot \cos t\cdot U^{-1}$, where $U:\ell^2\to L^2([0,\pi])$ is defined as $U\bm{c}=\sum_{n=1}^\infty c_n\sin nt$. The second one is not symmetric, but we can  study the matrix $\kappa^{-1}(A-\varepsilon\dot B)\kappa=A-\kappa^{-1}B\kappa$ where $\kappa$ is a diagonal matrix, we did in the proof for theorem 4. After the transformation, we have
\[\kappa^{-1}B\kappa=\begin{bmatrix}
0 &\frac{1}{\sqrt{2}}\\
\frac{1}{\sqrt{2}} & 0 &\frac{1}{2}\\
& \frac{1}{2} & 0 &\frac{1}{2}\\
& & \ddots&\ddots\end{bmatrix},\]
which is also a real symmetric matrix with continuous spectrum $[-1,1]$. Without loss of generality, we just assume $B$ is symmetric in the follows.

In the following proof, we only consider $\ep\to+\infty$ case, and the other half is then immediate by symmetry. 

On the $k$-th curve, we know that  
\[\delta=\inf_{M_k}\sup_{\bm{c}\in M_k} \frac{\bm{c}^T(A+\varepsilon B)\bm{c}}{\bm{c}^T\bm{c}},\]
where the infimum takes over all the $k$ dimensional subspace $M_k$ of dom$A$. It is easy to see the following lower bound of $\delta$
\[\delta\geq \inf_{M_k}\sup_{x\in M_k} \frac{\bm{c}^TA\bm{c}}{\bm{c}^T\bm{c}}+\inf_{M_k}\sup_{\bm{c}\in M_k} \frac{\varepsilon \bm{c}^TB\bm{c}}{\bm{c}^T\bm{c}}=\inf_{M_k}\sup_{x\in M_k} \frac{\bm{c}^TA\bm{c}}{\bm{c}^T\bm{c}}-|\varepsilon|.\]
%One can check that $N_l\subset dom A$ since $c_j\sim j^{-3}$, decreasing faster than $\frac{1}{\lambda_j}$ by Proposition 5.3.  

Next, we give an upper bound. For any $l<1$, we can find a $k$ dimensional subspace $N_k\subset dom A$, such that $-1<B|_{N_k}<l$, which means 
\[-1<\frac{\bm{c}^TB\bm{c}}{\bm{c}^T\bm{c}}<l,\forall \bm{c}\in N_k.\]
In fact, observe that $B$ has continuous spectrum $[-1,1]$, we can always find a $k$ dimensional subspace $\tilde{N}_k\subset \ell^2$ such that $-1<B|_{\tilde{N}_k}<l$. But $\tilde{N}_k$ is not necessarily in $dom A$. So we pick a large $L$, and define
\[N_k=\{\bm{c}=(c_1,c_2,\cdots,c_L,0,0,\cdots):\text{There exists }\bm{\tilde{c}}\in \tilde{N}_k \text{ such that }c_j=\tilde{c_j},\forall 1\leq j\leq L\}.\]
Obviously, $N_k\subset dom A$, and $-1<B|_{N_k}<l$ as long as $L$ is large enough.

Then we get the following upper bound 
\[\delta\leq \sup_{\bm{c}\in N_k} \frac{\bm{c}^TA\bm{c}}{\bm{c}^T\bm{c}}+\sup_{\bm{c}\in N_k}\frac{\varepsilon \bm{c}^TB\bm{c}}{\bm{c}^T\bm{c}}<\sup_{\bm{c}\in N_k} \frac{\bm{c}^TA\bm{c}}{\bm{c}^T\bm{c}}+l\varepsilon\text{, for }\varepsilon>0.\]
Combining the upper and lower bounds of $\delta$,  we get
\[-1\leq\liminf_{\varepsilon\to+\infty}\frac{\delta}{\varepsilon}\leq \limsup_{\varepsilon\to+\infty}\frac{\delta}{\varepsilon}\leq l.\]
for any $l<1$. Thus $\lim_{\varepsilon\to +\infty} \frac{\delta}{\varepsilon}=-1$. By symmetry, we also get $\lim_{\varepsilon\to -\infty} \frac{\delta}{\varepsilon}=1$.\hfill$\square$

\subsection{Solutions of the Mathieu Differential Equation}
Now we investigate solutions to the fractal MDE.

From Corollary 3.3 we deduce that, for $\ep\neq0,$ a nontrivial solution $\textbf{c}=(c_1,c_2,...)\in\ell^2$ (or $\bm{c}=(c_0,c_1,c_2,...)\in\ell^2$) to any of the four systems $M_ix=0$ $(i=1,2,3,4)$ in section 6.2 will decay with asymptotic behavior $\frac{c_j}{c_{j-1}}\sim \frac{\beta_{j-1}\varepsilon}{\lambda_j-\delta}.$ Hence, a solution $$u(x)=\sum_{j=0}^\infty c_j\varphi_j(x)$$ to any of the four versions of the fractal MDE may be well-approximated by the sum of just the first few terms in the series, and such an approximation will yield only a small error to the actual infinite-series solution.

Recall that in Section 4 we investigated, for the MDE on the line, how properties of solutions change if one fixes a transition curve and considers periodic solutions corresponding to various $(\delta,\ep)$ points along that transition curve. We will consider the same question for the $SG_\infty$ transition curves for the fractal MDE as well.

We investigate this question by first studying the solutions plotted in Figures 7.10-7.29. In constructing the solutions we choose eigenfunctions of generation of birth 2 and of initial eigenvalue $5$ or $6$,  arising from either of the initial functions shown in Figure 7.9,  Figures 7.10-7.15 correspond to solutions of Version 1, Figures 7.16-7.21 correspond to solutions of Version 2, Figures 7.22-7.25 correspond to solutions of Version 3, and Figures 7.26-7.29 correspond to solutions of Version 4. Each Figure includes two images. The one on the left is a plot of three different solutions corresponding to different $(\delta,\ep)$ pairs along the first transition curve; the set of three solutions in each Figure corresponds to either 5-series or 6-series, and each solution been normalized so that the maximum value it attains is 1. The image on the right in each Figure shows an overlook of the three solutions, where the color at each point indicates which solution takes on the biggest value among the three. From the plots, it is evident that these normalized solutions seem to converge as $\ep$ tends to infinity along the first transition curve.

The second way in which we investigate the question is by observing the behavior of the location of relative maxima of solutions as $\ep$ traverses along a transition curve. Here we still choose to study $(\delta,\ep)$ pairs and corresponding solutions on the first transition curve for each version. See Figure \ref{peaks version 1, 5series}-\ref{peaks version 4, 6series}. This is analogous to our study of the relative maxima of solutions to the MDE on the line in Section 4.2. Each of the Figures contains two images indicating the position of relative maxima: the image on the left shows how the positions of the maximas move as $\ep$ becomes larger, and the image on the right is an overlook of the left one where we can clearly see the locations of the peaks. We make several observations: 
 %\st{the image on the left shows vertical blue lines emanating from the locations of the relative maxima, and the image on the right shows an overlook in which the locations the peaks are indicated by blue dots.}
\begin{itemize}
   \item The first thing we can see is that the movement of the peaks are not large, and we do not observe the peaks converging to the boundary.  
    
\item    The second observation is that the $5$ and $6$ series behave quite differently. We can observe jumps of peaks in all the versions for $5$ series, but the number of peaks is usually 3. A special case is figure $\ref{peaks version 4, 5series}$, where many peaks occur when $\ep$ is quite large. However, for $6$ series, we can observe many more peaks when $\ep$ is very small, but most of them do not appear when $\ep$ is large. Until now, we do not have explanation for these behaviors.
\end{itemize}

\begin{comment}
\st{In Figure 7.10 we present three solutions to Version 1 of the fractal MDE, where $(\delta,\ep)$ lies on the first transition curve and the eigenvalues arise from 5-series; the three solutions correspond to $\ep$-coordinate 1000, 2000, and 3000. The solutions have been normalized so that each takes on a maximum value of 1. The left figure of Figure 7.10 shows the solutions as three-dimensional graphs; the right figure shows an overlook of the three solutions, where the color indicates which value takes the biggest value among the three.

In our experiment, we choose eigenfunctions of generation of birth 2, with initial functions shown in Figure \ref{initial_eigenfunction}.}

\st{The first thing we can see is that the movement of the peaks are not large, and we do not observe the peaks converging to the boundary. The second observation is that the $5$ and $6$ series behave quite differently. We can observe jumps of peaks in all the versions for $5$ series, but the number of peaks is usually 3. A special case is figure $\ref{peaks version 4, 5series}$, where many peaks occur when $\ep$ is quite large. However, for $6$ series, we can observe many more peaks when $\ep$ is very small, but most of them do not appear when $\ep$ is large. Until now, we do not have explanation for these behaviors. }
\end{comment}

\begin{figure}[h]
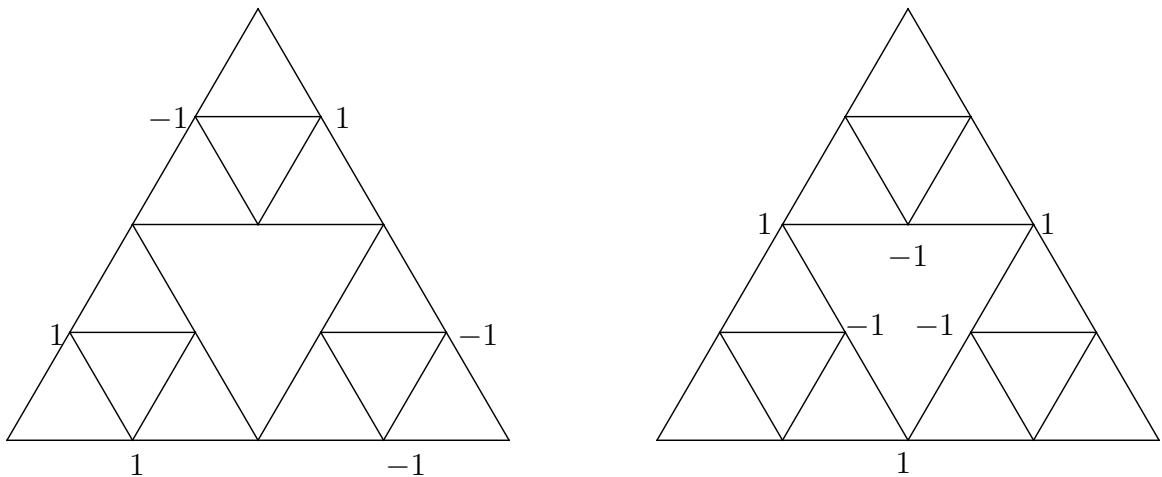

    \centering
    \includegraphics[height=6cm]{gamma2.pdf}\qquad\qquad
    \includegraphics[height=6cm]{gamma2.pdf}
    \begin{picture}(0,0)
    \put(-394,-9){$1$}
    \put(-298,-9){$-1$}
    \put(-424,40){$1$}
    \put(-387,122){$-1$} 
    \put(-317,122){$1$}
    \put(-271,40){$-1$}
%hello shiping    %hello anthony let's review introduction %good idea. Should I call you on skype? 
    \put(-53,82){$1$}
    \put(-159,82){$1$}
    \put(-107,-8){$1$}
    \put(-126,44){$-1$}
    \put(-100,44){$-1$}
    \put(-110,70){$-1$}
    \end{picture}
    \caption{Initial values of the eigenfunctions with birth of generation $2$ and initial eigenvalues $5$ (left) and $6$(right).}
    \label{initial_eigenfunction}
\end{figure}

%\st{In Figure \ref{sg-m1-5series-1000} we have three solutions, corresponding to the points along the first transition curve for the version 1 equation with $\varepsilon$-coordinate 1000, 2000, and 3000, respectively. Similarly, we have Figure \ref{sg-m1-5series-10000} for $\varepsilon$-coordinate 10000, 11000, and 12000, respectively, and figure \ref{sg-m1-5series-20000} $\varepsilon$-coordinate 20000, 25000, and 30000, respectively. The solutions have been normalized so that each takes on a value 1 at its highest peak (relative maximum).  On the right side of each graph, we put a overlook of the three solutions, where the color implies which functions takes the biggest value among the three. From the plots, is evident that these normalized solutions seem to converge as epsilon tends to infinity along the first transition curve.}

\begin{figure}[!h]
    \centering
    \includegraphics[height=6cm]{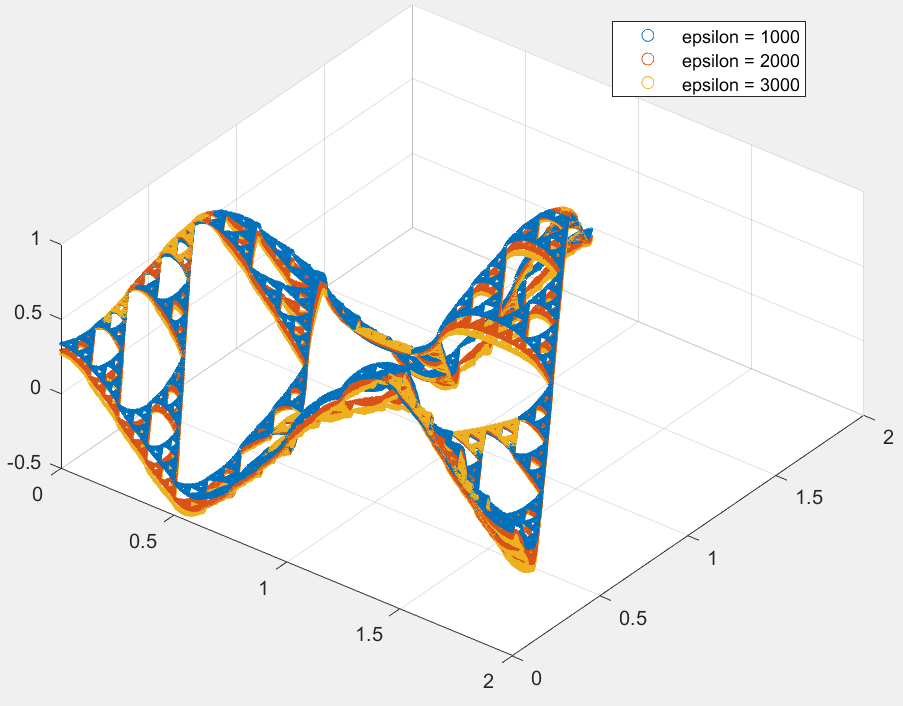}
    \includegraphics[height=6cm]{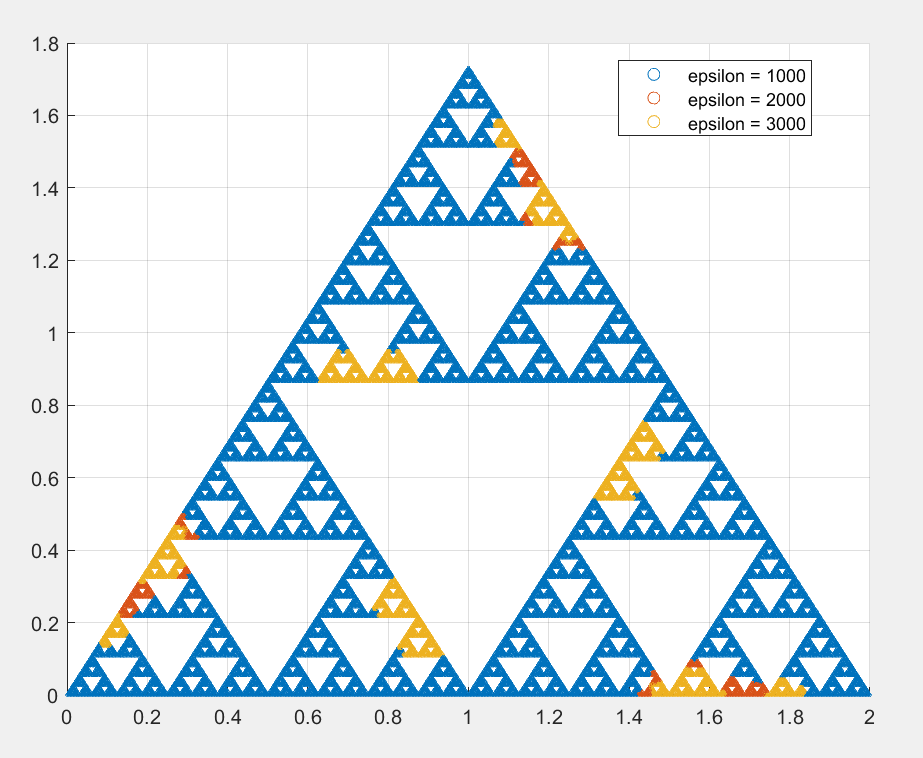}
    \caption{Solutions on the first curve for choice 1, initial eigenvalue $5$,  with $\ep=1000,2000,3000$.}
    \label{sg-m1-5series-1000}
\end{figure}

\begin{figure}[!h]
    \centering
    \includegraphics[height=6cm]{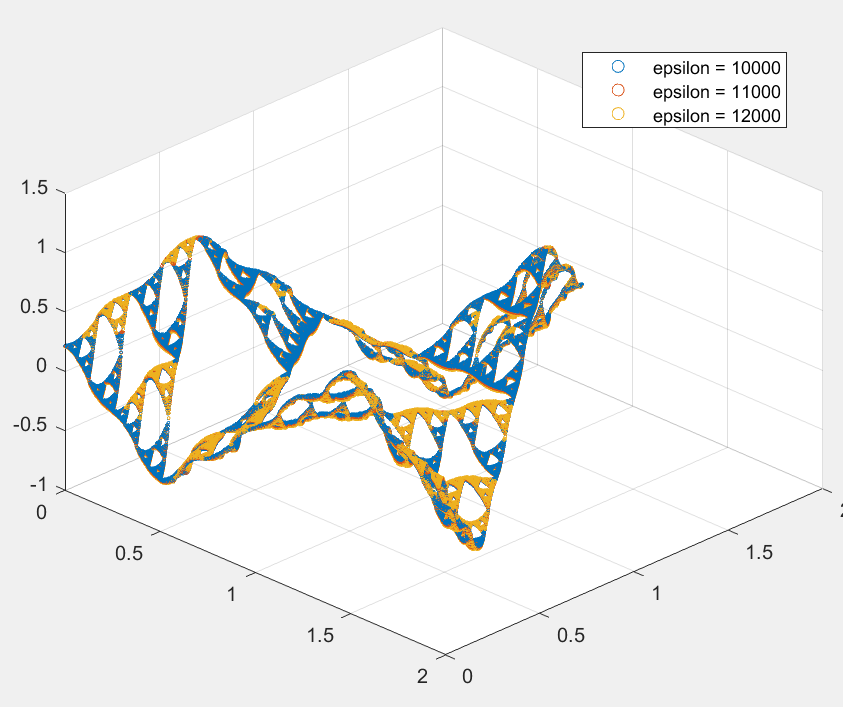}
    \includegraphics[height=6cm]{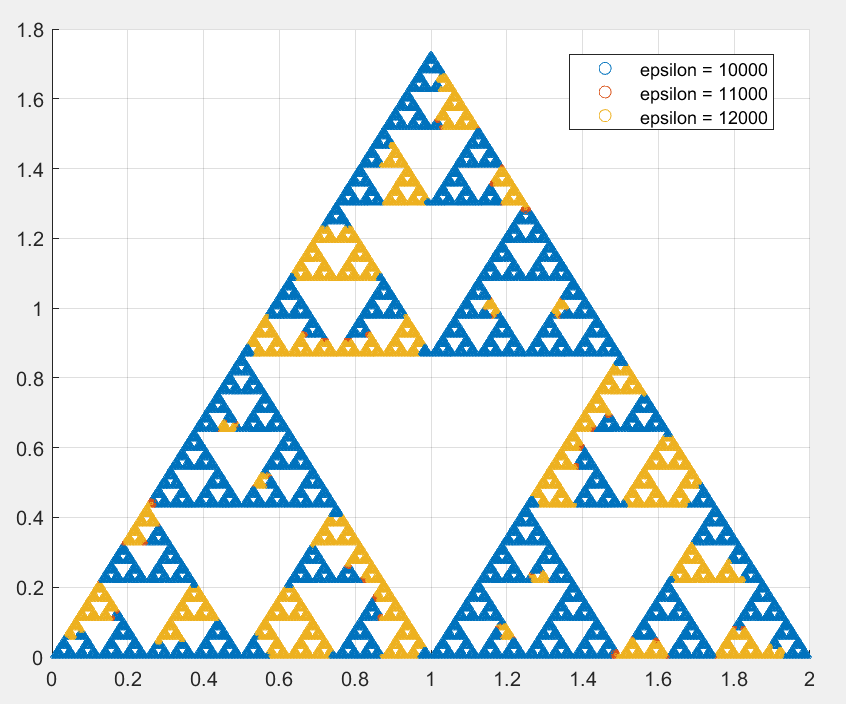}
    \caption{Solutions on the first curve for choice 1, initial eigenvalue $5$,  with $\ep=10000,11000,12000$.}
    \label{sg-m1-5series-10000}
\end{figure}

\begin{figure}[!h]
    \centering
    \includegraphics[height=6cm]{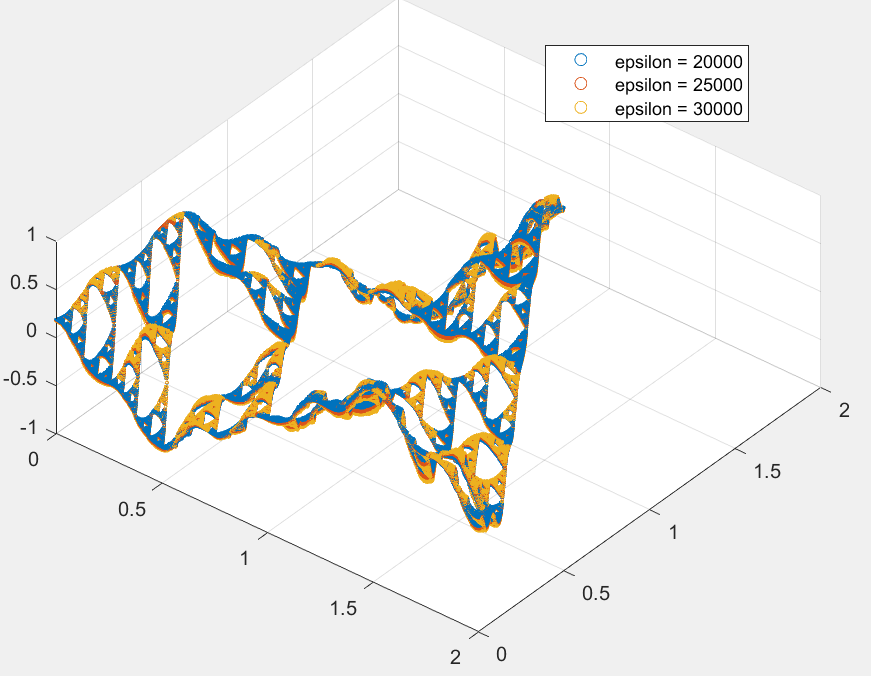}
    \includegraphics[height=6cm]{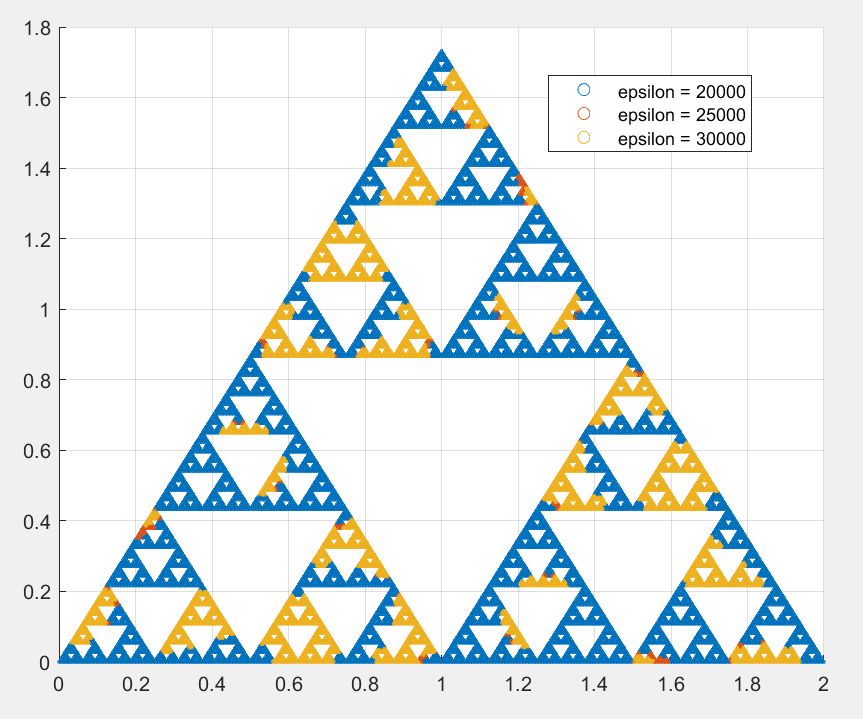}
    \caption{Solutions on the first curve for choice 1, initial eigenvalue $5$,  with $\ep=20000,25000,30000$.}
    \label{sg-m1-5series-20000}
\end{figure}
\clearpage{}

\begin{figure}[!h]
    \centering
    \includegraphics[height=6cm]{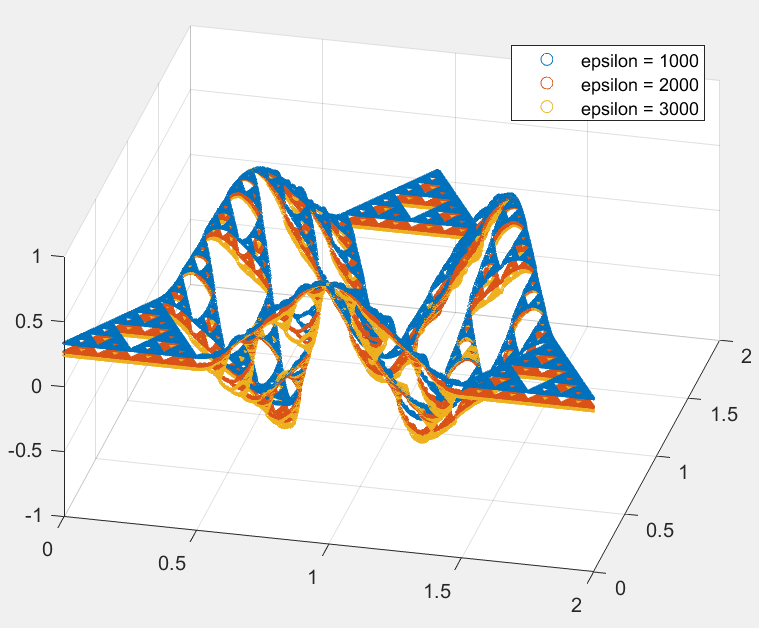}
    \includegraphics[height=6cm]{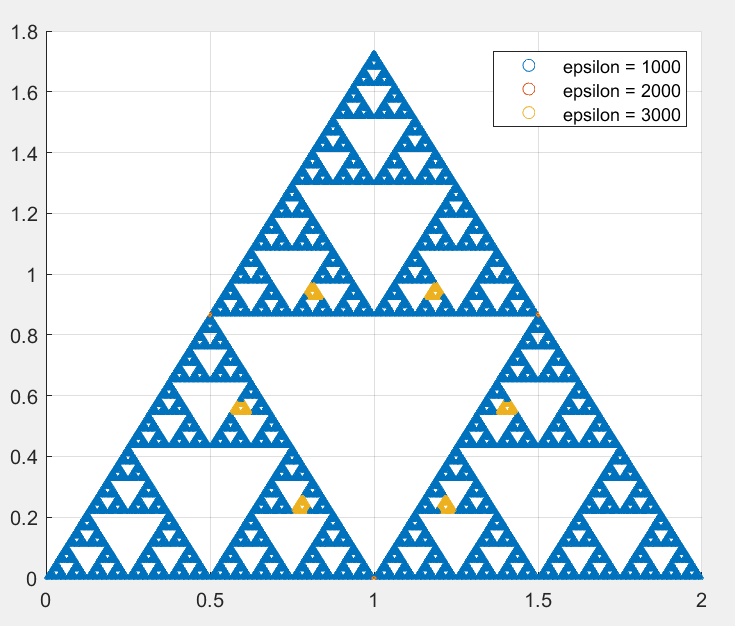}
    \caption{Solutions on the first curve for choice 1, initial eigenvalue $6$,  initial eigenvalue $6$, with $\ep=1000,2000,3000$.}
    \label{sg-m1-6series-1000}
\end{figure}

\begin{figure}[!h]
    \centering
    \includegraphics[height=6cm]{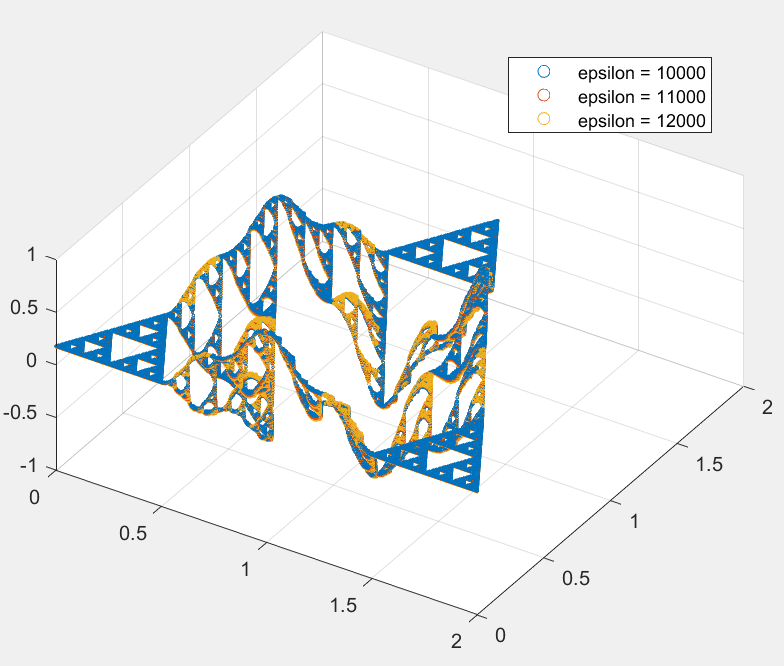}
    \includegraphics[height=6cm]{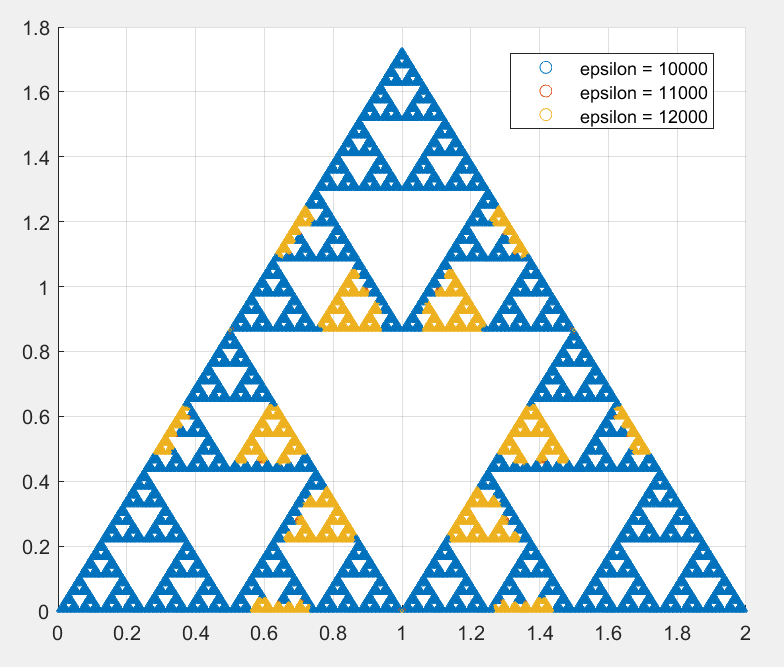}
    \caption{Solutions on the first curve for choice 1, initial eigenvalue $6$, initial eigenvalue $6$, with, with $\ep=10000,11000,12000$.}
    \label{sg-m1-6series-10000}
\end{figure}

\begin{figure}[!h]
    \centering
    \includegraphics[height=6cm]{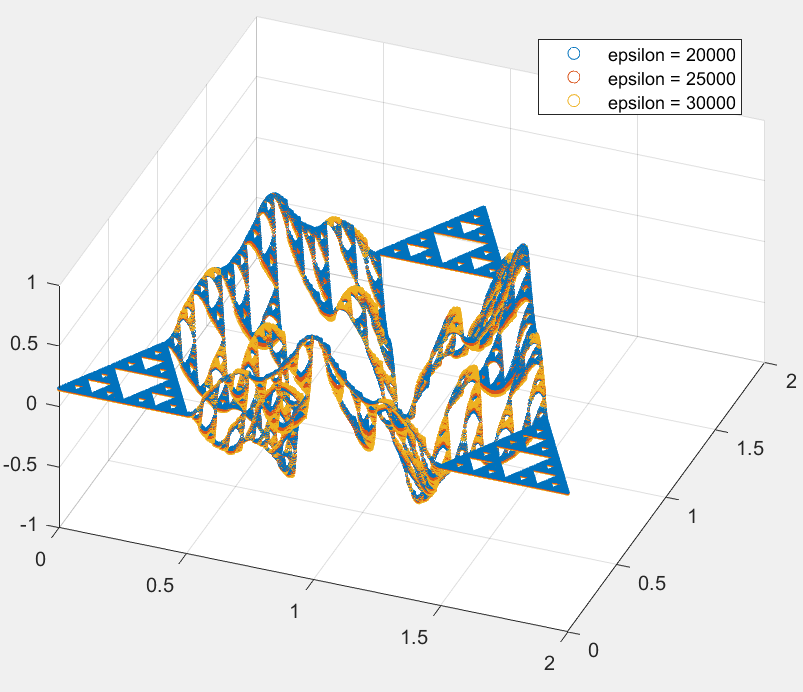}
    \includegraphics[height=6cm]{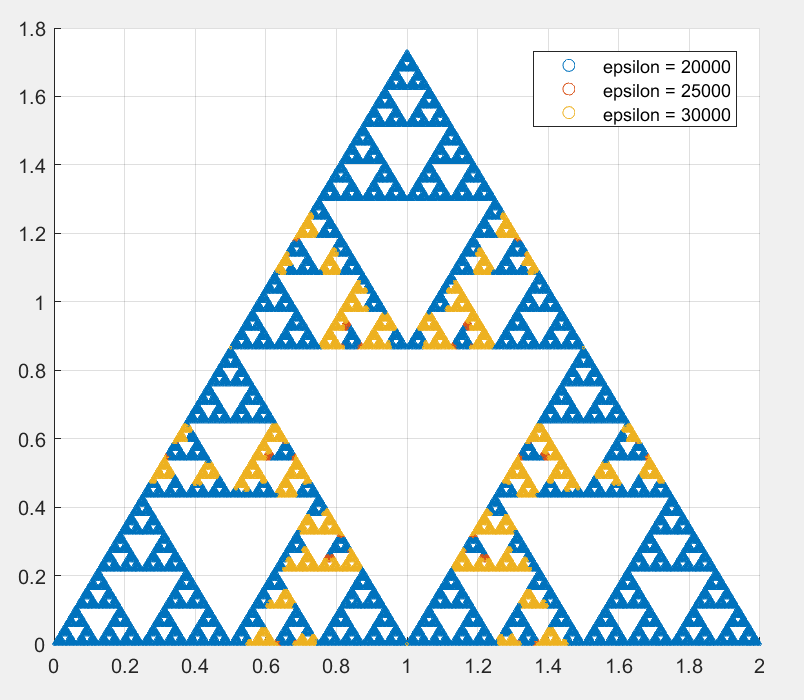}
    \caption{Solutions on the first curve for choice 1, initial eigenvalue $6$, initial eigenvalue $6$, with $\ep=20000,25000,30000$.}
    \label{sg-m1-6series-20000}
\end{figure}
\clearpage{}

%\st{We also do the same experiment on version 2, with the same class of eigenfunctions, see figure \ref{sg-m2-5series-1000},figure \ref{sg-m2-5series-10000} and figure \ref{sg-m2-5series-20000}.}

\begin{figure}[!h]
    \centering
    \includegraphics[height=6cm]{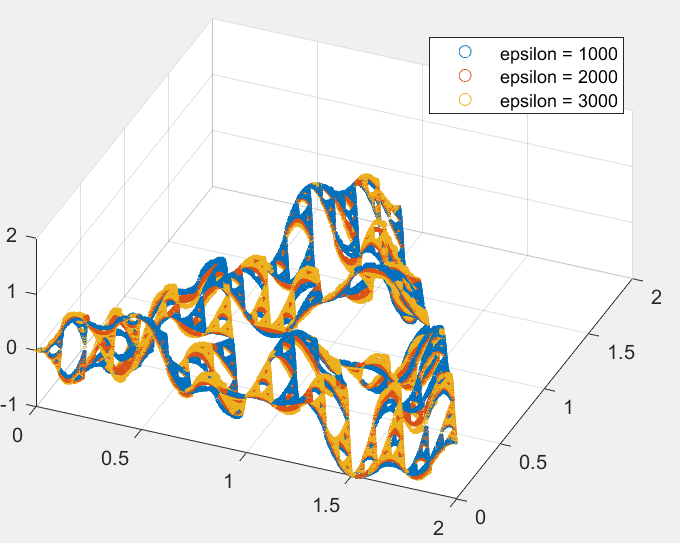}
    \includegraphics[height=6cm]{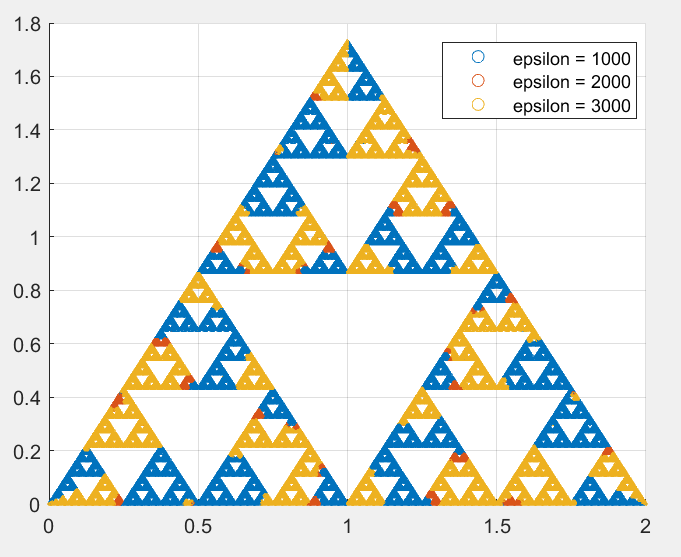}
    \caption{Solutions on the first curve for choice 2, initial eigenvalue $5$, with $\ep=1000,2000,3000$}
    \label{sg-m2-5series-1000}
\end{figure}

\begin{figure}[!h]
    \centering
    \includegraphics[height=6cm]{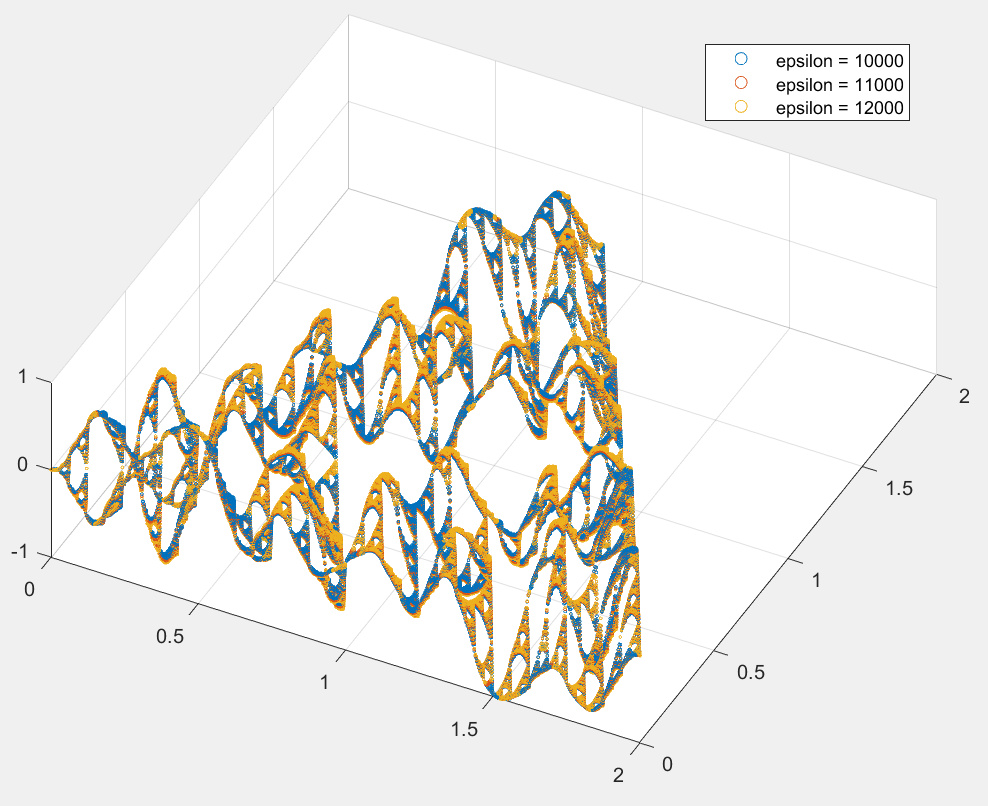}
    \includegraphics[height=6cm]{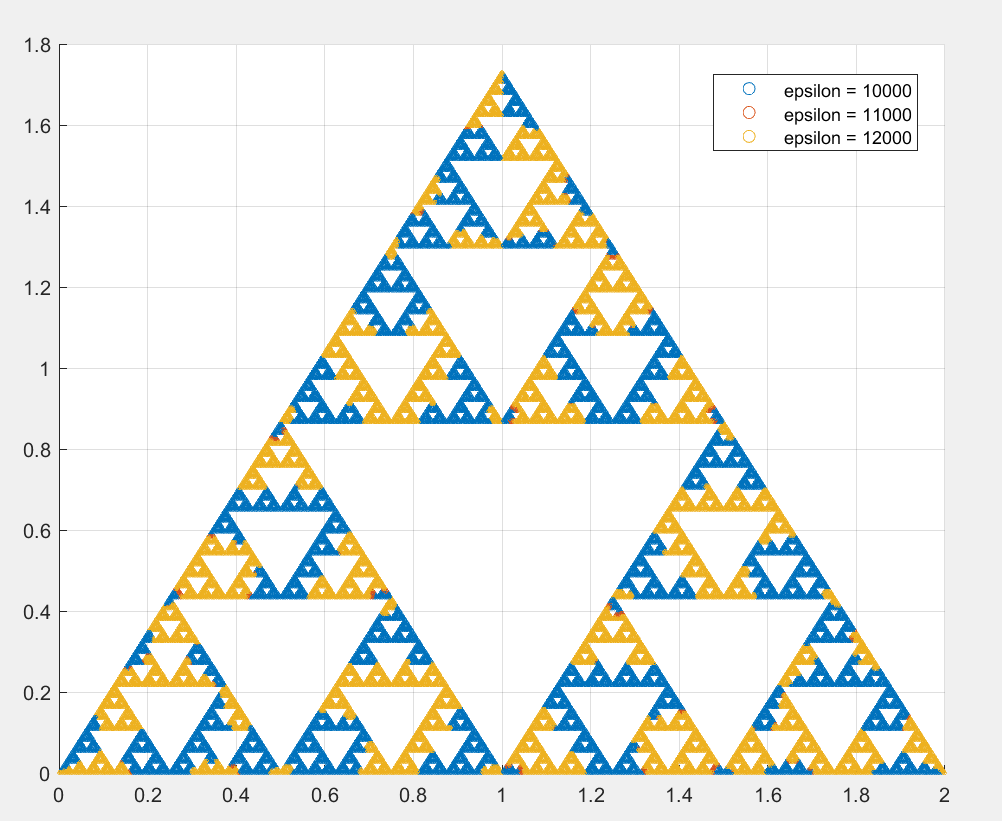}
    \caption{Solutions on the first curve for choice 2, initial eigenvalue $5$,  with $\ep=10000,11000,12000$}
    \label{sg-m2-5series-10000}
\end{figure}

\begin{figure}[!h]
    \centering
    \includegraphics[height=6cm]{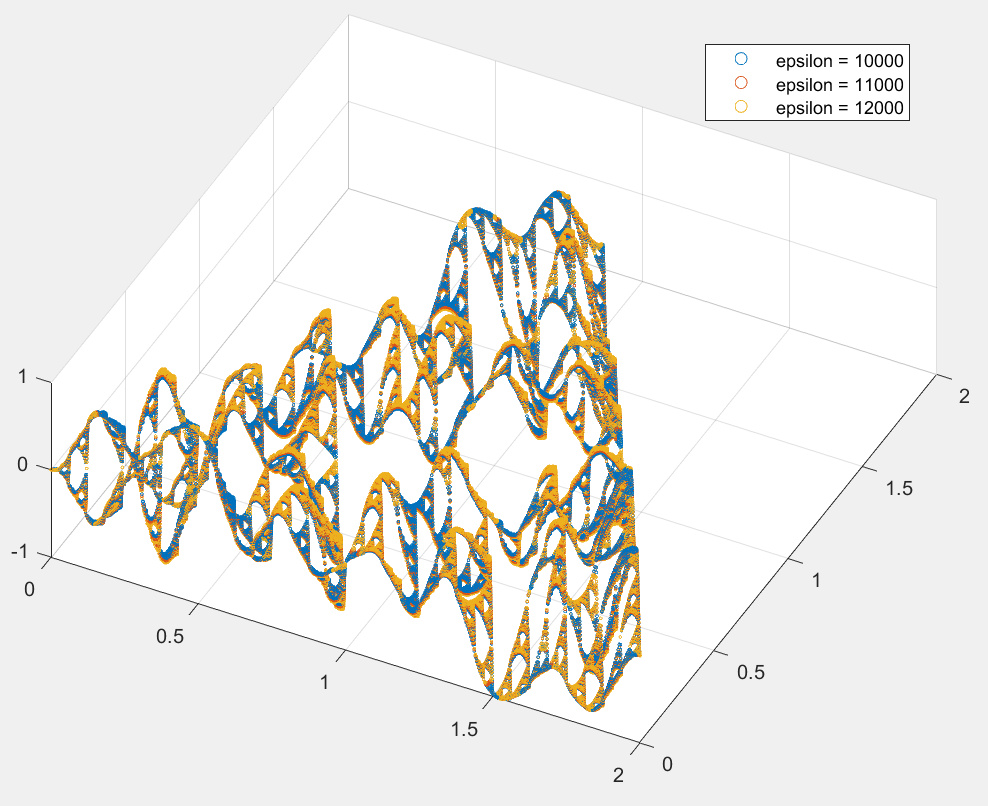}
    \includegraphics[height=6cm]{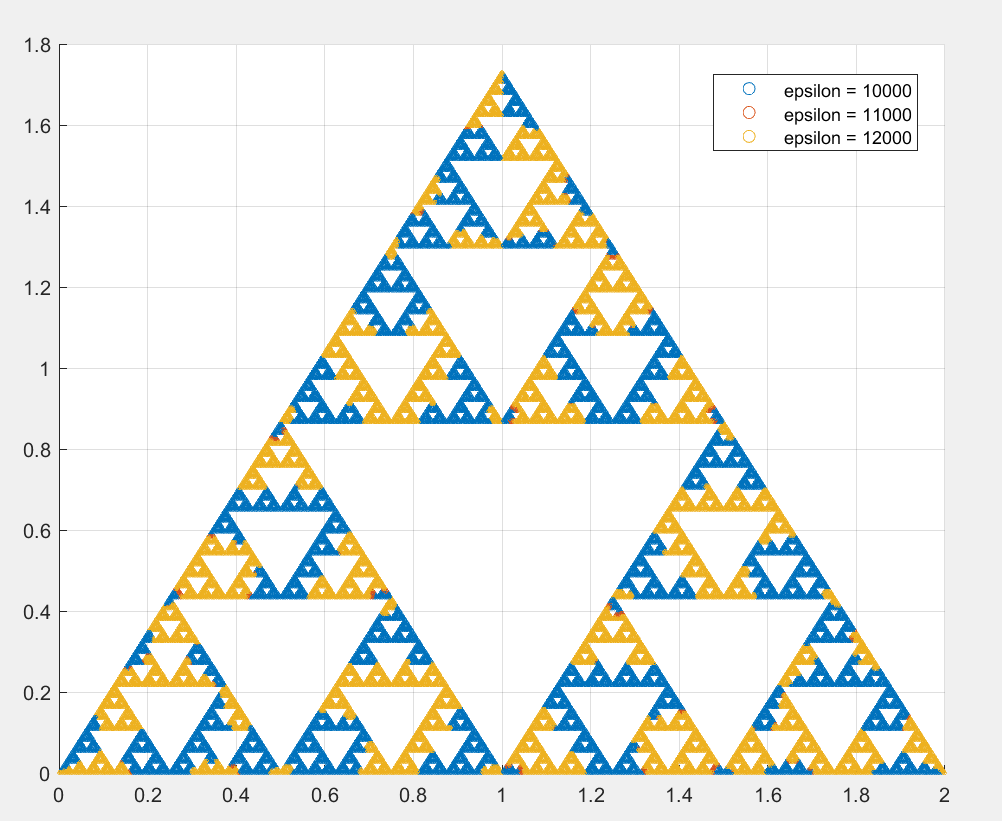}
    \caption{Solutions on the first curve for choice 2, initial eigenvalue $5$,  with $\ep=20000,25000,30000$}
    \label{sg-m2-5series-20000}
\end{figure}
\clearpage{}

\begin{figure}[!h]
    \centering
    \includegraphics[height=6cm]{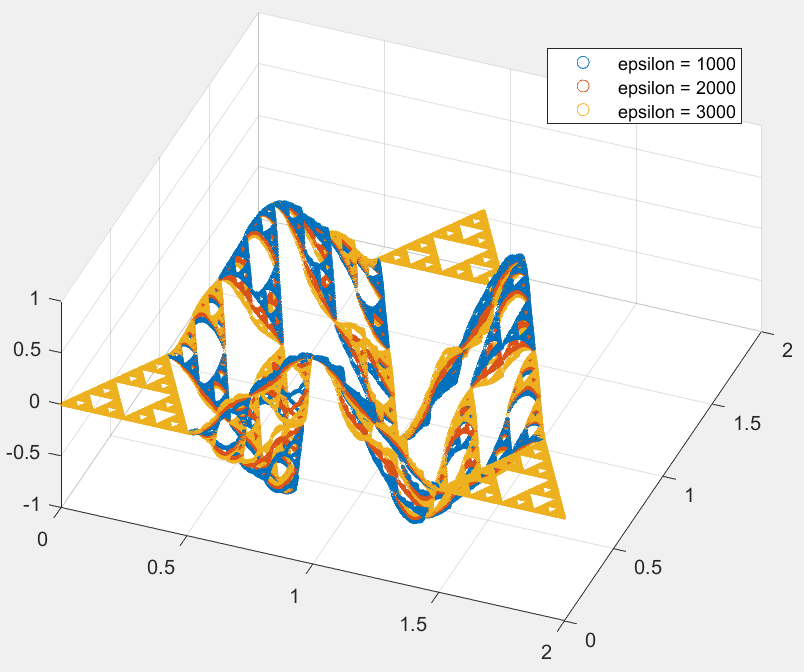}
    \includegraphics[height=6cm]{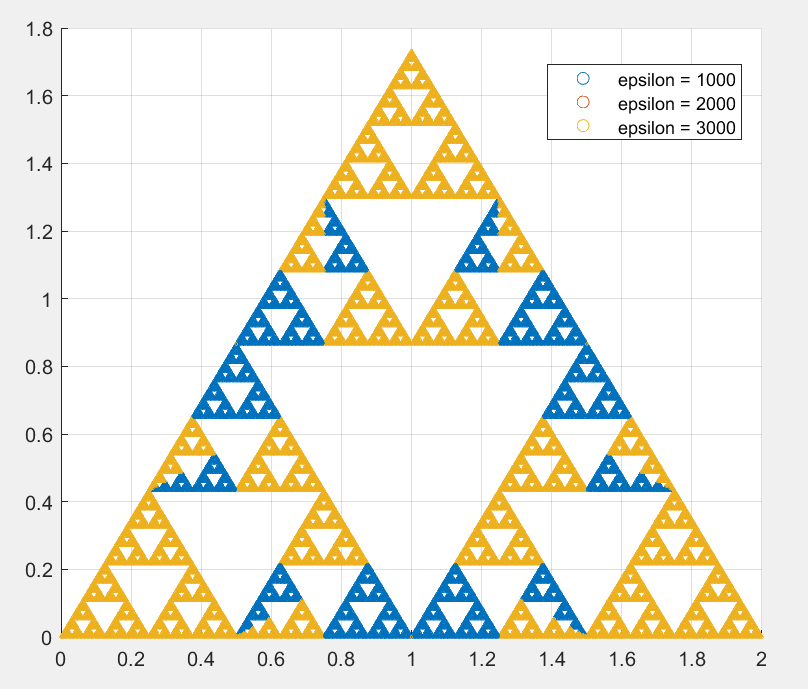}
    \caption{Solutions on the first curve for choice 2, initial eigenvalue $6$, with $\ep=1000,2000,3000$.}
    \label{sg-m2-6series-1000}
\end{figure}

\begin{figure}[!h]
    \centering
    \includegraphics[height=6cm]{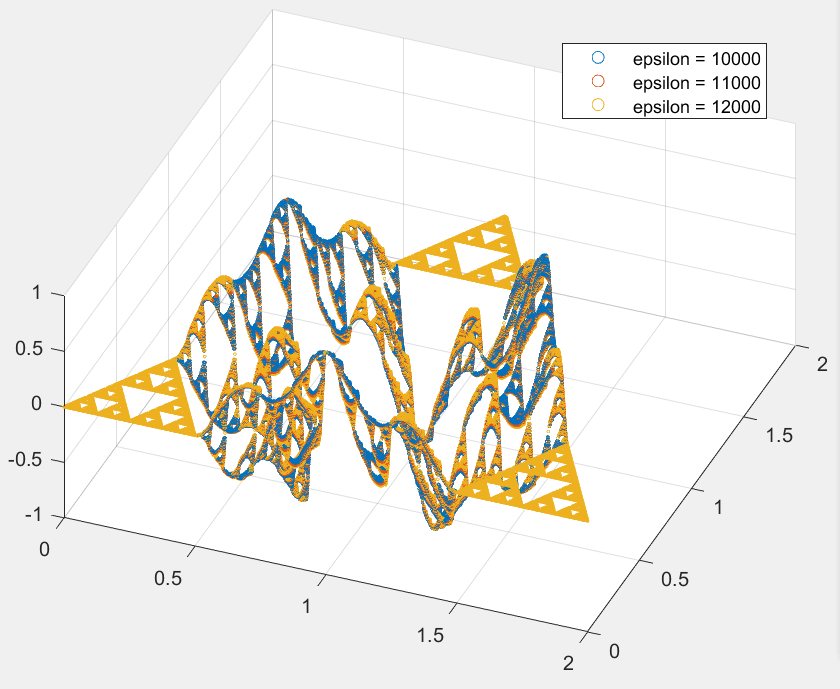}
    \includegraphics[height=6cm]{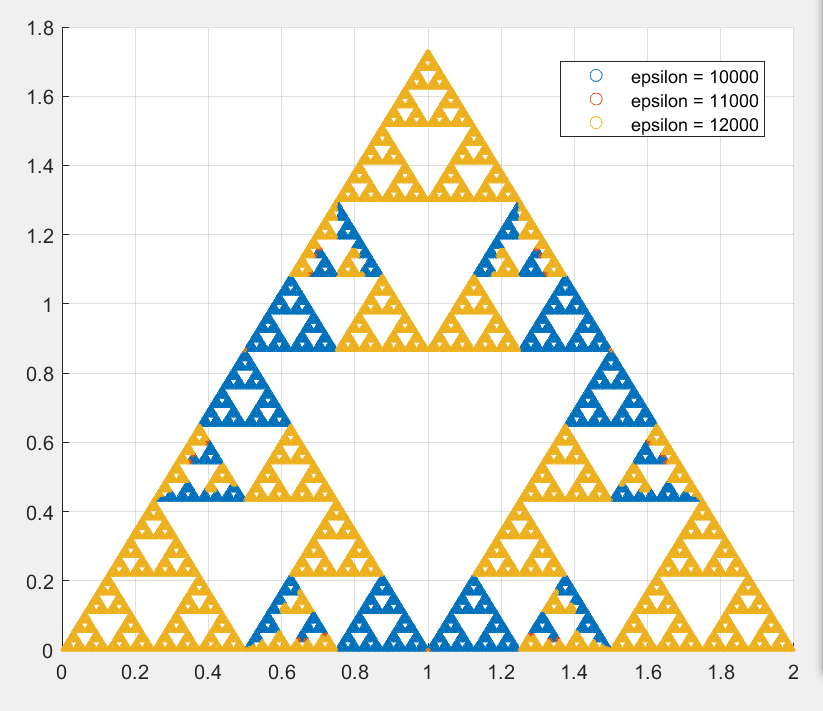}
    \caption{Solutions on the first curve for choice 2, initial eigenvalue $6$, initial eigenvalue $6$, with, with $\ep=10000,11000,12000$.}
    \label{sg-m2-6series-10000}
\end{figure}

\begin{figure}[!h]
    \centering    \includegraphics[height=6cm]{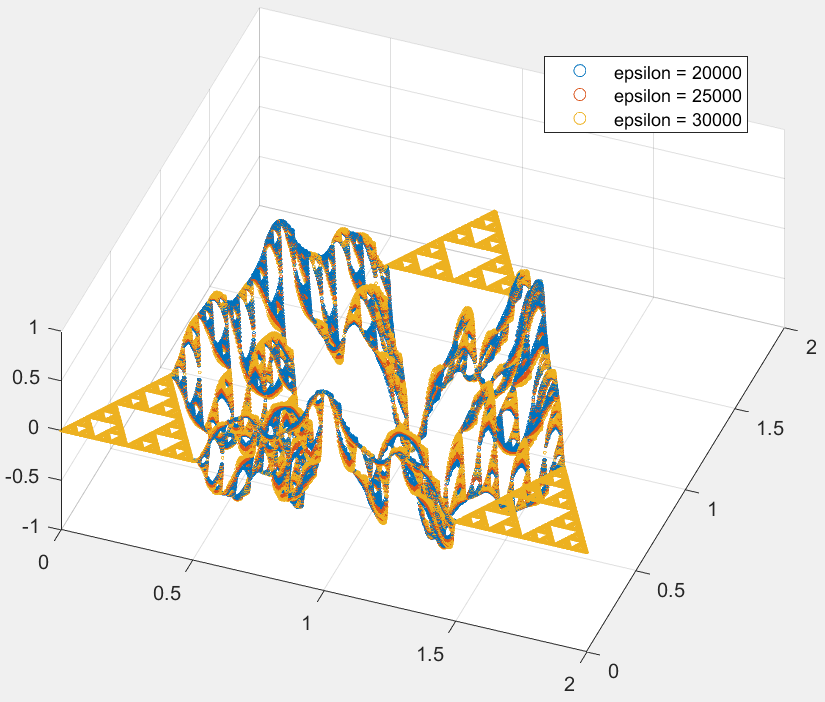}
    \includegraphics[height=6cm]{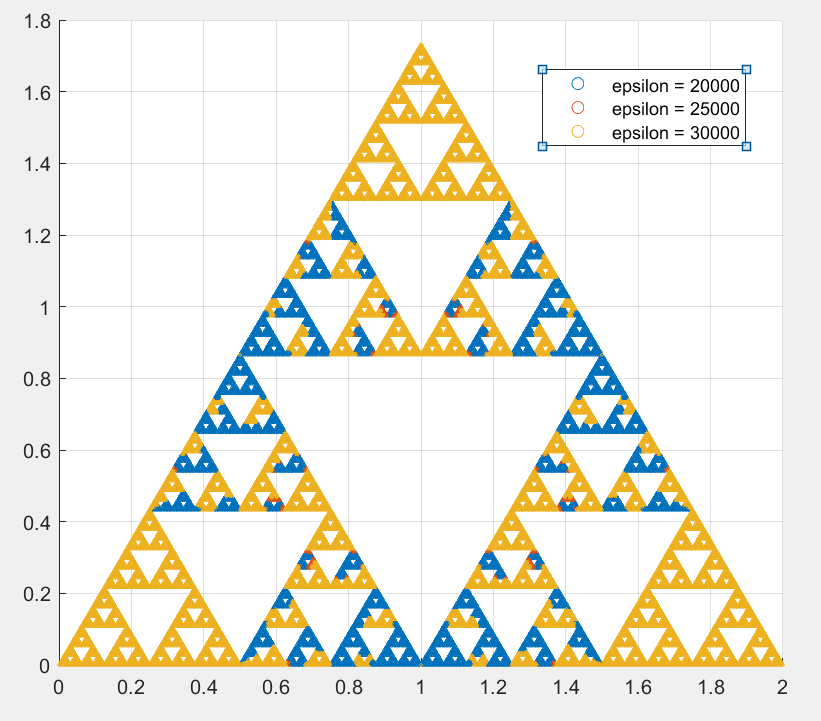}
    \caption{Solutions on the first curve for choice 2, initial eigenvalue $6$, initial eigenvalue $6$, with $\ep=20000,25000,30000$.}
    \label{sg-m2-6series-20000}
\end{figure}
\clearpage{}

%\st{We also do experiments on version 3 and version 4. The same classes of eigenfunctions are considered. }

\begin{figure}[!h]
    \centering
    \includegraphics[height=6cm]{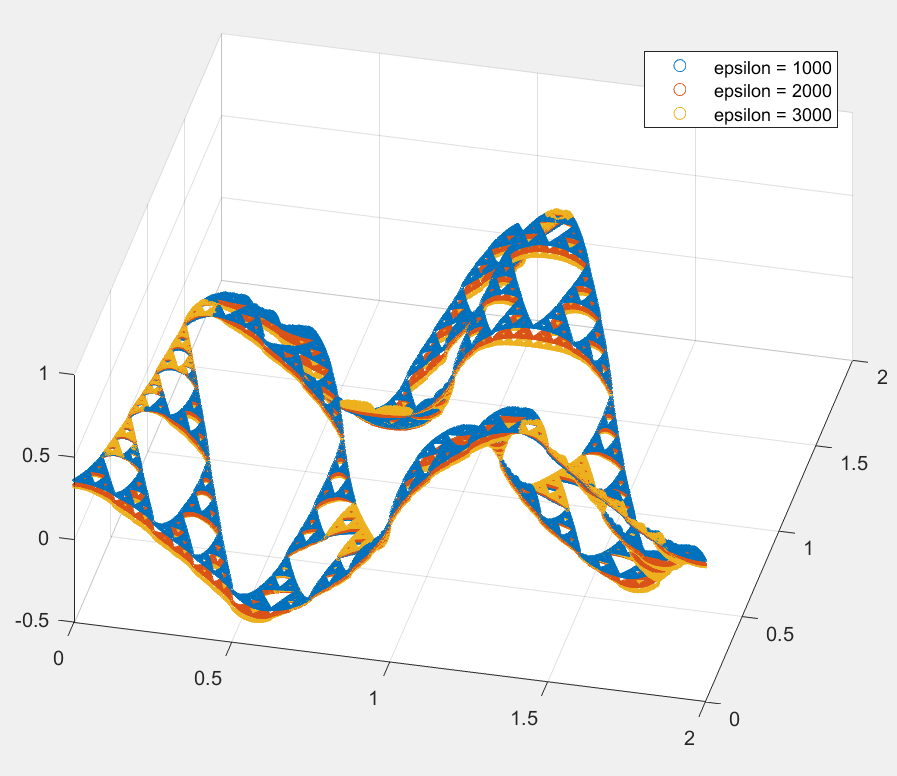}
    \includegraphics[height=6cm]{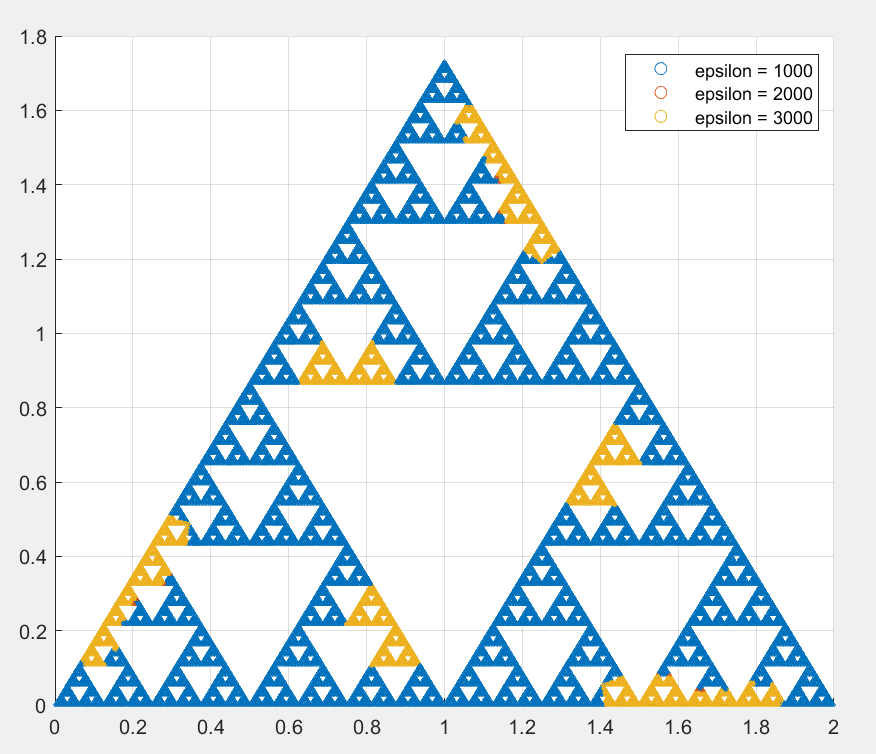}
    \caption{Solutions on the first curve for choice 3, initial eigenvalue $5$, with $\ep=1000,2000,3000$.}
    \label{sg-m3-5series-1000}
\end{figure}

\begin{figure}[!h]
    \centering
    \includegraphics[height=6cm]{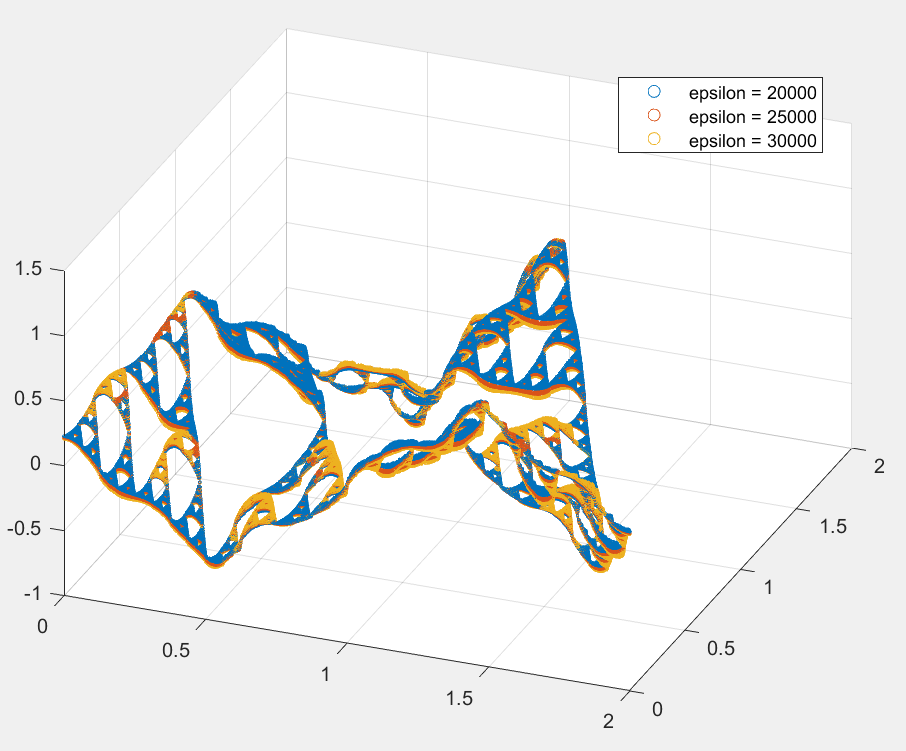}
    \includegraphics[height=6cm]{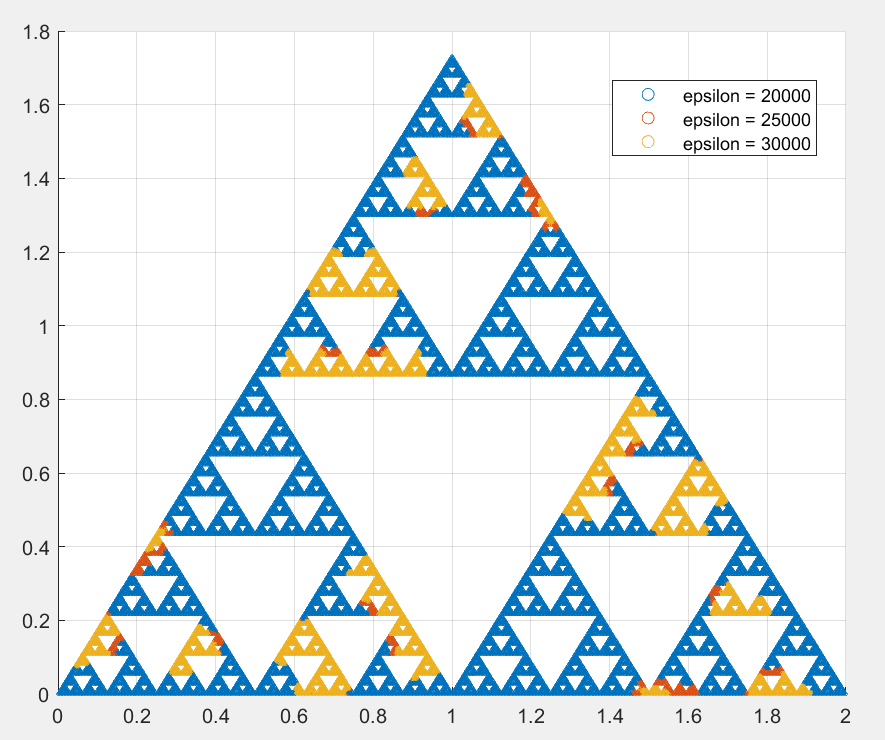}
    \caption{Solutions on the first curve for choice 3, initial eigenvalue $5$, with $\ep=20000,25000,30000$.}
    \label{sg-m3-5series-20000}
\end{figure}

\begin{figure}[!h]
    \centering
    \includegraphics[height=6cm]{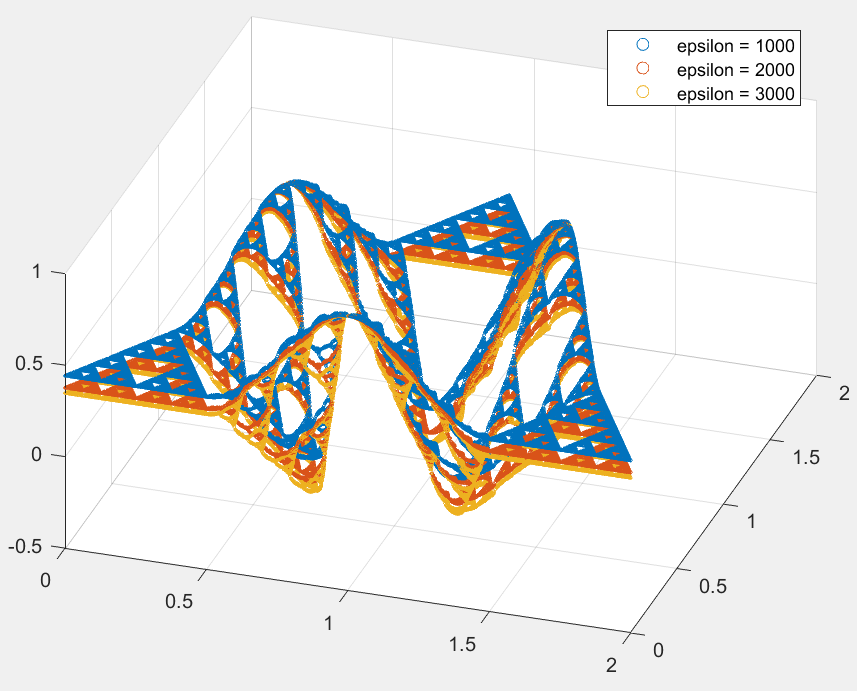}
    \includegraphics[height=6cm]{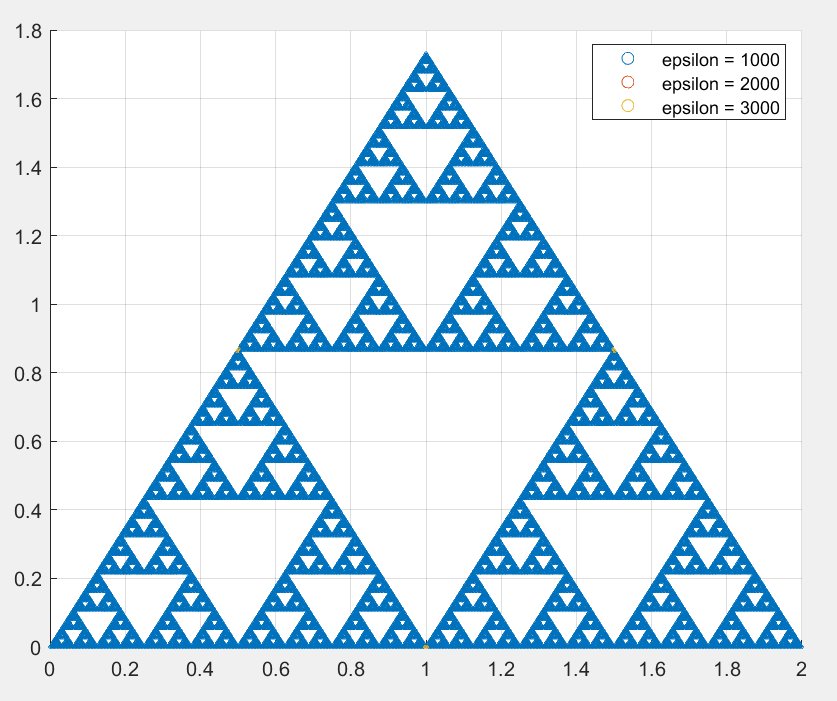}
    \caption{Solutions on the first curve for choice 3, initial eigenvalue $6$, with $\ep=1000,2000,3000$.}
    \label{sg-m3-6series-1000}
\end{figure}
\clearpage{}

\begin{figure}[!h]
    \centering
    \includegraphics[height=6cm]{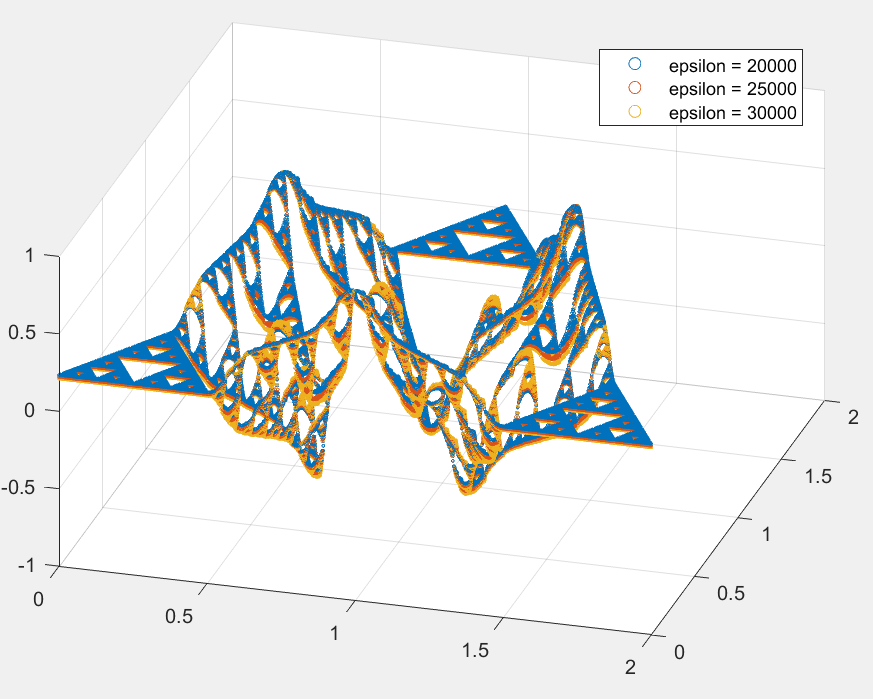}
    \includegraphics[height=6cm]{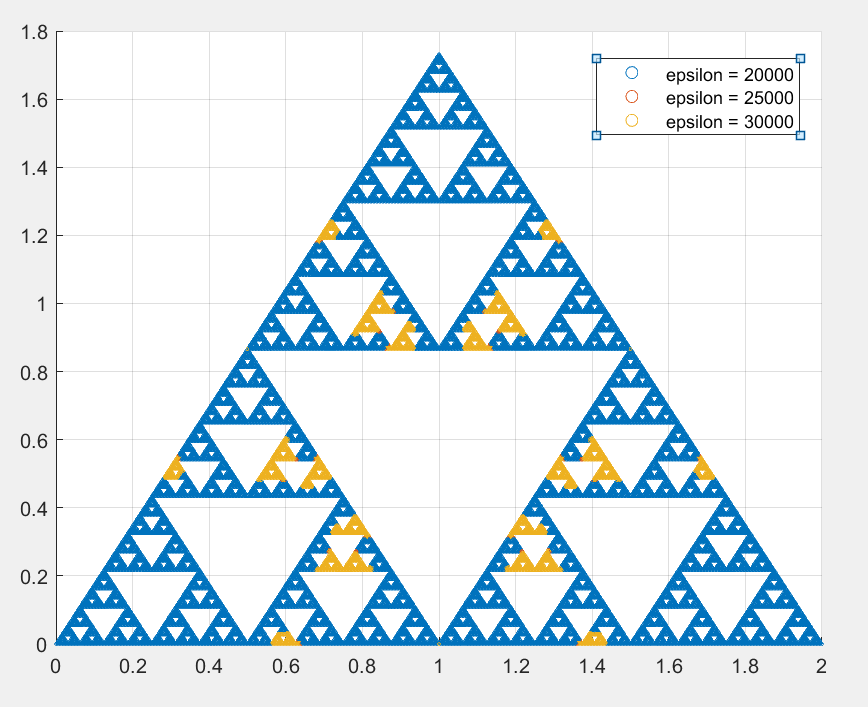}
    \caption{Solutions on the first curve for choice 3, initial eigenvalue $6$, with $\ep=20000,25000,30000$.}
    \label{sg-m3-6series-20000}
\end{figure}

\begin{figure}[!h]
    \centering
    \includegraphics[height=6cm]{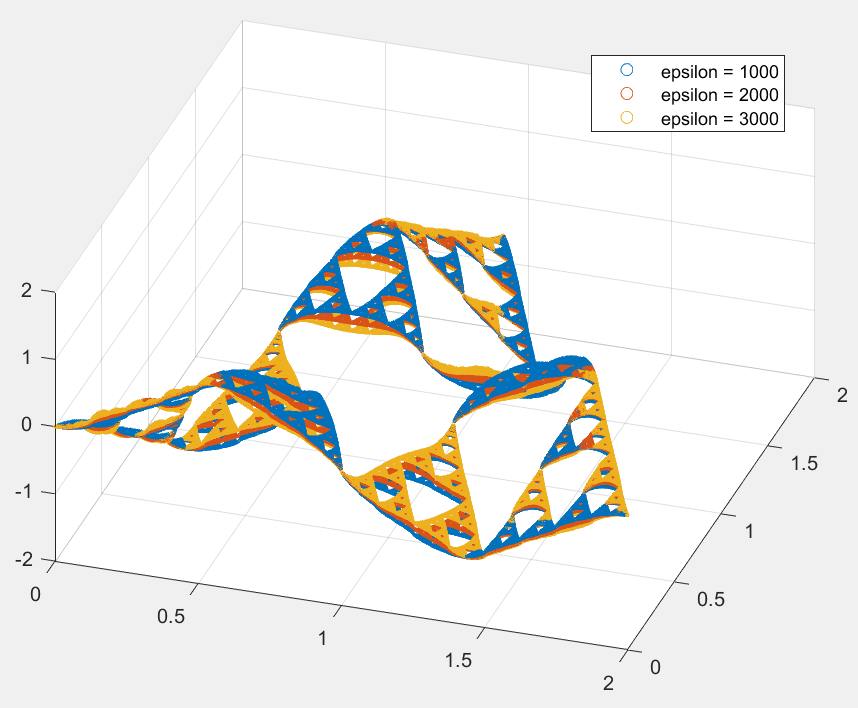}
    \includegraphics[height=6cm]{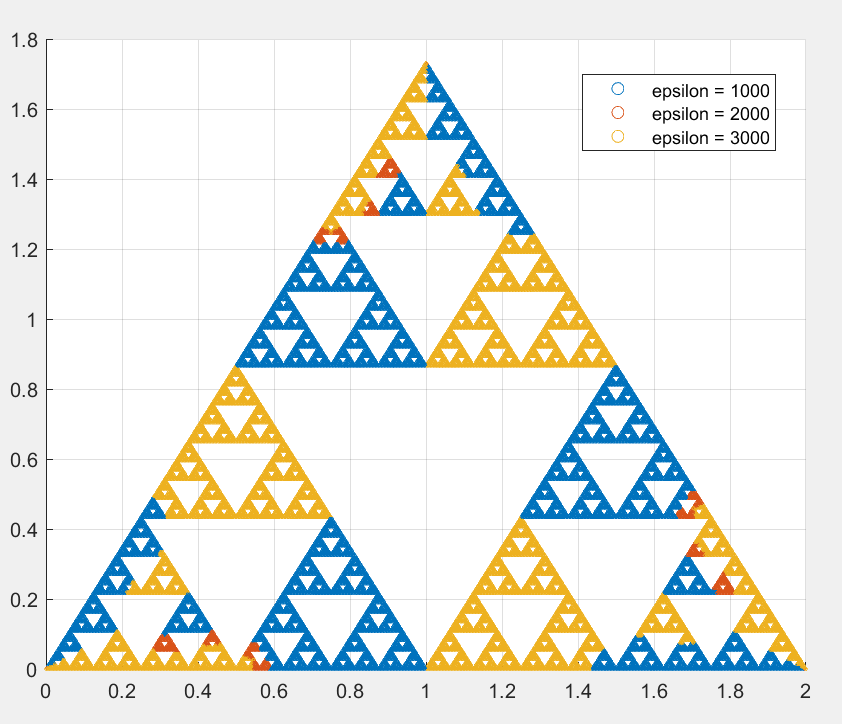}
    \caption{Solutions on the first curve for choice 4, initial eigenvalue $5$, with $\ep=1000,2000,3000$.}
    \label{sg-m4-5series-1000}
\end{figure}

\begin{figure}[!h]
    \centering
    \includegraphics[height=6cm]{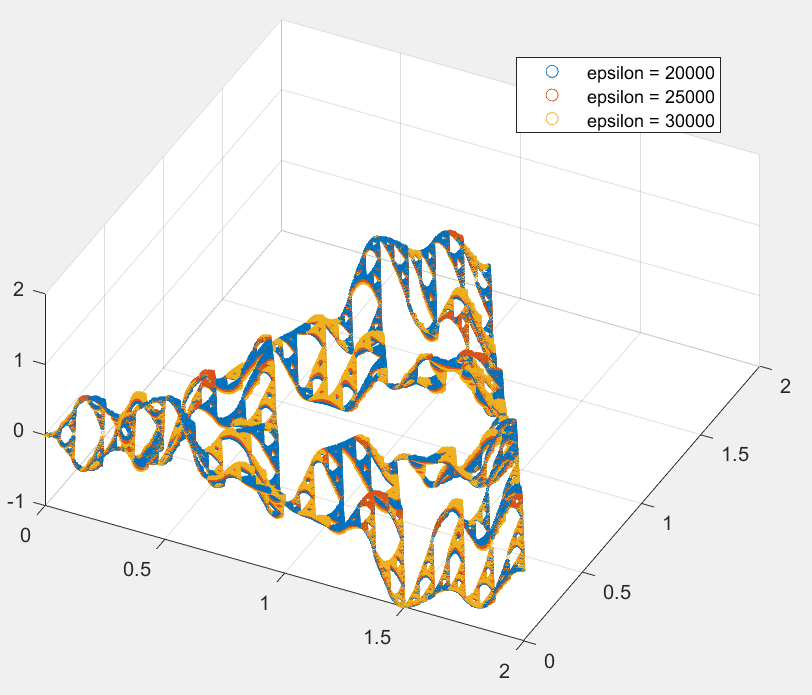}
    \includegraphics[height=6cm]{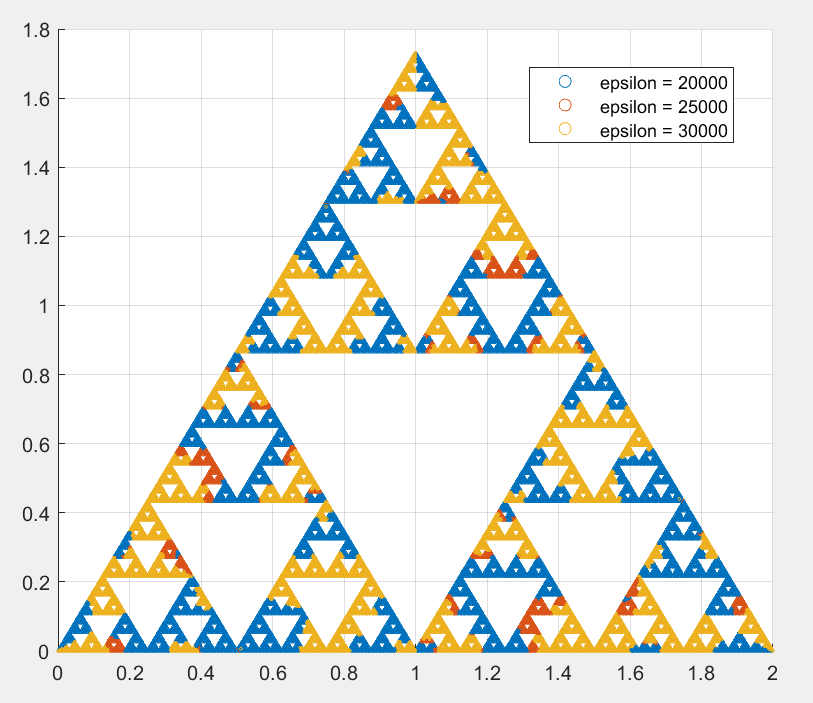}
    \caption{Solutions on the first curve for choice 4, initial eigenvalue $5$, with $\ep=20000,25000,30000$.}
    \label{sg-m4-5series-20000}
\end{figure}
\clearpage{}

\begin{figure}[!h]
    \centering
    \includegraphics[height=6cm]{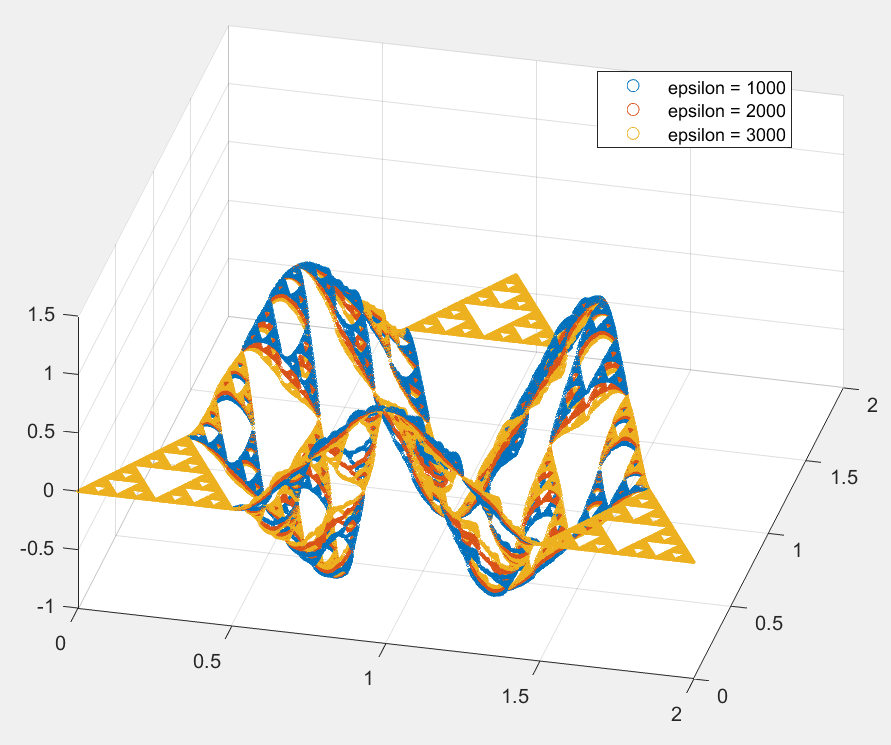}
    \includegraphics[height=6cm]{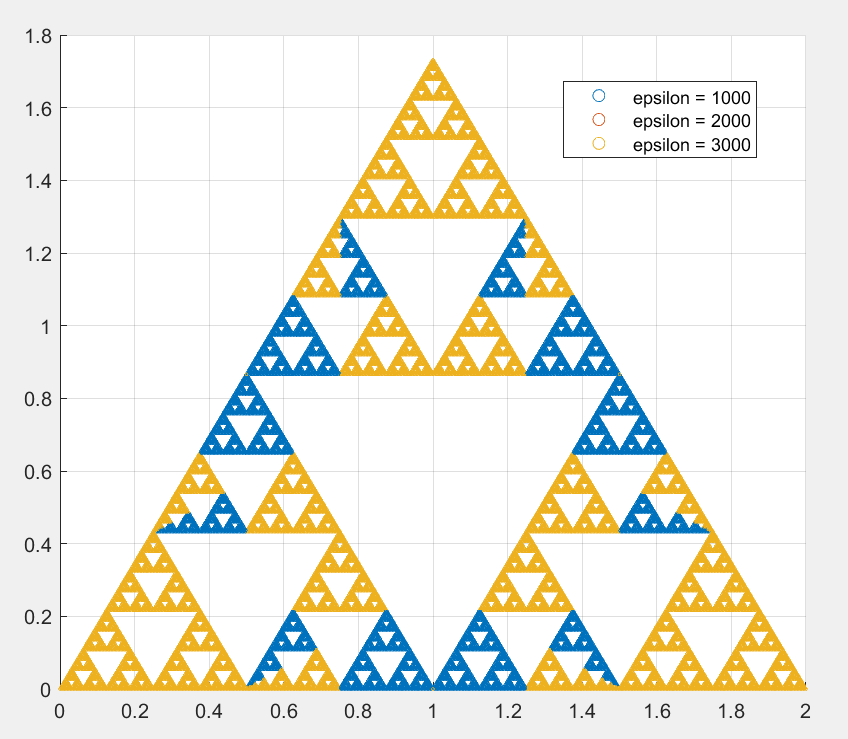}
    \caption{Solutions on the first curve for choice 4, initial eigenvalue $6$, with $\ep=1000,2000,3000$.}
    \label{sg-m4-6series-1000}
\end{figure}

\begin{figure}[!h]
    \centering
    \includegraphics[height=6cm]{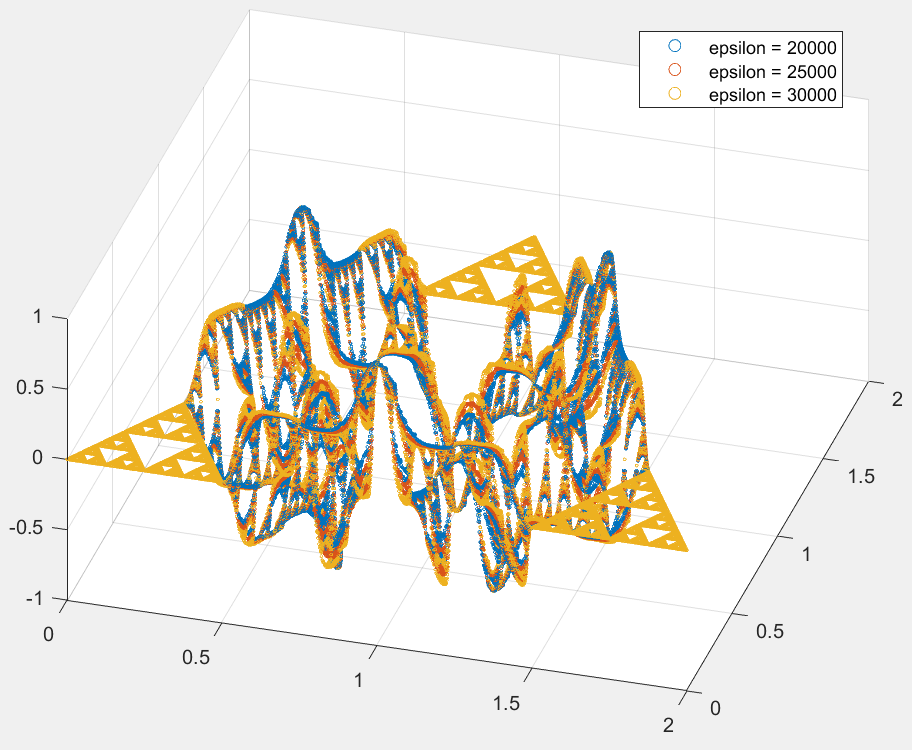}
    \includegraphics[height=6cm]{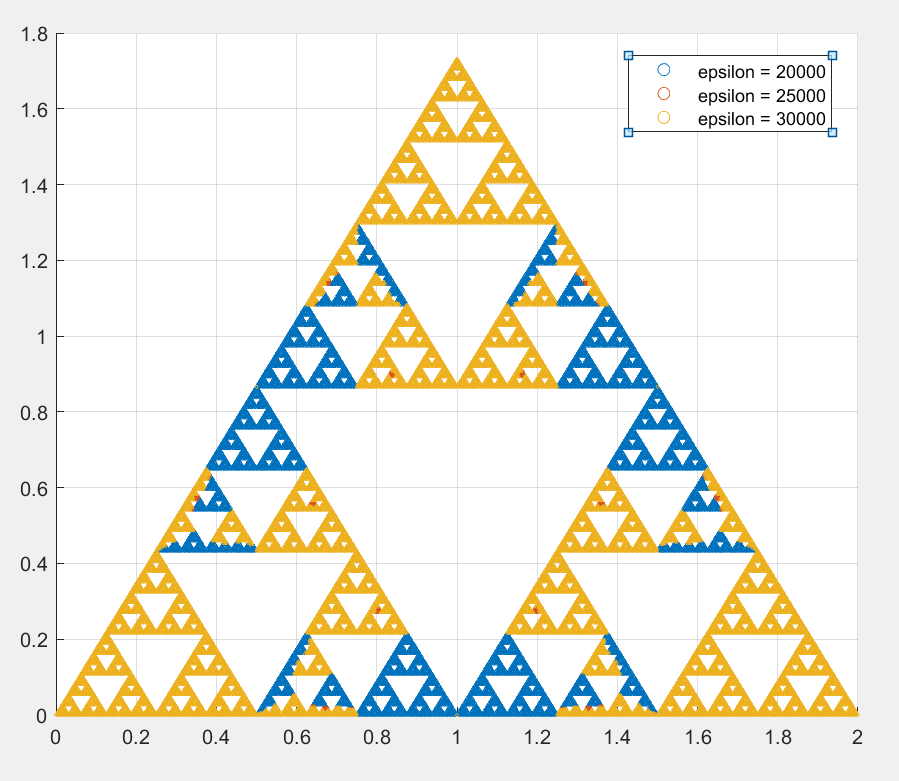}
    \caption{Solutions on the first curve for choice 4, initial eigenvalue $6$, with $\ep=20000,25000,30000$.}
    \label{sg-m4-6series-20000}
\end{figure}

\clearpage{}

\begin{comment}
\textcolor{green}{As shown in section 2, the transition curves can be organized into four different classes, made up of points $\{p(\sin kt,\ep)\},\{p(cos kt,\ep)\},\{p(\sin \frac{2k+1}{2}t,\ep)\},\{p(cos\frac{2k+1}{2}t,\ep)\}$ with our new notation. It is natural to study them separately.}

\textcolor{green}{First consider the class $\{p(\sin kt,\ep)\}$. In particular, we let $k=1,2,3$ and $\varepsilon=5,10,20,40,80,160$}(by the symmetry of the transition curves across the $\delta$-axis, it suffices to only consider positive values of $j$, not negative ones). The normalized solution corresponding to each of the first few terms of this sequence is plotted in Figure \textcolor{green}{\ref{linesolsconv1}-\ref{linesolsconv3}, where $15\times 15$ truncated matrices are used for computation}. By `normalized solution', we mean the plot of the equation of the solution multiplied by a scale factor so as to cause the relative maxima closest to the $\ep$-axis to have $\ep$-coordinate equal to 1. The solid black curve corresponds to $\ep=0,$ the red curve corresponds to $\ep=5$, the orange curve corresponds to $\ep=10$, the green curve corresponds to $\ep=20$, the blue curve corresponds to $\ep=40$, the purple curve corresponds to $\ep=80$. The horizontal axis is the $t$-axis, and the vertical axis is the $u$-axis.

\end{comment}

%\pagebreak

\begin{figure}[!h]
    \centering
    \includegraphics[height=5.5cm]{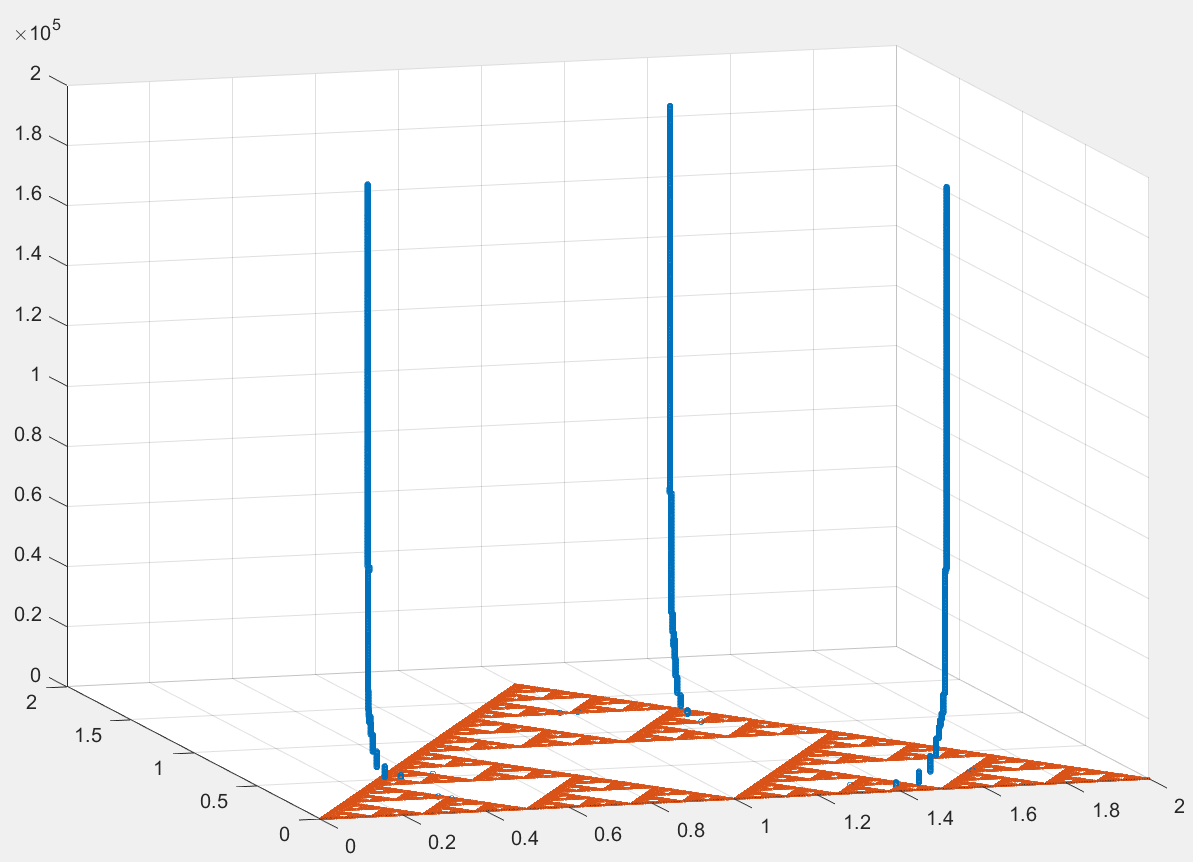}
    \includegraphics[height=5.5cm]{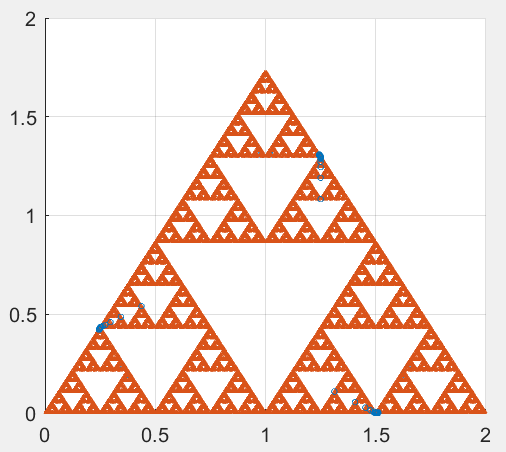}
    \caption{The position of peaks for version 1, 5 series.}
    \label{peaks version 1, 5series}
\end{figure}

\begin{figure}[!h]
    \centering
    \includegraphics[height=5.5cm]{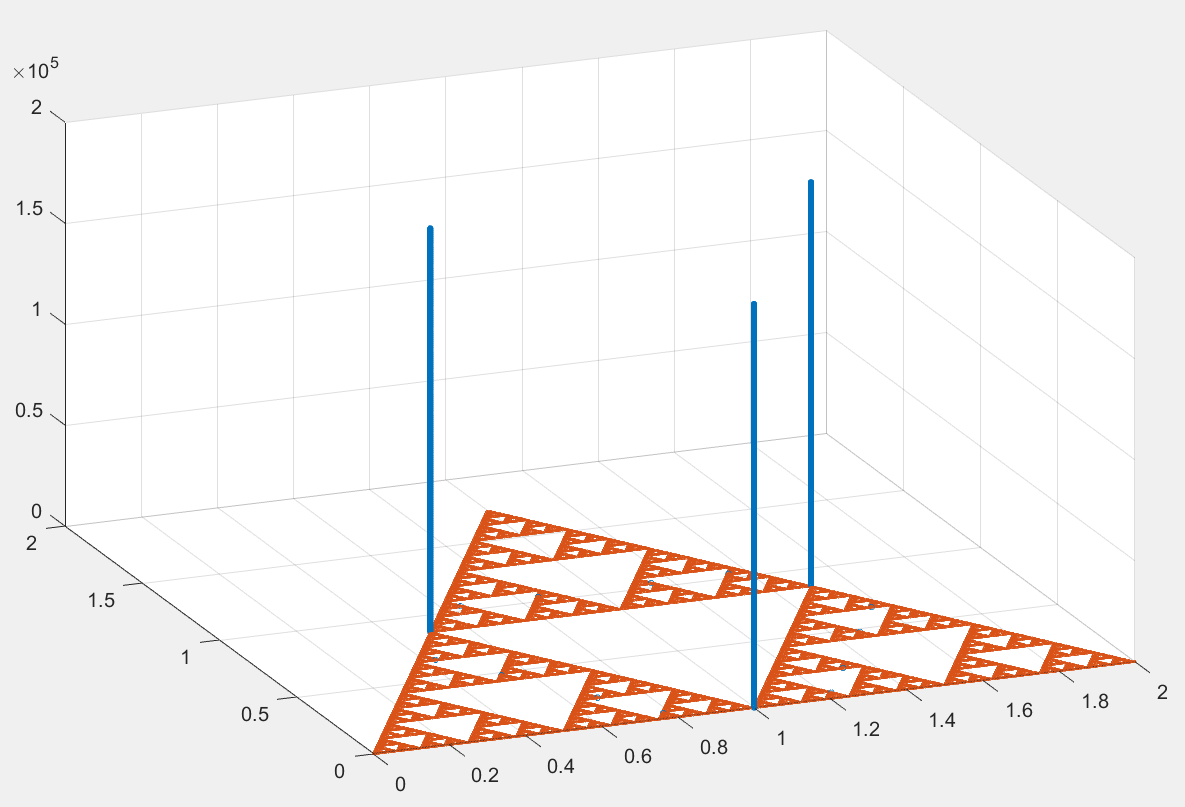}
    \includegraphics[height=5.5cm]{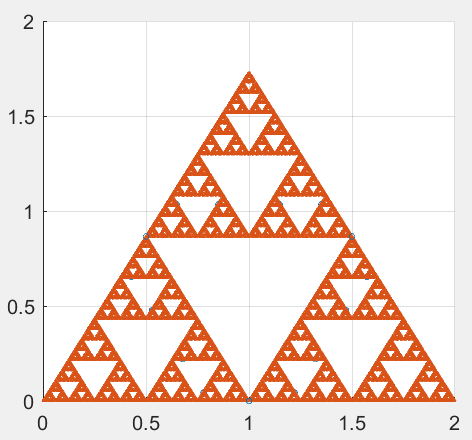}
    \caption{The position of peaks for version 1, 6 series.}
    \label{peaks version 1, 6series}
\end{figure}

\begin{figure}[!h]
    \centering
    \includegraphics[height=5.5cm]{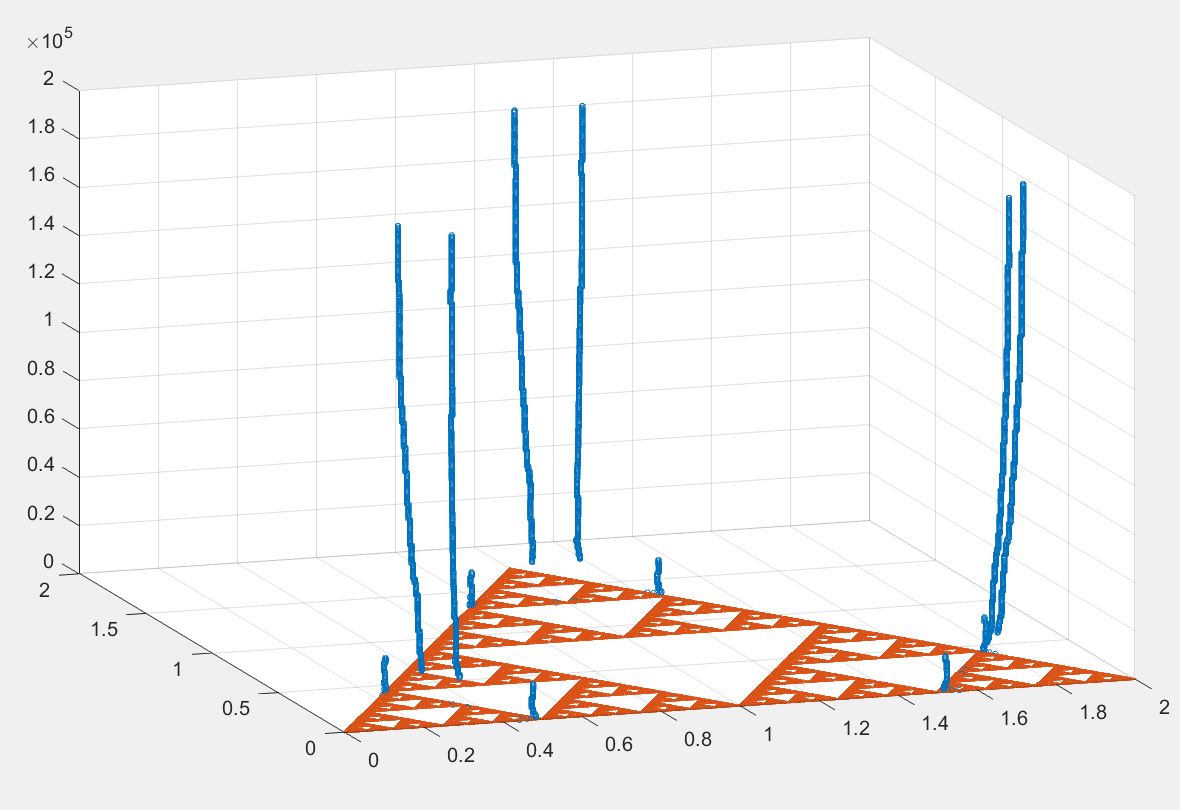}
    \includegraphics[height=5.5cm]{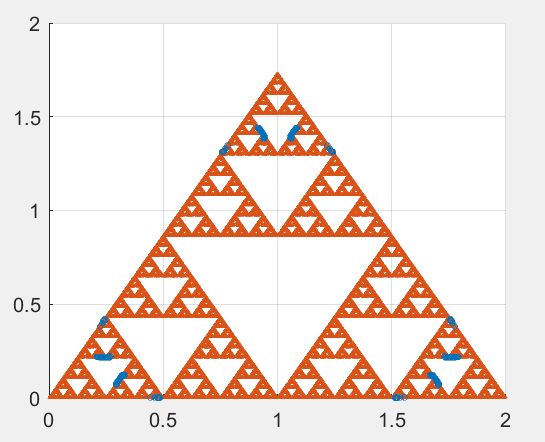}
    \caption{The position of peaks for version 2, 5 series.}
    \label{peaks version 2, 5series}
\end{figure}

\clearpage{}

\begin{figure}[!h]
    \centering
    \includegraphics[height=6cm]{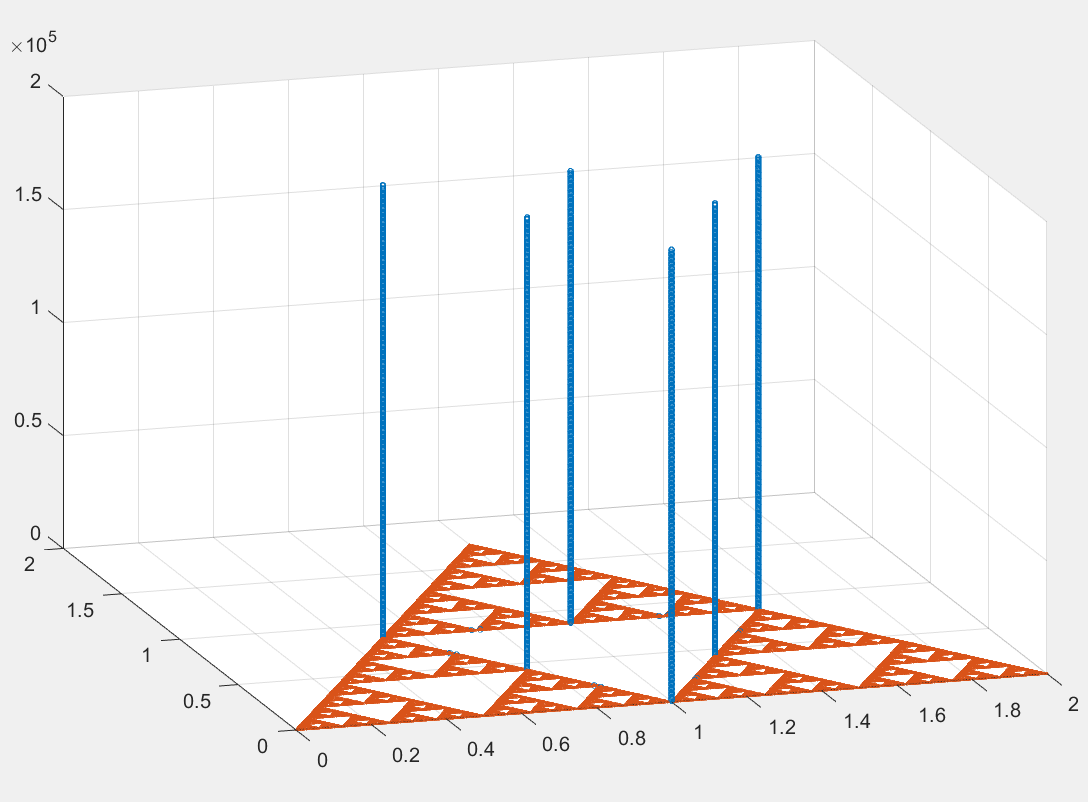}
    \includegraphics[height=6cm]{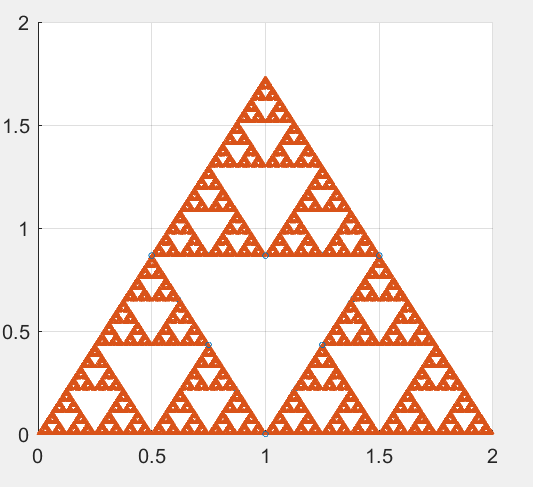}
    \caption{The position of peaks for version 2, 6 series.}
    \label{peaks version 2, 6series}
\end{figure}

%\st{Lastly, let's look at version 3 and 4. They are more complicated. }

\begin{figure}[!h]
    \centering
    \includegraphics[height=6cm]{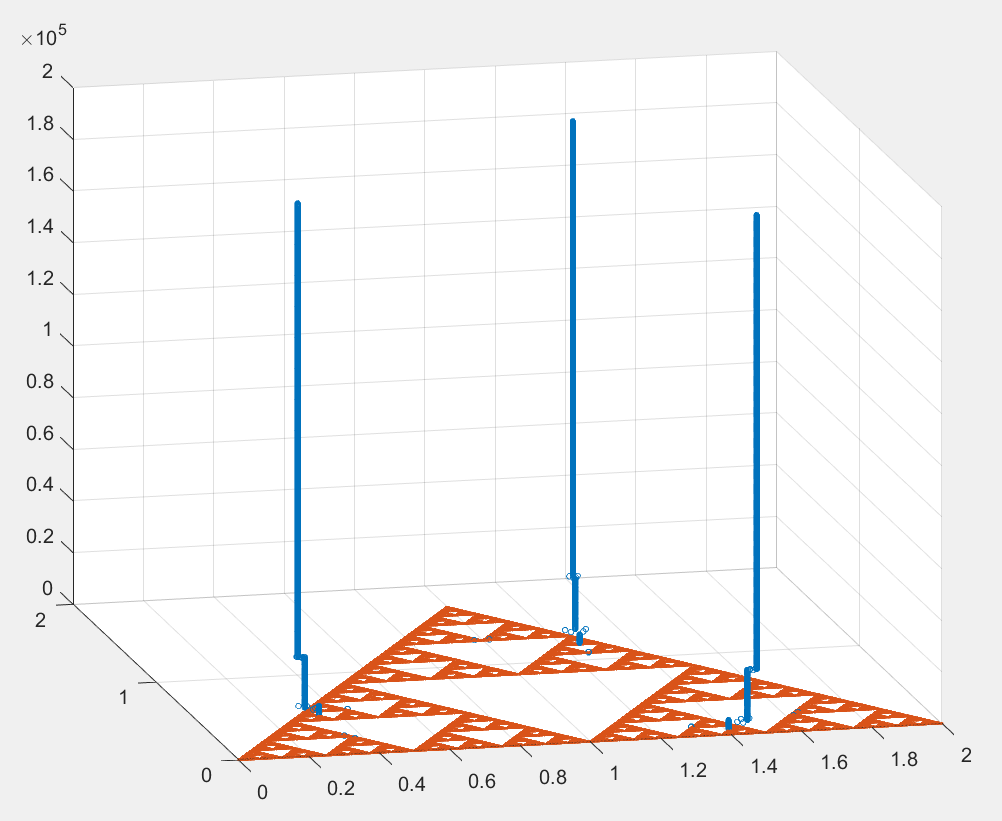}-
    \includegraphics[height=6cm]{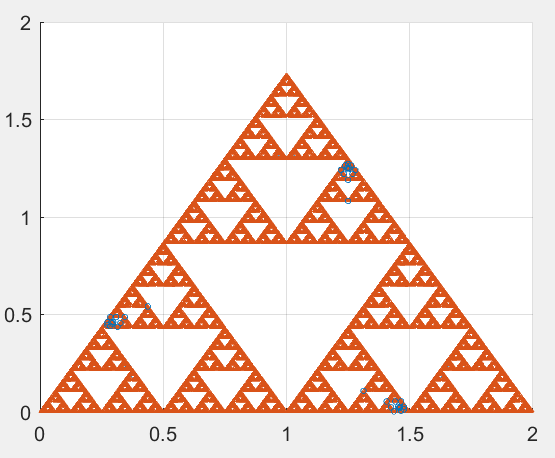}
    \caption{The position of peaks for version 3, 5 series.}
    \label{peaks version 3, 5series}
\end{figure}

\begin{figure}[!h]
    \centering
    \includegraphics[height=6cm]{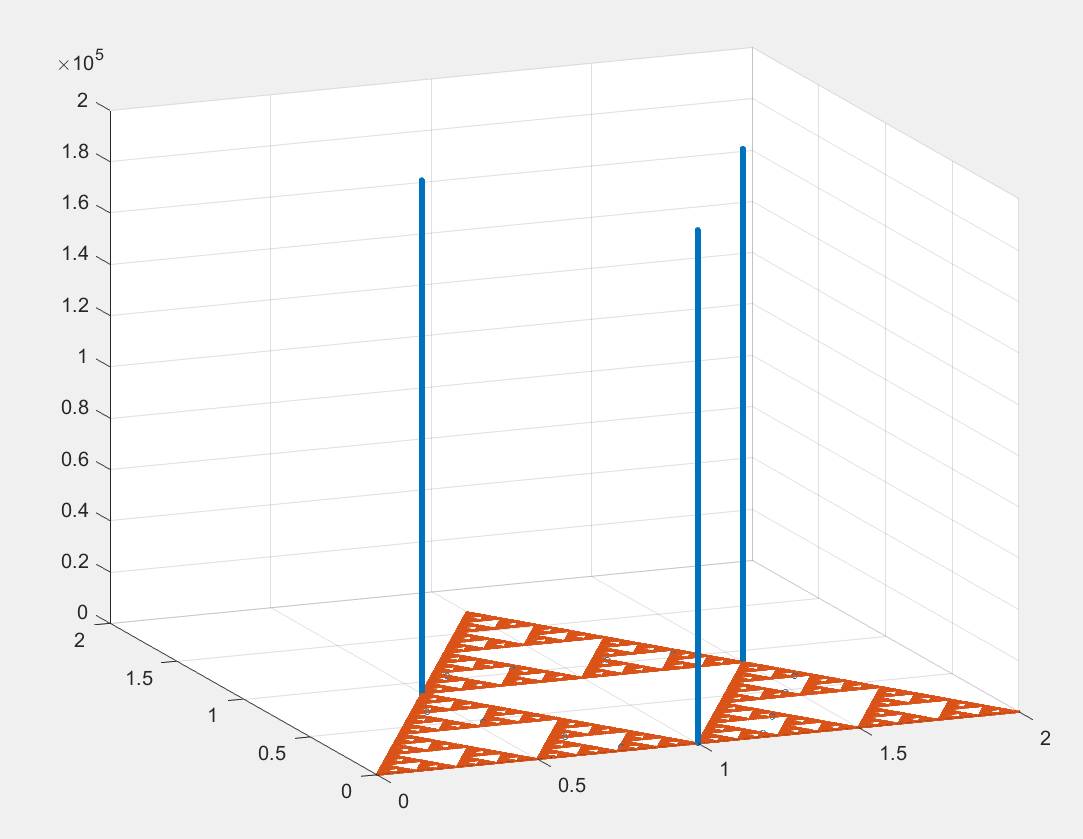}
    \includegraphics[height=6cm]{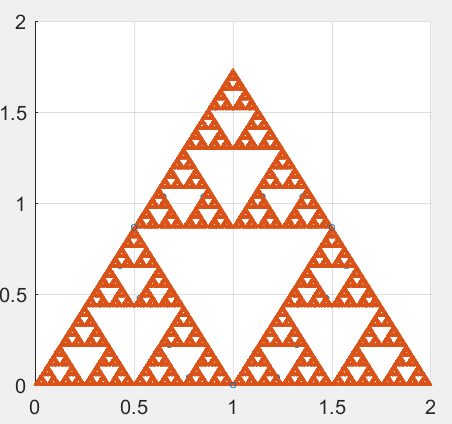}
    \caption{The position of peaks for version 3, 6 series.}
    \label{peaks version 3, 6series}
\end{figure}

\clearpage{}

\begin{figure}[!h]
    \centering
    \includegraphics[height=6cm]{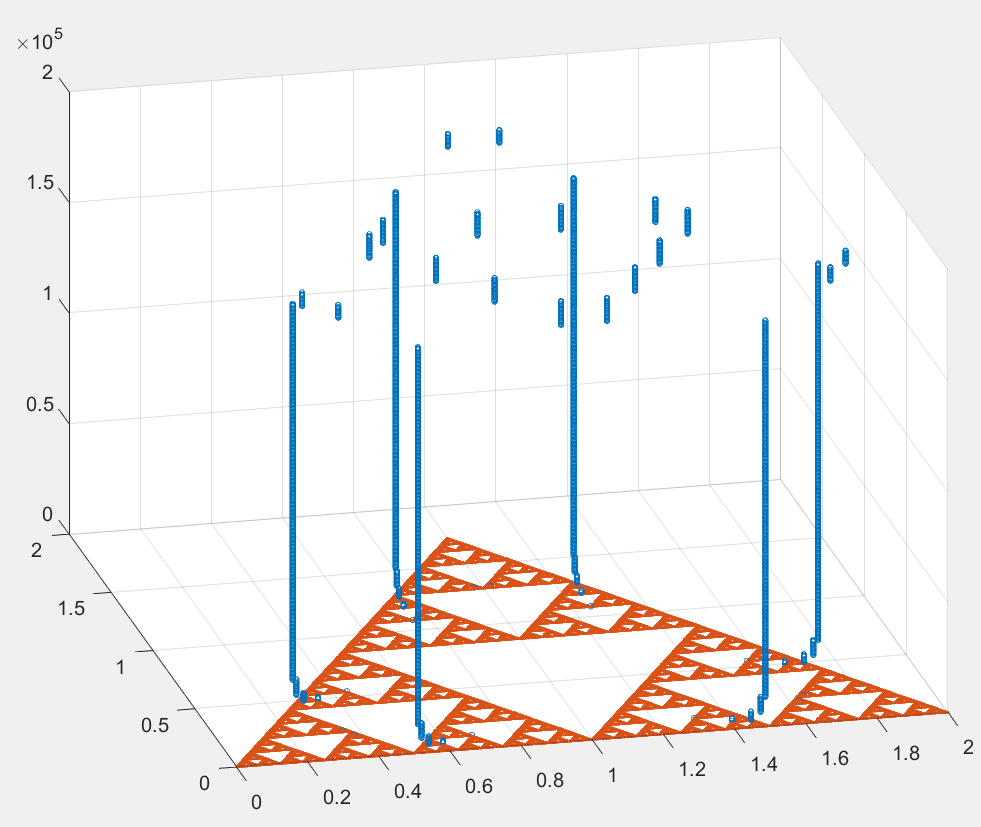}-
    \includegraphics[height=6cm]{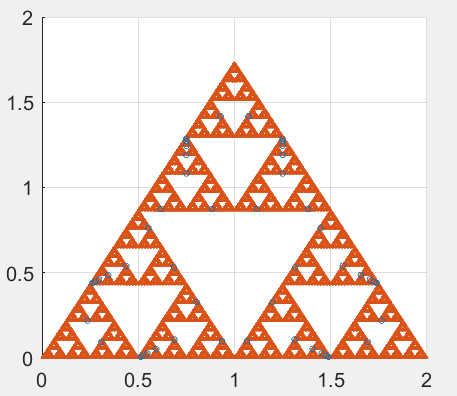}
    \caption{The position of peaks for version 4, 5 series.}
    \label{peaks version 4, 5series}
\end{figure}

\begin{figure}[!h]
    \centering
    \includegraphics[height=6cm]{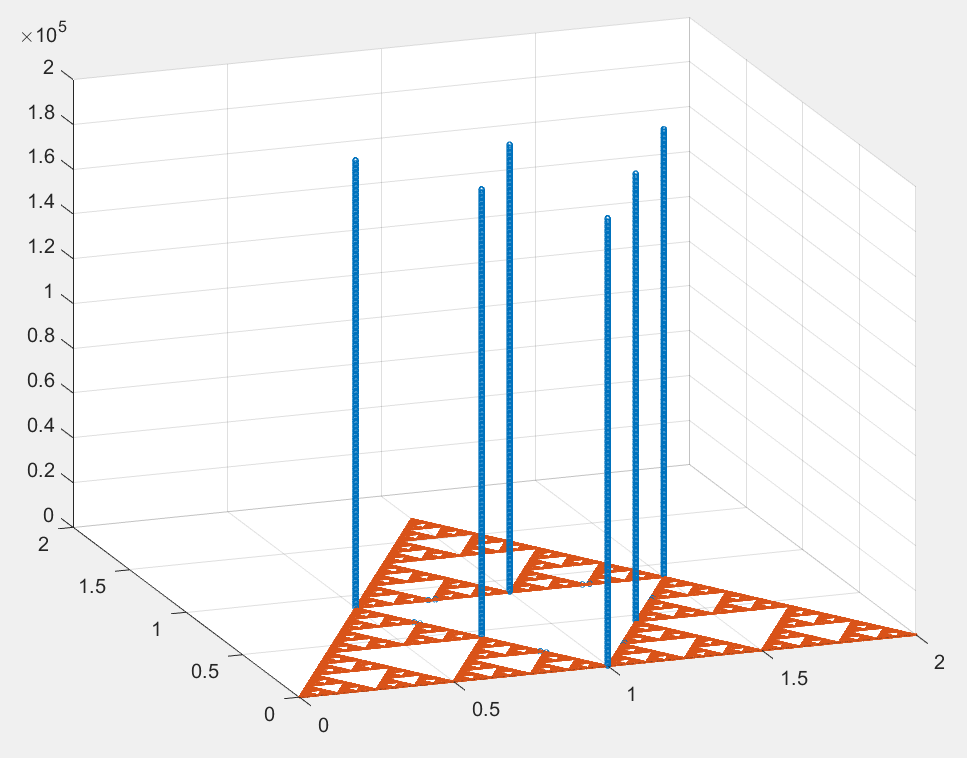}
    \includegraphics[height=6cm]{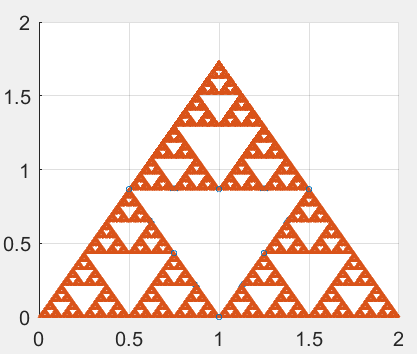}
    \caption{The position of peaks for version 4, 6 series.}
    \label{peaks version 4, 6series}
\end{figure}

\section{Further Research}
\label{furtherresearch}

Now we discuss further research questions that can be investigated in this area. 

\begin{comment}
\subsection{Probabilistic considerations}

One could investigate, both for the line case and for the fractal case, what the \textit{probability} is of choosing, at random, a pair $\left(\delta_0,\varepsilon_0\right)$ in the subset $\{(\delta,\varepsilon):x>|y|\}$ of the $\delta$-$\varepsilon$ plane such that $(\delta_0,\varepsilon_0)$ lies in a stable region.
\end{comment}
\begin{enumerate}

\begin{comment}Aside from the Mathieu differential equation, another well-known form of the Hill equation is the Meissner equation, given by

$$\frac{d^2u}{dx^2}+\left(\alpha^2+\omega^2\sgn\cos(x)\right)u=0.$$

A fractal analog which could be investigated in future research is one in which, on the blow-up of $SG$, $\sgn\cos(x)\cdot u$ could be defined to be $\pm u,$ where choice of plus or minus varies depending on each `copy' of $SG$ in the blow-up.
\end{comment}

\item \textbf{Asymptotic behavior of $SG_\infty$ transition curves}

One could further investigate the asymptotic behavior of the $SG_\infty$ transition curves. For the line case, we have Theorem 4.2 due to W.S. Loud, stating that $$\delta=-\ep+\left(k-\frac{1}{2}\right)\sqrt{2\ep}+O(\ep^{1/2})$$ as $\ep\to +\infty$ on the $k-$th transition curve. We can investigate whether estimates of a similar form hold on $SG_\infty$. This would include extending Proposition \ref{prop71} to Version 3 and Version 4 the MDE on $SG_\infty$.

\begin{comment}
\item \textbf{Convergence?}
We conjecture that the solutions in 
\end{comment}

\begin{comment}
\item \textbf{Section 7.2 for other transition curves}
\end{comment}

\item \textbf{Other modifications of MDE matrix}

Further research could investigate  different matrix versions of the fractal MDE, aside from Versions 1-4 presented in Section 6.

\item \textbf{Other Fractal Domains}

One could investigate how the Mathieu differential equation could be extended to other infinite fractafolds. Such infinite fractafolds may or may not be based on $SG$.

\item\textbf{The Hill Equation}

The Mathieu differential equation is actually a special case of the so-called \textit{Hill differential equation}. The Hill differential equation is given by
$$\frac{d^2u}{dt^2}+f(t)u=0,$$
where $f(t)$ is an arbitrary periodic function. Readers can read \cite{stoker,levykeller,loud,verticalconvergence} for details. Further research could investigate ways in which to extend, for various choices of $f(x)$, the Hill differential equation to be defined on fractal domains.

\end{enumerate}

\section{Appendix}
\label{Appendix}
In this appendix we derive equations in matrix form for solutions to the MDE on the line which have period $2N\pi$, $N\in\mathbb{N}$. We will use a similar procedure as in Section 2 and will obtain tridiagonal matrices. 

Let $u$ be a solution with period $2N\pi$. We can write the Fourier expansion of $u$ as 
\begin{equation}
   u(t)=\sum_{j=0}^\infty a_j\cos\left(\frac{j}{N}t\right)+ \sum_{j=1}^\infty b_j\sin\left(\frac{j}{N}t\right).
\end{equation}

Plugging in this Fourier series for $u$ into the Mathieu differential equation, we obtain two infinite systems of linear homogeneous equations for the cosine and sine coefficients, respectively, as follows:

\[
\text{cosine coefficients} 
\begin{cases}
    \delta a_0+\frac{\varepsilon}{2}a_N=0,\\
    (\delta-j^2) a_j+ \frac{\varepsilon}{2}(a_{N-j}+a_{N+j})=0&\text{ if } 1\leq j\leq N-1,\\
    (\delta-N^2)a_N+\frac{\varepsilon}{2}(2a_{0}+a_{2N})=0,\\
    (\delta-j^2)a_j+\frac{\varepsilon}{2}(a_{j-N}+a_{j+N})=0&\text{ if }j\geq N+1,
\end{cases}\]
and
\[
\text{sine coefficients} 
\begin{cases}
    (\delta-j^2)b_j+\frac{\varepsilon}{2}(b_{N+j}-b_{N-j})=0&\text{ if } 1\leq j\leq N-1,\\
    (\delta-N^2)b_N+\frac{\varepsilon}{2}b_{2N}=0,\\
    (\delta-j^2)b_j+\frac{\varepsilon}{2}(b_{j-N}+b_{j+N})=0&\text{ if }j\geq N+1.
 \end{cases}
\]

We can immediately write the equations in the matrix form, but the result is unwieldy. So, we make further classifications of the coefficients, just as we did in Section 2.2, where coefficients are separated into two classes---those with even indices and those with and odd indexes.

A natural criterion is to have coefficients that appear in the same equation be in a same class. For example, $a_0$ and $a_N$ should be in a same class. Based on this idea, we say that $a_j$ and $a_{j'}$ (or $b_j$ and $b_{j'}$) are in the `same class' if and only if there is a finite sequence
\[a_j=a_{j_0},a_{j_1},a_{j_2},\cdots,a_{j_L}=a_{j'},\]
such that $a_{j_l}$ and $a_{j_{l+1}}$ appear in a same equation for all $0\leq l\leq L-1$. For example, $a_0$ and $a_{3N}$ are in the same class, since we $a_0$ and $a_{N}$ are in the same equation, $a_{N}$ and $a_{2N}$ are in the same equation, and $a_{2N}$ and $a_{3N}$ are in the same equation. It is easy to check that the the property of being in the `same class' is an equivalence relation on the set $\{a_j:j\geq0\}$ and on the set $\{b_j:j\geq1\}$. This partitions $\{a_j:j\geq0\}$ into $[\frac{N}{2}]+1$ different equivalence classes and partitions $\{b_j:j\geq0\}$ into $[\frac{N}{2}]+1$ different equivalence classes. So we can write equations for different equivalence classes separately. Below, we discuss all the possible cases. First we give a discussion for the cosine coefficients $\{a_j\}$. As we will see, the case for the sine coefficients $\{b_j\}$ is similar.

There are three possible forms that the matrix corresponding to any particular equivalence class can take:

\begin{comment}
\st{A natural criterion is to have coefficients that appear in a same equation be in a same class. For example, $a_0$ and $a_N$ should be in a same class. Based on this idea, we say that $a_j$ and $a_{j'}$ (or $b_j$ and $b_{j'}$) are in \textit{a same class} if and only if there is a finite sequence
$a_j=a_{j_0},a_{j_1},a_{j_2},\cdots,a_{j_L}=a_{j'},$
such that $a_{j_l}$ and $a_{j_{l+1}}$ appear in a same equation. For example, $a_0$ and $a_{3N}$ are in a same class, since we have $a_0,a_{2N}$ in a same equation, and $a_{2N},a_{3N}$ in a same equation. As a result, all the coefficients in a class are related by equations, while coefficients in different class are not related. So we can write equations for different class separately. Below, we discuss all possible cases.}

\st{Let's first consider the equations for $\{a_j\}$. There are $[\frac{N}{2}]+1$ different classes of coefficients, with three possible forms of matrices.}
\end{comment}

1). The first possible equivalence class is $\{a_0,a_N,a_{2N},\cdots\}$, with the corresponding matrix equation

\[
\begin{pmatrix}
    \delta       & \frac{\varepsilon}{2} &  & \\ 
    \varepsilon & \delta-1 & \frac{\varepsilon}{2}   &   &  \\
           & \frac{\varepsilon}{2} & \delta-2^2 & \frac{\varepsilon}{2} &  &  \\
           & & \ddots & \ddots & \ddots \\ \\ 
\end{pmatrix}
\begin{pmatrix}
    a_0 \\
    a_N \\
    a_{2N} \\
    \vdots \\
    \vdots \\
\end{pmatrix}
=\begin{pmatrix}
    0 \\
    0 \\
    0 \\
    \vdots \\
    \vdots \\
\end{pmatrix}.
\] 
Note that this is just the equation $Ax=0$ in Section 2.2. If there is a nontrivial solution to the above equation, then the MDE has a $2\pi$ periodic solution, since $u(t)=\sum_{l=0}^\infty a_{lN}\cos (lt)$ solves the MDE.
%This matrix is actually for solutions of period $2\pi$, since we have the solution $u(t)=\sum_{l=0}^\infty a_{lN}\cos (lt)$ in this case.\textcolor{blue}{\st{[I think this only applies if the coefficients not in this equivalence class are all zero.]}}\textcolor{green}{[We have a nontrivial solution if and only if one of these matrix equations has a nontrivial solutions. The basis of the solution space is given by these forms, ... if other coefficients are nonzero, then we can write it as a combination of serveral solutions in this form. Delete this after you read.]} \st{So, we do not get something new here.}\textcolor{blue}{[?]}\textcolor{green}{[I mean we already have this equation in section 2.]}

2). The second form that an equivalence class can take is $$\{a_{|k-lN|}\}_{l=-\infty}^\infty=\{a_k,a_{k+N},a_{k_{2N},\cdots}\}\cup \{a_{N-k},a_{2N-k},\cdots\},$$ where $1\leq k<\frac{N}{2}$. The corresponding equations can be written in the form 
\[
\begin{pmatrix}
\quad\ddots & \ddots & \ddots\\
    & & \delta-\left(\frac{2N-k}{N}\right)^2       & \frac{\varepsilon}{2} &  & \\ 
    & & \frac{\varepsilon}{2} & \delta-\left(\frac{N-k}{N}\right)^2 & \frac{\varepsilon}{2}   &   &  \\
           & & & \frac{\varepsilon}{2} & \delta-\left(\frac{k}{N}\right)^2 & \frac{\varepsilon}{2} &  &  \\
         & & & & \frac{\varepsilon}{2} & \delta-\left(\frac{N+k}{N}\right)^2 & \frac{\varepsilon}{2} &  &  \\
           & & & & & \ddots & \ddots & \ddots \\
\end{pmatrix}
\begin{pmatrix}
    \vdots \\
    a_{2N-k} \\
    a_{N-k} \\
    a_k \\
    a_{k+N} \\
    \vdots \\
\end{pmatrix}
=\begin{pmatrix}
    \vdots \\
    0 \\
    0 \\
    0 \\
    0 \\
    0 \\
    \vdots \\
\end{pmatrix}.
\]
If there is a nontrivial solution to the above equation, then $u(t)=\sum_{l=-\infty}^\infty a_{k+lN}\cos (\frac{k+lN}{N}t)$ is a $\frac{2N\pi}{\text{gcd}(N,K)}$-periodic solution to the MDE, where $\text{gcd}(N,k)$ denotes the greatest common divisor of the pair $(N,k)$.

In particular, we get a $2N\pi$-periodic solution to the MDE from the equation above if and only if $k$ and $N$ are coprime. Thus, the matrices of the second form with $k$ and $N$ coprime yield $2N\pi$-periodic solutions.

3) The third form that an equivalence class can take is  $\{a_{\frac{N}{2}},a_{\frac{3}{2}N},a_{\frac{5}{2}N},\cdots\}.$ This form can only occur if $N$ is even. The corresponding equations can be written in the form

\[
\begin{pmatrix}
    \delta-\frac{1}{4}+\frac{\varepsilon}{2}       & \frac{\varepsilon}{2} &  & \\ 
    \frac{\varepsilon}{2} & \delta-\frac{3^2}{4} & \frac{\varepsilon}{2}   &   &  \\
           & \frac{\varepsilon}{2} & \delta-\frac{5^2}{4} & \frac{\varepsilon}{2} &  &  \\
           & & \ddots & \ddots & \ddots \\ \\ 
\end{pmatrix}
\begin{pmatrix}
    a_{\frac{N}{2}} \\
    a_{\frac{3N}{2}} \\
    a_{\frac{5N}{2}} \\
    \vdots \\
    \vdots \\
\end{pmatrix}
=\begin{pmatrix}
    0 \\
    0 \\
    0 \\
    \vdots \\
    \vdots \\
\end{pmatrix}.
\]
The above matrix is exactly the matrix $C$ in Section 2.2, which yield $4\pi$-periodic solutions.\\

\begin{comment}
\textit{To conclude, for solutions of period $2N\pi, N\geq 3$, the equations for cosine coefficents have the following form}

\[
\begin{bmatrix}
\quad\ddots & \ddots & \ddots\\
    & & \delta-(\frac{2N-k}{N})^2       & \frac{\varepsilon}{2} &  & \\ 
    & & \frac{\varepsilon}{2} & \delta-(\frac{N-k}{N})^2 & \frac{\varepsilon}{2}   &   &  \\
           & & & \frac{\varepsilon}{2} & \delta-(\frac{k}{N})^2 & \frac{\varepsilon}{2} &  &  \\
         & & & & \frac{\varepsilon}{2} & \delta-(\frac{N+k}{N})^2 & \frac{\varepsilon}{2} &  &  \\
           & & & & & \ddots & \ddots & \ddots \\
\end{bmatrix}
\begin{bmatrix}
    \vdots \\
    a_{2N-k} \\
    a_{N-k} \\
    a_k \\
    a_{k+N} \\
    \vdots \\
\end{bmatrix}
=\begin{bmatrix}
    \vdots \\
    0 \\
    0 \\
    0 \\
    0 \\
    0 \\
    \vdots \\
\end{bmatrix},
\]
\textit{where $k$ can by any integer such that $k$ and $N$ are coprime.}
\end{comment}

The matrix forms for sine coefficients are quite similar to those for the cosine coefficients. For solutions of period $2N\pi, N\geq 3$,  the equations for sine coefficents have the following form
\[
\begin{bmatrix}
\quad\ddots & \ddots & \ddots\\
    & & \delta-(\frac{2N-k}{N})^2       & \frac{\varepsilon}{2} &  & \\ 
    & & \frac{\varepsilon}{2} & \delta-(\frac{N-k}{N})^2 & -\frac{\varepsilon}{2}   &   &  \\
           & & & -\frac{\varepsilon}{2} & \delta-(\frac{k}{N})^2 & \frac{\varepsilon}{2} &  &  \\
         & & & & \frac{\varepsilon}{2} & \delta-(\frac{N+k}{N})^2 & \frac{\varepsilon}{2} &  &  \\
           & & & & & \ddots & \ddots & \ddots \\
\end{bmatrix}
\begin{bmatrix}
    \vdots \\
    b_{2N-k} \\
    b_{N-k} \\
    b_k \\
    b_{k+N} \\
    \vdots \\
\end{bmatrix}
=\begin{bmatrix}
    \vdots \\
    0 \\
    0 \\
    0 \\
    0 \\
    0 \\
    \vdots \\
\end{bmatrix},
\]
where $k$ can by any integer such that $k$ and $N$ are coprime.

\begin{comment}
Assume that $N$ is even (the case for odd $N$ is treated similarly). For $1\leq m\leq N,$ let $A_m:=\{a_k:k\equiv m (\text{mod} N)\}.$ We partition the set of coefficients $\{a_k\}_{k=0}^\infty$ into the following $\frac{N}{2}+1$ equivalence classes:
\begin{itemize}
    \item $A_N$
    \item $A_1\cup A_{N-1}$
    \item $A_2\cup A_{N-2}$
    \item $A_3\cup A_{N-3}$
    \item $\quad\quad\vdots$
    \item $A_{N/2-1}\cup A_{N/2+1}$
    \item $A_{N/2}$
\end{itemize}

%    \item $A_{N/2-1}\cup A_{N/2+1}$ if $N$ is even, $A_{(N-1)/2-1}\cup A_{(N+1)/2+1}$ if $N$ is odd
%    \item $A_{N/2}$ if $N$ is even, or $A_{(N-1)/2}\cup A_{(N+1)/2}$ if $N$ is odd

We can create a matrix for each of these $\ceil*{\frac{N}{2}}+1$ equivalence classes of coefficients.

In order to show how this works, it is best to present an example. We show the case for $N=4$ (which corresponds to a period of $2\pi N=8\pi$). We have
\begin{itemize}
    \item $A_4=\{a_0,a_4,a_8,...\}$\\
    \item $A_1=\{a_1,a_5,a_9,...\}$\\
    \item $A_2=\{a_2,a_6,a_{10},...\}$\\
    \item $A_3=\{a_3,a_7,a_{11},...\}$
\end{itemize}

For this case, we have $\ceil
*{\frac{N}{2}}+1=\ceil*{\frac{3}{2}}+1=3$ equivalence classes:
\begin{itemize}
    \item $A_4=\{a_0,a_4,a_8,...\}$\\
    \item $A_1\cup A_3=\{a_1,a_3,a_5,a_7,a_9...\}$\\
    \item $A_2=\{a_2,a_6,a_{10},...\}$\\
\end{itemize}

\end{comment}

\section*{Acknowledgements}
This research was hosted by the Cornell University Department of Mathematics through its 2018 Summer Program for Undergraduate Research. Anthony Coniglio's participation in this research was partly supported by Indiana University Bloomington. Xueyan Niu's participation in this research was partly supported by the Overseas Research Fellowship (ORF) of Faculty of Science, University of Hong Kong.

\nocite{*}
\printbibliography

\end{document}